\documentclass[11pt]{amsart}
\usepackage[margin=1in,includeheadfoot]{geometry}
\usepackage{enumerate}
\usepackage{comment}
\usepackage{bm}
\usepackage{tikz-cd}
\usepackage{stmaryrd}
\usepackage{amsfonts}
\usepackage{amssymb}
\usepackage{amsmath,amssymb,amsthm}
\usepackage{latexsym}
\usepackage{amscd}
\usepackage{esint}
\usepackage{mathrsfs}
\usepackage{dsfont}
\usepackage{indentfirst}
\usepackage{enumitem}
\newlist{Aenumerate}{enumerate}{1}
\setlist[Aenumerate]{label=(S\arabic*)}

\usepackage{enumerate}
\usepackage{slashed}
\usepackage{fancyhdr}
\usepackage{pdfsync}
\usepackage{setspace}
\usepackage{color}
\usepackage{xcolor}
\usepackage{graphicx}
\definecolor{revisionblue}{RGB}{0,0,0}
\definecolor{repairblue}{RGB}{0,0,0}
\newcommand{\repair}[1]{{\color{repairblue}#1}}

\usepackage{latexsym}
\usepackage{float}
\usepackage{needspace}
\definecolor{citeblue}{RGB}{0,0,0}
 \usepackage[colorlinks = true,
            linkcolor = citeblue,
            urlcolor  =  citeblue,
            citecolor =  citeblue,
            anchorcolor = citeblue]{hyperref}
            
\theoremstyle{plain}
\newtheorem{theorem}{Theorem}[section]
\newtheorem{coro}[theorem]{Corollary}
\newtheorem{lemma}[theorem]{Lemma}
\newtheorem{prop}[theorem]{Proposition}
\theoremstyle{definition}
\newtheorem{definition}[theorem]{Definition}

\theoremstyle{definition}

\newtheorem{rmk}[theorem]{Remark}
\numberwithin{equation}{section}
\numberwithin{theorem}{section}

\newcommand{\reg}{\mathrm{reg}}
\newcommand{\loc}{\mathrm{loc}}
\newcommand{\diam}{\operatorname{diam}}

\newcommand{\vol}{\operatorname{vol}}

\newcounter{stepnum}[section]

\newcounter{claimnum}[subsection]

\newcounter{casenum}[section]

\usepackage{etoolbox}

\makeatletter
\newcommand{\reduceoperator}[1]{%
  \@for\next:=#1\do{\expandafter\reduceoperator@\expandafter{\next}}%
}
\newcommand{\reduceoperator@}[1]{%
  \csletcs{normal@#1@}{#1@}%
  \csedef{#1@}{\noexpand\reduceoperator@@\csname normal@#1@\endcsname}%
}
\newcommand{\reduceoperator@@}[1]{%
  \mathop{\mathpalette\reduceoperator@@@{#1}}%
}
\newcommand{\reduceoperator@@@}[2]{%
  \ifx#1\displaystyle\textstyle\fi#2%
}
\makeatother

\reduceoperator{bigcup,bigcap,bigotimes}

\def\R{\mathbb R}
\def\S{\mathbb S}

\selectfont
\def\contraction{\mathop{\hbox{\vrule height 7pt width .5pt 
depth 0pt\vrule height .5pt width 6pt depth 0pt}}\nolimits}

\newcommand{\C}{\mathbb{C}}
\newcommand{\N}{\mathbb{N}}
\newcommand{\dd}{\,\mathrm{d}}
\newcommand{\cA}{\mathcal{A}}
\newcommand{\cB}{\mathscr{B}}

\newcommand{\cK}{\mathcal{K}}

\newcommand{\cN}{\mathcal N}
\newcommand{\dist}{\operatorname{dist}}

\allowdisplaybreaks[4]
\makeatletter
\def\@tocline#1#2#3#4#5#6#7{\relax
	\ifnum #1>\c@tocdepth % then omit
	\else
	\par \addpenalty\@secpenalty\addvspace{#2}%
	\begingroup \hyphenpenalty\@M
	\@ifempty{#4}{%
		\@tempdima\csname r@tocindent\number#1\endcsname\relax
	}{%
		\@tempdima#4\relax
	}%
	\parindent\z@ \leftskip#3\relax \advance\leftskip\@tempdima\relax
	\rightskip\@pnumwidth plus4em \parfillskip-\@pnumwidth
	#5\leavevmode\hskip-\@tempdima
	\ifcase #1
	\or\or \hskip 2em \vspace{0.5mm}\or \hskip 3em \else \hskip 3em \fi%
	#6\nobreak\relax
	\hfill\hbox to\@pnumwidth{\@tocpagenum{#7}}\par% <- \dotfill -> \hfill
	\nobreak
	\endgroup
	\fi}
\makeatother

\begin{document}
\title[Compactness of branched conformal immersions with free boundary]{Compactness of branched conformal immersions with free boundary}

    \author[Gao]{Rui Gao}
    \address{School of Mathematics and Statistics, Shaanxi Normal University, 620 West Chang'an Avenue, Xi'an, Shaanxi, 710119, P. R. China}
	\email{rui.gao@snnu.edu.cn}

    \author[Wang]{Zijun Wang}
    \address{Institute for Theoretical Sciences, Westlake University\\600 Dunyu Road  \\ Hangzhou, 310030 \\P. R. China}
    \email{wangzijun79@westlake.edu.cn}

	\author[Zhu]{Miaomiao Zhu}
	\address{School of Mathematical Sciences, Shanghai Jiao Tong University\\ 800 Dongchuan Road \\ Shanghai, 200240 \\ P. R. China}
	\email{mizhu@sjtu.edu.cn}
	
	 \thanks{}
	 \date{\today}

\begin{abstract}
We prove bubble tree compactness for branched conformal immersions with free boundary under uniform Dirichlet and Willmore energy bounds, with images in a fixed compact subset of $\R^K$ and possibly degenerating conformal structures. In particular, we prove that the convergence satisfies the Dirichlet energy identity and the no-neck property. For smooth Lagrangian immersions in K\"ahler surfaces, we characterize the natural free boundary conditions for stationarity of the area functional under exact infinitesimal variations on the domain. We then prove bubble tree compactness for Hamiltonian stationary Lagrangian (HSL) immersions in a complete K\"ahler surface under the orthogonal free boundary condition and the zero normal trace boundary condition for the associated mean curvature $1$-form. Finally, as applications of the bubble tree convergence results, we obtain Willmore energy minimizers for short free boundary annuli under a Steklov constraint and concentration of free boundary disk minimizers near maxima of the support mean curvature under a small area constraint. We also derive a Willmore energy lower bound under a calibration constraint and an intrinsic Willmore energy decomposition for free boundary HSL sequences in K\"ahler-Einstein surfaces.
\end{abstract}
\maketitle

{\tableofcontents}
\section{Introduction}

Let $(N,\bar g)$ be a Riemannian manifold, and let $\cK\subset N$ be a smooth support submanifold.
Let $(\Sigma,h)$ be a compact Riemannian surface with boundary.
A smooth immersion $u:\Sigma\to N$ is a free boundary minimal surface if it is a critical point of the area functional $\cA$ under variations satisfying $u_t(\partial\Sigma)\subset\cK$.
The first variation formula gives
\begin{equation*}
    \mathbf H_{u,N}=0 \quad\text{in }\Sigma,\qquad
    u(\partial\Sigma)\subset\cK,
    \qquad
    du(\nu)\perp T_u\cK \quad\text{on }\partial\Sigma ,
\end{equation*}
Here $\mathbf H_{u,N}$ is the trace mean curvature vector of $u$ in $N$, and $\nu$ is the outward unit conormal along $\partial\Sigma$ with respect to the induced metric $u^*\bar g$.
The orthogonality condition $du(\nu) \perp T_u\cK$ follows from the boundary term in the first variation formula.

The conformal approach to surfaces with free boundary goes back to the classical theory of free boundary minimal surfaces.
Building on the conformal treatment of Plateau's problem by Douglas \cite{douglas1931mapping} and Rad\'{o} \cite{rado1930plateau,rado1932contributions}, Courant \cite{Courant1950book} and Courant--Davids \cite{Courant-Davids} studied minimal surfaces whose boundaries are allowed to move on a prescribed support.
Lewy \cite{Lewy1951CPAM,Lewy1951PNAS} and Hildebrandt--Nitsche \cite{Hildebrandt-Nitsche1979} studied boundary regularity and the behavior of the free boundary.
Struwe \cite{Struwe1984invention} obtained a nonconstant conformal critical point of the Dirichlet energy among maps from the disk with boundary on a smooth embedded $2$-sphere in $\R^3$, yielding a minimal disk that meets the supporting $2$-sphere orthogonally.
For other early existence results, including embedded disks in convex bodies and area minimizing surfaces under topological assumptions, see \cite{Smyth1984,Gruter-Jost,Ye1991}.

% In the unit ball, free boundary minimal surfaces are related to Steklov eigenvalue and rigidity problems \cite{Fraser-Schoen-adv2011,Fraser-Schoen-2016,Fraser-Schoen2019,Fraser-Schoen2020}.
% Variational and min-max constructions appear in \cite{Marin-CPAM-2015,DeLellis-Ramic2018,Martin-ZhoujdgI,Martin-ZhoujdgII,haslhofer2023freeboundaryminimaldisks}.
% Compactness, curvature concentration, and bubbling for embedded free boundary minimal surfaces and hypersurfaces are studied in \cite{fraser2014,sharp2018,AmbrozioBuzanoCarlottoSharp2019}.
% See \cite{Carlotto2019,Martin-servey2017} for broader surveys.

For a conformal parametrization, the Dirichlet energy agrees with area, and harmonicity is equivalent to minimality away from branch points.
The free boundary conditions take the form $u(\partial\Sigma)\subset\cK$ and $du(\nu)\perp T_u\cK$.
This formulation turns the free boundary minimal surface problem into a variational problem for maps and applies in higher codimension as well.
Fraser \cite{FraserCPAM} developed a partial Morse theory using the Sacks-Uhlenbeck approximation \cite{sacks1981existence} for harmonic maps from the disk with free boundary on a fixed compact support.
For a complete homogeneously regular target $N$ and compact support $\cK$, she proved that $\pi_k(N,\cK)\neq0$, $k\geq2$, implies the existence of either a nonconstant free boundary minimal disk or a minimal $2$-sphere with Morse index at most $k-2$.
The proof relies on blow-up analysis of critical points of the Sacks-Uhlenbeck $\alpha$-energy as the exponent $\alpha>1$ tends to $1$. Both interior sphere bubbles and boundary disk bubbles can occur.
If the weak limit is constant, rescaling at a concentration point produces a nonconstant minimal sphere or disk. Localized variations then give the Morse index bound.
Fraser \cite{Fraser2002} later used this existence and index theory to study the topology of two-convex domains.

Bubble tree convergence also plays a central role in the min-max construction of Lin--Sun--Zhou \cite{Lin-Sun-Zhou2020}; see also Laurain--Petrides \cite{Laurain-Petrides}.
Using the harmonic replacement method of Colding--Minicozzi \cite{Colding-Minicozzi08}, they constructed free boundary minimal disks and minimal spheres from suitable min-max sequences. They proved that the total area of the limiting surfaces equals the min-max value.
The bubble tree accounts for concentration at different scales, while the energy identity excludes energy loss in the connecting necks.
For bordered surfaces of more general topology, Chen--Fraser--Pang \cite{Chen-Fraser-Pang2015} minimized over maps and conformal structures to obtain branched minimal immersions with free boundary under incompressibility assumptions.
In their argument, incompressibility prevents degeneration of the bordered domains and thus ensures compactness of the conformal structures.

Concentration compactness with bubbling also describes singularity formation in evolution problems.
Jost--Liu--Zhu \cite{Jost-Liu-Zhu2019} proved the energy identity and the no-neck property for approximate harmonic maps from a fixed bordered surface into a compact target, under free boundary conditions, a uniform energy bound, and an $L^2$ bound on the tension fields.
They applied these results to the harmonic map heat flow with free boundary, obtaining energy identities at finite and infinite singular times.
When the domain conformal structures degenerate, further concentration may occur along the degenerating collars. Liu--Song--Zhu \cite{Liu-Song-Zhu2022} established a generalized energy identity for harmonic maps from degenerating bordered Riemann surfaces that also accounts for this collar energy.

These variational constructions and evolution problems motivate a compactness theory for conformally parametrized surfaces with free boundary under geometric energy bounds.
We study sequences $u_n:(\Sigma,h_n)\to\R^K$ of branched conformal immersions with free boundary on a fixed smooth compact support $\cK\subset\R^K$ and uniformly bounded Dirichlet and Willmore energies. The conformal structures of $(\Sigma,h_n)$ are allowed to degenerate. The maps $u_n$ need not satisfy the minimal surface or harmonic map equations. This setting includes branched minimal immersions with free boundary, for which $\mathbf H_{u_n}=0$, as well as other surfaces with finite total curvature and Hamiltonian stationary Lagrangian surfaces. More precisely, in Theorem~\ref{thm:intro-main-compactness}, we prove bubble tree compactness for $\{u_n\}$ under these energy bounds when all images $u_n(\Sigma)$ lie in a fixed compact set, including when the bordered conformal structures degenerate. Unless the images $u_n(\Sigma)$ collapse to a point, we obtain a subsequence converging to a continuous map on a finite connected bordered stratified surface made up of principal components, sphere bubbles, and free boundary disk bubbles. In particular, we establish the Dirichlet energy identity and the no-neck property for $\{u_n\}$ and prove lower semicontinuity of the Willmore energy. 

Two principal difficulties arise in the free boundary analysis.
First, for closed surfaces in Euclidean space, the smallness condition in the local regularity estimates involves only the Willmore energy \cite{chen2014} and the corresponding positive Willmore energy gap rules out nonconstant closed minimal surfaces. For free boundary surfaces, an equatorial disk in a Euclidean ball already provides a nonconstant minimal surface with zero Willmore energy, so Willmore smallness alone does not exclude nonconstant components or boundary energy concentration. We therefore establish local boundary estimates under smallness of the sum of Dirichlet and Willmore energies and prove a gap for this combined energy.
Second, although geodesic reflection is an important tool in earlier free boundary analyses \cite{FraserCPAM,Lin-Sun-Zhou2020,Laurain-Petrides,Jost-Liu-Zhu2019,Liu-Song-Zhu2022}, reflecting $u_n$ across a curved support submanifold generally destroys conformality and adds a term quadratic in $du_n$.
We overcome this difficulty by using conformality on the original domain to derive boundary identities and absorbing the reflection errors under the combined small energy condition.

For Hamiltonian stationary Lagrangian (HSL) surfaces, the boundary constraint first requires a variational analysis of its own.
In Theorem~\ref{thm:intro-hsl-boundary}, we characterize smooth stationarity by a coupled condition for the conormal and the mean curvature 1-form.
The homogeneous case consists of orthogonal contact and zero normal trace.
In Theorem~\ref{main theorem}, we apply the compactness theory to this class with homogeneous boundary conditions in a complete K\"ahler surface and prove that every nonconstant limiting component is a smooth branched HSL immersion with the same boundary conditions.
The map equation involves the mean curvature 1-form $\alpha_n$, while the right-hand side of the Hodge system for $\alpha_n$ depends on $du_n$.
Since the Willmore bound gives only $L^2$ control of $\alpha_n$, the boundary estimates for the map and the 1-form must be obtained together.
We reflect this coupled homogeneous system and establish higher boundary estimates under small Dirichlet energy, without requiring small Willmore energy.

In Theorems~\ref{app:small-annuli}-\ref{app:calibration-barrier}, we apply Theorem~\ref{thm:intro-main-compactness} to Willmore energy minimization problems and to the curvature cost of contact with a support disjoint from a calibrated surface.
The boundary length measures in Lemma~\ref{lem:app-contact-measures} retain a prescribed boundary length, and a normalized Steklov lower bound prevents an annulus from concentrating into disk components.
For disks of small prescribed area, the area identity gives global weak minimizers; reflection and quantitative sphere rigidity then locate every minimizer near the maximum set of the support mean curvature.
In Theorem~\ref{app:hsl-energy}, we express the limiting Willmore energy as the sum of the Willmore energies of nonconstant principal components and additional contributions determined by squared $L^2$ norms of limiting harmonic $1$-forms on constant principal domains and of constant angular Fourier components on degenerating collars.

\subsection{Bubble tree convergence for branched conformal immersions with free boundary}
\label{subsec:fb-bubble tree-compactness}
We begin with the bubble tree convergence theorem for branched conformal immersions into Euclidean space with free boundary. Let $\Sigma$ be a connected compact oriented surface with nonempty smooth boundary, and let $\cK\subset\R^K$ be a compact smooth embedded submanifold. For a branched conformal immersion $u:(\Sigma,h)\to\R^K$ with $u\in W^{2,2}(\Sigma,\R^K)$, let $\cB_u$ be its finite branch set. On $\overline\Sigma\setminus\cB_u$, let $g_u:=u^*g_{\R^K}$ be the induced metric, $A_u$ the second fundamental form, and $\mathbf H_u=\operatorname{tr}_{g_u}A_u$ the trace mean curvature vector. We denote by $\mu_u$ the induced area measure, assigning zero mass to $\cB_u$. We use the normalizations
\begin{equation}\label{eq:intro-energy-normalization}
 E_h(u):=\frac12\int_\Sigma |du|_h^2\,\dd\vol_h
 =\cA(u),
 \qquad
 W_{\R^K}(u):=\frac14\int_{\Sigma\setminus\cB_u} |\mathbf H_u|^2\,\dd\mu_u.
\end{equation}
The free boundary conditions on a smooth supporting submanifold $\cK \subset \R^K$ are
\begin{equation*}
 u(\partial\Sigma)\subset\cK,
 \qquad du(\nu_h)\perp T_u\cK,
\end{equation*}
where $\nu_h$ is the outward unit conormal with respect to the metric $h$ and the orthogonality condition for $du(\nu_h)$ is understood in the Sobolev trace sense away from branch points.
Definition~\ref{def:fb-branched-conformal} specifies the class of weak branched conformal immersions used below.
By fixing an isometric embedding $N\hookrightarrow\R^K$, we also obtain compactness for maps into a compact Riemannian target $N$; see Remark~\ref{rem:intro-euclidean-reduction}.

Let $h_n$ be a sequence of smooth metrics on $\Sigma$, and denote by $j_n=[h_n]$ the associated conformal structures.
Assume that, after passing to a subsequence,
\begin{equation}\label{eq:intro-bordered-DM-convergence}
    (\Sigma,j_n)
    \xrightarrow{\mathrm{DM}}
    (\Sigma_0,j_0)
\end{equation}
in the bordered Deligne-Mumford sense introduced in Subsection~\ref{subsec:bordered-DM}. We denote by $\cN_0$ the nodal set of $\Sigma_0$, by $\Sigma_0^{\reg}:=\Sigma_0\setminus\cN_0$ the complement of its nodes, and by $\widetilde\Sigma_0$ its normalization.
For the convergence in \eqref{eq:intro-bordered-DM-convergence}, let $\Psi_n:V_n\to\Sigma$ be the smooth boundary preserving embeddings, where $V_n\Subset\widetilde\Sigma_0\setminus\widetilde{\cN}_0$ exhaust the complement of the node preimages $\widetilde{\cN}_0$, and $\Psi_n^*j_n\to\widetilde j_0$ smoothly on compact subsets; here $\widetilde j_0$ is the conformal structure on the normalization. We fix these embeddings and the conformal coordinates on the corresponding interior and boundary collars, as described in Subsection~\ref{subsec:bordered-DM}. For the notion of bordered stratified Riemann surfaces, see Subsection~\ref{subsec:bordered-stratified}, and for free boundary bubble tree convergence, see Subsection~\ref{subsec:bt-convergence-revised}.
We use throughout the notation and conventions introduced in Subsection~\ref{sec:fbbci}. 

With this terminology in place, we can now state our main bubble tree convergence results for branched conformal immersions with a free boundary constraint.

\begin{theorem}
\label{thm:intro-main-compactness}
\label{thm:fb-chen-li}
Let $Q\subset\R^K$ be compact, and let
\[
    u_n\in
    W^{2,2}_{b,c,\mathrm{fb}}
    \bigl((\Sigma,h_n),\R^K,\cK\bigr)
\]
satisfy $u_n(\Sigma)\subset Q$ for every $n$, and
\begin{equation}\label{eq:intro-main-bound}
    \sup_n
    \left(
        E_{h_n}(u_n)+W_{\R^K}(u_n)
    \right)
    <\infty.
\end{equation}
After passing to a subsequence, exactly one of the following alternatives holds.

\noindent
\textup{(a).} There exists a point $q\in Q$ such that $d_{\mathrm{Haus}}\bigl(u_n(\Sigma),\{q\}\bigr) \rightarrow 0.$

\noindent
\textup{(b).} There exist a finite connected bordered stratified Riemann surface $\Sigma_\infty$ and component maps $u_\infty^i: \Sigma_\infty^i\rightarrow Q$, $i\in\mathfrak I$, with at least one nonconstant component.
The surface $\Sigma_\infty$ is obtained from the compactified components of $\widetilde\Sigma_0$ by attaching finitely many sphere and free boundary disk components at energy concentration points and at the points lying over the nodes of $\Sigma_0$. The index set $\mathfrak I$ has the decomposition
\[
    \mathfrak I
    =
    \mathfrak I_{\mathrm{pr}}
    \sqcup
    \mathfrak I_{S^2}
    \sqcup
    \mathfrak I_{D},
\]
where $D=\{z\in\C:|z|<1\}$ is the unit disk, $\mathfrak I_{\mathrm{pr}}$ indexes the principal components, $\mathfrak I_{S^2}$ the sphere bubbles, and $\mathfrak I_{D}$ the free boundary disk bubbles.
Set $\mathfrak I_*:= \left\{ i\in\mathfrak I: u_\infty^i\text{ is nonconstant} \right\}.$ Then the following statements hold.
\begin{enumerate}[label=\textup{(\roman*)},leftmargin=2.4em]
\item For every $i\in\mathfrak I_*$, the map $u_\infty^i$ is a $W^{2,2}$ branched conformal immersion on $\Sigma_\infty^i$.
If $\partial\Sigma_\infty^i\neq\varnothing$, then $u_\infty^i$ satisfies
\begin{equation*}
        u_\infty^i
        (\partial\Sigma_\infty^i)
        \subset\cK,
        \qquad
        du_\infty^i(\nu_i)
        \perp
        T_{u_\infty^i}\cK
    \end{equation*}
for $\mathcal H^1$-a.e. boundary point away from the branch points of $u_\infty^i$.
Here $\nu_i$ denotes the outward unit conormal with respect to any smooth metric on $\Sigma_\infty^i$ compatible with its conformal structure.

\item For every $i\in\mathfrak I_*$ and $U\Subset\Sigma_\infty^i$ avoiding the branch points and the nodal or attaching points of that component, there are reparametrizations $\Phi_n^i:U\longrightarrow\Sigma$ such that
\begin{equation*}
        u_n\circ\Phi_n^i
        \rightharpoonup
        u_\infty^i
        \quad\text{weakly in }W^{2,2}(U)
    \end{equation*}
and
\begin{equation*}
        u_n\circ\Phi_n^i
        \longrightarrow
        u_\infty^i
        \quad\text{strongly in }W^{1,p}(U)
        \quad\text{for every }p<\infty.
    \end{equation*}
\item The Dirichlet energy identity holds:
\begin{equation*}        \lim_{n\to\infty}E_{h_n}(u_n)
        =
        \sum_{i\in\mathfrak I_*}
        E(u_\infty^i).
    \end{equation*}
\item The Willmore energy is lower semicontinuous:
\begin{equation*}        
        \sum_{i\in\mathfrak I_*}
        W_{\R^K}(u_\infty^i)
        \le
        \liminf_{n\to\infty}W_{\R^K}(u_n).
    \end{equation*}
\item For $u_\infty: \Sigma_\infty\rightarrow Q$, the image $u_\infty(\Sigma_\infty)$ is connected, and
\begin{equation*}
        d_{\mathrm{Haus}}
        \bigl(
            u_n(\Sigma),
            u_\infty(\Sigma_\infty)
        \bigr)
        \longrightarrow0.
    \end{equation*}
\item The pushforward area measures converge weakly as finite Radon measures on $\R^K$:
\begin{equation*}
        (u_n)_\#\mu_{u_n}
        \rightharpoonup
        \sum_{i\in\mathfrak I_*}
        (u_\infty^i)_\#
        \mu_{u_\infty^i}.
    \end{equation*}
\end{enumerate}
\end{theorem}

Now let $(N,\bar g)$ be a  compact Riemannian manifold, and let $\cK\subset N$ be a compact embedded submanifold.
Fix an isometric embedding $\iota:N\hookrightarrow\R^K$ and set $F_n=\iota\circ u_n$.
By Remark~\ref{rem:intro-euclidean-reduction}, an intrinsic Dirichlet energy and Willmore energy bound of the form
\[
    \sup_n
    \bigl(E_{h_n}(u_n)+W_N(u_n)\bigr)<\infty
\]
implies the Euclidean energy bound \eqref{eq:intro-main-bound} in Theorem~\ref{thm:intro-main-compactness}.
The inclusion $F_n(\Sigma)\subset\iota(N)$ provides the required compact set containing the images. Pulling the limit back by $\iota^{-1}$ gives the corresponding bubble tree convergence theorem for maps into $N$ with free boundary on $\cK$.

\subsection{Hamiltonian stationary Lagrangian immersions} 
We next apply Theorem~\ref{thm:intro-main-compactness} to Hamiltonian stationary Lagrangian immersions. Before stating the resulting compactness theorem, we recall the background of Hamiltonian stationary Lagrangian geometry and discuss the boundary conditions relevant to the free boundary problem.

Let $(M^4,\omega,J,\bar g)$ be a K\"ahler surface, where $\bar g(\cdot,\cdot)=\omega(\cdot,J\cdot)$.
An immersion $u:\Sigma\to M$ is Lagrangian if $u^*\omega=0$.
For the induced metric $g_u=u^*\bar g$ and the trace mean curvature vector $\mathbf H_{u,M}=\operatorname{tr}_{g_u} A_{u,M}$, define the mean curvature, or Maslov, 1-form by
\begin{equation}\label{eq:intro-geometric-mean-curvature-form}
\alpha(X):=\omega(\mathbf H_{u,M},du(X))
=\bar g(J\mathbf H_{u,M},du(X)).
\end{equation}
For a branched immersion, the geometric form $\alpha$ is initially defined on $\overline\Sigma\setminus\cB_u$; its extension and weak equations across $\cB_u$ are specified in Subsection~\ref{subsec:hsl-system}.
Hamiltonian stationarity means criticality of area under the admissible Hamiltonian tests, and its interior equation is $d_{g_u}^*\alpha=0$.
In a Calabi-Yau surface, a local real Lagrangian angle $\theta$ satisfies
\[
\mathbf H_{u,M}=Jdu(\nabla_{g_u}\theta),\qquad \alpha=-d\theta.
\]
For an HSL immersion the angle is locally harmonic, $\Delta_{g_u}\theta=0$.
Minimal Lagrangian surfaces are the special case $\mathbf H_{u,M}=0$.

Oh \cite{Oh1993,Oh1994} developed the variational theory of HSL submanifolds, studying criticality and stability of the volume functional under Hamiltonian deformations. Dazord \cite{Dazord1981} related the mean curvature 1-form $\alpha$ to the ambient Ricci form, providing the geometric identities used in the analytic formulation of Hamiltonian stationarity. In complex dimension two, H\'elein--Romon \cite{HeleinRomon2002} formulated the HSL equations in $\C^2$ as a completely integrable system and constructed a Weierstrass type representation. Schoen--Wolfson \cite{schoen2001} developed the existence and regularity theory for area minimizing Lagrangian maps into K\"ahler surfaces, including the analysis of singularities of the minimizers.

For closed Riemann surfaces, Chen--Ma \cite{chen2021} established the bubble tree compactness for smooth branched HSL immersions into a complete K\"ahler surface under uniform area and Willmore energy bounds, with all images in a fixed compact set.  The proof uses the bubble tree convergence theory for branched conformal immersions developed by Chen--Li \cite{chen2014}. In general dimension, Chen--Warren \cite{Chen-Warren} established compactness for closed HSL submanifolds $L^n\subset\C^n$ under uniform bounds on volume and total extrinsic curvature $\int_L|A|^n\,d\mu$, where $A$ is the second fundamental form. Chen--Ma \cite{Chen-Ma-2024} extended this compactness theory to compact symplectic manifolds with compatible metrics and proved convergence of the total volumes. Apart from collapse to a point, both compactness theorems give smooth local graphical convergence away from a finite set and obtain a Hamiltonian stationary varifold limit.

For the corresponding free boundary problem, the free boundary constraint $u(\partial\Sigma)\subset\cK$ restricts the choices of admissible Hamiltonian variations and produces additional terms in the first variation. Schoen \cite{Schoen2005} formulated relative Lagrangian variational problems with special Lagrangian supports or complex hypersurfaces in Calabi-Yau manifolds, relating these geometric constraints to Dirichlet or Neumann conditions for the Lagrangian angle. Li--Wang--Weng \cite{LiWangWeng2021} studied Legendrian boundary conditions with constant contact angle and proved rigidity of minimal Lagrangian disks in the unit Euclidean ball. Luo--Sun \cite{LuoSun2021} obtained rigidity for minimal Lagrangian surfaces with Legendrian free boundary and classified minimal Lagrangian annuli with Legendrian capillary boundary in the unit  ball. Gaia \cite{Gaia2025} proved minimality and, in the unit ball, equatorial disk rigidity for weakly conformal branched HSL immersions of disks with free boundary in $\C^2$ under suitable regularity and boundary assumptions. His characterization of the orthogonal free boundary condition uses a localized formulation of stationarity for weakly conformal branched HSL immersions of disks with free boundary and an additional tangential boundary assumption \cite[Theorem~7]{Gaia2025}.

% To our knowledge, no analogous bubble tree compactness theorem has been established for branched free boundary HSL immersions with uniformly bounded area and Willmore energy and images in a fixed compact set, while allowing the bordered conformal structures to degenerate. We use the interior estimates of \cite{chen2014,chen2021} and adapt the bubble tree construction in \cite[Section~4]{chen2021}. The boundary estimates and the analysis of degenerating boundary collars are developed below.

To our knowledge, the natural free boundary conditions for weakly conformal branched HSL surfaces whose boundaries are constrained to a general real supporting submanifold $\cK$ are not yet well established.
We therefore first give a complete variational characterization for smooth Lagrangian immersions, with the admissible infinitesimal variations specified in the following theorem.
For a smooth Lagrangian immersion $u:\Sigma\to M$ satisfying $u(\partial\Sigma)\subset\cK$, use the induced metric $g_u$ and let $\nu$ and $\tau$ be the outward unit conormal and unit boundary tangent vectors with respect to the metric $g_u$, chosen so that $(\nu,\tau)$ is positively oriented.
An admissible exact infinitesimal variation on the source $\Sigma$ is a smooth vector field $V\in\Gamma(u^*TM)$ with $V|_{\partial\Sigma}\in T_u\cK$ for which there exists a globally defined $\varphi\in C^\infty(\Sigma,\R)$ such that $\omega(V,du(X))=-d\varphi(X)$ for every $X\in T\Sigma$.
The exactness condition gives
\[
\partial_\tau\varphi=\langle V,Jdu(\tau)\rangle,
\]
linking the boundary trace of the potential $\varphi$ to the admissible boundary velocity $V$. Consequently, the boundary terms $\varphi\alpha(\nu)$ and $\langle V,du(\nu)\rangle$ in the first variation of area cannot in general be tested independently. Accounting for this coupling between the boundary geometry of $u$ relative to $\cK$ and the conormal component $\alpha(\nu)$ of its mean curvature $1$-form yields the following characterization, proved in Theorem~\ref{thm:1st_vari_natural_support}.

\begin{theorem}\label{thm:intro-hsl-boundary}
Let $(M^4,\omega,J,\bar g)$ be a K\"ahler surface and $\cK\subset M$ be a smooth embedded supporting submanifold. Suppose $\Sigma$ is a compact Riemann surface with nonempty smooth boundary and  $u:\Sigma\to M$ is a smooth Lagrangian immersion satisfying $u(\partial\Sigma)\subset\cK$.
Use the induced metric $g_u$ and let $\alpha$ be the mean curvature $1$-form defined in \eqref{eq:intro-geometric-mean-curvature-form}.
Let $\nu$ and $\tau$ be the outward unit conormal and unit boundary tangent vectors with respect to the metric $g_u$, chosen so that $(\nu,\tau)$ is positively oriented.
Then the first variation of area vanishes for every smooth field $V\in\Gamma(u^*TM)$ admitting a globally defined potential $\varphi\in C^\infty(\Sigma,\R)$ with
\[
\omega(V,du(X))=-d\varphi(X)\quad \text{for any } X\in T\Sigma,\qquad
V|_{\partial\Sigma}\in T_u\cK
\]
if and only if there exists a single-valued $q\in C^\infty(\partial\Sigma,\R)$, periodic on every boundary circle, such that
\begin{equation}\label{eq:intro-hsl-coupled-boundary}
\begin{gathered}
d_{g_u}^*\alpha=0\quad\text{in }\Sigma,\\
du(\nu)+qJdu(\tau)\perp T_u\cK,\qquad
\alpha(\nu)+\partial_\tau q=0\quad\text{on }\partial\Sigma.
\end{gathered}
\end{equation}
If the entire boundary map $u|_{\partial\Sigma}$ is an embedding, \eqref{eq:intro-hsl-coupled-boundary} is also equivalent to ambient Hamiltonian stationarity, tested by every $f\in C_c^\infty(M)$ for which $J\nabla f$ is tangent to $\cK$ along $u(\partial\Sigma)$.

Moreover, the homogeneous conditions
\begin{equation}\label{eq:intro-hsl-homogeneous-boundary}
d_{g_u}^*\alpha=0\quad\text{in }\Sigma,\qquad
du(\nu)\perp T_u\cK,\qquad \alpha(\nu)=0
\quad\text{on }\partial\Sigma
\end{equation}
are sufficient for both source and ambient stationarity classes without any boundary embedding assumption. If the immersion $u$ satisfies $Jdu(T\partial\Sigma)\subset T_u\cK$, these homogeneous conditions are also necessary for exact source stationarity, and for ambient Hamiltonian stationarity when the entire boundary map $u|_{\partial\Sigma}$ is embedded.
\end{theorem}

Here, we note that source stationarity means that the first variation of area vanishes for every admissible exact infinitesimal vector field $V$ along $u$, as specified in Theorem~\ref{thm:intro-hsl-boundary}. These conditions are imposed on the initial vector field $V$, without requiring a family of Lagrangian immersions with initial velocity $V$ and boundary constraint on $\cK$. Every admissible ambient Hamiltonian vector field $J\nabla f$, with induced $V=(J\nabla f)\circ u$, yields the source potential $\varphi=f\circ u$, although a general source potential need not admit this representation. When the entire boundary map $u|_{\partial\Sigma}$ is embedded, the compatible boundary values and first derivatives
\[
f\circ u=\varphi,\qquad (\nabla f)\circ u=-JV
\quad\text{on }\partial\Sigma
\]
extend to an ambient potential $f\in C_c^\infty(M)$. If boundary sheets have coincident images, an ambient potential $f$ and its gradient have the same values on those sheets. Their boundary contributions are therefore tested only through their sum and may cancel, as the cut-torus example in Subsection~\ref{subsec:hsl-system} shows.

For the smooth immersion $u$ in Theorem~\ref{thm:intro-hsl-boundary}, let $h$ be any smooth metric on $\Sigma$ in the conformal class of the induced metric $g_u$.
Then $g_u=e^{2\lambda}h$, where $e^{2\lambda}=\tfrac12\operatorname{tr}_h g_u=\tfrac12|du|_h^2>0$ uniquely determines the smooth function $\lambda$.
The outward unit conormal and the positively oriented unit boundary tangent with respect to the metric $h$ are $\nu_h=e^\lambda\nu$ and $\tau_h=e^\lambda\tau$, respectively.
Let $\alpha_h^\sharp$ be the vector field dual to the mean curvature $1$-form $\alpha$ with respect to the metric $h$. With the convention $d\alpha=u^*\operatorname{Rc}$ for the ambient Ricci form, see \cite{Dazord1981,schoen2001}, by \eqref{eq:intro-hsl-coupled-boundary} the general coupled boundary value problem for the map $u$ and its mean curvature $1$-form $\alpha$ is
\begin{equation}\label{eq:intro-hsl-intrinsic-system}
\left\{
\begin{aligned}
\operatorname{tr}_h\nabla du&=-J(u)du(\alpha_h^\sharp),
&\quad u^*\omega&=0 &&\text{in }\Sigma,\\
d\alpha&=u^*\operatorname{Rc},
& d_h^*\alpha&=0 &&\text{in }\Sigma,\\
u(\partial\Sigma)&\subset\cK, & du(\nu_h)+qJ(u)du(\tau_h)&\perp T_u\cK
&&\text{on }\partial\Sigma,\\
\alpha(\nu_h)+\partial_{\tau_h}q&=0
&&&&\text{on }\partial\Sigma.
\end{aligned}
\right.
\end{equation}
On a general supporting submanifold $\cK$, the boundary condition for $u$ in \eqref{eq:intro-hsl-intrinsic-system} involves both $du(\nu_h)$ and $qJ(u)du(\tau_h)$, whereas the condition for the mean curvature form prescribes its normal trace $\alpha(\nu_h)=-\partial_{\tau_h}q$.
Since $q$ is unknown, this first-order condition for $u$ and normal trace condition for $\alpha$ do not give separate standard boundary value problems for the map equation and the Hodge system for $\alpha$. Therefore, in this paper we study the homogeneous boundary conditions \eqref{eq:intro-hsl-homogeneous-boundary} for branched conformal HSL immersions with the free boundary constraint $u(\partial \Sigma)\subset \cK$ in Theorem~\ref{thm:intro-hsl-boundary}. Taking $q\equiv 0$ gives the orthogonal free boundary condition for $u$ in \eqref{sys1} and the zero normal trace for $\alpha$ in \eqref{sys2}.

In Corollary~\ref{coro:hsl-support-geometries}, we characterize the boundary conditions when $\cK$ is a complex curve, a real symplectic surface, a Lagrangian surface or a real hypersurface.
For exact source stationarity, \eqref{eq:intro-hsl-homogeneous-boundary} follows if $u$ satisfies $Jdu(T\partial\Sigma)\subset T_u\cK$. This includes complex curve supports and complex tangent boundary curves on real hypersurfaces.
If $\cK$ is minimal Lagrangian and $\alpha$ is globally exact, the exact source stationarity gives \eqref{eq:intro-hsl-homogeneous-boundary} and forces $\mathbf H_{u,M}=0$. The same conclusions hold for ambient Hamiltonian stationarity when the entire boundary map $u|_{\partial\Sigma}$ is an embedding. Moreover, in Subsection~\ref{subsec:hsl-system}, we give two examples in which the homogeneous conditions fail. 

For a branched conformal HSL immersion with free boundary constraint $u(\partial \Sigma)\subset \cK$, the corresponding geometric equations are initially understood away from the finitely many branch points of $u$, and their weak extension is discussed in Subsection~\ref{subsec:hsl-system}. Using Theorem~\ref{thm:intro-main-compactness} and the homogeneous boundary conditions \eqref{eq:intro-hsl-homogeneous-boundary} of Theorem~\ref{thm:intro-hsl-boundary}, we establish the following bubble tree compactness result for branched conformal HSL immersions with free boundary.

\begin{theorem}\label{main theorem}
Let $(M, \omega, J, \bar g)$ be a complete K\"{a}hler surface, and let $\cK\subset M$ be a compact smooth supporting submanifold.
Let $\Sigma$ be a compact Riemann surface with nonempty boundary $\partial\Sigma$.
For each $n$, let $h_n$ be a Riemannian metric of constant curvature on $\Sigma$, and let $u_n:\left(\Sigma, h_n\right) \rightarrow(M, \bar g)$ be a smooth branched conformal HSL immersion satisfying $u_n(\partial\Sigma)\subset\cK$.
Let $\alpha_n$ be its mean curvature $1$-form, and let $\nu_{h_n}$ be the outward unit conormal along $\partial\Sigma$ with respect to the metric $h_n$. Assume that $d_{h_n}^*\alpha_n=0$ in $\Sigma$ and $\alpha_n(\nu_{h_n})=0$ on $\partial\Sigma$.
Fix $p_0\in M$ and assume that
\begin{equation}\label{eq:distance bound}
    \sup_{n\in\mathbb N}\sup_{x\in\Sigma}
    d_M\bigl(u_n(x),p_0\bigr)<\infty.
\end{equation}
Assume also that:
\begin{itemize}
\item the areas $\cA(u_n)$ are uniformly bounded;
\item the Willmore energies $W_M(u_n)$ are uniformly bounded;
\item $u_n$ meets the supporting submanifold $\cK$ orthogonally:
\begin{equation*}
        du_n(\nu_{h_n}) \perp T_{u_n} \cK \quad \text{on}\ \partial \Sigma.
    \end{equation*}
\end{itemize}

Then, after passing to a subsequence, either the images $u_n(\Sigma)$ converge in Hausdorff distance to a single point, or there exist a stratified bordered surface $\Sigma_{\infty}$ and a continuous map $u_{\infty}: \Sigma_{\infty} \rightarrow M$ with the following properties:
\begin{enumerate}
\item $u_n$ converges to $u_\infty$ in the sense of free boundary bubble tree convergence described in Definition \ref{def:componentwise-free-boundary-bubble tree-convergence};
\item on each nonconstant component of $\Sigma_\infty$, $u_\infty$ is a smooth branched conformal HSL immersion;
\item the following area identity holds:
\begin{equation*}
        \lim_{n \to \infty} \cA(u_n) = \cA(u_\infty);
    \end{equation*}
\item
     
Let $\mathfrak I_*$ index the finitely many nonconstant components of the normalization of $\Sigma_\infty$, and retain all sheet multiplicities in the parametrized integral varifold
\[
        V_\infty
        :=
        \sum_{i\in\mathfrak I_*}
        (u_\infty^i)_\#|\Sigma_\infty^i|.
    \]
The pushforward interior mean curvature measure, see \eqref{eq:mean curvature measure}, has a Radon-Nikodym density $\overline{\mathbf H}_\infty\in L^2(\|V_\infty\|)$, and there is a vector-valued boundary conormal measure $\mathcal B_\infty$, supported on $\cK$ and singular with respect to $\|V_\infty\|$, such that
\[
        \delta V_\infty(X)
        =
        -\int_M
        \langle\overline{\mathbf H}_\infty,X\rangle
        \,d\|V_\infty\|
        +\mathcal B_\infty(X)
    \]
for every compactly supported smooth vector field $X$ on $M$.
Moreover,
\begin{equation*}
        \int_M
        \left\langle
            \overline{\mathbf H}_\infty,J\nabla f
        \right\rangle
        \,d\|V_\infty\|
        =0
        \qquad
        \text{for every }f\in C_c^\infty(M).
    \end{equation*}
If $X$ is tangent to $\cK$ along $\cK$, then $\mathcal B_\infty(X)=0$.
In particular, the full first variation vanishes for every admissible Hamiltonian field $X=J\nabla f$ tangent to $\cK$.
     
\end{enumerate}

\end{theorem}

Completeness of $M$ and the distance bound \eqref{eq:distance bound} in Theorem~\ref{main theorem} imply that all images of $u_n$ lie in a fixed compact subset of $M$. By Theorem \ref{main theorem}, see also Proposition~\ref{prop:stationary} for more details, $\delta V_\infty(J\nabla f)=0$ for every $f\in C_c^\infty(M)$ satisfying $J\nabla f(u_\infty^i(x))\in T_{u_\infty^i(x)}\cK$ at every $x\in\partial\Sigma_\infty^i\setminus\cB_{u_\infty^i}$ of each nonconstant limiting component. The HSL system with homogeneous boundary conditions on each limit bubble tree component yields this identity even when the boundary maps have self-intersections or coincident images. 

% \begin{rmk}\label{rem:closed-surface-comparison}
% The corresponding compactness results for closed surfaces are due to Chen--Li and Chen--Ma.
% Chen--Li \cite{chen2014} proved bubble tree compactness for branched conformal immersions of closed surfaces into $\R^n$ under uniform area and Willmore energy bounds. They assumed that all images meet a fixed ball and allowed the conformal structures to degenerate.
% Chen--Ma \cite[Theorem~1.1]{chen2021} obtained a bubble tree limit of smooth branched HSL immersions, together with the area identity, for closed HSL surfaces in a complete K\"ahler surface. Their assumptions are uniform area and Willmore bounds and containment of all images in a fixed compact set.
% We use the interior estimates of \cite{chen2014,chen2021} and follow the bubble tree construction in \cite[Section~4]{chen2021}.
% The boundary estimates and the analysis of degenerating boundary collars are developed below.
% \end{rmk}

\subsection{Applications of bubble tree compactness}\label{subsec:applications}

We first apply Theorem~\ref{thm:intro-main-compactness} to Willmore energy minimization problems under a boundary length or area constraint. We then study the Willmore energy needed to reach a support disjoint from an area minimizing surface and use Theorem~\ref{main theorem} to determine the limiting intrinsic Willmore energy of HSL immersions.

Let $\Omega\subset\R^3$ be a smooth bounded domain with boundary $S$ and outward unit normal $n_S$. Let $\Sigma$ be a compact smooth annulus with variable conformal structure. For a smooth metric $h$ on $\Sigma$ representing the chosen conformal structure, let $u:(\Sigma,h)\to\R^3$ be a branched conformal immersion with free boundary on $S$ in the sense of Definition~\ref{def:fb-branched-conformal}. Let $\cB_u$ denote the finite branch set of $u$. We use the induced metric $g_u$ and Euclidean Willmore energy $W(u)$ from Subsection~\ref{sec:fbbci}. Let $\nu_h$ and $\tau_h$ be the outward unit conormal and unit boundary tangent vectors with respect to the metric $h$. We use the induced boundary arclength measure $\dd s_{g_u}$, outward unit image conormal $\eta_u$, and total boundary length $L(u)$ defined in \eqref{eq:induced-boundary-geometry}. The length $L(u)$ includes both boundary curves of $\Sigma$. By the variational formulation in \cite[Section~4.1]{appGKL},  we define the first positive Steklov eigenvalue for this boundary measure by
\begin{equation}\label{eq:app-annular-steklov}
 \sigma_1(u)=\inf_{\substack{f\in W^{1,2}(\Sigma),\ \int_{\partial \Sigma}f\dd s_{g_u}=0\\
                        \int_{\partial \Sigma}f^2\dd s_{g_u}>0}}
 \frac{\int_\Sigma|\dd f|_h^2\dd\vol_h}
      {\int_{\partial \Sigma}f^2\dd s_{g_u}}.
\end{equation}
Conformal invariance of the numerator makes \eqref{eq:app-annular-steklov} independent of the representative $h$ in its conformal class for fixed $u$.
Set
\begin{equation*}
 s_{\mathrm{Cl}}
 :=4\pi\sqrt2\left(\sqrt2\coth\pi-\sqrt{2\coth^2\pi-1}\right)>2\pi.
\end{equation*}
This is the length-normalized first Steklov eigenvalue of the upper half Clifford torus $X$; see \eqref{sp:comparison-gap} for more details.

For $s\in(2\pi,s_{\mathrm{Cl}})$ and  $\ell>0$, let $\mathcal E_\Omega(\ell;s)$ consist of immersions $u$ satisfying
\begin{equation}\label{eq:app-annular-class}
 \begin{gathered}
 u(\Sigma)\subset\overline\Omega,\qquad u(\partial\Sigma)\subset S,\qquad
 \eta_u=n_S(u)\quad\text{a.e. on }\partial\Sigma\setminus\cB_u,\\
 L(u)=\ell,\qquad \sigma_1(u)L(u)\geq s.
 \end{gathered}
\end{equation}
Here, the lower bound $2\pi$ for $s$  is sharp; see \eqref{sp:disk-bound}. Its role is to give a uniform gap above $2\pi$ while retaining a closed constraint under passage to the limit. In a nonconstant limiting bubble tree, the conormal condition in \eqref{eq:app-annular-class} preserves the prescribed length by \eqref{eq:contact-scalar}. The lower bound on the normalized Steklov eigenvalue $\sigma_1(u)L(u)\geq s>2\pi$ excludes degeneration of the conformal class by Medvedev's formula \cite{appMedvedev} and excludes a disk carrying all the limiting boundary length, see Lemma~\ref{sp:spectral-transplantation} and \eqref{sp:disk-bound}. Furthermore, with a uniform Willmore energy bound below $4\pi$, sphere type bubbles are also excluded and only the principal annulus remains nonconstant, as proved in Theorem~\ref{sp:annular-compactness}. 

Applying bubble tree compactness in Theorem \ref{thm:intro-main-compactness}, we obtain an annular minimizer in $\mathcal E_\Omega(\ell;s)$ that is a continuous embedding with interior image in $\Omega$ for every sufficiently small $\ell>0$.
The minimum equals $\pi^2+O(\ell)$ as $\ell\to0$, and every minimizing sequence has a subsequence whose only nonconstant bubble tree component is a minimizing annulus.

\begin{theorem}\label{app:small-annuli}
Let $\Omega\subset\R^3$ be a smooth bounded domain and fix $s\in(2\pi,s_{\mathrm{Cl}})$. Let $\mathcal E_\Omega(\ell;s)$ be the class defined in \eqref{eq:app-annular-class}. There are constants $\ell_0>0$ and $C<\infty$, depending only on $\Omega$ and $s$, such that for every $0<\ell<\ell_0$ the infimum
\[
 I_\Omega(\ell;s):=\inf_{u\in\mathcal E_\Omega(\ell;s)}W(u)
\]
is attained and satisfies
\begin{equation}\label{sp:universal-profile}
 \bigl|I_\Omega(\ell;s)-\pi^2\bigr|\leq C\ell,\qquad I_\Omega(\ell;s)<4\pi.
\end{equation}
Every minimizer $u_\infty \in \mathcal E_\Omega(\ell;s)$ is a $W^{2,2}$ conformal immersion with free boundary and a homeomorphism onto its image, that is, $u_\infty$ is a topological embedding. Furthermore, $u_{\infty}^{-1}(\partial\Omega)=\partial\Sigma$, and
\[
 \operatorname{diam}u_{\infty}(\Sigma)\leq C\ell,\qquad \cA(u_\infty)\leq C\ell^2.
\]
Every minimizing sequence $u_n$ has a subsequence whose only nonconstant bubble tree component is a minimizer $u_\infty$ on the principal annulus $\Sigma$, with convergence defined in Definition~\ref{def:componentwise-free-boundary-bubble tree-convergence}. The images $u_n(\Sigma)$ converge to $u_\infty(\Sigma)$ in Hausdorff distance, and
\[
 \cA(u_n)\longrightarrow\cA(u_\infty),\qquad
 (u_n)_\#\dd s_{g_{u_n}}\rightharpoonup (u_\infty)_\#\dd s_{g_{u_\infty}}
 \quad\text{as boundary measures on }S.
\]
\end{theorem}

Next, we consider the Willmore energy minimization problem for branched conformal immersions of disks with free boundary and prescribed area, also with support $S=\partial\Omega$. Let $\nabla^{\mathbb R^3}$ denote the Euclidean connection and $n_S$ the outward unit normal to $S$. For tangent vectors $v,w\in TS$, define
\begin{equation}\label{ga:support-curvature}
 h^S(v,w)=\langle \nabla^{\mathbb R^3}_vn_S,w\rangle,\qquad
 H_S=\operatorname{tr}_{TS}h^S,\qquad H_*=\max_{p\in S}H_S(p).
\end{equation}
Let $\mathcal C(S)$ denote the class of nonconstant maps $u\in W^{2,2}_{b,c,\mathrm{fb}}((\overline{D},h),\R^3,S)$ from Definition~\ref{def:fb-branched-conformal} satisfying the orthogonal boundary condition with prescribed conormal sign
\begin{equation}\label{ga:signed-contact}
 u(\partial{D})\subset S,\qquad du(\nu_h)=|du(\tau_h)|\,n_S(u)
 \quad\text{a.e. on }\partial{D}\setminus\cB_u.
\end{equation}
On $\overline D\setminus\cB_u$, write $A_u^\circ=A_u-\tfrac12 g_u\,\mathbf H_u$ for the trace-free second fundamental form, where $\mathbf H_u=\operatorname{tr}_{g_u}A_u$. For $a>0$, define
\begin{equation}\label{ga:beta-definition}
 \beta_S(a)=\inf\{W(u):u\in\mathcal C(S),\ \cA(u)=a\},\qquad
 \lambda(a)=\sqrt{\frac{a}{2\pi}},\qquad
 \mathcal D(u)=\int_{D\setminus\cB_u}|A_u^\circ|^2\dd\mu_u.
\end{equation}

Using a local perturbation argument near hemispheres, Alessandroni--Kuwert~\cite{appAK} constructed disks meeting $S$ orthogonally that are critical for the Willmore energy under the constraint $\cA(u)=a$ for all sufficiently small $a>0$. Their Willmore energy upper bound yields $\beta_S(a)<2\pi$ (see \eqref{ga:AK-comparison}), which we use in a bubble tree compactness argument for global minimization over the mapping class $\{u\in\mathcal C(S):\cA(u)=a\}$. Combining this bound with the Gauss-Bonnet formula \eqref{ga:disk-cost} and the area identity in Theorem~\ref{thm:intro-main-compactness}, we prove in Proposition~\ref{ga:weak-attainment} that a single nonconstant free boundary disk carries all the prescribed area in the bubble tree limit of a minimizing sequence. For all sufficiently small $a>0$, we prove that the infimum $\beta_S(a)$ over the full class $\mathcal C(S)$ with prescribed area $a>0$ is attained and satisfies $\beta_S(a)=2\pi-\pi H_*\lambda(a)+O(a)$. Every global minimizer $u$ is a continuous embedding with $u({D})\subset\Omega$, is quantitatively close to a hemisphere after rescaling, and concentrates near the maximum set of $H_S$ as $a\to0$.

\begin{theorem}\label{app:small-disks}
Let $\Omega\subset\R^3$ be a smooth bounded domain with boundary $S$, and let $\beta_S(a)$ be as in \eqref{ga:beta-definition}. There are constants $a_0>0$ and $C<\infty$, depending only on $S$, such that $\beta_S(a)$ is attained for every $0<a<a_0$. Every minimizer $u$ is a weak unbranched conformal immersion of class $W^{2,2}$ and a continuous embedding of $\overline{D}$ satisfying
\begin{equation}\label{ga:properness}
 u({D})\subset\Omega,\qquad u(\partial{D})\subset S.
\end{equation}
Moreover,
\begin{align}
 \left|\beta_S(a)-2\pi+\pi H_*\lambda(a)\right|&\le Ca,
 \label{ga:energy-expansion}\\
 \mathcal D(u)&\le Ca,
 \label{ga:tracefree-bound}\\
 0\le H_*-H_S(q)&\le C\sqrt a,
 \qquad q\in u(\partial{D}).
 \label{ga:location-deficit}
\end{align}
Let $\mathbb S^2_-=\{x\in\mathbb S^2:x_3\le0\}$ and $e_3=(0,0,1)$. There exist $p_a\in S$, $R_a\in O(3)$ with $R_an_S(p_a)=e_3$, and a $W^{2,2}$ parametrization $v_a:\mathbb S^2_-\to\R^3$ of the rescaled image such that
\begin{equation}\label{ga:quantitative-hemisphere}
 v_a(\mathbb S^2_-)=\lambda(a)^{-1}R_a\bigl(u(\overline{D})-p_a\bigr),\qquad
 \|v_a-\operatorname{id}_{\mathbb S^2_-}\|_{W^{2,2}(\mathbb S^2_-)}\le C\sqrt a.
\end{equation}
Here $\operatorname{id}_{\mathbb S^2_-}$ is the inclusion of $\S^2_{-}$ into $\R^3$. At the original scale, the Hausdorff distance between $u(\overline{D})$ and $p_a+\lambda(a)R_a^{-1}\mathbb S^2_-$ is at most $Ca$. If $a_j\to0$, $u_j$ minimizes $\beta_S(a_j)$, and $x_j\in u_j(\overline{D})$ converges to $x$, then $x\in S$ and $H_S(x)=H_*$.
\end{theorem}

We next establish lower bounds for the area and Willmore energy of
surfaces with free boundary under a calibration constraint. Let $(N,\bar g)$ be a closed Riemannian manifold, $\cK\subset N$ a smooth compact embedded supporting submanifold, and $S_0\subset N\setminus\cK$ a closed oriented embedded surface. The comass $\|\omega_x\|_{*,\bar g}$ of a $2$-form $\omega$ on $N$ is the supremum of $|\omega_x(v,w)|$ over orthonormal pairs in $T_xN$. A smooth closed $2$-form $\omega$ of pointwise comass at most one calibrates $S_0$ if $\omega|_{TS_0}=\dd\mu_{S_0}$ \cite[Chapter~II, Sections~3-4]{HarveyLawson1982}. 

For branched conformal immersions $u:(\Sigma,h)\to(N,\bar g)$ from connected compact oriented surfaces with nonempty boundary, satisfying the free boundary constraint on $\cK$ and the period constraint \eqref{bar:period}, the calibration inequality gives $\cA(u)\ge m:=\cA(S_0)$. We study how approaching this calibrated area lower bound forces the Willmore energy to increase. For a fixed domain $\Sigma$, bubble tree compactness explains the qualitative obstruction to approaching $m$ with bounded Willmore energy. Indeed, under \eqref{bar:calibration}, the strict comass bound of $\omega$ near $\cK$ excludes limiting bubble tree collapse to a point, while the area identity and no-neck property obtained in Theorem~\ref{thm:intro-main-compactness} would put the entire limiting bubble tree image in $S_0$, contradicting the free boundary constraint on $\cK$. 

In Theorem~\ref{app:calibration-barrier} below, we apply Miura's diameter estimate \cite[Theorem~1.1]{appMiura} to obtain a lower bound of order $(\cA(u)-m)^{-1}$ for the intrinsic Willmore energy $W_N(u)$, uniformly in the topology and conformal structure of $\Sigma$. We also give an example of branched conformal immersions of disks with free boundary whose areas tend to $m$ and whose Willmore energies are of order $(\cA(u)-m)^{-1}$, showing that the order in \eqref{bar:inverse} is optimal; see Subsection~\ref{bar:section} for more details.

\begin{theorem}\label{app:calibration-barrier}
Let $(N,\bar g)$ be a closed Riemannian manifold, $\cK\subset N$ a smooth compact embedded supporting submanifold, and $S_0\subset N\setminus\cK$ a closed oriented embedded surface. Suppose there is a smooth closed $2$-form $\omega \in C^\infty(\wedge^2 T^*N)$ satisfying
\begin{equation}\label{bar:calibration}
 \omega|_{T\cK}=0,\qquad
 \|\omega_x\|_{*,\bar g}\leq1\quad(x\in N),\qquad
 \|\omega_x\|_{*,\bar g}<1\quad(x\notin S_0),\qquad
 \omega|_{TS_0}=\dd\mu_{S_0},
\end{equation}
and put $m=\cA(S_0)>0$. There are constants $c,C>0$, depending only on $N$, $\cK$, $S_0$ and $\omega$, with the following property. If $\Sigma$ is a connected compact oriented surface with nonempty boundary and $u:(\Sigma,h)\to(N,\bar g)$ is a branched conformal immersion with free boundary on $\cK$ in the sense of Definition~\ref{def:fb-branched-conformal} satisfying
\begin{equation}\label{bar:period}
 \int_\Sigma u^*\omega=m,
\end{equation}
then
\begin{equation}\label{bar:inverse}
 \cA(u)>m,\qquad W_N(u)\geq\frac{c}{\cA(u)-m}-C.
\end{equation}
Here $W_N$ is the intrinsic Willmore energy detailed in Subsection~\ref{sec:fbbci}, and constants $c$ and $C$ are independent of the topology and conformal structure of $\Sigma$ and of the branch orders of $u$.

There exist a smooth Riemannian metric $\bar g$ on $N=\mathbb S^2\times\mathbb S^1$, a smooth embedded supporting sphere $\cK$ disjoint from $S_0=\mathbb S^2\times\{0\}$, and a smooth closed $2$-form $\omega$ satisfying \eqref{bar:calibration}, with $m=4\pi$. Let $\iota:S_0\hookrightarrow N$ be the inclusion and
$j_*:H_2(N;\mathbb Z)\to H_2(N,\cK;\mathbb Z)$
be the canonical homomorphism from absolute to relative homology.
Set $\Gamma:=j_*(\iota_*[S_0])$, where
$[S_0]\in H_2(S_0;\mathbb Z)$ is the oriented fundamental class. The infimum of area over smooth conformal immersions
$u:\overline{D}\to N$ with free boundary on $\cK$
and $u_*[\overline{D},\partial{D}]=\Gamma$
is $4\pi$. If $\beta(a)$ denotes the infimum of $W_N(u)$ over weak branched conformal immersions of disks with free boundary satisfying $\cA(u)=a$ and $u_*[\overline{D},\partial{D}]=\Gamma$, then for all sufficiently small $a-4\pi>0$,
\begin{equation}\label{bar:sharp-order}
 \frac{c}{a-4\pi}-C\leq\beta(a)\leq\frac{C}{a-4\pi}.
\end{equation}
In particular, every sequence of branched conformal immersions $u_j:\overline{D}\to N$ with free boundary on $\cK$, relative class $(u_j)_*[\overline{D},\partial{D}]=\Gamma$, and $\cA(u_j)\to4\pi$ satisfies $W_N(u_j)\to\infty$.
\end{theorem}

Finally, we apply Theorem~\ref{main theorem} to study the convergence of the intrinsic Willmore energy for HSL immersions into K\"ahler-Einstein surfaces. The main question is whether the Willmore energies of the nonconstant components account for the full limiting Willmore energy, and how to quantify any additional terms arising from constant principal domains or degenerating collars. When $\operatorname{Rc}=\lambda\omega$, the Lagrangian condition gives $d\alpha=u^*\operatorname{Rc}=\lambda u^*\omega=0$. Under \eqref{eq:intro-hsl-homogeneous-boundary}, the mean curvature $1$-form $\alpha$ is therefore closed and coclosed with zero normal trace, and $W_M(u)=\tfrac14\|\alpha\|_{L^2}^2$ by \eqref{hs:normalization}.

Before giving a detailed description of Theorem \ref{app:hsl-energy}, we introduce some essential notational conventions. For the sequence of  smooth branched conformal Lagrangian
immersions $u_n:(\Sigma,h_n)\to M$ considered in Theorem \ref{app:hsl-energy}, let $\alpha_n$ be the corresponding mean curvature $1$-form defined by \eqref{eq:intro-geometric-mean-curvature-form}. In the subsequential limit supplied by Theorem~\ref{main theorem}, write $\Sigma_\infty^i$, $i\in\mathfrak I_{\mathrm{pr}}$, for the compactified components of the normalization of the bordered Deligne-Mumford limit, and $u_\infty^i$ for the corresponding principal component maps. All principal domains are retained, including those with constant maps. In the single point collapse alternative, every $u_\infty^i$ is the constant map to the limiting point. Let $\Psi_n$ denote the embeddings defined in Definition~\ref{def:bordered-DM-convergence}, $\widehat\Sigma_n$ be the closed complex double of $(\Sigma,[h_n])$ from Subsection~\ref{subsec:uniformization-boundary}, and let $\widehat\alpha_n$ be the reflection-invariant harmonic extension of $\alpha_n$. Index the pinched curves of the double $\widehat\Sigma_n$ by $e\in E$ and choose a full conformal collar for each curve by
\[
 C_{e,n}=(-T_{e,n},T_{e,n})\times(\R/2\pi\mathbb Z),\qquad
 T_{e,n}\longrightarrow\infty.
\]
Writing $\widehat\alpha_n=A_{e,n}\,\dd s+B_{e,n}\,\dd\theta$ on $C_{e,n}$, define
\begin{equation}\label{eq:hsl-collar-coefficients}
a_{e,n}=\frac1{2\pi}\int_0^{2\pi}A_{e,n}(s,\theta)\,\dd\theta,\qquad b_{e,n}=\frac1{2\pi}\int_0^{2\pi}B_{e,n}(s,\theta)\,\dd\theta.
\end{equation}
Closedness and coclosedness of $\alpha_n$ make these coefficients independent of $s$, and \eqref{hs:zero-mode-norm} gives
\[
 \bigl\|a_{e,n}\,\dd s+b_{e,n}\,\dd\theta\bigr\|_{L^2(C_{e,n})}^2
 =4\pi T_{e,n}(a_{e,n}^2+b_{e,n}^2).
\]
For the Fourier expansion of harmonic forms on cylinders, see \cite[Section~3]{appBMMS}. By \eqref{hs:double-norm}, $\|\widehat\alpha_n\|_{L^2(\widehat\Sigma_n)}^2=2\|\alpha_n\|_{L^2(\Sigma,h_n)}^2=8W_M(u_n)$. 

Combining the local convergence in Lemma~\ref{lem:hsl-local-convergence} with the collar identity in Lemma~\ref{hs:collar-identity}, we obtain the following decomposition of the limiting Willmore energy. The formula identifies the additional $L^2$ contributions from constant principal domains and degenerating collars and gives a necessary and sufficient condition for the validity of the Willmore energy identity.

\begin{theorem}\label{app:hsl-energy}
Let $(M,\omega,J,\bar g)$ be a complete K\"ahler-Einstein surface with $\operatorname{Rc}=\lambda\omega$ and  $\cK\subset M$ be a compact smooth supporting submanifold. Let $\Sigma$ be a compact connected oriented surface with nonempty smooth boundary, and let $u_n:(\Sigma,h_n)\to M$ be smooth branched conformal Lagrangian immersions, where the smooth metrics $h_n$ on $\Sigma$ represent possibly varying conformal structures. Assume $u_n(\partial\Sigma)\subset\cK$ and that $(u_n,\alpha_n)$ satisfies the homogeneous HSL elliptic systems \eqref{eq:intro-hsl-homogeneous-boundary}, written with respect to $h_n$ and interpreted weakly at branch points as in Subsection~\ref{subsec:hsl-system}. Suppose that the images lie in a fixed compact subset of $M$ and that
\[
 \sup_n\bigl(\cA(u_n)+W_M(u_n)\bigr)<\infty.
\]
Then, after passing to a subsequence, the following conclusions hold with the notation above.
\begin{enumerate}[label=\textup{(\roman*)},leftmargin=2.4em]
\item On each principal domain $\Sigma_\infty^i$, the pullbacks $\Psi_n^*\alpha_n$ converge smoothly on compact subsets away from the preimages of nodes and the marked points to a $1$-form $\alpha_\infty^i$. The $1$-form $\alpha_\infty^i$ extends to a smooth harmonic $1$-form on $\Sigma_\infty^i$ with zero normal trace along $\partial\Sigma_\infty^i$. If $u_\infty^i$ is nonconstant, $\alpha_\infty^i$ is its mean curvature $1$-form on $\Sigma_\infty^i\setminus\cB_{u_\infty^i}$.
\item Every sphere or free boundary disk bubble is minimal.
\item All terms in the following formula have finite limits, and the following decomposition of the limiting Willmore energy holds:
\begin{align}\label{hs:complete-formula}
 \lim_n W_M(u_n)
 ={}&\sum_{\substack{i\in\mathfrak I_{\mathrm{pr}}\\u_\infty^i\text{ nonconstant}}}W_M(u_\infty^i)
 +\frac14\sum_{\substack{i\in\mathfrak I_{\mathrm{pr}}\\u_\infty^i\text{ constant}}}\|\alpha_\infty^i\|_{L^2(\Sigma_\infty^i)}^2\\\notag
 &+\frac18\lim_n\sum_{e\in E}4\pi T_{e,n}(a_{e,n}^2+b_{e,n}^2).
\end{align}
\end{enumerate}
Consequently, the Willmore energy identity
\[
 \lim_n W_M(u_n)=\sum_{\substack{i\in\mathfrak I_{\mathrm{pr}}\\u_\infty^i\text{ nonconstant}}}W_M(u_\infty^i)
\]
holds if and only if $\alpha_\infty^i=0$ on every principal domain with constant map and the collar contribution in \eqref{hs:complete-formula} vanishes.
The degenerating collar contribution vanishes if the reflected 1-forms $\widehat\alpha_n$ represent a fixed cohomology class under symmetric markings adapted to the degeneration, as in Corollary~\ref{hs:fixed-class-corollary}. Under this fixed cohomology assumption, if $\Sigma$ has genus zero and $\cK$ is a Lagrangian surface, then $\alpha_\infty^i=0$ on every constant principal domain and the Willmore energy identity holds, by Proposition~\ref{hs:phase-support} and Corollary~\ref{hs:planar-identity}.
\end{theorem}

\subsection{Organization of the paper}\label{subsec:organization}

The rest of the paper is organized as follows.

In Section~\ref{subsec:contributions}, we outline the proofs of the compactness theorems and their applications.

In Section~\ref{sec:fb-bubble tree-preliminaries}, we specify the immersion class and the notions of convergence used in the paper. Subsections~\ref{sec:fbbci}-\ref{subsec:uniformization-boundary} fix the geometric conventions and boundary conformal coordinates, while Subsections~\ref{subsec:bordered-DM}-\ref{subsec:bt-convergence-revised} define bordered Deligne-Mumford convergence, bordered stratified surfaces, and free boundary bubble tree convergence.

In Section~\ref{sec:fundamental-ingredients-bubbling}, we establish the boundary estimates for the general compactness theorem. Subsections~\ref{subsec:general-small-energy}-\ref{subsec:boundary-curvature-control} give small energy estimates, removal of boundary singularities, and bounds for total curvature and branch orders. Subsections~\ref{subsec:general-neck-decay}-\ref{subsec:conformal-factor-compactness} establish decay estimates for the Dirichlet energy and the oscillation of the maps on neck regions, together with local compactness of the immersions.

In Section~\ref{sec:general-compactness-proof}, we prove Theorem~\ref{thm:intro-main-compactness}. Subsections~\ref{subsec:generic-extraction}-\ref{subsec:finite-profile-selection} construct the principal components and bubbles. Subsections~\ref{subsec:vanishing-necks}-\ref{subsec:general-completion} establish the no-neck property, extend the limiting maps across punctures, and prove the area identity and convergence of pushforward area measures.

In Section~\ref{sec:HSL-immersions}, we study the free boundary HSL problem. Subsection~\ref{subsec:hsl-system} characterizes stationarity for smooth Lagrangian surfaces and treats special supporting submanifolds. In Subsections~\ref{sec:geodesic-reflection}-\ref{subsec:hsl-branch-points}, we reflect the HSL system under the homogeneous conditions \eqref{eq:intro-hsl-homogeneous-boundary}, establish boundary regularity and removability, and identify the integer orders of the branch points.

In Section~\ref{section5}, we prove Theorem~\ref{main theorem}. Subsection~\ref{subsec:hsl-limit} identifies the limiting HSL components, Subsection~\ref{subsec:hsl-first-variation} establishes the first variation formula for the limiting varifold, and Subsection~\ref{subsec:hsl-consequences} gives related compactness results.

In Section~\ref{sec:applications}, we prove the applications stated in Subsection~\ref{subsec:applications}. Subsections~\ref{sp:section}-\ref{ga:section} treat Willmore energy minimization problems for annuli and disks, Subsection~\ref{bar:section} gives lower bounds for area and Willmore energy under a calibration constraint, and Subsection~\ref{hs:section} determines the limiting intrinsic Willmore energy of free boundary HSL sequences.

\vspace{1cm}

\section{Outline of the proofs}\label{subsec:contributions}

In this section, we outline the proofs of the compactness results in Theorems~\ref{thm:intro-main-compactness} and~\ref{main theorem} and their applications in Theorems~\ref{app:small-annuli}-\ref{app:hsl-energy}.

\subsection{Bubble tree compactness with free boundary}\label{subsec:outline-general-compactness}

To prove Theorem~\ref{thm:intro-main-compactness}, let $u_n:(\Sigma,h_n)\to\R^K$ be a sequence of branched conformal immersions with free boundary on $\cK$ satisfying the hypotheses of that theorem. Their images lie in a fixed compact set $Q\subset\R^K$, and the energy bound \eqref{eq:intro-main-bound} holds. Write $j_n=[h_n]$ and, after passing to a subsequence, assume the bordered Deligne-Mumford convergence \eqref{eq:intro-bordered-DM-convergence}. We first obtain local concentration compactness up to the free boundary and then prove that the Dirichlet energy and the oscillation of $u_n$ vanish on the neck regions connecting the limiting bubble tree components. The interior estimates and the construction on closed domains follow from  \cite{chen2014,chen2021}. However, the Gauss-Bonnet formula in the free boundary case also contains the integral of the geodesic curvature of $\partial\Sigma$ with respect to the induced metric $g_{u_n}$. We must bound this integral using only the Dirichlet and Willmore energy bounds. Moreover, geodesic reflection across a curved support $\cK$ generally destroys conformality of the mapping $u_n$ and introduces an additional term quadratic in $du_n$ into the reflected map equation. When the bordered conformal structures degenerate, we must also obtain decay estimates for the Dirichlet energy and the oscillation of $u_n$ on the boundary necks  uniformly in their lengths and up to the free boundary. In particular, for one-sided collars at contracting boundary components, the boundary term on the retained end circle must be controlled using the orthogonal contact angle condition and the geometry of $\cK$.

We begin by recovering total curvature and branching bounds from the assumed bound on $E_{h_n}(u_n)+W_{\R^K}(u_n)$. For a closed surface, the Gauss-Bonnet formula and the Gauss equation give this control directly. With a free boundary constraint, the same calculation also contains $\int_{\partial\Sigma}k_{g_{u_n}}\,\dd s_{g_{u_n}}$, where $k_{g_{u_n}}$ is the geodesic curvature of $\partial\Sigma$ in the induced metric $g_{u_n}$. Orthogonal contact expresses $k_{g_{u_n}}$ in terms of the second fundamental form of $\cK$, so we first need to bound the boundary length $L(u_n)$. Using a regularized distance vector field to $\cK$ in the first variation formula of the area functional, we prove in Proposition~\ref{prop:free-boundary-length-estimate} that
\[
 L(u_n)\leq C_{\cK}\left(\cA(u_n)+\sqrt{\cA(u_n)W_{\R^K}(u_n)}\right).
\]
The weak bordered Gauss-Bonnet formula given in Lemma~\ref{lem:free-boundary-gauss-bonnet} then yields the curvature and weighted branching bounds in Proposition~\ref{prop:automatic-curvature-branching-bound}.

The next task is to establish local estimates on boundary half disks under the conformal coordinates. For these local estimates, $u$ denotes a member of the sequence written in conformal coordinates. Following the geodesic reflection method for free boundary harmonic maps in \cite{gruter-hildebrandt-nitsche1986,FraserCPAM,JostLiuZhu2019,Jost-Liu-Zhu2019,Gao-Wang-Zhucvpde}, we use the normal reflection across $\cK$ in a tubular neighborhood. The orthogonal free boundary condition makes the first derivatives of $u$ match across the reflection boundary line. However, curvature of $\cK$ produces an additional term quadratic in $du$, and the reflected map need not be conformal. In Theorem~\ref{intsma0}, we absorb this term under smallness of $E+W_{\R^K}$ and obtain boundary $W^{2,p}$ estimates for $1<p<2$. To control the conformal factors, we instead reflect the Gauss maps after using an orthogonal gauge transformation to identify the tangent and normal spaces of $\cK$ with fixed subspaces. The weak connection identity in Lemma~\ref{lem:weak-boundary-connection} records the resulting boundary flux. Combining this identity with the projective Gauss map estimates of Kuwert--Li~\cite{kuwert2012} and a Neumann boundary problem gives boundary removability in Lemma~\ref{lem:general-boundary-removability} and the alternative between collapse to a single point and convergence to a conformal immersion in Lemma~\ref{lem:boundary-conformal-factor-compactness}.

Equipped with these estimates, we construct the components on the degenerating domains $(\Sigma,j_n)$ in \eqref{eq:intro-bordered-DM-convergence}. Taking the quotient of the degenerating collar neighborhoods on the complex double by the antiholomorphic involution interchanging its two copies gives three collar models on $(\Sigma,j_n)$: interior cylinders, longitudinal half-cylinders, and type~II collars with an end circle inherited from $\partial\Sigma$, see Definition~\ref{def:bordered-DM-convergence} and \cite{Katz-Liu2001,Liu2004,Liu-Song-Zhu2022} for more details. Unlike the closed case, a shrinking boundary component remains part of the surface on which the free boundary conditions must be used. Moreover, the embeddings $\Psi_n$ need not be conformal. In Subsection~\ref{subsec:generic-extraction}, we therefore choose smoothly converging metrics on fixed cores of $\widetilde\Sigma_0$ compatible with the pulled-back conformal structures $\Psi_n^*j_n$, and construct conformal boundary charts before applying the local estimates.

We then select all energy concentration scales, including those at which a sphere bubble separates from $\partial\Sigma$. A rescaling that keeps $\partial\Sigma$ at finite distance produces a disk bubble, whereas a rescaling that sends this boundary to infinity permits a sphere bubble. Here a Willmore energy upper bound alone cannot count nonconstant bubble components, since a nonconstant minimal disk with free boundary has zero Willmore energy. We therefore use the combined quantity $E+W_{\R^K}$, for which Theorem~\ref{gapth} gives a positive combined energy gap for nonconstant free boundary branched conformal immersions. In Proposition~\ref{prop:finite-profile-selection}, we assign pairwise disjoint energy concentration regions of $\Sigma$ to the recursive extraction steps in constructing a bubble tree, with each region carrying the same fixed positive amount of $E+W_{\R^K}$. The uniform bound on $E+W_{\R^K}$ therefore bounds the total number of steps and extracted bubbles. 

It remains to prove that the Dirichlet energy and the oscillation of $u_n$ vanish on the connecting neck regions. Smallness of Dirichlet energy on each unit strip does not by itself control a cylinder whose length tends to infinity. In addition, reflecting a longitudinal half-cylinder loses the exact equality between its longitudinal and angular derivative norms. We use conformality on the original half-cylinder to bound this error after reflection, and then derive a differential inequality involving the energy on a slightly enlarged cylinder. Lemma~\ref{lem:delayed-energy-decay} gives the exponential decay claimed in Theorem~\ref{deca}. The estimate on a type~II collar is one-sided: it retains the end circle inherited from $\partial\Sigma$ and truncates only the artificial attaching end. Consequently, Lemma~\ref{lem:noneck} gives vanishing Dirichlet energy and vanishing oscillation of $u_n$ on every neck region. Proposition~\ref{prop:area} yields the Dirichlet energy identity and convergence of the pushforward area measures. Weak $L^2$ convergence of $e^{-\lambda_n}\Delta u_n$, where $\lambda_n$ are the local conformal factors on noncollapsed cores, gives lower semicontinuity of the Willmore energy and completes the proof of Theorem \ref{thm:intro-main-compactness}.

\subsection{Natural boundary conditions and HSL compactness}\label{subsec:outline-hsl-compactness}

Before applying this result to HSL immersions, we determine which free boundary conditions follow from Hamiltonian stationarity of the area functional. For a smooth Lagrangian surface with free boundary constraint, an exact infinitesimal source field $V$ with potential $\varphi$ satisfies
\[
 \partial_\tau\varphi=\langle V,Jdu(\tau)\rangle
 \quad\text{on }\partial\Sigma.
\]
Thus the boundary velocity and the potential trace cannot be varied independently. In Theorem~\ref{thm:1st_vari_natural_support}, we characterize the realizable boundary data by an integral constraint on each boundary circle. Independent constants in $\varphi$ first give zero flux of $\alpha(\nu)$ on each circle, allowing us to construct a periodic primitive $q$. The remaining boundary test vector fields, with an adjustment of the additive constant when needed, yield the coupled conditions in \eqref{eq:intro-hsl-coupled-boundary}. Corollary~\ref{coro:hsl-support-geometries} identifies the supporting submanifold  geometries for which the coupling system \eqref{eq:intro-hsl-intrinsic-system} reduces to the homogeneous boundary condition \eqref{eq:intro-hsl-homogeneous-boundary}; hence it gives the proof of Theorem \ref{thm:intro-hsl-boundary}.

To prove Theorem~\ref{main theorem}, let $u_n:(\Sigma,h_n)\to M$ be a free boundary HSL sequence satisfying its hypotheses and denote the mean curvature $1$-forms by $\alpha_n$. Our task is to preserve this homogeneous system and obtain smoothness of the limiting HSL components in the bubble tree limit. The bound $$\|\alpha_n\|_{L^2(\Sigma,h_n)}^2=4W_M(u_n)$$
controls the mean curvature $1$-forms only in the $L^2$ norm. Hence the term $-J(u_n)du_n((\alpha_n)_{h_n}^\sharp)$ in the mapping equation \eqref{eq:intro-hsl-intrinsic-system} of $u_n$ cannot be treated as having a bounded prescribed coefficient. We instead estimate the map equation of $u_n$ and the Hodge type system for $\alpha_n$ together. Proposition~\ref{extsys} reflects both systems without introducing a measure on the reflection line on which the zero normal trace of $\alpha$ is essential. Combining the estimates for $u_n$ and $\alpha_n$ improves the integrability of both $\alpha_n$ and $du_n$. The coefficient in the absorption step of Theorem~\ref{higher} depends on $\|du_n\|_{L^2}$, so small Dirichlet energy suffices and the Willmore energy need only be bounded. This argument uses \eqref{eq:intro-hsl-homogeneous-boundary} essentially, and the general system \eqref{eq:intro-hsl-intrinsic-system} does not give the same reflectable boundary value problem.

Combining these boundary estimates with the interior estimates for closed HSL surfaces in \cite{chen2021}, we obtain smooth local convergence of $(u_n,\alpha_n)$ in Lemma~\ref{lem:hsl-local-convergence}. Theorem~\ref{remo} removes boundary punctures, and Lemma~\ref{lem:hsl-branch-orders} identifies the integer branch orders. Every nonconstant bubble tree component is therefore a smooth branched conformal HSL immersion satisfying the same homogeneous boundary conditions. Finally, Proposition~\ref{prop:stationary} separates the pushforward interior mean curvature measure from the boundary conormal measure $\mathcal B_\infty$. This proves the weak version for the first variation of the area functional without requiring the limiting boundary maps to be embedded, as claimed in Theorem \ref{main theorem}.

\subsection{Annular minimization under a Steklov constraint}\label{subsec:outline-annuli}

For the first application, Theorem~\ref{app:small-annuli}, we seek a Willmore minimizer in $\mathcal E_\Omega(\ell;s)$ for each sufficiently small $\ell>0$. We first construct admissible annuli by scaling the upper half of the Clifford torus and mapping it into $\Omega$ using inward normal coordinates along $S=\partial\Omega$. Choosing the scale to give total boundary length $\ell$ preserves the orthogonal free boundary condition with $S$. The strict spectral gap \eqref{sp:comparison-gap} persists by \eqref{sp:tubular-spectrum}, so the Steklov eigenvalue constraint holds. Their Willmore energies are $\pi^2+O(\ell)$, and hence $\mathcal E_\Omega(\ell;s)$ is nonempty and $I_\Omega(\ell;s)\leq\pi^2+C\ell<4\pi$ for small $\ell$. 

Next, let $u_n:(\Sigma,h_n)\to\R^3$ be a minimizing sequence in $\mathcal E_\Omega(\ell;s)$, where $\Sigma$ is an annulus. After discarding finitely many terms, its Willmore energies are uniformly below $4\pi$. Lemma~\ref{sp:coercivity} supplies the area bound needed to apply Theorem~\ref{thm:intro-main-compactness}, while the first variation formula of the area functional and Proposition~\ref{prop:free-boundary-length-estimate} exclude single point collapse. The main difficulty is then to preserve the boundary length and annular topology. Since the area identity alone does not control boundary length, we use the common conormal sign $\eta_{u_n}=n_S(u_n)$ in the first variation formula. Lemma~\ref{lem:app-contact-measures} gives convergence of the pushforward boundary arclength measures, together with convergence of boundary integrals depending continuously on the image point and its unoriented tangent line. In particular, no boundary length is lost. Meanwhile, Medvedev's degeneration formula~\cite{appMedvedev} and $\sigma_1(u_n)L(u_n)\geq s>2\pi$ keep the conformal structures in a compact subset of annular moduli space. Lemma~\ref{sp:spectral-transplantation} then excludes multiple nonconstant bordered components, since functions constant on each component would force $\sigma_1(u_n)\to0$. If a single disk carried all the boundary length, the same lemma and \eqref{sp:disk-bound} would instead give $\limsup_n\sigma_1(u_n)L(u_n)\leq2\pi$. Finally, the Willmore energy bound below $4\pi$ excludes nonconstant closed components. Thus the only nonconstant component is a principal annular immersion $u_\infty$, which retains $L(u_\infty)=\ell$ and $\sigma_1(u_\infty)L(u_\infty)\geq s$. Lower semicontinuity of $W$ now gives an admissible minimizer; see Theorem~\ref{sp:annular-compactness} and Corollary~\ref{sp:willmore-minimizer} for more details. 

It remains to prove embeddedness and the Willmore energy lower bound for a minimizer $u$. The diameter estimate in Lemma~\ref{sp:coercivity} places $u(\Sigma)$ in a single boundary chart when $\ell$ is small. We flatten $S$ in this chart and reflect the annulus across the resulting plane. The doubled map has Willmore energy below $8\pi$, so Lemma~\ref{lem:app-local-multiplicity} and the Li--Yau multiplicity inequality exclude branch points and multiple preimages; see \cite[Theorem~6]{Li-Yau1982} for the classical result and \cite[Section~2.2.2, equation~(2.12)]{appMS} for the weak form used here. Thus $u$ is a continuous embedding and satisfies $u^{-1}(S)=\partial\Sigma$. Smooth approximation from Lemma~\ref{lem:app-weak-approximation} and the Li--Yau inequality give embedded smooth tori. Applying the Marques--Neves Willmore inequality \cite[Theorem~A]{Marques-Neves2014} via inverse stereographic projection and passing to the limit yields $W(u)\geq\pi^2-C\ell$, as proved in Proposition~\ref{prop:app-short-annular-geometry}. Combining this lower bound with the constructed upper bound proves \eqref{sp:universal-profile} and completes Theorem~\ref{app:small-annuli}. Moreover, applying the quantitative rigidity theorem of Bi--Zhou~\cite{BiZhou2025} in Corollary~\ref{coro:app-normalized-annular-double} gives a parametrization of the doubled torus that is $O(\sqrt\ell)$-close in $W^{2,2}$ to the standard parametrization of the Clifford torus after inverse stereographic projection and an ambient M\"obius transformation.

\subsection{Small-area disks and concentration on the supporting hypersurface}\label{subsec:outline-disks}

For the second application, Theorem~\ref{app:small-disks}, we minimize the Willmore energy over $u\in\mathcal C(S)$ with prescribed area $\cA(u)=a$. The local construction of Alessandroni--Kuwert~\cite{appAK}, obtained by perturbing hemispheres, gives the upper bound \eqref{ga:AK-comparison} and hence $\beta_S(a)<2\pi$ for small $a$. To obtain a global minimizer, we must preserve the prescribed area under bubble tree convergence, although the interiors of competitors are not confined to $\overline\Omega$. After discarding finitely many terms, Lemma~\ref{ga:boundary-distance} places the images of a minimizing sequence in a fixed compact set, and the first variation formula of area excludes collapse to a point at fixed $a>0$. Applying Theorem~\ref{thm:intro-main-compactness} to this minimizing sequence gives a noncollapsed bubble tree whose component areas sum to $a$. The boundary length estimate \eqref{eq:sec2-boundary-length-area-willmore} and the bound \eqref{ga:disk-cost} give Willmore energy greater than $\pi$ for each nonconstant limiting disk, whereas every nonconstant sphere has Willmore energy at least $4\pi$. Since the minimizing sequence has Willmore energies uniformly below $2\pi$, exactly one nonconstant disk remains. The area identity assigns this disk area $a$, and lower semicontinuity of Willmore energy gives attainment of $\beta_S(a)$. Proposition~\ref{ga:weak-attainment} also excludes branch points of the minimizer.

To determine the shape and location of a minimizer $u$, we next estimate the integral $\int_{\partial D}k_{g_u}\,\dd s_{g_u}$ of the geodesic curvature of $\partial D$ in the induced metric $g_u$, as it appears in \eqref{ga:disk-defect}. With $\lambda(a)$ and $\mathcal D(u)$ as in \eqref{ga:beta-definition}, the preliminary bounds \eqref{ga:preliminary-bounds} allow us to rescale about $q\in u(\partial D)$ by $\lambda(a)^{-1}$, flatten $S$ in normal coordinates, and reflect the disk across the resulting plane. The reflected double is a weak $W^{2,2}$ immersion of a sphere. By Lemma~\ref{ga:tracefree-transport}, its trace-free second fundamental form has small $L^2$ norm. Using the sphere rigidity estimates of De Lellis--M\"uller~\cite{appDLM05,appDLM06}, Lemmas~\ref{ga:weak-dlm} and~\ref{ga:equivariant-gauge} give a conformal parametrization of this sphere that is close in $W^{2,2}$ to the standard parametrization of a round sphere and sends the equator onto the boundary of the flattened disk. Restricting this parametrization to the equator and applying the trace theorem controls its tangential derivative, yielding the boundary integral estimate in Lemma~\ref{ga:boundary-moment}. Transferring this estimate back to $S$ gives Proposition~\ref{ga:boundary-expansion}. Combining the expansion of $\int_{\partial D}k_{g_u}\,\dd s_{g_u}$ from Proposition~\ref{ga:boundary-expansion} with \eqref{ga:disk-defect} and the upper bound \eqref{ga:strict-upper}, we obtain
\[
 \frac12\mathcal D(u)+\pi\bigl(H_*-H_S(q)\bigr)\lambda(a)
 \leq C\lambda(a)\sqrt{\mathcal D(u)}+C\lambda(a)^2,
 \qquad q\in u(\partial D).
\]

Young's inequality bounds the right-hand side mixed term by $\tfrac14\mathcal D(u)+C\lambda(a)^2$. Absorbing $\tfrac14\mathcal D(u)$ into the left-hand side and using $\mathcal D(u)\geq0$ and $H_*-H_S(q)\geq0$, we obtain $\mathcal D(u)=O(a)$ and $H_*-H_S(q)=O(\sqrt a)$, uniformly for $q\in u(\partial D)$. The diameter bound \eqref{ga:preliminary-bounds} extends this concentration near the maximum set of $H_S$ from the boundary to the entire image. Substituting these improved bounds into the curvature expansion and \eqref{ga:disk-defect} yields the Willmore energy expansion \eqref{ga:energy-expansion}. Moreover, injectivity of the reflected double in Lemma~\ref{ga:weak-dlm} makes $u$ a continuous embedding; together with the prescribed conormal sign, it also gives $u(D)\subset\Omega$. Finally, restricting the sphere parametrization to the appropriate hemisphere and undoing the flattening, with a suitable recentering and rotation, proves the quantitative approximation \eqref{ga:quantitative-hemisphere} and Theorem \ref{app:small-disks}.

\subsection{Contact with a separated calibrated surface}\label{subsec:outline-calibration}

For the third application in Theorem~\ref{app:calibration-barrier}, we establish a lower bound for $W_N(u)$ in terms of the area excess $\cA(u)-m$, under the period constraint \eqref{bar:period} and the free boundary constraint on $\cK$. We first explain the qualitative obstruction using bubble tree compactness. Suppose a sequence on a fixed domain has areas tending to $m$ and uniformly bounded Willmore energies. The strict comass bound of $\omega$ near $\cK$ excludes collapse to a point. The calibration inequality and the area identity then force every nonconstant limiting component to lie in $S_0$. By the no-neck property, the Hausdorff limit of the entire images is also contained in $S_0$, contradicting the boundary contact with $\cK$.

To obtain a quantitative estimate uniform in the topology of the domain, we use a smooth cut-off function equal to $0$ near $S_0$ and to $1$ near $\cK$ to select a region separating these sets. The period constraint \eqref{bar:period} and the strict comass gap \eqref{bar:calibration} bound the area in this region by $C(\cA(u)-m)$. After excising small neighborhoods of the branch points and smoothly approximating the remaining immersion, the coarea formula gives two regular level sets with controlled total length. Connectedness of $u$ then provides a subsurface whose image joins these level sets and has diameter bounded below by a fixed positive constant. Its area and boundary length are bounded by $C(\cA(u)-m)$, up to errors tending to zero. Applying Miura's \cite{appMiura} diameter estimate to this subsurface and passing first to the approximation limit and then shrinking the excised neighborhoods gives
\[
 1\leq C\sqrt{W_N(u)\bigl(\cA(u)-m\bigr)}+C\bigl(\cA(u)-m\bigr).
\]
Thus $\cA(u)>m$, and for small positive area excess we obtain $W_N(u)\geq c/(\cA(u)-m)$. Increasing the additive constant gives \eqref{bar:inverse} over the remaining range. Finally, in the example constructed in Subsection~\ref{bar:section}, replacing a small disk in $S_0$ by a thin tube ending orthogonally on $\cK$ produces smooth embedded disks with the matching upper bound \eqref{bar:sharp-order}. This completes the proof of Theorem \ref{app:calibration-barrier}.

\subsection{Intrinsic Willmore energy for HSL limits}\label{subsec:outline-hsl-energy}

Finally, let $u_n:(\Sigma,h_n)\to M$ satisfy the hypotheses of Theorem~\ref{app:hsl-energy}, and let $\alpha_n$ be their mean curvature $1$-forms. We determine whether the Willmore energies of the nonconstant limiting components account for $\lim_nW_M(u_n)$ and identify the additional terms when they do not. Since $(M,\omega,J,\bar g)$ is K\"ahler-Einstein and $u_n$ are Lagrangian, we have $d\alpha_n=u_n^*\operatorname{Rc}=0$. Since $\alpha_n$ is also coclosed and has zero normal trace on $\partial \Sigma$, reflection across $\partial\Sigma$ extends it to a harmonic $1$-form $\widehat\alpha_n$ on the complex double of $\Sigma$. This reduces the problem to decomposing the limiting squared $L^2$ norms of these harmonic forms.

We first identify the contributions from the principal domains.  Compactness and removability of isolated singularities of harmonic 1-forms give a smooth harmonic 1-form limit $\alpha_\infty^i$ on each compactified principal domain $\Sigma_\infty^i$. When $u_\infty^i$ is nonconstant, Lemma~\ref{lem:hsl-local-convergence} identifies $\alpha_\infty^i$ with its mean curvature $1$-form away from the branch points. When $u_\infty^i$ is constant, however, $\alpha_\infty^i$ need not vanish, and the term $\tfrac14\|\alpha_\infty^i\|_{L^2(\Sigma_\infty^i)}^2$ must be retained separately. Local estimates for harmonic 1-forms exclude concentration of the squared $L^2$ norms at bubbling points away from the degenerating nodal collars. Moreover, the harmonic mean curvature $1$-form of every sphere or free boundary disk bubble vanishes, since a sphere has no nonzero harmonic $1$-forms and a disk has none with zero normal trace. Thus all bubbles are minimal.

It remains to account for the degenerating collars, where local convergence of harmonic 1-forms alone does not determine the limiting squared $L^2$ norm. In Lemma~\ref{hs:collar-identity}, the Fourier expansion for harmonic forms on cylinders~\cite{appBMMS} separates the constant and nonconstant modes of $\widehat\alpha_n$ on each degenerating collar. The nonconstant modes decay exponentially away from the two ends, whereas the constant part has squared norm
\[
 \|a_{e,n}\,\dd s+b_{e,n}\,\dd\theta\|_{L^2(C_{e,n})}^2
 =4\pi T_{e,n}(a_{e,n}^2+b_{e,n}^2)
\]
by \eqref{hs:zero-mode-norm}. Taking $n\to\infty$ and then exhausting the normalized components separates these constant Fourier contributions from the component norms. Combining this decomposition with $\|\widehat\alpha_n\|_{L^2}^2=8W_M(u_n)$ from \eqref{hs:double-norm} proves \eqref{hs:complete-formula}.

We then examine when the additional terms in \eqref{hs:complete-formula} vanish. If the reflected harmonic 1-forms $\widehat{\alpha}_n$ represent a fixed cohomology class under symmetric markings adapted to the degeneration, Corollary~\ref{hs:fixed-class-corollary} uses the minimizing property of harmonic 1-form representatives to eliminate the degenerating collar term. For constant principal domains, the local phase identity in Proposition~\ref{hs:phase-support} shows that the periods of $\alpha_\infty^i$ lie in $2\pi\mathbb Z$. If $\cK$ is Lagrangian and $\partial\Sigma_\infty^i\ne\varnothing$, the boundary identity \eqref{eq:hsl-lagrangian-boundary-form} also gives zero tangential trace on an open boundary arc inherited from $\partial\Sigma$. Together with the zero normal trace boundary condition and unique continuation, this implies $\alpha_\infty^i=0$. Consequently, \eqref{hs:complete-formula} reduces to the Willmore energy identity precisely when the collar term vanishes and $\alpha_\infty^i=0$ on every constant principal domain, completing the proof of Theorem \ref{app:hsl-energy}.

\vskip1cm

\section{\texorpdfstring{Preliminaries}{Preliminaries}}
\label{sec:fb-bubble tree-preliminaries}

We specify the immersion class in Subsection~\ref{sec:fbbci} and the notions of domain and map convergence in Subsections~\ref{subsec:bordered-DM}--\ref{subsec:bt-convergence-revised}.
Throughout this section, $\Sigma$ is a compact connected oriented smooth surface with nonempty smooth boundary, $(N,\bar g)$ is a compact Riemannian manifold, and $\cK\subset N$ is a compact smooth embedded submanifold.
We denote the K\"ahler target in the HSL results by $M$ and use $\R^K$ for Euclidean embedding spaces.
A smooth metric on $\Sigma$ is denoted by $h$, and $g_u=u^*\bar g$ denotes the metric induced on $\Sigma$ away from the branch points of $u$.
The same local definitions apply with $N=\R^K$ and a compact smooth embedded support $\cK\subset\R^K$.
For a smooth metric $h$ on $\Sigma$, we write $j=[h]$ for the associated conformal structure.
Local Sobolev spaces and boundary traces are understood in boundary conformal charts, using the Sobolev and trace theorems in~\cite{AdamsFournier2003}.
We use $i$ for indices of bubble tree components and $\alpha$ for the mean curvature $1$-form in the HSL results.

\Needspace{6\baselineskip}
\subsection{\texorpdfstring{Branched conformal immersions with free boundary}{Branched conformal immersions with free boundary}}
\label{sec:fbbci}

Let $U\subset\overline\Sigma$ be relatively open.
A map $u\in W^{2,2}_{\loc}(U,N)$ is called a $W^{2,2}$ conformal immersion on $U$ if there exists $\lambda\in L^\infty_{\loc}(U)$ such that
\begin{equation}\label{eq:sec2-intrinsic-conformality}
    u^*\bar g=e^{2\lambda}h
    \qquad\text{a.e. on }U.
\end{equation}
Let $u\in W^{2,2}(\Sigma,N)$ be a conformal immersion on $\overline\Sigma\setminus\cB_u$, where $\cB_u\subset\overline\Sigma$ is finite.
On $\overline\Sigma\setminus\cB_u$, write $g_u:=u^*\bar g$ and $\dd\mu_u:=\dd\vol_{g_u}$.
For measurable $\Omega\subset\Sigma$, set $\cA(u;\Omega):=\int_{\Omega\setminus\cB_u}\dd\mu_u$ and $\cA(u):=\cA(u;\Sigma)$.
Since $|du|_h^2=2e^{2\lambda}$ and $\dd\mu_u=e^{2\lambda}\dd\vol_h$,
\begin{equation}\label{eq:sec2-energy-area}
 E_h(u;\Omega):=\frac12\int_\Omega|du|_h^2\,\dd\vol_h
 =\cA(u;\Omega).
\end{equation}
On $\overline\Sigma\setminus\cB_u$, let $A_{u,N}$ be the second fundamental form and $\mathbf H_{u,N}:=\operatorname{tr}_{g_u}A_{u,N}$ its trace mean curvature vector.
Curvature quantities and $g_u$-unit vectors use the inverse induced metric and are defined almost everywhere on this set.
We write
\begin{align*}
T_N(u;\Omega):=\int_{\Omega\setminus\cB_u}|A_{u,N}|_{\bar g}^2\,\dd\mu_u,
\end{align*}
\begin{align*}
W_N(u;\Omega):=\frac14\int_{\Omega\setminus\cB_u}|\mathbf H_{u,N}|_{\bar g}^2\,\dd\mu_u
 =\frac18\int_{\Omega\setminus\cB_u}|\mathbf H_{u,N}|_{\bar g}^2|du|_h^2\,\dd\vol_h.
 \end{align*}
Omitting $\Omega$ means integration over $\Sigma$; geometric integrals exclude the finite branch set. We write $E(u;\Omega)=E_h(u;\Omega)$ when the conformal structure is understood. For Euclidean immersions, we abbreviate $A_u=A_{u,\R^K}$ and $\mathbf H_u=\mathbf H_{u,\R^K}$, and write $W(u)=W_{\R^K}(u)$.
The tension field $\tau_h(u):=\operatorname{tr}_h\nabla du$ satisfies $\tau_h(u)=\frac12|du|_h^2\mathbf H_{u,N}$ at almost every point of $\overline\Sigma\setminus\cB_u$.
Thus a conformal immersion is harmonic as an $N$-valued map if and only if it is minimal in $N$.
For an ambient vector field $X$, $\operatorname{div}_uX$ denotes its tangential divergence along $u$, computed on $du(T\Sigma)$ away from branch points.

\begin{definition}
\label{def:fb-branched-conformal}
A map $u\in W^{2,2}(\Sigma,N)$ belongs to
\[
    W^{2,2}_{b,c,\mathrm{fb}}
    \bigl((\Sigma,h),N,\cK\bigr)
\]
if there exists a finite set $\cB_u\subset\overline\Sigma$ such that the following conditions hold
\begin{enumerate}[label=\textup{(\roman*)},leftmargin=2.4em]
\item The map is conformal on $\overline\Sigma\setminus\cB_u$, with
\[
        u^*\bar g=e^{2\lambda}h
        \qquad\text{a.e. on }
        \overline\Sigma\setminus\cB_u
    \]
for some $\lambda\in L^\infty_{\loc}(\overline\Sigma\setminus\cB_u)$.

\item The immersion has finite area and finite total curvature
\begin{equation}\label{eq:sec2-finite-total-curvature}
        \int_{\Sigma\setminus\cB_u}
            \bigl(1+|A_{u,N}|_{\bar g}^2\bigr)
            \,\dd\mu_u
        <\infty.
    \end{equation}
\item The continuous representative of $u$ satisfies $u(\partial\Sigma)\subset\cK$, and the Sobolev trace of its differential satisfies
\[
du(\nu_h)\perp_{\bar g}T_u\cK
\quad\text{for }\mathcal H^1\text{-a.e. point of }
\partial\Sigma\setminus\cB_u,
\]
where $\nu_h$ is the outward unit conormal with respect to the metric $h$ on $\Sigma$.
\end{enumerate}
\end{definition}

In the general compactness theorem and its Willmore applications, a branched conformal immersion with free boundary means a map in Definition~\ref{def:fb-branched-conformal}.
The terms disk and annulus specify the domain. Embeddedness and higher regularity are specified separately.

Let $\nu_h$ and $\tau_h$ be the outward unit conormal and the unit boundary tangent with respect to $h$, chosen so that $(\nu_h,\tau_h)$ is positively oriented. Along $\partial\Sigma\setminus\cB_u$, the induced arclength measure, outward unit image conormal, and unit image tangent are
\begin{equation}\label{eq:induced-boundary-geometry}
\begin{gathered}
 \dd s_{g_u}=|du(\tau_h)|_{\bar g}\,\dd s_h,\qquad
 L(u)=\int_{\partial\Sigma\setminus\cB_u}\dd s_{g_u},\\
 \eta_u=\frac{du(\nu_h)}{|du(\nu_h)|_{\bar g}},\qquad
 T=\frac{du(\tau_h)}{|du(\tau_h)|_{\bar g}}.
\end{gathered}
\end{equation}
The boundary values of $du$ are understood in the Sobolev trace sense, and these definitions are independent of the representative $h$ in the conformal class.

The set $\cB_u$ in Definition~\ref{def:fb-branched-conformal} is initially a finite exceptional set.
The boundary expansion in Lemma~\ref{lem:general-boundary-removability} and the interior expansion in \cite[Theorem~3.1]{kuwert2012} assign a nonnegative integer order to each such point.
We denote by $\cB_u$ the set of points with positive branch order.

The Sobolev embedding of $W^{2,2}(\Sigma,N)$ into $C^0(\overline\Sigma,N)$ gives these maps a continuous representative.
Image constraints and Hausdorff convergence always refer to this representative.

Conformal rescaling of $h$ preserves conformality and Dirichlet energy and multiplies the outward unit conormal by a positive scalar.
Thus Definition~\ref{def:fb-branched-conformal} depends only on $[h]$.
Since $|\mathbf H_{u,N}|^2\le 2|A_{u,N}|^2$, condition \eqref{eq:sec2-finite-total-curvature} implies that $\cA(u)$, $T_N(u)$ and $W_N(u)$ are finite.

\begin{rmk}
\label{rem:intro-euclidean-reduction}
Fix an isometric embedding $\iota:(N,\bar g)\hookrightarrow\R^K$ and set $F:=\iota\circ u$.
The induced metric, area measure, area and Dirichlet energy are unchanged.
After identifying $N$ with $\iota(N)$, let $A^N_{\R^K}$ denote the second fundamental form of $N\subset\R^K$.
On $\overline\Sigma\setminus\cB_u$,
\begin{equation}\label{eq:sec2-second-fundamental-form-splitting}
    A_{F,\R^K}(X,Y)
    =
    \iota_*A_{u,N}(X,Y)
    +
    A^N_{\R^K}
    \bigl(dF(X),dF(Y)\bigr).
\end{equation}
If $\{e_1,e_2\}$ is a local frame orthonormal with respect to the induced metric $g_u$, then
\begin{equation*}   
\mathbf H_{F,\R^K}    =
    \iota_*\mathbf H_{u,N}
    +
    \sum_{i=1}^2
    A^N_{\R^K}
    \bigl(dF(e_i),dF(e_i)\bigr).
\end{equation*}
The two summands in \eqref{eq:sec2-second-fundamental-form-splitting} are orthogonal.
Since $N$ is compact, there exists $C_\iota<\infty$ such that
\begin{align}
    T_N(u)
    &\le T_{\R^K}(F)
    \le T_N(u)+C_\iota\cA(u),
    \label{eq:sec2-total-curvature-comparison}
    \\
    W_N(u)
    &\le W_{\R^K}(F)
    \le W_N(u)+C_\iota\cA(u).
    \label{eq:sec2-willmore-comparison}
\end{align}
Moreover, the intrinsic free boundary condition is equivalent to
\[
    F(\partial\Sigma)\subset\iota(\cK),
    \qquad
    dF(\nu_h)\perp_{\R^K}T_F\iota(\cK).
\]
The result for compact targets therefore follows by applying the Euclidean theorem to $F$, whose image lies in the fixed compact set $\iota(N)$.
\end{rmk}

In a Euclidean conformal coordinate $z=x_1+\mathrm i x_2$, write $g_F=e^{2\lambda}(\dd x_1^2+\dd x_2^2)$.
Then $|\nabla F|^2=2e^{2\lambda}$ and $\Delta F=\frac12\mathbf H_{F,\R^K}|\nabla F|^2$.
Extending the trace mean curvature vector by zero where $|\nabla F|=0$, we write
\[
 \mathbf H_{F,\R^K}:=
 \begin{cases}
  2\Delta F/|\nabla F|^2,&|\nabla F|\ne0,\\
  0,&|\nabla F|=0.
 \end{cases}
\]
This extension satisfies
\[
 W_{\R^K}(F;\Omega)=\frac18\int_\Omega|\mathbf H_{F,\R^K}|^2|\nabla F|^2\,\dd x
\]
on every conformal coordinate domain $\Omega$.

\subsection{Boundary conformal coordinates and uniformization}
\label{subsec:uniformization-boundary}

For $r>0$, set
\begin{align}\label{denotion-ir}
    D_r:=\{z\in\C:|z|<r\},
    \qquad
    D_r^+:=D_r\cap\{\operatorname{Im}z>0\},
    \qquad
    I_r:=(-r,r)=D_r\cap\R.
\end{align}
We write $D:=D_1$, $D^+:=D_1^+$, and $I:=I_1$.
In local charts and collar decompositions, the boundary inherited from $\partial\Sigma$ is the part on which the free boundary conditions hold.
Boundaries introduced by cutting or truncating a domain are called artificial boundaries.

For each $p\in\partial\Sigma$, choose a conformal chart $\zeta:U\to D_r^+\cup I_r$ with $\zeta(p)=0$, mapping $U\cap\Sigma^\circ$ onto $D_r^+$ and $U\cap\partial\Sigma$ onto $I_r$.
If $z=x+\mathrm i y$ and $h=e^{2\varphi}(\dd x^2+\dd y^2)$, the outward unit conormal on $I_r$ with respect to the metric $h$ is $\nu_h=-e^{-\varphi}\partial_y$.
Thus the free boundary condition becomes
\[
 u(I_r)\subset\cK,\qquad \partial_yu\perp T_u\cK\quad\text{on }I_r
\]
in the Sobolev trace sense.

For a surface of genus $g$ with $k$ boundary components, the closed complex double $\widehat\Sigma$ has genus $2g+k-1$; see~\cite{Katz-Liu2001,Liu2004}.
It carries an antiholomorphic involution $\sigma:\widehat\Sigma\to\widehat\Sigma$ whose fixed set projects to $\partial\Sigma$, with $\Sigma\simeq\widehat\Sigma/\langle\sigma\rangle$.
If $\chi(\Sigma)<0$, then each conformal class on $\Sigma$ contains a unique hyperbolic metric with totally geodesic boundary.
Equivalently, this metric is the restriction of the unique $\sigma$-invariant hyperbolic metric on the complex double; see~\cite{hummel1997} for general hyperbolic geometry and~\cite[Section~4.2.2]{Liu2004} for this bordered uniformization.
We use these canonical representatives to describe bordered Deligne--Mumford degeneration. For the disk, we use the round hemisphere. For the annulus, we use $\mathfrak A_\ell=[0,\ell]\times(\R/2\pi\mathbb Z)$ with the unit-area flat metric
\[
    h_\ell=\frac{dt^2+d\theta^2}{2\pi\ell}.
\]
As $\ell\to0$, the transverse geodesic arcs have length $\sqrt{\ell/(2\pi)}\to0$; as $\ell\to\infty$, the boundary circles have length $\sqrt{2\pi/\ell}\to0$.
These give type~I and type~II collars, respectively; a finite positive limiting value of $\ell$ gives no degeneration.

\Needspace{6\baselineskip}
\subsection{\texorpdfstring{Bordered Deligne--Mumford convergence}{Bordered Deligne--Mumford convergence}}
\label{subsec:bordered-DM}

We use the nodal bordered surfaces and boundary node types described in \cite[Definition~3.4 and Section~3.2]{Liu2004}.
\begin{definition}
\label{def:nodal-bordered-surface}
A nodal bordered Riemann surface is a quotient
\[
    \Sigma_0
    =
    \widehat\Sigma_0/\langle\sigma_0\rangle,
\]
where $\widehat\Sigma_0$ is a nodal closed Riemann surface and $\sigma_0$ is an antiholomorphic involution preserving its nodal set, with an orientation chosen on the quotient away from the images of the nodes.
This is the orientation inherited from $\Sigma$ when the quotient arises by degeneration of its complex double.
The nodes of $\Sigma_0$ fall into the following three types:
\begin{enumerate}[label=\textup{(\roman*)},leftmargin=2.4em]
\item An \emph{interior node} is the image of a pair of distinct nodes of $\widehat\Sigma_0$ exchanged by $\sigma_0$;

\item A \emph{type~I boundary node} is the image of a node fixed by $\sigma_0$ such that the two local branches are preserved separately by the involution;

\item A \emph{type~II boundary node} is the image of a node fixed by $\sigma_0$ such that the two local branches are exchanged by the involution.
\end{enumerate}
We denote the sets of nodes of these three types by $\cN_0^{\mathrm{int}}$, $\cN_0^{\mathrm I}$, and $\cN_0^{\mathrm{II}}$, respectively.
The finite nodal set is denoted by
\[
    \cN_0
    =
    \cN_0^{\mathrm{int}}
    \sqcup
    \cN_0^{\mathrm I}
    \sqcup
    \cN_0^{\mathrm{II}}.
\]
The complement of the nodes in $\Sigma_0$ is denoted by
\[
    \Sigma_0^{\mathrm{reg}}:=\Sigma_0\setminus\cN_0.
\]
\end{definition}

Let $\pi_0:\widetilde\Sigma_0\to\Sigma_0$ be the normalization.
The inverse image of an interior node consists of two interior points, and the inverse image of a type~I boundary node consists of two boundary points.
A type~II boundary node has a single preimage under normalization.
This interior special point records the contracted boundary component.

\begin{figure}[H]
\centering
\includegraphics[width=0.6\textwidth]{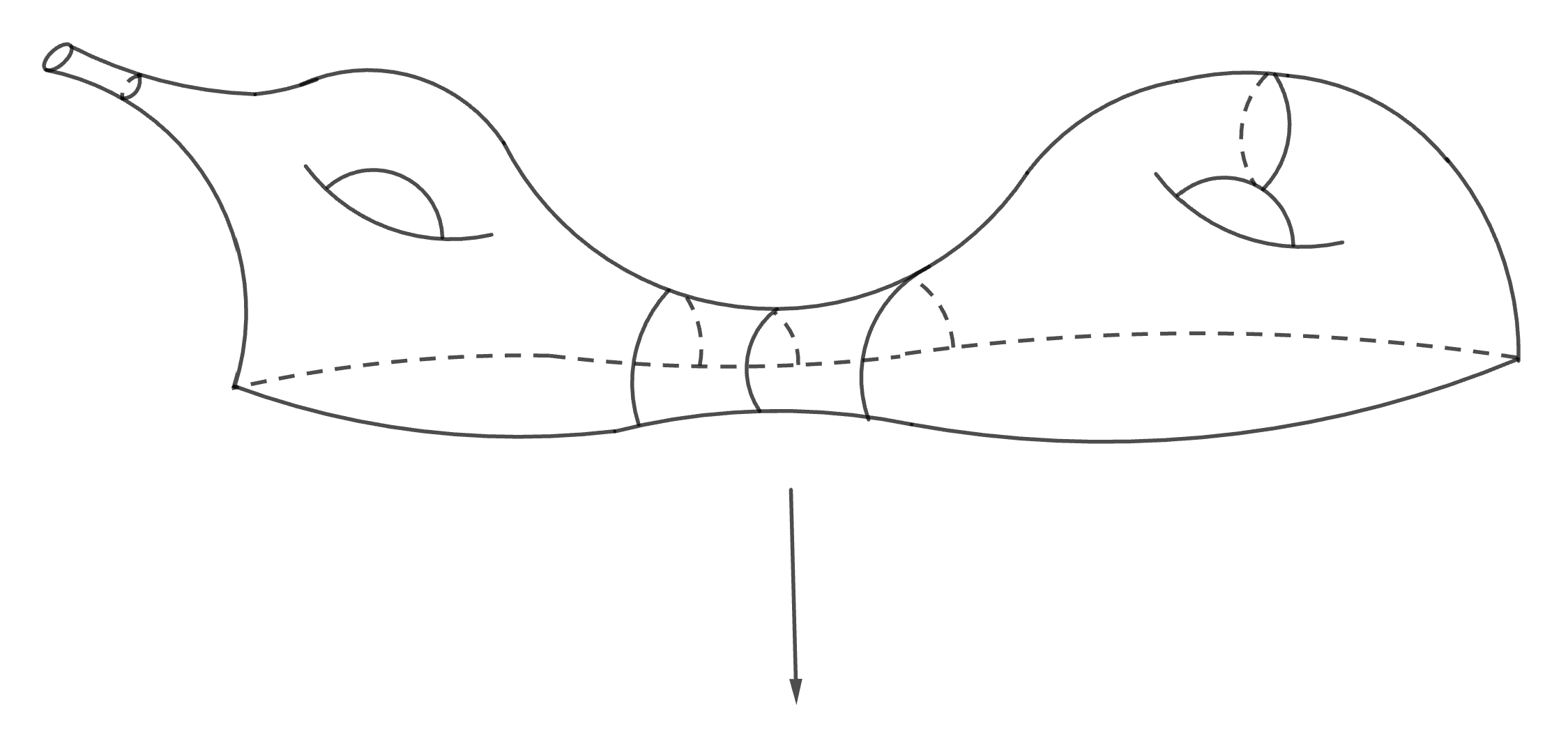}
\centering
\includegraphics[width=0.6\textwidth]{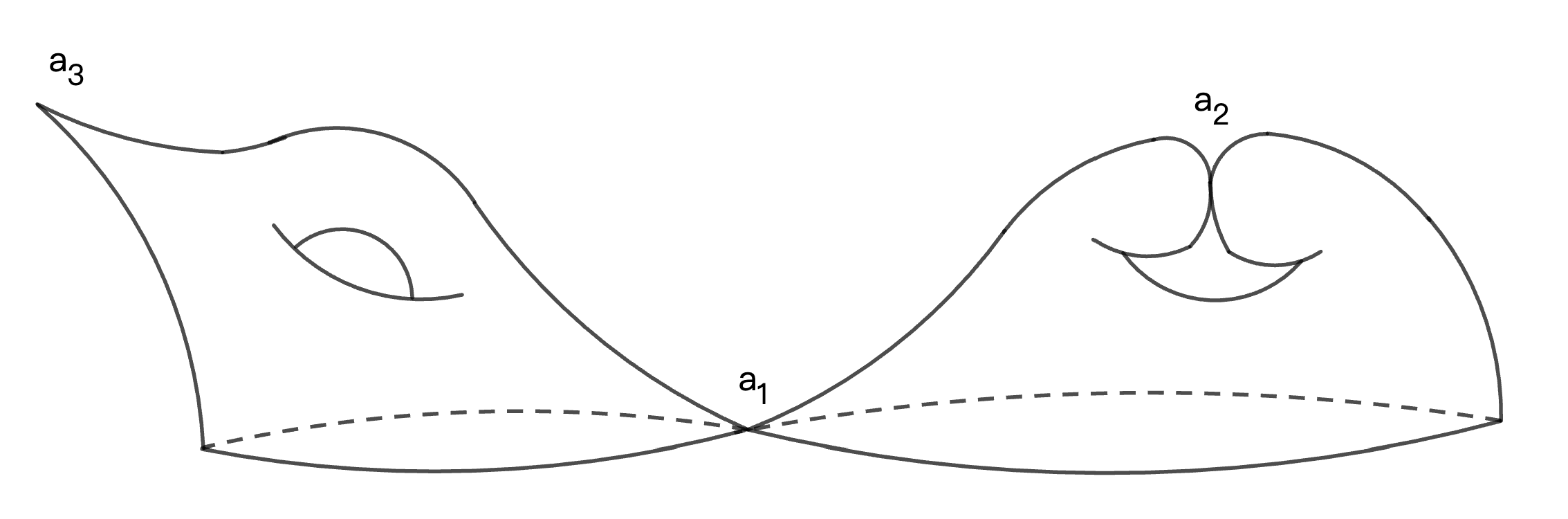}
\caption{Degeneration to a nodal surface.}
\label{fig:nodal-degeneration}
\end{figure}

Figure~\ref{fig:nodal-degeneration} illustrates the three types of nodes. Pinching an arc with endpoints on the boundary produces the type~I boundary node $a_1$, while pinching an interior simple closed curve produces the interior node $a_2$. Contracting a boundary component produces the type~II boundary node $a_3$.

\begin{definition}
\label{def:bordered-DM-convergence}
Let $j_n$ be conformal structures on the fixed smooth bordered surface $\Sigma$.
We say that
\[
    (\Sigma,j_n)
    \longrightarrow
    (\Sigma_0,j_0)
\]
in the bordered Deligne--Mumford sense if, after passing to a subsequence, the following conditions hold.

Choose the constant curvature representatives with geodesic boundary from Subsection~\ref{subsec:uniformization-boundary}, using the normalized models there for the disk and annulus.
There are finitely many pairwise disjoint degenerating geodesics
\[
    \Gamma_n = \Gamma_n^{\mathrm{int}} \sqcup \Gamma_n^{\mathrm I}\sqcup\Gamma_n^{\mathrm{II}},
\]
whose lengths tend to zero, where
\begin{enumerate}[label=\textup{(\alph*)},leftmargin=2.4em]
\item $\Gamma_n^{\mathrm{int}}$ consists of interior simple closed geodesics;

\item $\Gamma_n^{\mathrm I}$ consists of properly embedded geodesic arcs meeting $\partial\Sigma$ orthogonally;

\item $\Gamma_n^{\mathrm{II}}$ consists of boundary components.
\end{enumerate}
The corresponding collar neighborhoods are conformally equivalent, respectively, to
\begin{align}
    C_n^{\mathrm{int}}
    &=(-T_n,T_n)\times S^1,
    \label{eq:sec2-interior-collar-model}
    \\C_n^{\mathrm I}
    &=(-T_n,T_n)\times[0,\pi],
    \label{eq:sec2-type-I-collar-model}\\
    C_n^{\mathrm{II}}
    &=[0,T_n)\times S^1,
    \label{eq:sec2-type-II-collar-model}
\end{align}
where the corresponding numbers $T_n$ tend to infinity.
In the type~I model, the boundary portions are
\[
    (-T_n,T_n)\times\{0,\pi\},
\]
and the collapsing arc is represented by $\{0\}\times[0,\pi]$.
In the type~II model, the shrinking boundary component is represented by
\[
    \{0\}\times S^1.
\]
Let
\[
    \widetilde{\mathcal N}_0:=\pi_0^{-1}(\cN_0)
    \subset\widetilde\Sigma_0.
\]
There exist an exhaustion
\[
    V_1\Subset V_2\Subset\cdots
    \Subset\widetilde\Sigma_0\setminus\widetilde{\mathcal N}_0,
    \qquad
    \bigcup_{m=1}^\infty V_m
    =
    \widetilde\Sigma_0\setminus\widetilde{\mathcal N}_0,
\]
and, after passing to a diagonal subsequence, smooth embeddings that preserve the boundary
\[
    \Psi_n:V_n\longrightarrow\Sigma
\]
whose images are the complements of suitable truncations of the collars, such that
\begin{equation*}    \Psi_n^*j_n
    \longrightarrow
    \widetilde j_0
    \qquad\text{in  }C^\infty_{\loc}\bigl(\widetilde\Sigma_0\setminus\widetilde{\mathcal N}_0\bigr).
\end{equation*}
Here $\widetilde j_0$ is the conformal structure on the normalization, and the convergence is understood up to boundary points away from the nodes, using boundary conformal charts.
\end{definition}

When $\chi(\Sigma)<0$, subsequential compactness in Definition~\ref{def:bordered-DM-convergence} follows from bordered Deligne--Mumford compactness; see \cite[Theorem~4.9]{Liu2004}.
Equivalently, one applies compactness equivariantly to the complex doubles and takes the quotient by the antiholomorphic involution.
Formulations using complex doubles and bordered stable curves appear in \cite{Katz-Liu2001,Liu2004,Liu-Song-Zhu2022}.
The three collar models \eqref{eq:sec2-interior-collar-model}-\eqref{eq:sec2-type-II-collar-model} are exhaustive; type~II must be included unless boundary component lengths are bounded away from zero.

\Needspace{6\baselineskip}
\subsection{\texorpdfstring{Bordered stratified surfaces and bubble tree domains}{Bordered stratified surfaces and bubble tree domains}}
\label{subsec:bordered-stratified}

\begin{definition}
\label{def:bordered-stratified-surface}
A bordered stratified Riemann surface is a compact connected nodal bordered Riemann surface $\Sigma_\infty$ whose normalization is a finite disjoint union
\[
    \widetilde\Sigma_\infty
    =
    \bigsqcup_{i\in\mathfrak I}
        \Sigma_\infty^i
\]
of compact connected smooth closed or bordered Riemann surfaces.
The components $\Sigma_\infty^i$ are called the irreducible components.
We denote by
\[
    \pi_\infty:
    \widetilde\Sigma_\infty
    \longrightarrow
    \Sigma_\infty
\]
the normalization map and by $\cN_\infty$ the nodal set.
For each component, let
\[
    Z_i
    :=
    \pi_\infty^{-1}(\cN_\infty)
    \cap\Sigma_\infty^i
\]
be its finite set of special points.
At a type~II boundary node, $Z_i$ contains the unique interior preimage of that node.
\end{definition}

The regular part is denoted by $\Sigma_\infty^{\mathrm{reg}}:=\Sigma_\infty\setminus\cN_\infty$. This notation removes the nodes and is independent of any branch points of a map.
Its boundary is
\[
 \partial\Sigma_\infty^{\mathrm{reg}}
 :=\pi_\infty\left(\bigsqcup_{i\in\mathfrak I}
       (\partial\Sigma_\infty^i\setminus Z_i)\right).
\]
Ordinary interior nodes identify two interior special points, while ordinary boundary nodes identify two boundary special points.
Type~II nodes retain the one-branch local model specified in Definition~\ref{def:nodal-bordered-surface}.

\begin{definition}
\label{def:free-boundary-bubble tree-surface}
Let $\Sigma_0$ be a nodal bordered Riemann surface.
A \emph{free boundary bubble tree surface over $\Sigma_0$} is a bordered stratified Riemann surface $\Sigma_\infty$ together with a continuous stabilization map
\[
    \mathfrak s:\Sigma_\infty\longrightarrow\Sigma_0
\]
with the following properties.

The index set of the normalization has a decomposition
\[
    \mathfrak I
    =
    \mathfrak I_{\mathrm{pr}}
    \sqcup
    \mathfrak I_{S^2}
    \sqcup
    \mathfrak I_{D},
\]
where
\begin{enumerate}[label=\textup{(\roman*)},leftmargin=2.4em]
\item the components indexed by $\mathfrak I_{\mathrm{pr}}$ are the \emph{principal components} and are naturally identified with the components of the normalization $\widetilde\Sigma_0$;

\item $\Sigma_\infty^i\simeq S^2$ for $i\in\mathfrak I_{S^2}$;

\item $\Sigma_\infty^i\simeq\overline D$ for $i\in\mathfrak I_{D}$;

\item on each principal component, $\mathfrak s\circ\pi_\infty$ agrees with the corresponding restriction of the normalization map $\pi_0:\widetilde\Sigma_0\to\Sigma_0$, and $\mathfrak s\circ\pi_\infty$ is constant on every sphere or disk component;

\item For each $x\in\Sigma_0$, define the \emph{local fiber graph} as follows.
Its vertices are the individual points of $\pi_0^{-1}(x)$, counted separately even when two such points lie on the same principal component, together with the sphere and disk components on which $\mathfrak s$ is constant with value $x$.
Its edges correspond to the ordinary nodes over $x$ and join the vertices incident to their two branches.
A retained type~II one-branch node is represented by a marking at its unique incident vertex, rather than by a loop edge.
This graph is connected and contains no cycle.
The global dual graph need not be a tree.
Cycles inherited from the original nodal surface, including those created by nonseparating pinches, are retained.
In particular, an ordinary node with no attached bubble joins its two distinct branch vertices even if both preimages lie on the same principal component. For $x\in\Sigma_0\setminus\cN_0$ with no bubbles, the graph has one vertex.
\end{enumerate}
When component maps are specified, a component on which the map is constant may be retained to record the nodal incidence relations.
We call such a component a \emph{ghost component}.
\end{definition}

Sphere bubbles attach at interior points of principal, sphere, or disk components.
At a boundary concentration point, a disk component, possibly a ghost component, attaches through its boundary; spheres attach at interior points of that disk.
Thus no boundary point is joined directly to a sphere interior.
A type~II end either remains a type~II node or is resolved by disk and sphere components at its unique normalized special point.
All nodes have the types in Definition~\ref{def:nodal-bordered-surface}.

A collection $\mathbf u_\infty=\{u_\infty^i:\Sigma_\infty^i\to N\}_{i\in\mathfrak I}$ is a \emph{free boundary bubble tree map} if every nonconstant component is a $W^{2,2}$ branched conformal immersion and each bordered component satisfies
\begin{equation}\label{eq:sec2-component-free-boundary-condition}
    u_\infty^i(\partial\Sigma_\infty^i)\subset\cK,
    \qquad
    du_\infty^i(\nu_i)
    \perp_{\bar g}T_{u_\infty^i}\cK
\end{equation}
for $\mathcal H^1$-a.e. point of $\partial\Sigma_\infty^i$, in the Sobolev trace sense.
On a nonconstant component, the Sobolev boundary trace of $du_\infty^i(\nu_i)$ is taken on $\partial\Sigma_\infty^i\setminus\cB_{u_\infty^i}$; on a constant component, $du_\infty^i=0$, so the orthogonality condition is automatic.
Here $\nu_i$ is the outward unit conormal with respect to any smooth metric on $\Sigma_\infty^i$ compatible with its conformal structure.
If the continuous representatives of the component maps agree at the two branches of every ordinary node, they induce a continuous map
\[
u_\infty:\Sigma_\infty\longrightarrow N.
\]
\Needspace{6\baselineskip}
\subsection{\texorpdfstring{Free boundary bubble tree convergence}{Free boundary bubble tree convergence}}
\label{subsec:bt-convergence-revised}

Let $h_n$ be smooth metrics on $\Sigma$, put $j_n=[h_n]$, and assume $(\Sigma,j_n)\to(\Sigma_0,j_0)$ as in Definition~\ref{def:bordered-DM-convergence}.
Let $\Sigma_\infty$ be a free boundary bubble tree surface over $\Sigma_0$, with normalization $\widetilde\Sigma_\infty=\bigsqcup_{i\in\mathfrak I}\Sigma_\infty^i$.
\begin{definition}
\label{def:componentwise-free-boundary-bubble tree-convergence}
Let
\[
    u_n\in
    W^{2,2}_{b,c,\mathrm{fb}}
    \bigl((\Sigma,h_n),N,\cK\bigr).
\]
We say that $u_n$ converges \emph{componentwise in the free boundary bubble tree sense} to a collection
\[
    \mathbf u_\infty
    =
    \{u_\infty^i:
      \Sigma_\infty^i\to N\}_{i\in\mathfrak I}
\]
if, after passing to a subsequence, the following conditions hold.

\begin{enumerate}[label=\textup{(\roman*)},leftmargin=2.4em]
\item For each $i\in\mathfrak I$, there are relatively compact open sets
\[
        V_n^i
        \Subset
        \Sigma_\infty^i\setminus Z_i,
        \qquad
        V_n^i\subset V_{n+1}^i,
        \qquad
        \bigcup_{n=1}^\infty V_n^i
        =
        \Sigma_\infty^i\setminus Z_i,
    \]
and smooth embeddings
\[
        \Phi_n^i:
        V_n^i\longrightarrow\Sigma.
    \]
For fixed $n$, the images of the maps $\Phi_n^i$ are pairwise disjoint.
If $\Sigma_\infty^i$ is bordered, then
\begin{equation*}        \Phi_n^i
        \bigl(V_n^i\cap\partial\Sigma_\infty^i\bigr)
        \subset\partial\Sigma.
    \end{equation*}
\item If $j_i$ denotes the conformal structure on $\Sigma_\infty^i$, then for every $L\Subset\Sigma_\infty^i\setminus Z_i$,
\begin{equation*}        (\Phi_n^i)^*j_n
        \longrightarrow
        j_i
        \qquad\text{in }C^\infty(L).
    \end{equation*}
On a principal component, the maps $\Phi_n^i$ are induced by the Deligne--Mumford comparison maps.
On a sphere or disk component, they are the conformal blow-up reparametrizations arising in the bubble extraction.

\item There is a finite exceptional set $\cB_\infty^i\subset\Sigma_\infty^i$ such that, for every
\[
        L\Subset
        \Sigma_\infty^i
        \setminus
        \bigl(Z_i\cup\cB_\infty^i\bigr),
    \]
one has
\begin{align*}
        u_n\circ\Phi_n^i
        &\rightharpoonup
        u_\infty^i
        &&\text{weakly in }W^{2,2}(L),
        \\
        u_n\circ\Phi_n^i
        &\longrightarrow
        u_\infty^i
        &&\text{strongly in }W^{1,p}(L)
        \quad\text{for every }p<\infty.
    \end{align*}
If $L$ meets $\partial\Sigma_\infty^i$, the convergence is understood in boundary conformal charts.

\item Every nonconstant $u_\infty^i$ is a $W^{2,2}$ branched conformal immersion.
If $\Sigma_\infty^i$ is bordered, then $u_\infty^i$ satisfies the free boundary condition \eqref{eq:sec2-component-free-boundary-condition}.

\item The complement
\begin{equation*}        \Sigma\setminus
        \bigcup_{i\in\mathfrak I}
            \Phi_n^i(V_n^i)
    \end{equation*}
is a finite union of annular necks, half-annular necks, and truncated subregions of the three degenerating collar models \eqref{eq:sec2-interior-collar-model}-\eqref{eq:sec2-type-II-collar-model}.
\end{enumerate}
\end{definition}

For each ordinary node, bubble attachment, or type~II end, let $\cN_{n,e}(R)$ denote the associated neck after removing subregions of conformal length $R$ adjacent to its artificial attaching ends.
The portions inherited from $\partial\Sigma$ are retained.
In particular, this convention never removes the circle $\{0\}\times S^1$ in a type~II collar.
We fix these neck parametrizations and their artificial attaching ends.
Uniformly bounded longitudinal shifts at the artificial ends leave both limits in Definition~\ref{def:no-neck-properties} unchanged, provided the boundary inherited from $\partial\Sigma$ is kept fixed.

For a map $u$ into a Riemannian manifold $N$ and a subset $\mathfrak A$ of its domain, we write
\begin{equation}\label{eq:oscillation-definition}
    \operatorname{Osc}^N_\mathfrak A u
    :=\sup_{x,y\in \mathfrak A}d_N\bigl(u(x),u(y)\bigr).
\end{equation}
We set $\operatorname{Osc}^N_{\varnothing}u=0$.
Here $d_N$ is the Riemannian distance; for $N=\R^K$, it is the Euclidean distance. We omit the superscript when the target is understood.

\begin{definition}
\label{def:no-neck-properties}
A componentwise bubble tree convergent sequence has the \emph{no-neck property} if, for every neck label $e$,
\begin{equation*}    \lim_{R\to\infty}
    \limsup_{n\to\infty}
    \operatorname{Osc}_{\cN_{n,e}(R)}u_n
    =0.
\end{equation*}
It has the \emph{property of vanishing neck energy} if
\begin{equation*}    \lim_{R\to\infty}
    \limsup_{n\to\infty}
    E_{h_n}
    \bigl(u_n;\cN_{n,e}(R)\bigr)
    =0
\end{equation*}
for every neck label $e$.
\end{definition}

The no-neck property ensures that the component values agree at every ordinary attaching pair.
At a type~II end, suppose the oscillation of $u_n$ on the entire unresolved chain from the shrinking component of $\partial\Sigma$ to the unique normalized special point tends to zero.
The shrinking boundary component and the normalized special point then have the same limiting value, which belongs to $\cK$ because $\cK$ is closed.

If only the terminal profile containing the shrinking boundary circle collapses while a nonconstant disk or sphere chain remains, only the value of that collapsed profile is forced to lie in $\cK$.
The value at the principal normalized point need not lie in $\cK$.
If a nonconstant terminal disk resolves the end, its boundary circle maps into $\cK$, but its interior attaching value need not belong to $\cK$.

In the latter two cases, the no-neck property ensures that the values agree at each ordinary attachment along the intervening disk and sphere chain.
In every case, $\mathbf u_\infty$ induces a continuous map
\[
    u_\infty:\Sigma_\infty\longrightarrow N.
\]
For a bubble tree map, let $\mathfrak I_*:=\{i\in\mathfrak I:u_\infty^i\text{ is nonconstant}\}$.
The quantities $E(u_\infty)$, $\cA(u_\infty)$, and $W_N(u_\infty)$ are defined by summing the corresponding Dirichlet energies, areas, and Willmore energies over the nonconstant component maps.
The Dirichlet energy of $u_\infty^i$ is computed with any smooth metric on $\Sigma_\infty^i$ compatible with its conformal structure.
Ghost components contribute zero to these sums.

\vspace{1cm}

\section{\texorpdfstring{Boundary estimates for conformal immersions}{Boundary estimates for conformal immersions}}
\label{sec:fundamental-ingredients-bubbling}

We prove the boundary estimates used in Theorem~\ref{thm:intro-main-compactness} for the class in Definition~\ref{def:fb-branched-conformal}.

\Needspace{6\baselineskip}
\subsection{\texorpdfstring{Small energy estimates}{Small energy estimates}}
\label{subsec:general-small-energy}

\begin{theorem}\label{intsma0}

Fix a smooth compact embedded submanifold $\cK\subset\R^K$.
For any $1<p<2$ and $0<r<R$, there are constants $\varepsilon_0>0$ and $C<\infty$, depending only on $p,K,\cK$, and $r/R$, with the following property.
Let
\[
u\in W^{2,2}( D_R^+,\R^K)
\]
be weakly conformal.
Suppose its trace satisfies
\[
u(I_R)\subset\cK,
\qquad
\partial_yu\perp T_u\cK
\quad\text{a.e. on }I_R,
\]
and that there is a measurable vector field $\mathbf H:D_R^+\to\R^K$ such that
\[
\Delta u=\frac12\mathbf H|\nabla u|^2
\quad\text{a.e. on }D_R^+.
\]
If
\begin{align*}E(u; D_R^+)+W_{\R^K}(u; D_R^+)<\varepsilon_0^2,
\end{align*}
where
\[
E(u; D_R^+)=\frac12\int_{ D_R^+}|\nabla u|^2,
\qquad
W_{\R^K}(u; D_R^+)
=\frac18\int_{ D_R^+}|\mathbf H|^2|\nabla u|^2,
\]
then
\[
R^{1-\frac2p}\|\nabla u\|_{L^p( D_r^+)}
+
R^{2-\frac2p}\|D^2u\|_{L^p( D_r^+)}
\le
C\|\nabla u\|_{L^2( D_R^+)}.
\]
\end{theorem}

\begin{proof}

We first consider $R=1$ and $r=1/2$.
Choose a tubular neighborhood of $\cK$, and let $\pi_{\cK}$ denote its nearest-point projection.
The normal reflection
\[
S_0=2\pi_{\cK}-\operatorname{Id}
\]
is smooth on this neighborhood.
Extend it to a smooth map $S:\R^K\to\R^K$ with bounded first and second derivatives.
On $\cK$, $S$ is the identity, while $DS$ is the identity on $T\cK$ and minus the identity on $(T\cK)^\perp$.

Set $\rho(x,y)=(x,-y)$, and define
\[
U(x,y)=
\begin{cases}
u(x,y),&y\geq0,\\
S(u(x,-y)),&y<0.
\end{cases}
\]
The trace of $u_x$ is tangent to $\cK$.
The trace of $u_y$ is normal to $\cK$.
Hence, on the interface,
\[
DS(u)u_x=u_x,
\qquad
-DS(u)u_y=u_y.
\]
Thus the traces of $U,U_x,U_y$ match across the interface.
These matching traces permit Sobolev gluing.
The two-dimensional Sobolev embedding~\cite[Chapter~4]{AdamsFournier2003} makes the first derivatives integrable to the fourth power.
The bounded first and second derivatives of the reflection therefore justify the $W^{2,2}$ chain rule and give $U\in W^{2,2}( D_1,\R^K)$.
Moreover, $\Delta U$ contains no measure supported on the reflection line $y=0$, and
\[
\|\nabla U\|_{L^2( D_1)}
\leq C_{\cK}\|\nabla u\|_{L^2( D_1^+)}.
\]
On the lower half disk, write $v=u\circ\rho$.
The chain rule and the conformal immersion equation give
\[
\Delta U
=
DS(v)\Delta v+D^2S(v)[\nabla v,\nabla v].
\]
Consequently,
\[
|\Delta U|
\leq
C_{\cK}(1+|\mathbf H\circ\rho|)|\nabla v|^2.
\]
Also,
\[
\|(1+|\mathbf H|)|\nabla u|\|_{L^2( D_1^+)}
\leq
C\bigl(E(u; D_1^+)^{1/2}
      +W_{\R^K}(u; D_1^+)^{1/2}\bigr),
\]
since
$$\||\mathbf H||\nabla u|\|_{L^2( D_1^+)}^2 =8W_{\R^K}(u; D_1^+).$$

Choose a radial cutoff $\eta\in C_c^\infty( D_1)$ equal to one on $ D_{1/2}$.
Let $c$ be the average of $U$ over $D_1$, and set $w=\eta(U-c)$, $p^*=2p/(2-p)$.
Apply H\"older's inequality, reflect the integrals over the lower half, and use the identity $\eta\nabla u=\nabla w-(u-c)\nabla\eta$ on the upper half to obtain
\[
\|\eta\Delta U\|_{L^p}
\leq C_{\cK}(E^{1/2}+W^{1/2})
       \|\eta\nabla u\|_{L^{p^*}},
\qquad
\|\eta\nabla u\|_{L^{p^*}}
\leq C_p\bigl(\|w\|_{W^{2,p}}+\|\nabla U\|_{L^2}\bigr).
\]
Here both energies are evaluated on $ D_1^+$.
The second estimate follows from Sobolev embedding and Poincar\'e's inequality.
Since
\[\Delta w=\eta\Delta U+2\nabla\eta\cdot\nabla U+(U-c)\Delta\eta,\]
the zero-boundary Calder\'on--Zygmund estimate~\cite{GilbargTrudinger2001} yields
\[
\|w\|_{W^{2,p}}
\leq C_p\|\nabla u\|_{L^2( D_1^+)}
+C_{p,\cK}(E^{1/2}+W^{1/2})\|w\|_{W^{2,p}}.
\]
Choose $\varepsilon_0$ so that $2C_{p,\cK}\varepsilon_0<1$ and absorb the last term.
The same argument with $\eta=1$ on $ D_{r/R}$ proves the estimate for any fixed $r/R<1$, with constants depending on this ratio.
For $u_R(z)=u(Rz)$, the $L^p$-norms of $\nabla u_R$ and $D^2u_R$ scale by $R^{1-2/p}$ and $R^{2-2/p}$, respectively, while $E$ and $W$ are unchanged.
This gives the stated estimate.
Poincar\'e's inequality, with $\bar u_R=| D_R^+|^{-1}\int_{ D_R^+}u$, also gives the translation-invariant form
\begin{align*}&R^{-\frac2p}\|u-\bar u_R\|_{L^p( D_r^+)}
+R^{1-\frac2p}\|\nabla u\|_{L^p( D_r^+)}
+R^{2-\frac2p}\|D^2u\|_{L^p( D_r^+)}
\leq C\|\nabla u\|_{L^2( D_R^+)}.
\end{align*}
\end{proof}

For the projective Gauss maps below, let $g_{\mathrm{FS}}$ be the Fubini--Study metric on $\mathbb{CP}^{K-1}$, normalized so that the Hopf map from the unit sphere is a Riemannian submersion.
Let $\omega_{\mathrm{FS}}$ be the associated K\"ahler form.
We use $\Omega_{\mathbb{CP}}=2\omega_{\mathrm{FS}}$, as in \cite[Theorem~2.1]{kuwert2012}.
All projective differential norms and area Jacobians are computed with $g_{\mathrm{FS}}$, and
\[
E(G;U)=\frac12\int_U|dG|_{g_{\mathrm{FS}}}^2.
\]

\begin{lemma}\label{lem:weak-boundary-connection}
Let $u\in W^{2,2}(D_r^+,\R^K)$ be an unbranched conformal immersion with
\[
u^*g_{\R^K}=e^{2\lambda}(dx^2+dy^2),
\qquad \lambda\in L^\infty(D_r^+).
\]
Suppose that $u(I_r)\subset\cK$ and $u_y\perp T_u\cK$ in the trace sense.
On a neighborhood of the image, choose a smooth orthogonal matrix field $O$ and a fixed subspace $V\subset\R^K$ such that
\[
O(a)T_a\cK=V,
\qquad O(a)(T_a\cK)^\perp=V^\perp
\quad\text{for }a\in\cK\text{ in this neighborhood}.
\]
Set
\[
\tau_1=e^{-\lambda}u_x,\qquad
\tau_2=e^{-\lambda}u_y,\qquad
e_j=O(u)\tau_j,
\]
and write $\beta=\beta_1\,dx+\beta_2\,dy$, where
\[
\beta_j=\langle\partial_j(O(u))\tau_1,e_2\rangle,
\qquad
F_\beta=(-\beta_2,\beta_1).
\]
Let $\widehat G$ be the extension to $D_r$ of $G=[(e_1+i e_2)/\sqrt2]$ by reflection under the antiholomorphic isometry $[z]\mapsto[R_V\overline z]$, where $R_V$ is reflection across $V$.
If an even function $\psi\in W^{1,2}_{\loc}(D_r)$ satisfies
\[
-\Delta\psi=*\widehat G^*\Omega_{\mathbb{CP}}
\quad\text{in }\mathcal D'(D_r),
\]
then $h:=\lambda-\psi|_{D_r^+}$ satisfies
\begin{equation}\label{eq:weak-boundary-conformal-factor}
\int_{D_r^+}\nabla h\cdot\nabla\varphi
=\int_{D_r^+}F_\beta\cdot\nabla\varphi
\end{equation}
for every smooth test function whose support is compact in $D_r^+\cup I_r$.
Thus $\Delta h=\operatorname{div}F_\beta$ with the natural Neumann condition on $I_r$, both in the weak sense.
On every smaller half-chart, the weak curvature identity is
\begin{equation}\label{eq:weak-boundary-curvature}
\int_{D_r^+}\nabla\lambda\cdot\nabla\varphi
=\int_{D_r^+}K_{g_u} e^{2\lambda}\varphi
 -\int_{I_r}\beta_1\varphi\,dx,
\end{equation}
with $g_u$ as in Subsection~\ref{sec:fbbci}, and
\begin{equation}\label{eq:weak-boundary-support-curvature}
-\beta_1
=e^\lambda
\langle A^{\cK}_{\R^K}(\tau_1,\tau_1),\tau_2\rangle
\quad\text{a.e. on }I_r.
\end{equation}
Thus $-\beta_1\,dx$ is the geodesic curvature measure of the boundary interval $I_r$ with respect to the induced metric $g_u$, computed with the inward unit conormal.
The assertion also holds locally on a punctured half disk for test functions supported away from the puncture.
\end{lemma}

\begin{proof}
Applying the chain rule to $\lambda=\log|u_x|$ gives $\lambda\in W^{1,2}$ on smaller closed half disks.
Consequently $\tau_1,\tau_2,e_1,e_2\in W^{1,2}\cap L^\infty$ there.
The trace of $u_x$ is tangent to $\cK$.
Since the trace of $u_y$ is normal to $\cK$, the traces of $e_1$ and $e_2$ belong to $V$ and $V^\perp$, respectively.
Reflect these frames by
\[
\widehat e_1(x,y)=R_Ve_1(x,-y),\qquad
\widehat e_2(x,y)=-R_Ve_2(x,-y)
\quad\text{for }y<0,
\]
and retain $e_1,e_2$ for $y\geq0$.
The matching traces give an orthonormal $W^{1,2}$ frame on every smaller full disk.
This frame represents $\widehat G$.

Let $a_j=\langle\partial_j e_1,e_2\rangle$.
The first derivatives of $u$ are bounded and belong to $W^{1,2}$, so we may differentiate its conformality identities to obtain
\[
a_1=-\lambda_y+\beta_1,
\qquad
a_2=\lambda_x+\beta_2.
\]
The connection $\widehat a_j=\langle\partial_j\widehat e_1,\widehat e_2\rangle$ belongs to $L^2$ on the full disk.
Its first component is the odd reflection of $a_1$, and its second component is the even reflection of $a_2$.
Writing $\lambda^{\mathrm e}$ for the even reflection of $\lambda$ and reflecting $\beta_1$ oddly and $\beta_2$ evenly, we therefore have
\[
\widehat a
=(-\partial_y\lambda^{\mathrm e},\partial_x\lambda^{\mathrm e})
+\widehat\beta
\quad\text{a.e. on }D_r.
\]
The projective lift satisfies the distributional curvature identity
\begin{equation}\label{eq:reflected-connection-curvature}
*\widehat G^*\Omega_{\mathbb{CP}}
=-\operatorname{curl}\widehat a
=-\Delta\lambda^{\mathrm e}
 -\operatorname{curl}\widehat\beta.
\end{equation}
Indeed, the $W^{1,2}$ product rule gives
\[
-\operatorname{curl}\widehat a
=\sum_{\ell=1}^K
\bigl(
\partial_x\widehat e_1^\ell\,\partial_y\widehat e_2^\ell
-\partial_y\widehat e_1^\ell\,\partial_x\widehat e_2^\ell
\bigr)
\]
in distributions.
To justify this identity, mollify the bounded frame components on a smaller full disk.
The connection products converge in $L^2$, and the displayed Jacobians converge in $L^1$.
Passing to the limit in the smooth identity gives the claim.
The right-hand side is also the a.e. pullback of $\Omega_{\mathbb{CP}}$, because the unit lift $(\widehat e_1+i\widehat e_2)/\sqrt2$ has connection $-\widehat a$.
In particular, the reflected curvature identity contains no interface measure.

Subtracting the equation for $\psi$ from \eqref{eq:reflected-connection-curvature} yields
\[
\Delta(\lambda^{\mathrm e}-\psi)
=-\operatorname{curl}\widehat\beta
=\operatorname{div}(-\widehat\beta_2,\widehat\beta_1)
\quad\text{in }\mathcal D'(D_r).
\]
Both sides are divergences of $L^2$ fields.
By density, the weak equation therefore holds for compactly supported $W^{1,2}$ test functions.
Use the even extension of $\varphi$ as a test function.
Since $\lambda^{\mathrm e}-\psi$ is even and the components of $\widehat\beta$ have the stated parities, each integral is twice the corresponding integral over the upper half.
Dividing by $2$ proves \eqref{eq:weak-boundary-conformal-factor}.

We next identify the boundary flux without assuming that $\lambda$ has a classical normal derivative.
On a smaller half-chart, $\nabla u\in L^\infty$ and the normalized frames belong to $W^{1,2}\cap L^\infty$.
The definition of $\beta$ gives $\beta\in W^{1,2}$ there.
Its boundary trace is therefore well defined.
The same frame-curvature calculation on the upper half disk gives
\[
*G^*\Omega_{\mathbb{CP}}=K_{g_u} e^{2\lambda}-\operatorname{curl}\beta.
\]
Choose the even zero-Dirichlet solution of $-\Delta\psi=*\widehat G^*\Omega_{\mathbb{CP}}$ on a smaller full disk.
The right-hand side is $-\operatorname{curl}\widehat a\in W^{-1,2}$, and uniqueness gives evenness.
For smooth test functions as above, even reflection gives
\[
\int_{D_r^+}\nabla\psi\cdot\nabla\varphi
=\int_{D_r^+}
\bigl(K_{g_u} e^{2\lambda}-\operatorname{curl}\beta\bigr)\varphi.
\]
Adding this identity to \eqref{eq:weak-boundary-conformal-factor} and integrating the $W^{1,2}$ field $\beta$ by parts cancels its interior curl.
Since the outward normal to the diameter is $-\partial_y$, the remaining boundary term is $-\beta_1\varphi$.
This proves \eqref{eq:weak-boundary-curvature}.

To identify the boundary curvature density, fix $a\in\cK$, a unit tangent vector $v\in T_a\cK$, and a unit normal vector $n\in(T_a\cK)^\perp$.
Choose a smooth tangent vector field $X$ on $\cK$ with $X(a)=v$.
The derivative along $v$ of the vector field $O X\in V$ is orthogonal to $O(a)n\in V^\perp$.
Consequently,
\[
\langle DO(a)[v]v,O(a)n\rangle
=-\langle A^{\cK}_{\R^K}(v,v),n\rangle.
\]
Apply this pointwise identity to the boundary traces $v=\tau_1$, $n=\tau_2$, and use $u_x=e^\lambda\tau_1$.
This proves \eqref{eq:weak-boundary-support-curvature}.
Since the construction is local, the same argument applies away from a boundary puncture.
\end{proof}

\Needspace{6\baselineskip}
\subsection{\texorpdfstring{Boundary removability and branch orders}{Boundary removability and branch orders}}
\label{subsec:general-boundary-removability}

Lemma~\ref{lem:general-boundary-removability} proves boundary removability and integer branch order.
The scalar estimate in Lemma~\ref{lem:delayed-energy-decay} is also used in Theorem~\ref{deca}.

\begin{lemma}\label{lem:delayed-energy-decay}
Let $L>2$, $\lambda>0$, and $\delta\geq0$, and set $\beta=\delta e^{2\lambda}$.
Suppose that $\beta<\lambda$ and that $\phi:[0,L]\to[0,\infty)$ is nondecreasing and absolutely continuous, with
\[
 \lambda\phi(t)\leq\phi'(t)+\delta\phi(t+2)
 \quad\text{for a.e. }0\leq t\leq L-2.
\]
Then, for every $0\leq t\leq L$,
\[
 \phi(t)\leq
 \left(1+\frac{\beta}{\lambda}\right)e^{2\lambda}
 \phi(L)e^{-(\lambda-\beta)(L-t)}.
\]
\end{lemma}

\begin{proof}
Multiplying by $e^{-\lambda t}$ and integrating from $t$ to $L-2$ gives
\[
e^{-\lambda t}\phi(t)
\leq
e^{-\lambda(L-2)}\phi(L)
+
\beta\int_{t+2}^{L}e^{-\lambda r}\phi(r)\,dr.
\]
With
\[
F(t)=\int_t^{L-2}e^{-\lambda r}\phi(r)\,dr,
\qquad
A=
\left(1+\frac{\beta}{\lambda}\right)
e^{-\lambda(L-2)}\phi(L),
\]
monotonicity of $\phi$ gives
\[
e^{-\lambda t}\phi(t)\leq A+\beta F(t),
\qquad
-F'(t)\leq A+\beta F(t).
\]
Since $F(L-2)=0$,
\[
F(t)
\leq
\frac{A}{\beta}
\left(e^{\beta(L-2-t)}-1\right)
\]
when $\beta>0$; the case $\beta=0$ follows directly.
Substituting this bound gives
\[
\phi(t)
\leq
\left(1+\frac{\beta}{\lambda}\right)e^{2\lambda}
\phi(L)e^{-(\lambda-\beta)(L-t)}
\]
for $0\leq t\leq L-2$.
For $L-2\leq t\leq L$, the same estimate follows from monotonicity.
\end{proof}

\begin{lemma}\label{lem:boundary-layer-length}
Let $u\in W^{2,2}_{\loc}(D_r^+\cup I_r,\R^K)$ be an unbranched conformal immersion with $\lambda\in L^\infty_{\loc}(D_r^+\cup I_r)$, and suppose that $u(I_r)\subset\cK$ and $u_y\perp T_u\cK$ in the trace sense.
Write $v=\operatorname{dist}(u,\cK)$, let $\tau>0$ be a tubular radius of $\cK$, and choose a nonnegative smooth function $\chi$ on $[0,\infty)$ supported in $[0,\tau)$ and equal to $1$ near $0$.
For every nonnegative continuous function $\varphi$ compactly supported in $D_r^+\cup I_r$, we have
\begin{equation}\label{eq:boundary-layer-length}
 \int_{I_r}\varphi(x,0)|u_x(x,0)|\,dx
 \leq
 \liminf_{\varepsilon\downarrow0}
 \int_{D_r^+}
 \varphi\,\chi(v)
 \frac{\varepsilon^2}{(v^2+\varepsilon^2)^{3/2}}
 |\nabla v|^2\,dx\,dy.
\end{equation}
The same conclusion holds on a punctured half disk, with the boundary integral understood as the supremum over compact subarcs avoiding the puncture.
\end{lemma}

\begin{proof}
Fix a compact subinterval $J\Subset I_r$.
Choose $\delta>0$ so that the rectangle $J\times(0,\delta)$ lies in the coordinate domain and its image lies in the chosen tubular neighborhood of $\cK$.
Let $\pi_{\cK}$ be the nearest-point projection there.
By Fubini's theorem, for almost every $x\in J$ the function $y\mapsto u(x,y)$ belongs to $W^{2,2}(0,\delta)$, and its first derivative is continuous up to $y=0$.
Its value and derivative at $0$ coincide with the Sobolev traces.
Set $a(x)=|u_y(x,0)|$.
Conformality and the locally bounded conformal factor give
\[
 a(x)=|u_x(x,0)|>0
 \quad\text{for a.e. }x\in J.
\]
For such an $x$, the normal displacement
\[
 u(x,y)-\pi_{\cK}(u(x,y))
 =y u_y(x,0)+o(y)
\]
follows from $u(x,0)\in\cK$ and the boundary condition that $u_y$ is normal to $\cK$.
Differentiating along this restriction gives
\[
 v(x,y)=a(x)y+o(y),\qquad
 v_y(x,y)\longrightarrow a(x)
 \quad\text{as }y\downarrow0.
\]
In particular, $v(x,y)>0$ for all sufficiently small positive $y$, which justifies differentiating the distance here.
Bound the nonnegative integral from below using $|v_y|^2\leq|\nabla v|^2$, and substitute $y=\varepsilon s$.
For every fixed $s>0$, the rescaled integrand converges for almost every $x\in J$ to
\[
 \varphi(x,0)\,
 \frac{a(x)^2}{(1+a(x)^2s^2)^{3/2}}.
\]
Fatou's lemma in $(x,s)\in J\times(0,\infty)$, with the rescaled integrand set equal to zero when $\varepsilon s\geq\delta$, yields
\[
\begin{split}
 \liminf_{\varepsilon\downarrow0}
 \int_{J\times(0,\delta)}
 \varphi\,\chi(v)
 \frac{\varepsilon^2|\nabla v|^2}
      {(v^2+\varepsilon^2)^{3/2}}\,dx\,dy
 \ &\geq
 \int_J\varphi(x,0)
 \left(\int_0^\infty
       \frac{a(x)^2}{(1+a(x)^2s^2)^{3/2}}\,ds\right)dx \\
  &=\int_J\varphi(x,0)a(x)\,dx.
\end{split}
\]
Exhausting $I_r$ by compact intervals proves \eqref{eq:boundary-layer-length}.
For a punctured half disk, take compact intervals on either side of the puncture and exhaust each side separately.
This argument does not assume that the full boundary length is finite.
\end{proof}

The following lemma applies to punctured limits with only local $W^{2,2}$ regularity away from the puncture.

\begin{lemma}\label{lem:general-boundary-removability}
Let $u\in W^{1,2}(D^+,\R^K)\cap\allowbreak W^{2,2}_{\loc}(\overline{D^+}\setminus\{0\},\R^K)$ be a conformal immersion away from $0$, with
\[
    u(I)\subset\cK,\qquad
    u_y\perp T_u\cK
    \quad\text{a.e. on }I\setminus\{0\}.
\]
Assume that
\[
    \mu_u(D^+)+\int_{D^+}|A_u|^2\,d\mu_u<\infty.
\]
For every $0<s<1$, the map $u$ extends to $W^{2,2}(D_s^+,\R^K)$ and continuously to $\overline{D_s^+}$.
Its conformal factor satisfies
\[
\lambda(z)=m\log|z|+\omega(z),\qquad
m\in\mathbb N\cup\{0\},\qquad
\omega\in C^0(\overline{D_s^+})\cap W^{1,2}(D_s^+).
\]
The extension preserves the trace constraint in $\cK$ and the almost-everywhere conormal condition on $I_s$.
\end{lemma}

\begin{proof}
Put $\rho(q)=\operatorname{dist}(q,\cK)$.
Choose a tubular radius $\tau>0$ on which $\rho^2/2$ is smooth, and choose a smooth function $\chi:[0,\infty)\to[0,1]$ supported in $[0,\tau)$ and equal to one near $0$.
For $\epsilon>0$, set
\[
    X_\epsilon
    =
    \chi(\rho)
    \frac{\nabla(\rho^2/2)}
         {(\rho^2+\epsilon^2)^{1/2}}.
\]
This vector field is uniformly bounded and vanishes on $\cK$.
Let $\varphi$ be a smooth cutoff on the parameter domain that equals one on $D_{1/2}^+$ and vanishes near the outer semicircle.
Apply the first variation identity on $(D^+\cup I)\setminus\{0\}$ to the field $z\mapsto\varphi(z)X_\epsilon(u(z))$ along $u$.
Remove the puncture cutoff using logarithmic cutoffs $\eta_j$ around $0$.
The cutoff error is bounded by
\[
 C\,
 \mu_u(\operatorname{supp}\nabla\eta_j)^{1/2}
 \left(\int|\nabla\eta_j|^2\,dx\,dy\right)^{1/2},
\]
which tends to zero.
The first factor tends to zero, and the logarithmic cutoffs can be chosen with Dirichlet energy tending to zero.
Passing to the limit gives
\[
 \int_{D^+}\varphi\,
       \operatorname{div}_uX_\epsilon\,d\mu_u =
 -\int_{D^+}\varphi\langle \mathbf H_u,X_\epsilon\rangle\,d\mu_u 
 -\int_{D^+}
   \left\langle
      X_\epsilon,du(\nabla_{g_u}\varphi)
   \right\rangle d\mu_u.
\]
Both terms on the right are bounded independently of $\epsilon$.
Indeed, $|\mathbf H_u|^2\leq2|A_u|^2$, and the cutoff term is controlled by Cauchy-Schwarz and the conformal invariance of the two-dimensional Dirichlet energy of $\varphi$.

The Hessian formula in the tubular neighborhood gives
\[
 \operatorname{div}_uX_\epsilon
 \geq
 \chi(\rho)
 \frac{\epsilon^2}
      {(\rho^2+\epsilon^2)^{3/2}}
 |\nabla_{g_u}(\rho\circ u)|^2
 -C_{\cK}.
\]
Take $0\leq\varphi\leq1$.
The preceding estimates bound
\[
 \int_{D^+}\varphi\,\chi(\rho\circ u)
 \frac{\epsilon^2}
      {((\rho\circ u)^2+\epsilon^2)^{3/2}}
 |\nabla_{g_u}(\rho\circ u)|^2\,d\mu_u
\]
independently of $\epsilon$.
Since $|\nabla_{g_u}(\rho\circ u)|^2\,d\mu_u=|\nabla(\rho\circ u)|^2\,dx\,dy$, Lemma~\ref{lem:boundary-layer-length} bounds the induced length of compact diameter subarcs avoiding $0$.
Exhaustion gives finite length of $I_{1/2}$ with respect to $g_u$.
The boundary trace has finite total variation on each side of $0$, so it has one-sided limits $p_+,p_-\in\cK$.

The trace of $u\in W^{1,2}(D^+)$ belongs to $W^{1/2,2}(I)$.
If $p_+\neq p_-$, then for sufficiently small $s,t>0$,
\[
    |u(t,0)-u(-s,0)|\geq\frac12|p_+-p_-|.
\]
Pairs on opposite sides of the origin would then contribute at least the following amount to its $W^{1/2,2}$ seminorm:
\[
 \frac{|p_+-p_-|^2}{4}
 \int_0^\delta\int_0^\delta\frac{ds\,dt}{(s+t)^2}
 =\infty.
\]
Thus $p_+=p_-=:p\in\cK$.

For $0<r<1/2$, let
\[
    \ell(r)=\int_0^\pi e^{\lambda(re^{i\theta})}r\,d\theta
\]
be the image length of the parameter semicircle.
Cauchy-Schwarz gives
\[
    \int_0^{1/2}\frac{\ell(r)^2}{r}\,dr
    \leq\pi\mu_u(D_{1/2}^+)<\infty.
\]
Hence there are decreasing radii $r_j\to0$ with $\ell(r_j)\to0$.
Finite boundary variation and the common one-sided limit imply that the image $u(I_{r_j})$ of the straight boundary interval converges uniformly to $p$.
The oscillation of $u$ on the semicircle is at most $\ell(r_j)$, so the entire boundary image converges uniformly to $p$.

Moreover,
\[
 \mu_u(D_{r_j}^+)\longrightarrow0,
 \qquad
 \int_{D_{r_j}^+}|\mathbf H_u|^2\,d\mu_u\longrightarrow0.
\]
\repair{Let $V_j$ be the integral varifold induced by $u|_{D_{r_j}^+\setminus\{0\}}$, and set
\[
    B_j=\overline{u(\partial D_{r_j}^+\setminus\{0\})}.
\]
The boundary limits show that $B_j$ is compact and that $B_j$ converges to $\{p\}$.
For any compactly supported smooth ambient field $X$, the same logarithmic cutoffs give a vanishing error bounded by the preceding estimate times $\|X\|_\infty$.
Finite boundary length allows us to pass to the limit in the boundary integral.
The puncture therefore contributes no additional first variation term.
Away from $B_j$, the generalized mean curvature of $V_j$ is the multiplicity-weighted average of the sheetwise vectors $\mathbf H_u$.
By Jensen's inequality, its squared $L^2$ norm is at most $\int_{D_{r_j}^+}|\mathbf H_u|^2\,d\mu_u$.

For almost every $q\in D_{r_j}^+$, put $d=\operatorname{dist}(u(q),B_j)$.
If $d>0$, apply the monotonicity identity in $B_d(u(q))$, where $V_j$ has no boundary term.
Since the density is at least one, we obtain
\[
 1\leq C\left(
       \frac{\mu_u(D_{r_j}^+)}{d^2}
       +\int_{D_{r_j}^+}|\mathbf H_u|^2\,d\mu_u
     \right).
\]
Continuity extends this distance estimate to every point of the punctured half disk.} For large $j$, the curvature term is less than $1/(2C)$, and hence
\[
 \operatorname{dist}
 \bigl(u(q),\repair{B_j}\bigr)
 \leq C\mu_u(D_{r_j}^+)^{1/2}
 \quad\text{for every }q\in D_{r_j}^+.
\]
It follows that
\[
    \sup_{D_{r_j}^+}|u-p|\longrightarrow0.
\]
Thus $u$ extends continuously to $0$, with $u(0)=p$.

We next derive decay estimates from the conformal equation.
Shrink the half disk so that its image lies in a tubular neighborhood of $\cK$.
Let $\sigma$ be the tubular reflection and set
\[
 F(z)=
 \begin{cases}
 u(z),&\operatorname{Im}z\geq0,\\
 \sigma(u(\overline z)),&\operatorname{Im}z<0.
 \end{cases}
\]
The tangential derivative of $u$ lies in $T\cK$, and its normal derivative lies in $(T\cK)^\perp$.
On these spaces, $D\sigma$ equals $+\operatorname{Id}$ and $-\operatorname{Id}$, respectively, so the first derivatives of $F$ have matching traces.
Thus $F\in W^{2,2}_{\loc}(D\setminus\{0\})$, with no interface measure in $\Delta F$.

On the upper half disk,
\[
    \Delta u=\frac12\mathbf H_u|\nabla u|^2.
\]
The chain rule on the lower half disk gives
\begin{equation}\label{eq:general-reflected-puncture-equation}
    |\Delta F|
    \leq
    C\bigl(|\nabla F|^2+b|\nabla F|\bigr),
\end{equation}
where the reflected scalar field $b$ satisfies
\[
    \int_D b^2
    \leq C\int_{D^+}|\mathbf H_u|^2\,d\mu_u<\infty.
\]
Write $z=e^{-t+i\theta}$ and $\widetilde F(t,\theta)=F(z)$.
In cylindrical coordinates, \eqref{eq:general-reflected-puncture-equation} has the same form with $F$ replaced by $\widetilde F$ and $b$ replaced by
\[
    b_{\mathrm{cyl}}(t,\theta):=e^{-t}b(e^{-t+i\theta}).
\]
Moreover,
\[
    \int_{(T,\infty)\times S^1}|b_{\mathrm{cyl}}|^2\,dt\,d\theta
    =\int_{D_{e^{-T}}}|b|^2\,dx\,dy.
\]
In the cylindrical argument below, we write $b$ for $b_{\mathrm{cyl}}$.
The Dirichlet integral is conformally invariant.
Both it and the $L^2$ integral of the cylindrical coefficient tend to zero on cylinder tails as $T\to\infty$.

The upper half of $\widetilde F$ is conformal.
On the lower half-cylinder, the matrix of $2\widetilde F^*g_{\R^K}/(|\widetilde F_t|^2+|\widetilde F_\theta|^2)$ in the coordinates $(t,\theta)$ differs from the identity by at most $C\operatorname{dist}(u,\cK)$.
Fix a sufficiently small $\epsilon>0$.
Choose $T_0$ so large that
\[
\sup_{\substack{t\geq T_0-2,\,\,0\leq\theta\leq\pi}}
\operatorname{dist}\bigl(u(e^{-t+i\theta}),\cK\bigr)\leq\epsilon,
\qquad
\int_{(T_0-2,\infty)\times S^1}
\bigl(|\nabla\widetilde F|^2+b^2\bigr)\leq\epsilon^2.
\]
The first bound follows from continuity at the puncture, and the second from the integrability established above.
We work on $(T_0,\infty)\times S^1$, so the enlarged strips below satisfy these bounds.
Every subcylinder $Q\subset(T_0,\infty)\times S^1$ then satisfies
\begin{equation}\label{eq:general-puncture-near-conformality}
\int_Q\left||\widetilde F_t|^2-|\widetilde F_\theta|^2\right|
\leq C\epsilon\int_Q|\nabla\widetilde F|^2.
\end{equation}
Cover the tail by unit strips $C_j$, and let $C_j^*$ denote their fixed enlargements by two units.
Applying Theorem~\ref{intsma0} on the original upper half annuli, together with its interior cutoff argument, gives
\[
    \operatorname{Osc}_{C_j}\widetilde F
    \leq C E_j^{1/2},
    \qquad
    E_j:=\int_{C_j^*}|\nabla\widetilde F|^2.
\]
The same estimate holds on the reflected lower half because $\sigma$ and its inverse have bounded derivatives.
If $B_j=\int_{C_j^*}b^2$, then \eqref{eq:general-reflected-puncture-equation} gives
\[
\begin{split}
 \left|
  \int_{C_j}
    \Delta\widetilde F\cdot
    (\widetilde F-\widetilde F^*)\right|\leq
 C E_j^{1/2}
 \left(E_j+B_j^{1/2}E_j^{1/2}\right)
 \leq C\epsilon E_j,
\end{split}
\]
where $\widetilde F^*(t)$ is the angular average.
The strips have finite overlap.
Summing the strip estimates gives the same bound in terms of the energy on a region enlarged by two units.

Fix $T>2$ and a center $t_c$ such that $t_c-T\geq T_0$.
Translate the cylinder coordinate by $t_c$, retain the notation $\widetilde F$, and set
\[
\Phi(s)=\int_{[-s,s]\times S^1}|\nabla\widetilde F|^2,
\qquad 0\leq s\leq T.
\]
Jensen's inequality and \eqref{eq:general-puncture-near-conformality} give
\[
\int_{[-s,s]\times S^1}
\nabla(\widetilde F-\widetilde F^*)\cdot\nabla\widetilde F
\geq\frac12(1-C\epsilon)\Phi(s).
\]
Integration by parts, the circle Poincar\'e inequality, and the strip error estimate therefore yield
\[
(1-C\epsilon)\Phi(s)\leq\Phi'(s)+C\epsilon\Phi(s+2)
\quad\text{for a.e. }0<s<T-2.
\]
Apply Lemma~\ref{lem:delayed-energy-decay} with $L=T$, $\lambda=1-C\epsilon$, and $\delta=C\epsilon$; its hypothesis $\delta e^{2\lambda}<\lambda$ holds after decreasing $\epsilon$.
The lemma gives
\[
    \Phi(s)\leq C\Phi(T)e^{-(1-C_0\epsilon)(T-s)}.
\]
Set $\delta=1-C_0\epsilon\in(0,1)$.
For $t>T_0+3$, apply the preceding estimate with center $t$, outer half-length $T=t-T_0$, and $s=1$.
It bounds the energy on $[t-1,t+1]\times S^1$ by $CE_{\mathrm{tail}}e^{-\delta(t-T_0-1)}$, where $E_{\mathrm{tail}}$ is the finite energy of the fixed tail.
Summing over consecutive unit strips and returning to $r=e^{-t}$ gives $r_0>0$ and $C<\infty$ such that
\begin{equation*}
    \int_{D_r}|\nabla F|^2\leq Cr^\delta,
    \qquad 0<r<r_0.
\end{equation*}
Set $\mathfrak A_r^+=D_r^+\setminus\overline{D_{r/2}^+}$ and fix
\[
2<q<\frac{4}{2-\delta}.
\]
Covering these half-annuli by boundary half disks and interior disks, the scaled estimate of Theorem~\ref{intsma0} and its interior proof give
\[
\|\nabla u\|_{L^q(\mathfrak A_r^+)}
\leq C r^{2/q-1+\delta/2}.
\]
The $q$-th powers are summable over dyadic radii because $2-q+q\delta/2>0$; hence
\begin{equation}\label{eq:general-boundary-puncture-Lq}
\nabla u\in L^q(D_{r_0}^+)
\quad\text{for some }q>2.
\end{equation}
We now identify the singular part of the conformal factor.
Choose a smooth orthogonal gauge $O$ that sends $T_a\cK$ to a fixed subspace $V$ and $(T_a\cK)^\perp$ to $V^\perp$.
Put
\[
 \tau_1=e^{-\lambda}u_x,\qquad
 \tau_2=e^{-\lambda}u_y,
\]
and represent the gauged Gauss map by
\[
 \widetilde G
 =
 \left[
   \frac{O(u)\tau_1+iO(u)\tau_2}{\sqrt2}
 \right]\in\mathbb{CP}^{K-1}.
\]
Its energy is bounded by
\[
    \int|\nabla\widetilde G|^2
    \leq C_{\cK}
       \left(\int|A_u|^2\,d\mu_u+E(u)\right).
\]
Reflection by $\mathfrak C_V([z])=[R_V\overline z]$, where $R_V$ is the orthogonal reflection across $V$, doubles $\widetilde G$ to a $W^{1,2}$ projective map on a full disk.
After shrinking the disk, its energy is arbitrarily small.
Extend it to the plane by radial inversion.
The exterior copy has opposite orientation, so the total pulled-back projective K\"ahler form is zero.
After further shrinking, the area Jacobian integral of the extension is less than $2\pi$, since inversion doubles the Dirichlet integral.
The projective estimate for the curvature potential \cite[Theorem~2.1]{kuwert2012} then supplies a bounded potential $v\in W^{1,2}_{\loc}(\R^2)$ with $\nabla v\in L^2(\R^2)$ and $v(z)\to0$ as $|z|\to\infty$, such that
\[
    -\Delta v=*\widehat G^*\Omega_{\mathbb{CP}}.
\]
The reflected right-hand side is even; uniqueness therefore implies that $v$ is even.
The boundary flux will be used only in the weak form of Lemma~\ref{lem:weak-boundary-connection}.

We next show that the potential is continuous.
Rescale the disk to $D=D_1$.
View $\widehat G$ in a fixed Euclidean embedding of $\mathbb{CP}^{K-1}$ and extend it across the circle by radial inversion.
For a standard mollifier $\eta_\delta$, the Poincar\'e inequality gives
\[
\sup_{x\in\overline D}
\operatorname{dist}\bigl((\eta_\delta*\widehat G)(x),\mathbb{CP}^{K-1}\bigr)^2
\leq C\sup_{x\in\overline D}
\int_{B_\delta(x)}|\nabla\widehat G|^2
\longrightarrow0.
\]
The nearest-point projection is therefore defined for all sufficiently small $\delta$, and projecting these mollifications gives smooth projective maps converging strongly in $W^{1,2}(D)$.

Let $G_j$ be such smooth approximations on $\overline D$, and extend each map to the plane by $G_j(z)=G_j(1/\overline z)$ for $|z|>1$.
Then
\[
E(G_j;\R^2)=2E(G_j;D),\qquad
\int_{\R^2}G_j^*\Omega_{\mathbb{CP}}=0.
\]
Choose the disk and the approximations so that $\int_{\R^2}JG_j\leq\gamma<2\pi$ uniformly in $j$, where $JG_j$ is the nonnegative area Jacobian.
Let $v_j$ be the normalized potential in \cite[Theorem~2.1]{kuwert2012}, and let $C_\gamma$ be twice the constant in that theorem.
Thus
\[
\|dv_j\|_{L^2(\R^2)}\leq C_\gamma E(G_j;\R^2).
\]
Let $w_j\in W^{1,2}_0(D)$ solve
\[
-\Delta w_j=*G_j^*\Omega_{\mathbb{CP}}.
\]
The difference $v_j|_D-w_j$ is harmonic.
Testing its equation by $w_j$ gives
\[
\int_D|dw_j|^2\leq\int_D|dv_j|^2
\leq C_\gamma^2 E(G_j;\R^2)^2.
\]
Choose a smooth unit Hopf lift $\eta_j$ of $G_j|_{\overline D}$.
Solving the compatible Neumann problem for its phase gives a lift $\zeta_j$ whose connection
\[
a_j=\operatorname{Im}\sum_{\ell=1}^K
\overline{\zeta_{j,\ell}}\,d\zeta_{j,\ell}
\]
is divergence free and has zero normal trace on $\partial D$.
Since $da_j=G_j^*\Omega_{\mathbb{CP}}$, Hodge uniqueness on the disk gives $a_j=-*dw_j$.
The horizontal and vertical parts of the Hopf metric therefore give
\[
\begin{split}
\int_D|d\zeta_j|^2
&=2E(G_j;D)+\|dw_j\|_{L^2(D)}^2 =E(G_j;\R^2)+\|dw_j\|_{L^2(D)}^2\\
&\leq E(G_j;\R^2)+C_\gamma^2E(G_j;\R^2)^2.
\end{split}
\]
Passing to a weakly convergent subsequence gives a unit lift $\zeta\in W^{1,2}$ of $\widehat G$, with $a=-*dw$.
Here $w\in W^{1,2}_0(D)$ solves the limiting zero-boundary problem.

Writing $\zeta_\ell=X_\ell+iY_\ell$, the connection identity gives
\[
 *\widehat G^*\Omega_{\mathbb{CP}}
 =
 2\sum_\ell J(X_\ell,Y_\ell)
\]
distributionally.
Approximate the fixed functions $X_\ell,Y_\ell$ strongly in $W^{1,2}$.
By the Wente difference estimate~\cite[Theorem~3.1.2]{helein}, the corresponding zero-boundary solutions form a Cauchy sequence both uniformly and in $W^{1,2}$.
Their limit is $w$, so $w$ is continuous on the closed disk.
Finally, $v-w$ is harmonic on the disk.
Hence $v$ is continuous there as well.

Use the connection form $\beta$ and vector field $F_\beta$ from Lemma~\ref{lem:weak-boundary-connection}, with $e_2=O(u)\tau_2$.
By \eqref{eq:general-boundary-puncture-Lq}, $\beta\in L^q$.
By Lemma~\ref{lem:weak-boundary-connection}, the function $h:=\lambda-v$ satisfies
\[
 \int\nabla h\cdot\nabla\varphi
 =\int F_\beta\cdot\nabla\varphi
\]
for every smooth test function supported away from the puncture and the artificial semicircle.
These test functions need not vanish on the diameter representing $\partial\Sigma$.
Choose a smaller smooth half-domain $\Omega^+$ containing a relative neighborhood of $0$ in the closed upper half-plane.
Let $q_0\in W^{1,2}(\Omega^+)$ be the mean-zero solution defined by
\[
 \int_{\Omega^+}\nabla q_0\cdot\nabla\varphi
 =\int_{\Omega^+}F_\beta\cdot\nabla\varphi
 \quad\text{for every }\varphi\in W^{1,2}(\Omega^+).
\]
The right-hand side vanishes on constants, so the weak Neumann problem is compatible.
The $L^q$ Neumann estimate gives $q_0\in W^{1,q}$.
Hence $q_0$ is continuous.
The equation for $q_0$ has the same right-hand side as \eqref{eq:weak-boundary-conformal-factor}.
Thus $H=h-q_0$ satisfies $\int\nabla H\cdot\nabla\varphi=0$ for all test functions supported away from the puncture and the artificial boundary, including those that meet the diameter.
Its even reflection is weakly harmonic on the punctured full disk.
Elliptic regularity then shows that it is harmonic there.

Since $v$ and $q_0$ are bounded and $u$ has finite area, $e^{2H}$ is integrable.
The subharmonic mean-value inequality gives
\[
    H(z)\leq C+\log\frac1{|z|}.
\]
After increasing $C$, the function $C+\log(1/|z|)-H$ is nonnegative and harmonic on the punctured disk.
The two-dimensional B\^ocher theorem therefore gives
\[
    H(z)=a\log|z|+h_0(z),
\]
where $h_0$ is harmonic across $0$.
Since $e^{2H}$ is integrable, we must have $a>-1$.
We have proved
\begin{equation*}    \lambda(z)=a\log|z|+\omega(z),
    \qquad
    a>-1,\qquad
    \omega:=v+q_0+h_0\in C^0\cap W^{1,2}.
\end{equation*}
It remains to show that $a$ is a nonnegative integer and that the singularity is removable in $W^{2,2}$.
For $r>0$, define the target rescaling
\[
 u_r(z)=\frac{u(rz)-p}{r^{a+1}},
 \qquad
 \lambda_r(z)=\lambda(rz)-a\log r.
\]
On every compact upper half-annulus $\mathfrak A$,
\[
 \lambda_r(z)
 =
 a\log|z|+\omega(rz)
 \longrightarrow
 a\log|z|+\omega(0)
\]
uniformly, while
\[
 \int_\mathfrak A|\nabla(\omega(rz))|^2\,dz
 =\int_{r\mathfrak A}|\nabla\omega|^2\,dz
 \longrightarrow0.
\]
The conformal Hessian decomposition gives
\[
 |D^2u_r|^2
 \leq
 C e^{2\lambda_r}|\nabla\lambda_r|^2
 +
 C e^{4\lambda_r}|A_{u_r}|_{g_{u_r}}^2.
\]
The first term is uniformly integrable on $\mathfrak A$.
For the second term,
\[
\begin{split}
 \int_\mathfrak Ae^{4\lambda_r}|A_{u_r}|_{g_{u_r}}^2\,dz =
 \int_\mathfrak Ae^{2\lambda_r}|A_{u_r}|_{g_{u_r}}^2\,d\mu_{u_r} \leq
 C\int_{r\mathfrak A}|A_u|_{g_u}^2\,d\mu_u
 \longrightarrow0.
\end{split}
\]
On every fixed closed upper half-annulus, the conformal factor expansion gives $|\nabla u(rz)|\leq Cr^a$ for all sufficiently small $r>0$.
Since $a>-1$, radial integration from $u(0)=p$ bounds $u_r$ there.
After passing to a sequence $r_j\to0$, the maps $u_{r_j}$ converge weakly in $W^{2,2}$, strongly in $W^{1,p}$ for every finite $p$, and uniformly on every such compact half-annulus, including its intervals on the real axis.
The limit is a nonconstant conformal map $u_0$ with
\[
\lambda_0(z)=a\log|z|+\omega(0).
\]
On each compact upper half-annulus $\mathfrak A$, the conformal factors $e^{\lambda_r}$ are bounded above and below.
The oriented Gauss map $G_{u_r}$ satisfies
\[
\int_\mathfrak A|\nabla G_{u_r}|^2
\leq C\int_{r\mathfrak A}|A_u|^2\,d\mu_u\longrightarrow0.
\]
Thus the Gauss map of $u_0$ is constant, so $u_0$ has a constant tangent plane.
The vanishing normal Hessian term and the conformal equation also give $\Delta u_0=0$.

The rescaled supports
\[
    \mathcal K_r=\frac{\cK-p}{r^{a+1}}
\]
converge in $C^1$ on bounded sets to $T_p\cK$.
In smooth graph coordinates, the normal component of $\cK-p$ is bounded by the square of its tangential component.
Since $a+1>0$, the rescaled supports converge as claimed.
The local uniform convergence of $u_r$ and the $C^1$ convergence of $\mathcal K_r$ imply that $u_0(x,0)\in T_p\cK$.
The boundary derivatives converge weakly in their $W^{1/2,2}$ trace spaces, while the tangent projections of $\mathcal K_r$ along $u_r$ converge uniformly to the projection onto $T_p\cK$.
Passing to the limit in the rescaled conormal condition therefore gives
\[
(u_0)_y(x,0)\perp T_p\cK
\quad\text{a.e. on the boundary intervals away from }0.
\]
Let $R_p$ be orthogonal reflection across $T_p\cK$, and define
\[
 F_0(z)=
 \begin{cases}
 u_0(z),&\operatorname{Im}z\geq0,\\
 R_pu_0(\overline z),&\operatorname{Im}z<0.
 \end{cases}
\]
The values and first derivatives have matching traces.
Thus $F_0$ is distributionally harmonic and conformal on the punctured disk.
Since $a>-1$, it has finite Dirichlet energy near $0$.
In two dimensions, a point cannot support the Laplacian of a $W^{1,2}$ harmonic map.
Thus $F_0$ extends harmonically across $0$.

The vector-valued function
\[
    W=(F_0)_x-i(F_0)_y
\]
is holomorphic and not identically zero.
Therefore
\[
    W(z)=z^mW_0(z)
\]
for a unique $m\in\mathbb N\cup\{0\}$, with $W_0(0)\neq0$.
On the upper half disk,
\[
    |W(z)|
    =
    \sqrt2e^{\omega(0)}|z|^a,
\]
whereas holomorphic factorization gives $c|z|^m\leq|W(z)|\leq C|z|^m$ for some constants $c,\ C>0$ and all sufficiently small $|z|>0$.
Hence $a=m$.

Finally, $m\geq0$ and boundedness of $\omega$ imply that $e^{2\lambda}$ is bounded near $0$.
The unscaled Hessian decomposition yields
\[
 \int|D^2u|^2
 \leq
 C\int e^{2\lambda}|\nabla\lambda|^2
 +
 C\int e^{2\lambda}|A_u|_{g_u}^2\,d\mu_u.
\]
The first integral is finite: the logarithmic term is absent when $m=0$.
For $m\geq1$, its singular part is bounded by
\[
    C\int_0^\epsilon r^{2m-2}r\,dr<\infty,
\]
and the $\omega$-part is controlled by $\|\nabla\omega\|_{L^2}$.
The second integral is finite because $e^{2\lambda}$ is bounded and the total curvature is finite.
Thus $u\in W^{2,2}$ across $0$.
Continuity of the Sobolev trace preserves the constraint $u(I)\subset\cK$.
Changing the value at one boundary point does not affect the almost-everywhere conormal identity.
\end{proof}

We use the positive-order branch set $\cB_u$ specified in Subsection~\ref{sec:fbbci}, with the orders supplied by Lemma~\ref{lem:general-boundary-removability} at the boundary and \cite[Theorem~3.1]{kuwert2012} in the interior.

For a branch point $p\in\cB_u$, let $m_u(p)\in\mathbb N$ be the local branching order in the conformal factor expansion
\begin{equation*}
    \lambda(z)=m_u(p)\log|z|+\omega(z),
\end{equation*}
where $\omega\in C^0\cap W^{1,2}$ in an interior disk chart or a boundary half disk chart.
We define the Gauss--Bonnet weighted total branching order by
\begin{equation*}    \mathfrak b(u)
    :=
    \sum_{p\in\cB_u\cap\Sigma^\circ}m_u(p)
    +\frac12\sum_{p\in\cB_u\cap\partial\Sigma}m_u(p).
\end{equation*}
In particular, $\#\cB_u\leq2\mathfrak b(u)$, where $\#$ denotes cardinality.
An interior branch point contributes $2\pi m_u(p)$ to Gauss--Bonnet, whereas a boundary branch point contributes $\pi m_u(p)$.

\Needspace{6\baselineskip}
\subsection{\texorpdfstring{Boundary length, total curvature, and branching}{Boundary length, total curvature, and branching}}
\label{subsec:boundary-curvature-control}

We use the total boundary length $L(u)$ defined in \eqref{eq:induced-boundary-geometry}.
\begin{prop}
\label{prop:free-boundary-length-estimate}
Let $\cK\subset\R^K$ be compact and smoothly embedded.
There exists a constant $C_{\cK}<\infty$ such that every
\[
    u\in
    W^{2,2}_{b,c,\mathrm{fb}}
    \bigl((\Sigma,h),\R^K,\cK\bigr)
\]
satisfies
\begin{equation}\label{eq:sec2-boundary-length-estimate}
    L(u)
    \le
    C_{\cK}
    \left(
        \cA(u)
        +
        \int_\Sigma
            |\mathbf H_u|\,\dd\mu_u
    \right).
\end{equation}
In particular,
\begin{equation}\label{eq:sec2-boundary-length-area-willmore}
    L(u)
    \le
    C_{\cK}
    \left(
        \cA(u)
        +
        \sqrt{\cA(u)W_{\R^K}(u)}
    \right).
\end{equation}
\end{prop}

\begin{proof}
Use the distance field $X_\varepsilon$ from the proof of Lemma~\ref{lem:general-boundary-removability}, with $\rho=\operatorname{dist}(\cdot,\cK)$ and the same tubular cutoff $\chi$.
It is uniformly bounded and vanishes on $\cK$.
The branch expansions make the artificial boundary terms tend to zero, so the first variation formula gives
\begin{equation}\label{eq:sec2-first-variation-distance-field}
 \int_\Sigma\operatorname{div}_uX_\varepsilon\,\dd\mu_u
 =-\int_\Sigma\langle\mathbf H_u,X_\varepsilon\rangle\,\dd\mu_u.
\end{equation}
The tubular Hessian estimate from that proof is
\begin{equation}\label{eq:sec2-divergence-distance-field}
 \operatorname{div}_uX_\varepsilon
 \geq\chi(\rho)\frac{\varepsilon^2}{(\rho^2+\varepsilon^2)^{3/2}}
       |\nabla_{g_u}(\rho\circ u)|^2-C_{\cK}.
\end{equation}
Integrating \eqref{eq:sec2-divergence-distance-field} and using \eqref{eq:sec2-first-variation-distance-field} bounds the nonnegative boundary-layer integral by $C_{\cK}(\cA(u)+\int_\Sigma|\mathbf H_u|\,\dd\mu_u)$.
Apply Lemma~\ref{lem:boundary-layer-length} in boundary conformal charts with a nonnegative partition of unity.
This bounds the induced arclength of compact subsets of $\partial\Sigma\setminus\cB_u$ by the lower limit of the corresponding boundary-layer integrals.
Conformal invariance identifies the coordinate integrals with the intrinsic integrals, and the partition cutoffs sum to at most $1$.
Exhausting $\partial\Sigma\setminus\cB_u$ therefore gives the bound for $L(u)$ in \eqref{eq:sec2-boundary-length-estimate}.
Finally, Cauchy--Schwarz gives $\int_\Sigma|\mathbf H_u|\,\dd\mu_u\leq2\sqrt{\cA(u)W_{\R^K}(u)}$, which proves \eqref{eq:sec2-boundary-length-area-willmore}.
\end{proof}

For the induced metric $g_u$, let $\nu_u^{\mathrm{in}}$ be the inward unit conormal along $\partial\Sigma\setminus\cB_u$, and let $\tau_{g_u}$ be a unit tangent vector there.
For a smooth immersion, we use the convention
\[
    k_{g_u}=\left\langle\nabla^{g_u}_{\tau_{g_u}} \tau_{g_u},\nu_u^{\mathrm{in}}\right\rangle_{g_u}.
\]
For a $W^{2,2}$ immersion, we define the geodesic curvature measure through the boundary connection identity.
The trace formula \eqref{eq:sec2-boundary-geodesic-curvature-support} below defines its density and agrees with this convention in the smooth case.

\begin{lemma}
\label{lem:free-boundary-gauss-bonnet}
Let
\[
    u\in
    W^{2,2}_{b,c,\mathrm{fb}}
    \bigl((\Sigma,h),\R^K,\cK\bigr).
\]
Then
\begin{equation}\label{eq:sec2-branched-gauss-bonnet}
    \int_{\Sigma\setminus\cB_u}
        K_{g_u}\,\dd\mu_u
    +
    \int_{\partial\Sigma\setminus\cB_u}
        k_{g_u}\,\dd s_{g_u}
    =
    2\pi\bigl(\chi(\Sigma)+\mathfrak b(u)\bigr).
\end{equation}
If $A^{\cK}_{\R^K}$ denotes the second fundamental form of $\cK\subset\R^K$, then along $\partial\Sigma\setminus\cB_u$
\begin{equation}\label{eq:sec2-boundary-geodesic-curvature-support}
    k_{g_u}
    =
    \left\langle
        A^{\cK}_{\R^K}
        \bigl(du(\tau_{g_u}),du(\tau_{g_u})\bigr),
        du(\nu_u^{\mathrm{in}})
    \right\rangle_{\R^K}.
\end{equation}
Consequently,
\begin{equation}\label{eq:sec2-boundary-geodesic-curvature-bound}
    |k_{g_u}|
    \le
    \|A^{\cK}_{\R^K}\|_{L^\infty(\cK)},
\end{equation}
and hence
\begin{equation*}    \left|
        \int_{\partial\Sigma\setminus\cB_u}
            k_{g_u}\,\dd s_{g_u}
    \right|
    \le
    C_{\cK}L(u).
\end{equation*}
\end{lemma}

\begin{proof}
On a boundary conformal chart disjoint from $\cB_u$, write $g_u=e^{2\ell}|dz|^2$ and set $\tau_1=e^{-\ell}u_x$, $\tau_2=e^{-\ell}u_y$.
Lemma~\ref{lem:weak-boundary-connection}, specifically \eqref{eq:weak-boundary-curvature}-\eqref{eq:weak-boundary-support-curvature}, identifies the weak geodesic curvature measure of $\partial\Sigma$ with respect to $g_u$ on this boundary chart as
\[
 k_{g_u}\,ds_{g_u}
 =-\beta_1\,dx
 =e^\ell\left\langle
 A^{\cK}_{\R^K}(\tau_1,\tau_1),\tau_2
 \right\rangle dx.
\]
Choose $\tau_{g_u}=e^{-\ell}\partial_x$ and $\nu_u^{\mathrm{in}}=e^{-\ell}\partial_y$ on this chart.
Then $du(\tau_{g_u})=\tau_1$, $du(\nu_u^{\mathrm{in}})=\tau_2$, and $ds_{g_u}=e^\ell dx$, proving \eqref{eq:sec2-boundary-geodesic-curvature-support}.
Both trace vectors have unit length, so \eqref{eq:sec2-boundary-geodesic-curvature-bound} follows.
Proposition~\ref{prop:free-boundary-length-estimate} therefore shows that the boundary curvature measure is finite.

Fix a smooth metric $h$ on $\Sigma$ in its prescribed conformal class and write $g_u=e^{2\lambda}h$ on $\overline\Sigma\setminus\cB_u$.
Denote the Gauss curvature of $h$ by $K_h$ and the geodesic curvature of $\partial\Sigma$ with respect to $h$, computed with the inward unit conormal, by $k_h$.
The interior curvature equation and the preceding weak boundary identity give
\begin{equation}\label{eq:weak-gauss-bonnet-cutoff}
\begin{split}
 \int_\Sigma K_{g_u}\eta\,d\mu_u
 +\int_{\partial\Sigma}k_{g_u}\eta\,ds_{g_u}
 &=\int_\Sigma K_h\eta\,\dd\vol_h
   +\int_{\partial\Sigma}k_h\eta\,ds_h\\
&\quad+\int_\Sigma\langle\nabla_h\lambda,\nabla_h\eta\rangle_h\,\dd\vol_h
\end{split}
\end{equation}
for every smooth $\eta$ which vanishes near the branch set.
In a boundary chart, write $h=e^{2\phi}|dz|^2$ and $\ell=\lambda+\phi$ in the preceding identity.
Subtracting the smooth identity for $\phi$ gives the claim.
Interior charts follow from $-\Delta_h\lambda=K_{g_u}e^{2\lambda}-K_h$.
A partition of unity gives \eqref{eq:weak-gauss-bonnet-cutoff} on the whole surface.

Choose disjoint branch charts and nonnegative logarithmic cutoffs $\eta_\varepsilon$ which vanish for $|z|\leq\varepsilon^2$ and equal $1$ for $|z|\geq\varepsilon$ in each chart, and equal $1$ outside these charts.
We can choose them to be smooth and to satisfy
\[
 \int|\nabla\eta_\varepsilon|^2\leq\frac{C}{|\log\varepsilon|}.
\]
At a branch point of order $m$, the conformal factor expansion gives $\lambda=m\log|z|+\omega$ relative to $h$, after absorbing a smooth coordinate factor into $\omega\in C^0\cap W^{1,2}$.
Conformal invariance and radial integration imply
\[
 \int\langle\nabla\log|z|,\nabla\eta_\varepsilon\rangle\,dx\,dy
 =\begin{cases}
 2\pi,&\text{at an interior branch point},\\
 \pi,&\text{at a boundary branch point}.
 \end{cases}
\]
The remainder contributes zero in the limit, since
\[
 \left|\int\nabla\omega\cdot\nabla\eta_\varepsilon\right|
 \leq\|\nabla\omega\|_{L^2(\operatorname{supp}\nabla\eta_\varepsilon)}
       \|\nabla\eta_\varepsilon\|_{L^2}
 \longrightarrow0.
\]
Thus the last term of \eqref{eq:weak-gauss-bonnet-cutoff} tends to
\[
 2\pi\sum_{p\in\cB_u\cap\Sigma^\circ}m_u(p)
 +\pi\sum_{p\in\cB_u\cap\partial\Sigma}m_u(p).
\]
Finite total curvature makes the Gauss curvature measure integrable.
The boundary measure is finite by the preceding argument.
Dominated convergence therefore applies to the other terms.
The classical Gauss--Bonnet formula for the smooth metric $h$ gives
\[
 \int_\Sigma K_h\,\dd\vol_h+\int_{\partial\Sigma}k_h\,ds_h
 =2\pi\chi(\Sigma).
\]
Taking $\varepsilon\downarrow0$ in \eqref{eq:weak-gauss-bonnet-cutoff} proves \eqref{eq:sec2-branched-gauss-bonnet}.
The final estimate for the boundary integral follows by integrating \eqref{eq:sec2-boundary-geodesic-curvature-bound}.
\end{proof}

\begin{prop}
\label{prop:automatic-curvature-branching-bound}
There exists a constant $C=C(\Sigma,\cK)<\infty$ such that every
\[
    u\in
    W^{2,2}_{b,c,\mathrm{fb}}
    \bigl((\Sigma,h),\R^K,\cK\bigr)
\]
satisfies
\begin{equation}\label{eq:sec2-automatic-curvature-branching-bound}
    T_{\R^K}(u)+\mathfrak b(u)
    \le
    C\bigl(1+\cA(u)+W_{\R^K}(u)\bigr).
\end{equation}
In particular, if
\[
    \sup_n
    \bigl(E_{h_n}(u_n)+W_{\R^K}(u_n)\bigr)<\infty,
\]
then
\begin{equation}\label{eq:sec2-sequential-curvature-branching-bound}
    \sup_n
    \bigl(
        T_{\R^K}(u_n)+\mathfrak b(u_n)
    \bigr)<\infty.
\end{equation}
\end{prop}

\begin{proof}
The Euclidean Gauss equation and \eqref{eq:sec2-branched-gauss-bonnet} give
\[
 T_{\R^K}(u)+4\pi\mathfrak b(u)
 =4W_{\R^K}(u)-4\pi\chi(\Sigma)
  +2\int_{\partial\Sigma\setminus\cB_u}k_{g_u}\,\dd s_{g_u}.
\]
Proposition~\ref{prop:free-boundary-length-estimate}, Lemma~\ref{lem:free-boundary-gauss-bonnet}, and $2\sqrt{ab}\leq a+b$ bound the right-hand side by $C(\Sigma,\cK)(1+\cA(u)+W_{\R^K}(u))$.
Since $T_{\R^K}(u)$ and $\mathfrak b(u)$ are nonnegative, this proves \eqref{eq:sec2-automatic-curvature-branching-bound}.
The identity $E_{h_n}(u_n)=\cA(u_n)$ in \eqref{eq:sec2-energy-area} then gives \eqref{eq:sec2-sequential-curvature-branching-bound}.
\end{proof}

\begin{rmk}\label{rem:no-extra-total-curvature-hypothesis}
Proposition~\ref{prop:automatic-curvature-branching-bound} is the bordered analogue of the Gauss--Bonnet estimate used in Chen--Li~\cite{chen2014}.
Thus Theorem~\ref{thm:fb-chen-li} needs no separate bound on $T_{\R^K}(u_n)$ or $\mathfrak b(u_n)$.
For maps into a compact target $N$, the same conclusion follows from \eqref{eq:sec2-total-curvature-comparison} and \eqref{eq:sec2-willmore-comparison} after fixing the isometric embedding.
\end{rmk}

\Needspace{6\baselineskip}
\subsection{\texorpdfstring{Energy gap and neck decay}{Energy gap and neck decay}}
\label{subsec:general-neck-decay}

\begin{theorem}\label{gapth}
Fix an isometric embedding $N\hookrightarrow\R^K$, a compact subset $Q\subset N$, and a compact smooth supporting submanifold $\cK\subset N$.
There exists $\epsilon_*=\epsilon_*(K,Q,\cK)>0$ such that every nonconstant
\[
u\in
W^{2,2}_{b,c,\mathrm{fb}}
\bigl((\Sigma,h),N,\cK\bigr),
\qquad
u(\Sigma)\subset Q,
\]
on a compact connected surface with nonempty boundary satisfies
\[
\cA(u)+W_N(u)\geq\epsilon_*.
\]
For closed surfaces, there is a constant $c_K>0$ such that every nonconstant branched conformal immersion $u:\Sigma\to\R^K$ with finite area and finite total curvature satisfies
\[
\cA(u)+W_{\R^K}(u)\geq c_K.
\]
\end{theorem}

\begin{proof}
Compose $u$ with the fixed isometric embedding, retaining the notation $u$.
Lemma~\ref{lem:general-boundary-removability} and \cite[Theorem~3.1]{kuwert2012} give the conformal factor expansions with integer branch orders.
In particular, $|du|\leq Cr^m$ near a branch point of order $m$, so artificial boundary fluxes are $O(r^{m+1})$ and vanish as the radii of the excised neighborhoods tend to zero.
The boundary length is finite by Proposition~\ref{prop:free-boundary-length-estimate}.

Let $V_u$ be the integer rectifiable two-varifold parametrized by $u$, retaining all sheet multiplicities.
Its mass is $\cA(u)$ and its density is a positive integer almost everywhere.
Passing to the limit in the weak integration by parts formula outside small coordinate neighborhoods of the branch points gives
\[
 \delta V_u(X)
 =-\int_\Sigma\langle\mathbf H_u,X(u)\rangle\,d\mu_u
  -\int_{\partial\Sigma}\langle du(\nu_u^{\mathrm{in}}),X(u)\rangle\,ds_{g_u}
\]
for every compactly supported smooth ambient field $X$.
The isoperimetric inequality \cite[Section~VII, Proposition~11]{Allard1987}, the boundary length estimate, and Cauchy--Schwarz imply
\[
\begin{split}
 \sqrt{\cA(u)}
 &\leq C_K\|\delta V_u\|
 \leq C_K\left(\int_\Sigma|\mathbf H_u|\,d\mu_u+L(u)\right)\\
 &\leq C_{K,\cK}\left(\cA(u)+\sqrt{\cA(u)W_{\R^K}(u)}\right).
\end{split}
\]
A nonconstant branched conformal immersion has positive area.
Dividing by $\sqrt{\cA(u)}$ gives $1\leq C_{K,\cK}\left(\sqrt{\cA(u)}+\sqrt{W_{\R^K}(u)}\right)$, and hence $\cA(u)+W_{\R^K}(u)\geq c_{K,\cK}>0$.
For a closed surface, the same argument applies without the boundary term.
On the compact target set $Q$, the comparison $W_{\R^K}(u)\leq W_N(u)+C_Q\cA(u)$ gives the intrinsic lower bound after reducing the constant.
\end{proof}

\begin{theorem}\label{deca}
Fix a positive integer $K$, a compact set $Q\subset\R^K$, and a smooth compact embedded submanifold $\cK\subset\R^K$.
There are constants $\epsilon_*>0$, $C<\infty$, and $C_0<\infty$, depending only on $K,Q,\cK$, with the following property.

Let $T>0$, and let $u\in W^{1,2}$ have image in $Q$ on one of the flat domains
\par\noindent
\(
\begin{array}{@{}ll@{}}
\quad\textup{(a)} & C_T=(-T,T)\times S^1,\\
\quad\textup{(b)} & L_T=(-T,T)\times[0,\pi],\\
\quad\textup{(c)} & P_T=[0,T)\times S^1.
\end{array}
\)
\par
The free boundary consists of the longitudinal sides $(-T,T)\times\{0,\pi\}$ in case \textup{(b)} and the circle $\{0\}\times S^1$ in case \textup{(c)}; case \textup{(a)} has no such boundary.
Assume one of the following two regularity conditions.
In the case of finitely many branch points, $u\in W^{2,2}_{\loc}$ on the chosen domain, including its free boundary and branch points.
It is a conformal immersion away from the finite branch set and has finite induced area and finite Euclidean total curvature.
In the unbranched case, $u$ is a smooth unbranched conformal immersion in the interior and up to the free boundary, with finite total Dirichlet and Euclidean Willmore energies.
No total curvature hypothesis is imposed in this case.
In both cases, the branch set excludes the artificial open end faces.
In case \textup{(b)}, assume
\[
u(t,0),u(t,\pi)\in\cK,
\qquad
u_\theta(t,0),u_\theta(t,\pi)
\perp T_u\cK
\]
almost everywhere on the longitudinal boundary away from the branch points of $u$.
In case \textup{(c)}, assume
\[
u(0,\theta)\in\cK,
\qquad
u_t(0,\theta)\perp T_u\cK
\]
almost everywhere on $\{0\}\times S^1$ away from the branch points of $u$.
We use the Dirichlet and Willmore energies $E(u;\mathfrak A)$ and $W_{\R^K}(u;\mathfrak A)$ from Subsection~\ref{sec:fbbci}, computed with the flat metric on the cylinder, for each measurable subset $\mathfrak A$ of the domain.
For $I=S^1$ in cases \textup{(a)} and \textup{(c)}, and $I=[0,\pi]$ in case \textup{(b)}, write $S_a=(a,a+1)\times I$.
In case \textup{(c)}, the strip adjacent to the boundary inherited from $\partial\Sigma$ is understood as $S_0=[0,1)\times S^1$.
Suppose $0<\epsilon\leq\epsilon_*$ and
\[
E(u;S_a)+W_{\R^K}(u;S_a)\leq\epsilon^2
\]
for every such unit strip contained in the domain.
Put
\[
\gamma=1-C_0\epsilon.
\]
After decreasing $\epsilon_*$, we have $\gamma>0$.
For every $0\leq R\leq T$,
\[
E\bigl(u;(-T+R,T-R)\times I\bigr)
\leq
C E\bigl(u;(-T,T)\times I\bigr)e^{-\gamma R}
\]
in cases \textup{(a)} and \textup{(b)}, where $I=S^1$ or $[0,\pi]$, respectively.
In case \textup{(c)},
\[
E\bigl(u;[0,T-R]\times S^1\bigr)
\leq
C E\bigl(u;[0,T]\times S^1\bigr)e^{-\gamma R}.
\]
If $4\leq R\leq T$, then
\[
\operatorname{Osc}_{(-T+R,T-R)\times I}u
\leq
C E\bigl(u;(-T,T)\times I\bigr)^{1/2}
e^{-\gamma R/2}
\]
in cases \textup{(a)} and \textup{(b)}, while
\[
\operatorname{Osc}_{[0,T-R]\times S^1}u
\leq
C E\bigl(u;[0,T]\times S^1\bigr)^{1/2}
e^{-\gamma R/2}
\]
in case \textup{(c)}.
The oscillation estimates apply only to regions separated from uncontrolled artificial ends.

For an HSL immersion into a fixed compact subset of an isometrically embedded target $M$, the comparison
\[
W_{\R^K}(u;\mathfrak{A})
\leq W_M(u;\mathfrak{A})+C_QE(u;\mathfrak{A})
\]
shows that the same conclusions follow from small intrinsic Dirichlet and Willmore energy.
We use this result for the HSL necks below.

\end{theorem}

\begin{proof}

In the smooth unbranched case, we work first on compact truncations.
The proof uses the conformal equation, the small Dirichlet and Euclidean Willmore energies on each strip, the boundary conditions, and the bounded derivatives of the tubular reflection.
It needs no total curvature bound or branch point removal argument in this case.

When there are finitely many branch points, Lemma~\ref{lem:general-boundary-removability} gives the conformal factor expansion at each branch point on a free boundary edge,
\[
\lambda=m\log|z|+\omega,
\qquad
m\in\mathbb N\cup\{0\},
\qquad
\omega\in C^0\cap W^{1,2}.
\]
At an interior branch point, the same conclusion follows from \cite[Theorem~3.1]{kuwert2012}.
Since $|du|=\sqrt2e^\lambda\leq Cr^m$, the flux through an excised circle or semicircle is $O(r^{m+1})$.
Thus no point-supported term appears in the weak integrations by parts below.
Assign arbitrary values to $\mathbf H_u$ on the finite branch set.
The conformal equation then holds almost everywhere.

Near an artificial open end, truncate the domain at a distance $\delta>0$ from the end.
The assumed local $W^{2,2}$ regularity holds on each fixed truncation, so reflection, chain rules, traces, and integration by parts apply.
The estimates are uniform in $\delta$ on regions separated from the artificial end.
Let $\delta\downarrow0$.
Monotone convergence gives the energy estimates.
Exhaustion gives the oscillation estimates.

If $T\leq6$, monotonicity gives the energy estimates after increasing $C$.
For the oscillation estimate, $R\geq4$.
When $4\leq T\leq6$, cover the region away from the ends by finitely many charts and apply the local estimates below.
It remains to consider $T>6$.

Fix $4/3<p<2$.
Theorem~\ref{intsma0} gives the boundary $W^{2,p}$ estimate.
On a unit interior disk, apply the cutoff calculation in the proof of Theorem~\ref{intsma0} directly to $u$.
With $w=\eta(u-\overline u)$ and $p^*=2p/(2-p)$, the same Sobolev and Calder\'on--Zygmund estimates~\cite{AdamsFournier2003,GilbargTrudinger2001} give
\[
 \|w\|_{W^{2,p}}
 \leq C\|\nabla u\|_{L^2}
      +C W_{\R^K}(u)^{1/2}\|w\|_{W^{2,p}}.
\]
Small Willmore energy absorbs the last term, and rescaling gives the weighted interior estimate of Theorem~\ref{intsma0}.

Cover a closed unit strip $\mathfrak A$ by finitely many such disks and half disks.
Let $\mathfrak A^*$ be a fixed enlargement with a margin of $2$ units in the $t$-direction.
This cover may meet a genuine free boundary.
Each enlarged chart is contained in a fixed number of the unit strips from the hypothesis.
After decreasing $\epsilon_*$, the local estimates and Morrey's inequality give
\begin{equation}\label{decay-osc-corrected}
\operatorname{Osc}_\mathfrak A u
\leq
C E(u;\mathfrak A^*)^{1/2},
\qquad
\|\nabla u\|_{W^{1,p}(\mathfrak A)}
\leq
C\|\nabla u\|_{L^2(\mathfrak A^*)}.
\end{equation}
These estimates also hold after tubular reflection, since both the reflection and its inverse have bounded first and second derivatives.
Since $p>4/3$, the trace theorem also gives, in case \textup{(c)},
\[
\|u_t\|_{L^2(\{0\}\times S^1)}
\leq
C E\bigl(u;[0,2]\times S^1\bigr)^{1/2}.
\]
We use these estimates only on strips separated from uncontrolled artificial ends.
Each such strip has oscillation at most $C\epsilon$.

Consider case \textup{(b)} and put $T'=T-2$.
Every point of $(-T',T')\times[0,\pi]$ stays a fixed distance from the ends in the $t$-direction.
By \eqref{decay-osc-corrected}, the image of each local strip lies in a fixed tubular neighborhood of $\cK$.
Let $\pi_{\cK}$ be the nearest-point projection there and set
\[
\sigma(y)=2\pi_{\cK}(y)-y.
\]
For a sufficiently small tubular radius, $\sigma$ is a smooth involutive diffeomorphism with uniformly bounded first and second derivatives.
On $\cK$, its differential is the identity on $T\cK$ and minus the identity on the normal bundle.

Define a periodic map on the shortened cylinder by
\[
U(t,\theta)=
\begin{cases}
u(t,\theta),&0\leq\theta\leq\pi,\\
\sigma(u(t,-\theta)),&-\pi\leq\theta\leq0.
\end{cases}
\]
By the free boundary trace conditions, $U$ and its first derivatives agree at $\theta=0$ and across the identified seam $\theta=-\pi=\pi$.
Hence $U\in W^{2,2}_{\loc}$ on the periodic cylinder.
The reflected map $U$ need not be conformal.
Since $D\sigma$ is smooth and invertible,
\begin{equation*}C^{-1}|\nabla u|^2
\leq|\nabla U|^2
\leq C|\nabla u|^2
\end{equation*}
on corresponding reflected regions.

The pullback $\sigma^*g_{\R^K}$ of the Euclidean metric by the reflection $\sigma$ satisfies
\[
\left|(D\sigma_y)^TD\sigma_y-I\right|
\leq C\operatorname{dist}(y,\cK).
\]
On every strip under consideration, the oscillation estimate gives $\operatorname{dist}(u,\cK)\leq C\epsilon$.
Since $u$ is conformal on the original half-cylinder,
\begin{equation}\label{pohozaev-before-reflection}
\left|
|U_t|^2-|U_\theta|^2
\right|
\leq
C\epsilon|\nabla U|^2
\end{equation}
on the reflected half, while the left side vanishes on the original half.
Thus, for $0<s<T'$,
\[
\int_{[-s,s]\times S^1}
\left||U_t|^2-|U_\theta|^2\right|
\leq
C\epsilon
\int_{[-s,s]\times S^1}|\nabla U|^2.
\]
On every original conformal piece,
\[
\Delta u=\frac12\mathbf H_u|\nabla u|^2.
\]
On the reflected half, put $\rho(t,\theta)=(t,-\theta)$, $v=u\circ\rho$, and $\mathbf H_v=\mathbf H_u\circ\rho$.
The chain rule gives
\[
\Delta U
=
D\sigma(v)\Delta v
+
D^2\sigma(v)[\nabla v,\nabla v],
\]
and hence
\[
|\Delta U|\leq C(1+|\mathbf H_v|)|\nabla v|^2.
\]
Let $\mathfrak A$ be a unit strip in the doubled cylinder and put $\mathfrak A_+=\mathfrak A\cap\{0\leq\theta\leq\pi\}$.
Let $\mathfrak A_+^*$ be its fixed enlargement in the original half-cylinder, chosen away from the artificial end faces.
Also set
\[
U^*(t)=\frac1{2\pi}\int_{S^1}U(t,\theta)\,d\theta.
\]
The equations on the two halves, the oscillation estimate, and Cauchy--Schwarz give
\[
\begin{split}
\left|\int_\mathfrak A\Delta U\cdot(U-U^*)\right|
&\leq C\operatorname{Osc}_\mathfrak A U
\left(E(u;\mathfrak A_+)+W_{\R^K}(u;\mathfrak A_+)^{1/2}E(u;\mathfrak A_+)^{1/2}\right)\\
&\leq C\epsilon E(u;\mathfrak A_+^*).
\end{split}
\]
The same estimate holds without reflection in case \textup{(a)}.
The strips have bounded overlap.
Summing their estimates gives, for almost every $0<s<T'-2$ in case \textup{(b)},
\begin{equation}\label{interior-error-summed}
\left|
\int_{[-s,s]\times S^1}
\Delta U\cdot(U-U^*)
\right|
\leq
C\epsilon
\int_{[-s-2,s+2]\times S^1}|\nabla U|^2.
\end{equation}
Put
\[
\Phi(s)=\int_{[-s,s]\times S^1}|\nabla U|^2.
\]
Jensen's inequality and \eqref{pohozaev-before-reflection} imply
\begin{equation}\label{lower-bound-no-conformality}
\int_{[-s,s]\times S^1}
\nabla(U-U^*)\cdot\nabla U
\geq
\int_{[-s,s]\times S^1}|U_\theta|^2
\geq
\frac12(1-C\epsilon)\Phi(s).
\end{equation}
Integration by parts gives a boundary term at $t=\pm s$.
Since $U-U^*$ has zero angular average, the circle Poincar\'e inequality yields
\begin{align}
\left|
\int_{\partial([-s,s]\times S^1)}
\partial_\nu U\cdot(U-U^*)
\right|
&\leq
\left(\int_{\partial([-s,s]\times S^1)}|U_t|^2\right)^{1/2}
\left(\int_{\partial([-s,s]\times S^1)}|U_\theta|^2\right)^{1/2}
\notag\\
&\leq\frac12\Phi'(s).\label{boundary-term-cylinder}
\end{align}
Combining \eqref{interior-error-summed}-\eqref{boundary-term-cylinder} gives
\begin{equation}\label{ode-decay-cylinder}
(1-a\epsilon)\Phi(s)
\leq
\Phi'(s)+b\epsilon\Phi(s+2)
\end{equation}
for fixed constants $a,b$ and almost every $0<s<T'-2$.

To apply Lemma~\ref{lem:delayed-energy-decay} to $\Phi$ in \eqref{ode-decay-cylinder}, choose $\lambda=1-a\epsilon$ and $\delta=b\epsilon$.
For sufficiently small $\epsilon$, these parameters satisfy
\[
\lambda-\delta e^{2\lambda}
\geq1-(a+be^2)\epsilon>0.
\]
Lemma~\ref{lem:delayed-energy-decay} therefore gives exponential decay of $\Phi$ with rate $1-C_0\epsilon$.

In case \textup{(b)}, the preceding calculation with $L=T'$ bounds the Dirichlet energy on the reflected cylinder.
Comparing this energy with that of $u$ on the original half-cylinder gives
\[
E\bigl(u;(-s,s)\times[0,\pi]\bigr)
\leq
C E\bigl(u;(-T',T')\times[0,\pi]\bigr)
e^{-\gamma(T'-s)}.
\]
Since $T'=T-2$, the factor $e^{2\gamma}$ can be absorbed into $C$.
Energy monotonicity covers the remaining range $T'<s<T$.
Setting $s=T-R$ then proves case \textup{(b)}.
In case \textup{(a)}, take $U=u$ and $L=T$.
We only need the error estimate for $0<s<T-2$. In this range, the enlarged cylinder $[-s-2,s+2]\times S^1$ lies inside $C_T$. The same calculation and energy monotonicity near the ends prove case \textup{(a)}.

In case \textup{(c)}, let $u^*(t)$ be the circle average and set
\[
F(s)=\int_{[0,s]\times S^1}|\nabla u|^2.
\]
Conformality and Jensen's inequality give
\[
\int_{[0,s]\times S^1}
\nabla(u-u^*)\cdot\nabla u
\geq\frac12F(s).
\]
Integration by parts, the circle Poincar\'e inequality on $\{s\}\times S^1$, and the equation error estimate away from the artificial end give
\[
F(s)
\leq
F'(s)+C\epsilon F(s+2)+2|B_0|
\]
for almost every $0<s<T-2$, where
\[
B_0
=
\int_{S^1}
\left\langle
u_t(0,\theta),u(0,\theta)-u^*(0)
\right\rangle\,d\theta.
\]
For sufficiently close $x,y\in\cK$ and every $n\perp T_x\cK$, the uniform Taylor estimate for $\cK$ gives
\[
|\langle n,x-y\rangle|
\leq C|n||x-y|^2.
\]
Using
\[
u(0,\theta)-u^*(0)
=
\frac1{2\pi}\int_{S^1}
\bigl(u(0,\theta)-u(0,\varphi)\bigr)\,d\varphi
\]
and $u_t(0,\theta)\perp T_{u(0,\theta)}\cK$, we obtain
\[
|B_0|
\leq
C\bigl(\operatorname{Osc}_{S^1}u(0,\cdot)\bigr)^2
\|u_t(0,\cdot)\|_{L^1(S^1)}.
\]
The oscillation estimate \eqref{decay-osc-corrected} and the $L^2$ trace bound for $u_t(0,\cdot)$ give
\[
|B_0|
\leq
C E\bigl(u;[0,2]\times S^1\bigr)^{3/2}
\leq
C\epsilon E\bigl(u;[0,2]\times S^1\bigr)
\leq
C\epsilon F(s+2).
\]
Thus
\[
(1-a\epsilon)F(s)
\leq
F'(s)+b\epsilon F(s+2).
\]
Applying the same scalar estimate gives
\[
F(s)\leq C F(T)e^{-\gamma(T-s)}
\]
for $0\leq s\leq T$, using monotonicity near the end.
Taking $s=T-R$ proves case \textup{(c)}.

Assume $R\geq4$, and let $E_{\mathrm{tot}}$ denote the total Dirichlet energy $E(u;C_T)$, $E(u;L_T)$, or $E(u;P_T)$ in cases \textup{(a)}, \textup{(b)}, and \textup{(c)} in Theorem \ref{deca}, respectively.
Every unit strip meeting the truncated region stays a fixed positive distance from each uncontrolled artificial end.
In cases \textup{(a)} and \textup{(b)}, suppose such a strip has distance $d$ from the nearer artificial end.
Applying the energy estimate to its fixed enlargement gives
\[
E(u;\mathfrak A^*)
\leq
C E_{\mathrm{tot}}e^{-\gamma(d-2)}.
\]
In case \textup{(c)}, a strip ending near coordinate $t$ satisfies
\[
E(u;\mathfrak A^*)
\leq
C E_{\mathrm{tot}}e^{-\gamma(T-t-2)}.
\]
By \eqref{decay-osc-corrected}, the oscillation of $u$ on each strip is bounded by the square root of the corresponding energy bound.
Join two points in the truncated region by a chain of consecutive unit strips.
Summing from the nearest artificial end gives
\[
\sum_{j\geq R-2}
C E_{\mathrm{tot}}^{1/2}e^{-\gamma j/2}
\leq
C E_{\mathrm{tot}}^{1/2}e^{-\gamma R/2}.
\]
Choose $\epsilon_*$ small enough that $\gamma$ has a fixed positive lower bound.
The constant in the geometric series is then uniform.
This proves the oscillation estimates in cases \textup{(a)}-\textup{(c)}.

\end{proof}

\Needspace{6\baselineskip}
\subsection{\texorpdfstring{Local compactness and conformal factors}{Local compactness and conformal factors}}
\label{subsec:conformal-factor-compactness}

Theorem~\ref{sub} recalls H\'elein's interior compactness theorem in the sharp form of \cite[Theorem~5.1]{kuwert2012} and \cite{chen2014}.
Lemma~\ref{lem:boundary-conformal-factor-compactness} proves the boundary alternative between collapse and convergence to an unbranched conformal immersion.

\begin{theorem}\label{sub}

Let $u_k\in W^{2,2}( D_1,\R^K)$ be conformal immersions with induced metrics
\[
(g_{u_k})_{ij}=e^{2\lambda_k}\delta_{ij}.
\]
Assume
\[
\int_{ D_1}|A_{u_k}|^2\,d\mu_{u_k}
\leq\gamma<\gamma_K,
\qquad
\gamma_K=
\begin{cases}
8\pi,&K=3,\\
4\pi,&K\geq4,
\end{cases}
\]
and $\mu_{u_k}( D_1)\leq C$.
Suppose also that $|u_k(0)|\leq C_0$; equivalently, one may translate the maps so that $u_k(0)=0$.
Then $u_k$ is bounded in $W^{2,2}_{\mathrm{loc}}( D_1,\R^K)$.
After passing to a subsequence, exactly one of the following alternatives holds
\begin{enumerate}
\item The functions $\lambda_k$ are locally bounded above and below, and $u_k$ converges weakly in $W^{2,2}_{\mathrm{loc}}( D_1,\R^K)$ to a conformal immersion.
\item The functions $\lambda_k$ tend to $-\infty$ locally uniformly, and $u_k$ converges weakly in $W^{2,2}_{\mathrm{loc}}( D_1,\R^K)$ and locally uniformly to a constant.
\end{enumerate}

For $K=2$, apply the $K=3$ case by composing with the linear isometric embedding $\R^2\hookrightarrow\R^3$.
Composing the maps with this embedding preserves their induced metrics, Dirichlet and Willmore energies, and total curvature.

\end{theorem}

\begin{lemma}\label{lem:boundary-conformal-factor-compactness}
Let $\cK\subset\R^K$ be a smooth compact embedded submanifold and let $Q\subset\R^K$ be compact.
There exists $\epsilon_B>0$, depending only on $K$ and the fixed geometry of $\cK$, with the following property.
Let $u_k: D_1^+\longrightarrow Q$ be a sequence of $W^{2,2}$ conformal immersions with
\[
(u_k)^*g_{\R^K}
=
e^{2\lambda_k}(dx^2+dy^2).
\]
Assume
\[
u_k(I_1)\subset\cK,
\qquad (u_k)_y\perp T_{u_k}\cK
\quad\text{a.e. on }I_1,
\]
and
\[
E(u_k; D_1^+)+W_{\R^K}(u_k; D_1^+)+T_{\R^K}(u_k; D_1^+)\leq\epsilon_B,
\]
with the total curvature $T_{\R^K}$ defined in Subsection~\ref{sec:fbbci}.
After passing to a subsequence, exactly one of the following alternatives holds on every $ D_s^+$, $0<s<1/2$.

\begin{enumerate}
\item[\textup{(a)}] There is $q\in\cK$ such that
\[
u_k\longrightarrow q
\quad\text{uniformly and strongly in }W^{2,2}( D_s^+)
\]
and strongly in $W^{1,r}( D_s^+)$ for every finite $r$.

\item[\textup{(b)}] The functions $\lambda_k$ are bounded in $L^\infty( D_s^+)\cap W^{1,2}( D_s^+)$, the maps $u_k$ are bounded in $W^{2,2}( D_s^+)$, and
\[
u_k\rightharpoonup u_\infty
\quad\text{in }W^{2,2}( D_s^+),
\qquad
u_k\longrightarrow u_\infty
\quad\text{in }W^{1,r}( D_s^+)
\]
for every finite $r$.
The limit is an unbranched conformal immersion and satisfies
\[
u_\infty(I_s)\subset\cK,
\qquad
(u_\infty)_y\perp T_{u_\infty}\cK
\quad\text{a.e. on }I_s.
\]
\end{enumerate}
\end{lemma}

\begin{proof}
Apply Theorem~\ref{intsma0} with $p=3/2$ to the maps $u_k$ and decrease $\epsilon_B$ if necessary to obtain
\[
\|\nabla u_k\|_{W^{1,3/2}( D_{3/4}^+)}
\leq C\|\nabla u_k\|_{L^2( D_1^+)}.
\]
Thus $\nabla u_k$ is bounded in $L^6( D_{3/4}^+)$, and $u_k$ is uniformly H\"older continuous there.
Fix the boundary point $0\in I_1$ and pass to a subsequence such that $u_k(0)\to a\in\cK$.
The preceding estimate and Morrey's inequality give $\operatorname{Osc}_{D_{3/4}^+}u_k\leq C\epsilon_B^{1/2}$.
By compactness of $\cK$, we may choose $\epsilon_B$ uniformly small enough that, for all large $k$, the images $u_k(D_{3/4}^+)$ lie in one tubular coordinate neighborhood on which $T\cK$ can be sent to a fixed subspace by a smooth orthogonal gauge.

Let $G_k$ be the oriented Gauss map of $u_k$.
Choose the orthogonal gauge $O$ on the common tubular neighborhood and set $\widetilde G_k=O(u_k)G_k$, as in the proof of Lemma~\ref{lem:general-boundary-removability}.
The Gauss map identity gives
\begin{equation}\label{eq:boundary-gauge-gauss-energy}
 \int_{D_{3/4}^+}|\nabla\widetilde G_k|^2
 \leq C_{\cK}\left(T_{\R^K}(u_k;D_{3/4}^+)+E(u_k;D_{3/4}^+)\right).
\end{equation}
Use the same projective reflection $[z]\mapsto[R_V\overline z]$ to obtain $\widehat G_k$ on $D_{3/4}$, and extend it to the plane by radial inversion after rescaling.
The trace matching and the opposite orientation of the exterior copy give
\[
 \int_{\R^2}\widehat G_k^*\Omega_{\mathbb{CP}}=0.
\]
Here $\Omega_{\mathbb{CP}}$ has the normalization in \cite[Theorem~2.1]{kuwert2012}.
Reflection and inversion multiply the energy by a fixed factor.
Thus \eqref{eq:boundary-gauge-gauss-energy} allows us to choose $\epsilon_B$ so that the full-plane absolute Jacobian integral is at most a fixed $\gamma<2\pi$.
The normalized projective potential from \cite[Theorem~2.1]{kuwert2012} satisfies
\[
 -\Delta\psi_k=*\widehat G_k^*\Omega_{\mathbb{CP}},
 \qquad
 \|\psi_k\|_{L^\infty(\R^2)}+\|\nabla\psi_k\|_{L^2(\R^2)}
 \leq C(\gamma)\int_{\R^2}|\nabla\widehat G_k|^2.
\]
Its normalization is $\psi_k(z)\to0$ at infinity.
The right-hand side is even.
Uniqueness therefore makes $\psi_k$ even.

Use the connection form $\beta_k$ and vector field $F_{\beta_k}$ from Lemma~\ref{lem:weak-boundary-connection}, applied to $u_k$.
Since $|\beta_k|\leq C_{\cK}|\nabla u_k|$, the fields $F_{\beta_k}$ are bounded in $L^6$.
Lemma~\ref{lem:weak-boundary-connection} gives
\[
\int_{D_{3/4}^+}\nabla(\lambda_k-\psi_k)\cdot\nabla\varphi
=\int_{D_{3/4}^+}F_{\beta_k}\cdot\nabla\varphi
\]
for every smooth $\varphi$ compactly supported in $D_{3/4}^+\cup I_{3/4}$.
In particular, these test functions may be nonzero on the diameter $I_{3/4}$.
Choose a fixed smooth half-domain $\Omega^+$ with
\[
 D_{2/3}^+\cup I_{2/3}\Subset_{\overline{\mathbb H}}\Omega^+
 \Subset_{\overline{\mathbb H}}D_{3/4}^+\cup I_{3/4},
\]
flat along $I_{2/3}$.
As in the proof of Lemma~\ref{lem:general-boundary-removability}, solve the mean-zero Neumann problem
\[
 \int_{\Omega^+}\nabla q_k\cdot\nabla\varphi
 =\int_{\Omega^+}F_{\beta_k}\cdot\nabla\varphi,
 \quad\varphi\in W^{1,2}(\Omega^+).
\]
The right-hand side vanishes on constants.
The $L^6$ estimate gives
\[
 \|q_k\|_{W^{1,6}(\Omega^+)}\leq C\|F_{\beta_k}\|_{L^6(\Omega^+)}.
\]
Thus $q_k$ is uniformly bounded and H\"older continuous.
Subtracting the weak equations shows that $H_k:=\lambda_k-\psi_k-q_k$ is harmonic with zero weak normal flux on the diameter.
Its even reflection is weakly harmonic on a full neighborhood of $D_{1/2}$, and elliptic regularity makes it harmonic there.

Since
\[
\int_{ D_1^+}e^{2\lambda_k}
=
E(u_k; D_1^+)
\]
and $\psi_k,q_k$ are uniformly bounded,
\[
\int_{ D_{1/2}^+}e^{2H_k}\leq C.
\]
The mean-value inequality for the positive subharmonic function $e^{2H_k}$, applied after even reflection, bounds $H_k$ uniformly from above on smaller disks.
Apply Harnack's inequality~\cite{GilbargTrudinger2001} to a fixed upper bound minus $H_k$ and pass to a subsequence.
Then either
\[
H_k\longrightarrow-\infty
\]
locally uniformly or $H_k$ is locally bounded above and below, and all its derivatives are locally bounded.

In the first case, the bounds for $\psi_k$ and $q_k$ give $\lambda_k\to-\infty$ locally uniformly.
Since $|\nabla u_k|=\sqrt2e^{\lambda_k}$, the gradients converge uniformly to zero on every smaller half disk.
By compactness of $Q$, Poincar\'e's inequality, and the boundary trace condition, a subsequence satisfies $u_k\longrightarrow q\in\cK$ uniformly and strongly in $W^{1,r}$ for every finite $r$.

To complete the collapsing case, we estimate the second derivatives.
Let $U_k$ be the reflected map constructed in the proof of Theorem~\ref{intsma0}.
The same trace matching gives $U_k\in W^{2,2}_{\loc}$ with no interface measure in its Laplacian.
For $0<s<t<1/2$, the conformal equation and the chain rule on the reflected half give
\[
\|\Delta U_k\|_{L^2( D_t)}
\leq
C\|\nabla u_k\|_{L^\infty( D_t^+)}
\left(
 \|\nabla u_k\|_{L^2( D_t^+)}
 +W_{\R^K}(u_k; D_t^+)^{1/2}
\right).
\]
The right-hand side tends to zero because the energies are uniformly bounded and $\|\nabla u_k\|_{L^\infty( D_t^+)}\to0$.
An interior cutoff Calder\'on--Zygmund estimate~\cite{GilbargTrudinger2001} therefore gives
\[
\|D^2U_k\|_{L^2( D_s)}
\leq
C_{s,t}\left(
 \|\Delta U_k\|_{L^2( D_t)}
 +\|\nabla U_k\|_{L^2( D_t)}
\right)
\longrightarrow0.
\]
Restricting to the upper half disk proves strong $W^{2,2}( D_s^+)$ convergence to $q$.
This is alternative \textup{(a)}.

In the second case, $\psi_k$ is bounded in $L^\infty\cap W^{1,2}$, $q_k$ is bounded in $W^{1,6}$, and $H_k$ has local derivative bounds.
Thus $\lambda_k$ is locally bounded in $L^\infty\cap W^{1,2}$.
The conformal Hessian estimate
\[
|D^2u_k|^2
\leq
Ce^{2\lambda_k}
\bigl(
|\nabla G_k|^2+|\nabla\lambda_k|^2
\bigr),
\]
together with
\[
\int|\nabla G_k|^2
\leq CT_{\R^K}(u_k),
\]
gives a uniform $W^{2,2}$ bound on smaller half disks.
Weak compactness in $W^{2,2}$ and the compact embedding $W^{2,2}\hookrightarrow W^{1,r}$ for finite $r$ give the claimed convergence.
The upper and lower bounds for $\lambda_k$ ensure that the limit is an unbranched conformal immersion.
Uniform convergence preserves the boundary condition in $\cK$.
The boundary derivatives converge weakly in their $W^{1/2,2}$ trace spaces.
The tangent projections of $\cK$ converge uniformly along the maps, so the orthogonality condition also passes to the limit.
This proves alternative \textup{(b)}.
\end{proof}

For a sequence $u_k$ on a fixed parameter domain, take a weak limit $\nu$ of the uniformly bounded measures
\[
\nu_k=
\left(1+\frac14|\mathbf H_{u_k}|^2+|A_{u_k}|^2\right)d\mu_{u_k}.
\]
The set $S=\{x:\nu(\{x\})\geq\epsilon_B\}$ is finite.
For the sequences considered below, Proposition~\ref{prop:automatic-curvature-branching-bound} also bounds the number of branch points.
After a further subsequence, their accumulation points form a finite set $B$.
On compact subsets away from $S\cup B$, sufficiently small charts satisfy the hypotheses of Lemma~\ref{lem:boundary-conformal-factor-compactness}.

\vspace{1cm}

\section{\texorpdfstring{Proof of the general compactness theorem}{Proof of the general compactness theorem}}
\label{sec:general-compactness-proof}

We prove Theorem~\ref{thm:intro-main-compactness} using the construction of \cite[Section~4]{chen2021}, the interior estimates of \cite[Sections~2-3]{chen2014}, and the boundary estimates of Section~\ref{sec:fundamental-ingredients-bubbling}.

\Needspace{6\baselineskip}
\subsection{\texorpdfstring{Domains and local convergence}{Domains and local convergence}}
\label{subsec:generic-extraction}

We construct conformal coordinates on fixed cores of the degenerating domains and prove local convergence away from the specified exceptional points.

The bordered Fenchel-Nielsen convergence in Definition~\ref{def:bordered-DM-convergence} follows from \cite[Theorem~4.9]{Liu2004} when $\chi(\Sigma)<0$.
The node classification is given in \cite[Definition~3.4 and Section~3.2]{Liu2004}.

On the symmetric doubles, exchanged geodesics, invariant geodesics with two fixed longitudinal arcs, and fixed geodesics give the three quotient collars \eqref{eq:sec2-interior-collar-model}-\eqref{eq:sec2-type-II-collar-model}, respectively.
The normalized disk does not degenerate.
For the annuli $\mathfrak A_{\ell_n}$ from Subsection~\ref{subsec:uniformization-boundary}, pass to a subsequence on which $\ell_n$ tends to a finite positive value, zero, or infinity.
In the first case, the conformal structures converge smoothly.
If $\ell_n\to0$, cut along a transverse arc and use the coordinates $(-\pi\theta/\ell_n,\pi t/\ell_n)$ to obtain a longitudinal half-cylinder with cross-section $[0,\pi]$.
If $\ell_n\to\infty$, the coordinate $z=e^{t-\ell_n+i\theta}$ identifies the annulus with $\{e^{-\ell_n}\leq|z|\leq1\}$, whose inner boundary gives a type~II end and whose remaining compact subsets exhaust the punctured disk.

We call a local limiting map in principal, rescaled, or translated coordinates a profile.
A compact core is a compact region of its parameter domain obtained by removing neighborhoods of the points at which additional selected profiles attach and truncating the artificial ends.
A core may still contain low-mass curvature concentration points, as explained below.

On a fixed principal core in $\widetilde\Sigma_0$, choose a smooth metric $b$ on that core and set
\[
    j_n':=\Psi_n^*j_n,\qquad
    \widetilde h_n(X,Y):=\frac12\left(b(X,Y)+b(j_n'X,j_n'Y)\right).
\]
Define $\widetilde h_0$ by the same formula, using $\widetilde j_0$ instead of $j_n'$.
Definition~\ref{def:bordered-DM-convergence} gives $j_n'\to\widetilde j_0$ smoothly, and hence $\widetilde h_n\to \widetilde h_0$ smoothly.
Each metric $\widetilde h_n$ is compatible with the conformal structure $j_n'$.
In dimension two, it is therefore conformally equivalent to $\Psi_n^*h_n$.
Replacing $\Psi_n^*h_n$ by $\widetilde h_n$ preserves conformality and Dirichlet energy, while the induced area is independent of this choice.

At an interior point, \cite{chen2014} gives disk maps $\theta_n$ conformal with respect to the metric $\widetilde h_n$ and converging smoothly to a map $\theta_0$ conformal with respect to $\widetilde h_0$.
For a point on the boundary inherited from $\partial\Sigma$, start with an oriented conformal coordinate $x_0+iy_0$ for the metric $\widetilde h_0$, where $y_0=0$ on this boundary.
Choose a simply connected half-neighborhood $\Omega$, and let
\[
    \Delta_{\widetilde h_n}y_n=0,
    \qquad
    y_n|_{\partial\Omega}=y_0|_{\partial\Omega}.
\]
The difference $y_n-y_0$ has zero Dirichlet data and satisfies
\[
 \Delta_{\widetilde h_n}(y_n-y_0)
 =(\Delta_{\widetilde h_0}-\Delta_{\widetilde h_n})y_0.
\]
The right-hand side tends smoothly to zero.
The maximum principle and uniform local boundary elliptic estimates~\cite{GilbargTrudinger2001} give $y_n\to y_0$ smoothly on smaller closed half-neighborhoods away from the artificial corners.
After shrinking once more, $dy_n\neq0$.
Use the oriented Hodge star convention $*_{\widetilde h_n}dx_n=dy_n$, equivalently
\[
    dx_n=-*_{\widetilde h_n}dy_n.
\]
The right-hand side is closed because $y_n$ is harmonic.
Normalize $x_n$ by requiring $x_n=x_0$ at the base point.
Then $x_n+iy_n$ is an oriented conformal boundary coordinate for the metric $\widetilde h_n$ that maps this boundary to the real axis and converges smoothly to $x_0+iy_0$.
Inverting on a common smaller half disk gives the required maps $\theta_n$.

Since $u_n\circ\Psi_n\circ\theta_n$ is conformal, we apply the local estimates to this composition; $\Psi_n$ itself need not be exactly conformal.
Smooth convergence of $\theta_n$ transfers Sobolev convergence and mean curvature lower semicontinuity to fixed cores.
This also applies to the normalized disk and annulus models.

Proposition~\ref{prop:automatic-curvature-branching-bound} gives uniform bounds for total curvature and weighted branching. We use separate measures for bubble selection and conformal factor compactness.
In a conformal chart, set
\[
\begin{split}
\mu_n(\mathfrak{A})
:=E(u_n;\mathfrak{A})+W_{\R^K}(u_n;\mathfrak{A}),\qquad
\kappa_n(\mathfrak{A})
:=\int_\mathfrak{A}|A_{u_n}|^2\,d\mu_{u_n},\qquad
\nu_n:=\mu_n+\kappa_n.
\end{split}
\]
These measures agree on overlapping conformal charts and have uniformly bounded total mass.
Fix $\Lambda<\infty$ such that $\mu_n(\Sigma)\leq\Lambda$ for every $n$.
The coordinate changes below preserve this bound.

Choose $\vartheta>0$ below every $E+W$ threshold used in the interior and boundary small energy estimates and in Theorem~\ref{deca}.
Require also that $\vartheta<\pi$, that $\sqrt{\vartheta}\leq\epsilon_*$ and $1-C_0\sqrt{\vartheta}>0$ for the constants in Theorem~\ref{deca}, and that $\vartheta$ lie below both the joint $E+W+T$ threshold in Lemma~\ref{lem:boundary-conformal-factor-compactness} and the curvature threshold in Theorem~\ref{sub}.
Choose an integer $m_0$ exceeding the number of disks, half disks, or unit strips needed to cover any of the fixed enlargements in these estimates.
Fix
\begin{equation}\label{eq:selection-delta-margin}
    0<\delta
    \quad\text{so that}\quad
    64m_0\delta<\vartheta.
\end{equation}
On every principal, rescaled, transition, or collar profile, take weak limits of the pulled-back measures $\mu_n$ and $\nu_n$.
Each atom of the $\mu$-limit with mass at least $4\delta$ initiates a recursive extraction step at a concentration point.
An atom of the $\nu$-limit with mass at least $\delta$ is marked for conformal factor compactness.
Repeat both selections on each new profile.
An extraction step is an immediate successor of another step if it is initiated at a concentration point on a profile produced by that step, with no intervening recursive extraction step.
The step producing that profile is the immediate predecessor of the new extraction step.
This predecessor relation differs from adjacency between profiles.

A $\nu$-mark whose $\mu$-mass is less than $4\delta$ is kept inside the profile core.
Choose a neighborhood whose $\mu_n$-mass is less than $8\delta$ for all sufficiently large indices.
The interior or boundary small energy estimate, with a fixed exponent $p>1$, gives a uniform $W^{2,p}$ bound there.
Compactness of $W^{2,p}\hookrightarrow W^{1,2}\cap C^0$ on smaller charts gives strong $W^{1,2}$ convergence and uniform convergence across the mark.
Thus the Dirichlet energy and the oscillation of $u_n$ on shrinking neighborhoods tend to zero.
The $\nu$-atom may still carry positive concentrated curvature or Willmore energy.
Weak $W^{2,2}$ convergence has so far been proved only away from these marks.
For a nonconstant component, Subsection~\ref{subsec:compactification} extends this convergence across each such mark at which the extended limit is unbranched.
A mark where the limit is branched remains in the allowed branch set.

Every compact subset avoiding the finitely many $\nu$-marks and branch points admits a finite cover on which the joint mass lies below the conformal factor threshold.
On each interior chart or chart along the inherited boundary, Theorem~\ref{sub} and Lemma~\ref{lem:boundary-conformal-factor-compactness} give two alternatives: the maps either collapse to a point or converge to a conformal immersion. These alternatives agree on overlaps.
In the noncollapsed case, the logarithmic conformal factors are locally bounded above and below.
In the collapsed case, the gradients tend uniformly to zero on compact subsets.

\begin{lemma}\label{lem:local-convergence}
Let $L$ be a compact subset of the parameter domain of a selected profile, away from its nodal and attaching points, the concentration points of the pulled-back measures $\nu_n$, and the limits of the branch points of the pulled-back maps.
After passing to a subsequence, the pulled-back maps converge weakly in $W^{2,2}(L)$ and strongly in $W^{1,r}(L)$ for every $r > 1$.
The limiting profile is either an unbranched conformal immersion or a constant.
In the constant case, convergence is uniform on compact subsets of the connected core where the local alternatives have been established.
The same statements hold up to the boundary inherited from $\partial\Sigma$.
\end{lemma}

\begin{proof}
The preceding argument shows that the local alternatives agree on overlaps.
In an interior chart, Theorem~\ref{sub} gives the local $W^{2,2}$ bound in both alternatives.
In a chart along the boundary inherited from $\partial\Sigma$, the noncollapsed bound and the collapsing $W^{2,2}$ estimate are the two alternatives of Lemma~\ref{lem:boundary-conformal-factor-compactness}.
\end{proof}

\Needspace{6\baselineskip}
\subsection{\texorpdfstring{Bubble extraction and neck decomposition}{Bubble extraction and neck decomposition}}
\label{subsec:finite-profile-selection}

The following proposition gives finitely many profiles with disjoint embedded cores and decomposes their complement into neck regions.

\begin{prop}\label{prop:finite-profile-selection}
After passing to a subsequence, we obtain a finite family of principal, disk, sphere, transition, and collar profiles with the following properties.
We retain all constant profiles in this family before contraction.
\begin{enumerate}[label=\textup{(\roman*)},leftmargin=2.4em]
\item Any two source regions belonging to recursive extraction steps are disjoint or nested.
Source regions for distinct immediate successors of the same extraction step have disjoint closures and avoid the region assigned to that step in \textup{(ii)}.

\item For an extraction step $c$, let $U_{n,c}(\rho)$ be the image in the original surface of a radius-$\rho$ neighborhood of the atom initiating $c$, in the profile coordinates in which that atom was selected, and set
\[
    m(c):=\lim_{\rho\downarrow0}\lim_{n\to\infty}
    \mu_n\bigl(U_{n,c}(\rho)\bigr),
\]
where $\rho$ ranges over continuity radii.
To each step $c$, we assign a measurable region $D_{n,c}$ in the original surface such that
\[
        \mu_n(D_{n,c})=2\delta.
    \]
All the sets $D_{n,c}$ are pairwise disjoint.
If $\operatorname{ch}(c)$ is the set of the immediate successor extraction steps of $c$, then
\begin{equation}\label{eq:recursive-call-budget}
        2\delta+
        \sum_{c'\in\operatorname{ch}(c)}m(c')
        \leq m(c).
    \end{equation}
Thus both the number and the depth of the recursive extraction steps are finite.
Translated window profiles need not be assigned separate regions in the original surface.
The mass bound and the window selection threshold bound their number within each extraction step.

\item In the logarithmic coordinates of every complementary annulus, half-annulus, or degenerating collar, each residual unit strip, together with the fixed enlargement used in Theorem~\ref{deca}, has $\mu_n$-mass strictly less than $\vartheta$ for all sufficiently large $n$.
Fixed coordinate neighborhoods of profile centers and artificial attaching ends are included in the adjacent cores.
The portions inherited from $\partial\Sigma$, including the circle at $t=0$ of a type~II collar, are retained.

\item Choose compact exhaustions $K_{v,j}$ of the finitely many profile domains with their nodal and attaching points removed, before contraction.
For every fixed $j$ and all sufficiently large $n$, there are smooth embeddings
\[
        \Phi_{n,v,j}:K_{v,j}\longrightarrow\Sigma
    \]
with pairwise disjoint images.
On a principal component these are restrictions of the Deligne--Mumford comparison maps, and the pulled-back conformal structures converge smoothly.
On disk, sphere, transition, and collar profiles, the embeddings are conformal in the selected profile coordinates.
Up to boundaries of measure zero, the complement is the union of the corresponding truncated residual annuli, half-annuli, and collars from \textup{(iii)}.
Subthreshold marks arising only from curvature concentration remain inside the profile cores.
Convergence across them is strong in $W^{1,2}$ and uniform.
They are excluded only from the statement of $W^{2,2}$ convergence.

\item The incidence graph before contraction is finite.
Before principal branch vertices are identified in the global nodal surface, the local graph over each point is connected and acyclic.
Every boundary circle of a type~II collar is retained.
\end{enumerate}
\end{prop}

\begin{proof}
We retain all constant profiles.
To prove finiteness, we assign pairwise disjoint regions of fixed $\mu_n$-mass in the original surface.

A central profile is obtained by rescaling at a selected center and radius $(z_n,r_n)$.
If the boundary stays at bounded rescaled distance, we center the rescaling at the boundary projection of $z_n$.
If $r_n\ll d_n\to0$, where $d_n$ is the distance of $z_n$ from the boundary inherited from $\partial\Sigma$ in the chart in which the center and radius were selected, we also retain the intermediate rescaling by $d_n$.
The interior profile attaches at a distinguished interior point of the transition profile.
A neighborhood of this attaching point is omitted from the transition core; the interior profile is neither selected again nor assigned another region in the original surface.
Every remaining annular region is written in logarithmic coordinates, and a translated window profile is obtained by shifting its longitudinal coordinate by a selected sequence $t_n$.
In a fixed annular or collar coordinate, centers at bounded longitudinal distance belong to the same translated or endpoint profile.
We select centers separately only when their mutual distances tend to infinity.
Compact cores are chosen after removing neighborhoods of the points at which additional selected profiles attach and fixed coordinate neighborhoods of the artificial ends.

Consider an atom $a$ of a pulled-back $\mu$-limit with mass $m\geq4\delta$.
There are only finitely many atoms of mass at least $\delta$ in a relatively compact coordinate chart.
Choose relatively open coordinate neighborhoods $ a\in V$, $\overline V\subset U$, where the closure is taken relative to the closed half-plane in a chart along $\partial\Sigma$.
Take $\overline V$ compact, with $a$ in the relative interior of $V$, and ensure that no other $\mu$-atom of mass at least $\delta$ lies in $U$.
Choose boundaries of zero limiting measure, with source images disjoint from those of the other atom neighborhoods on the same profile.
At a distinguished attaching end, first pass to the adjacent annular or collar coordinate.
Exclude the end itself from the atom list.

Write $B_r^\diamond(z)=B_r(z)$ in an interior chart and $B_r^\diamond(z)=B_r(z)\cap\mathbb H$ in a chart along $\partial\Sigma$.
Choose $r_0>0$ so that the closure relative to the chart of every $B_r^\diamond(z)$, with $z\in\overline V$ and $0\leq r\leq r_0$, lies in $U$.
The parameter set $\overline V\times[0,r_0]$ is compact.
Absolute continuity of $\mu_n$ makes
\[
(z,r)\longmapsto\mu_n(B_r^\diamond(z))
\]
continuous.

For every sufficiently small continuity radius $\rho$, weak convergence gives $\mu_n(B_\rho^\diamond(a))>2\delta$ for large $n$.
Continuity from radius zero therefore gives disks centered at $a$ with mass $2\delta$ and radii tending to zero.
Among the pairs in $\overline V\times[0,r_0]$ whose disks have mass $2\delta$, choose one with minimum radius and set
\[
    D_{n,c}=B_{r_{n,c}}^\diamond(z_{n,c}).
\]
The disks centered at $a$ show that $r_{n,c}\to0$.

If $z_{n,c}\to b$ along a subsequence, then $D_{n,c}\subset B_\rho^\diamond(b)$ for every fixed continuity radius $\rho$ and all large $n$.
Passing to the limit and then sending $\rho\downarrow0$ gives $\mu(\{b\})\geq2\delta$.
The choice of $U$ forces $b=a$.
Thus $z_{n,c}\to a$, and the assigned region lies in every fixed neighborhood of $a$ for large $n$.

An admissible disk of radius less than $r_{n,c}$ has mass less than $2\delta$: otherwise shrinking it to mass $2\delta$ would contradict minimality.
In particular, every admissible disk of radius at most $r_{n,c}/2$ has mass less than $2\delta$.
Since $z_{n,c}\to a$, all disks used to cover a fixed compact subset of the rescaled chart are admissible for large $n$.
A finite cover by disks of rescaled radius at most $1/2$, together with \eqref{eq:selection-delta-margin}, puts every fixed enlargement below the $E+W$ regularity threshold.
The central rescaled limit may be constant.

We assign the region $D_{n,c}$ to this extraction step.
An atom of the rescaled limiting measure with mass at least $4\delta$ cannot occur at bounded rescaled distance inside the assigned disk: a sufficiently small rescaled disk about that atom would pull back to an admissible disk of radius at most $r_{n,c}/2$ and mass greater than $2\delta$, contrary to minimality.
Divide the remaining source neighborhood into the radial bands and transition regions described below.
Perform each new recursive extraction there, away from the distinguished attaching points.

Once profiles at bounded distance have been identified, fixed compact subsets of distinct central, transition, and window profiles pull back to disjoint radial or strip bands.
For any finite collection of proposed immediate successor extraction steps, first choose their atom neighborhoods and the compact profile subsets containing them.
For all sufficiently large $n$, their pullbacks are pairwise disjoint, lie outside $D_{n,c}$, and lie inside every fixed source neighborhood of the current extraction step.
This separation holds for each finite family, independently of finiteness of the full construction.
It also shows that source regions of steps obtained recursively from distinct immediate successors remain disjoint, and that each region assigned to a later step is disjoint from the regions assigned to all preceding steps in the same nested sequence.
Candidates converging to the same atom belong to a single extraction step.
Candidates converging to distinct atoms lie in disjoint isolating neighborhoods.

Let $c_1,\ldots,c_q$ be any finite collection of immediate successor extraction steps of $c$.
Choose their isolating neighborhoods $U_{n,c_i}(\rho_i)$ inside the source neighborhood $U_{n,c}(\rho)$.
By the separation argument above, these sets are pairwise disjoint and avoid $D_{n,c}$ for large $n$.
Finite additivity gives
\[
  2\delta+
  \sum_{i=1}^q
  \mu_n\bigl(U_{n,c_i}(\rho_i)\bigr)
  \leq
  \mu_n\bigl(U_{n,c}(\rho)\bigr).
\]
First let $n\to\infty$ along continuity radii.
Then let each $\rho_i\downarrow0$, and finally let $\rho\downarrow0$.
We obtain
\[
    2\delta+\sum_{i=1}^q m(c_i)\leq m(c).
\]
The collection was arbitrary, so this proves \eqref{eq:recursive-call-budget}.
Induction shows that the regions assigned to all extraction steps are pairwise disjoint in the original surface.
Hence
\[
    2\delta\,\#\mathscr C
    =\sum_{c\in\mathscr C}\mu_n(D_{n,c})
    \leq\mu_n(\Sigma)\leq\Lambda
\]
for every finite subcollection $\mathscr C$.
There are therefore at most $\Lambda/(2\delta)$ extraction steps, and no infinite nested chain is possible.

On every central, transition, and window profile, take weak limits of the pulled-back measures $\mu_n$ and $\nu_n$ and pass to a subsequence on which the branch points have limits in the compactified profile domain.
The $\nu$-atoms of mass at least $\delta$, together with the branch point limits lying in that domain away from its attaching ends, form a finite marked set.
A marked point with $\mu$-mass at least $4\delta$ starts a new extraction step.
The remaining marks stay inside the core, with strong $W^{1,2}$ and uniform convergence by the small-$\mu$ argument in Subsection~\ref{subsec:generic-extraction}.

Write every annular region in logarithmic coordinates and use the same coordinates on full collars, longitudinal half-collars, and one-sided type~II collars.
Let $S_n(t)$ be the unit strip centered at $t$, and let $S_n^*(t)$ be the fixed enlargement required in Theorem~\ref{deca}.
At a boundary circle inherited from $\partial\Sigma$, use a one-sided enlargement.
Set $\sigma:=16m_0\delta$.
Then $\sigma<{\vartheta}/{4}$ by \eqref{eq:selection-delta-margin}.

Within one extraction step or one initial collar, select a center sequence $t_n^1$ only when
\[
    \limsup_{n\to\infty}\mu_n(S_n^*(t_n^1))\geq\sigma
\]
and its distance from every existing profile center and every artificial or profile end tends to infinity.
Pass to a subsequence on which $\mu_n(S_n^*(t_n^1))$ converges to a value at least $\sigma$.
Inductively, select $t_n^{j+1}$ on the current subsequence under the same separation condition, and pass to a nested subsequence on which its strip mass converges to a value at least $\sigma$.
Each nested subsequence preserves the mass limits and separations already obtained.

For every finite selected family, the enlarged strips are pairwise disjoint for large $n$.
On the common subsequence for that family, finite additivity and passage to the limit give
\[
    J\sigma
    \leq
    \lim_{n\to\infty}
    \sum_{j=1}^J\mu_n(S_n^*(t_n^j))
    \leq\Lambda.
\]
Thus there are at most $\Lambda/\sigma$ selections, with no need to assign further fixed-mass regions in the original surface.
Retain the final nested subsequence.

Let $Z_n$ be the resulting finite set of selected centers and profile ends in the logarithmic coordinate of the region.
At the boundary end of a type~II collar, $t=0$ belongs to the profile containing the boundary circle and is included in $Z_n$; that circle is retained.
A center $t$ is admissible if the fixed enlargement $S_n^*(t)$ lies in the region, with the corresponding one-sided convention at the boundary inherited from $\partial\Sigma$.
Define
\[
 F(R)
 :=
 \limsup_{n\to\infty}
 \sup\left\{
   \mu_n(S_n^*(t)):
   t\text{ admissible},\ 
   \operatorname{dist}(t,Z_n)\geq R
 \right\},
\]
where the supremum of the empty set is zero.

We have
\[
    \lim_{R\to\infty}F(R)\leq\sigma.
\]
If the limit were larger than $\sigma$, we could choose $\epsilon>0$, radii $R_j\to\infty$, indices $n_j\to\infty$, and admissible centers $s_{n_j}$ such that
\[
 \operatorname{dist}(s_{n_j},Z_{n_j})\geq R_j,
 \qquad
 \mu_{n_j}(S_{n_j}^*(s_{n_j}))\geq\sigma+\epsilon.
\]
Along this subsequence, the centers move arbitrarily far from every selected center and profile end.
This contradicts the maximality of the selected family.

Choose $R_0$ so that $F(R_0)<2\sigma$.
By the definition of the limsup, for all sufficiently large $n$,
\[
 \sup\left\{
   \mu_n(S_n^*(t)):
   t\text{ admissible},\ 
   \operatorname{dist}(t,Z_n)\geq R_0
 \right\}
 <3\sigma<\vartheta.
\]
Include the longitudinal neighborhoods of width $R_0$, together with the fixed enlargements, in the adjacent endpoint or translated profile cores.
Every residual enlarged strip then has mass strictly less than $\vartheta$.
Fix these profiles and neighborhoods before taking $n$ large.

Apply the same $\mu$-atom and $\nu$-atom selection to every translated profile.
New recursive extraction steps are initiated at its selected concentration points.
Once neighborhoods of these points are removed, the remaining window core need not retain positive mass.
On a type~II collar, every center at bounded distance from $t=0$ belongs to the single profile associated with that boundary circle.
The circle $t=0$ remains a boundary component of its parameter domain.
Only centers moving arbitrarily far from this end can give additional profiles.

Let $z_n$ be the center of a selected disk, let $z_n'$ be its boundary projection, and put $d_n=|z_n-z_n'|$.
If $d_n/r_n$ is bounded, rescaling about $z_n'$ by $r_n$ gives a half-plane profile, which compactifies to a disk.
If $d_n/r_n\to\infty$ and $d_n$ does not tend to zero in the coordinates used to select the disk, rescaling about $z_n$ gives a plane profile, which compactifies to a sphere.

Suppose that $r_n/d_n\to0$ and $d_n\to0$.
Choose $\rho_n\downarrow0$ so that $d_n/\rho_n\to0$ and the half disk of radius $\rho_n$ remains inside the isolating source neighborhood.
In flat boundary coordinates take $z_n'=0$ and $z_n=(0,d_n)$.
Up to shared boundaries, partition the region between the selected scale and this half disk into
\[
%\begin{split}
 \mathfrak A_n^{\mathrm{in}}=B_{d_n/4}(z_n)\setminus
 \overline{B_{Rr_n}(z_n)},\quad
 T_n=
 \bigl(\mathbb H\cap B_{4d_n}(0)\bigr)
 \setminus\overline{B_{d_n/4}(z_n)},\quad
  \mathfrak A_n^{\mathrm{out}}
 = \bigl(\mathbb H\cap B_{\rho_n}(0)\bigr)
 \setminus\overline{B_{4d_n}(0)}.
%\end{split}
\]
Here $R$ is fixed and gives the inner truncation radius in units of the selected radius.
The first region is a full annulus, and the third is a half-annulus.
Scaling $T_n$ by $d_n$ gives a bounded half-plane transition core with the distinguished interior attaching point at $(0,1)$.
A compact exhaustion of this transition core and the two annular tails cover the entire intermediate region.

Apply the window construction separately on $\mathfrak A_n^{\mathrm{in}}$ and $\mathfrak A_n^{\mathrm{out}}$.
Windows at bounded distance from the transition core are absorbed into that core.
A new $\mu$-atom in the transition profile starts a new extraction step in a small source neighborhood disjoint from the neighborhood removed at the distinguished attaching point.
No additional region is assigned to the already selected interior rescaling.
Thus the assigned regions remain disjoint in the original domain.
The mass bound \eqref{eq:recursive-call-budget} remains valid.

The same construction applies in the three collar models.
A full collar has interior translated profiles.
A longitudinal half-collar has disk profiles that meet its longitudinal sides $\theta=0,\pi$ and sphere profiles whose rescaled boundary escapes to infinity.
A type~II collar is treated from $t=0$ toward its single artificial end.
A nonconstant profile at the boundary end $t=0$ compactifies to a disk whose boundary circle maps to $\cK$; its interior puncture attaches through the selected chain to the principal special point.
If the profile containing $t=0$ is constant, uniform convergence on every fixed compact core meeting this circle gives vanishing oscillation on that core.
Its limiting value belongs to the closed set $\cK$.

There are finitely many initial principal and collar pieces.
Each recursive extraction step creates only a bounded number of radial regions, at most one transition profile for each selected interior rescaling whose center approaches the inherited boundary at an intermediate scale, and at most $\Lambda/\sigma$ translated window profiles in each such region.
There are at most $\Lambda/(2\delta)$ extraction steps, so the profile family is finite before contraction.

The source regions are disjoint or nested: isolating neighborhoods for distinct immediate successors of the same extraction step are disjoint, all further extractions initiated in one such neighborhood remain inside it, and the radial bands between scales are disjoint.
Window profiles whose centers remain at bounded mutual distance have already been identified.
Hence no two profile cores partially overlap.
For each selected profile, choose an initial compact core containing all subthreshold $\nu$-marks and remove neighborhoods of its attaching points.
Choose the compact exhaustion $K_{v,j}$ by cutting each logarithmic end at length $j$ from a fixed cross-section of that end.
We use these same cross-sections for the neck truncations below and the same exhaustion index $j$ for all profiles.
On a nonprincipal profile, compose its conformal profile coordinate with the original chart.
On a principal profile, restrict the Deligne--Mumford comparison map.
In both cases, this gives $\Phi_{n,v,j}:K_{v,j}\longrightarrow\Sigma$.
For fixed $j$ and large $n$, these embeddings have pairwise disjoint images.
On principal cores, $(\Phi_{n,v,j})^*j_n$ converges smoothly.
On principal cores $K_{v,j}$, use the smoothly converging metrics $\widetilde h_n$ constructed above.
The other profile embeddings are conformal in their fixed profile coordinates.

The radial and collar partitions show that the complement of the core images consists exactly of the residual full annuli, free boundary half annuli, and collar intervals in part \textup{(iii)}.
Subthreshold marks arising only from curvature concentration remain in the cores.
We exclude them only from the local $W^{2,2}$ convergence statement.
On each fixed set in an exhaustion of a collapsed profile core, the pulled-back maps converge strongly in $W^{1,2}$ and uniformly to the constant limiting map.

Keep the two normalized branches of every ordinary node as distinct root vertices.
Ordered collar profiles subdivide the edge between these roots.
A type~II collar has one principal root and, when present, one terminal disk at the boundary circle $t=0$.
Recursively selected profiles attach in rooted trees.
The local incidence graph is connected and acyclic.
We identify the distinct branch roots when forming the global nodal surface.
This preserves the cycles created by nonseparating pinches.
This proves the proposition.
\end{proof}

\Needspace{6\baselineskip}
\subsection{\texorpdfstring{Vanishing Dirichlet energy and oscillation on necks}{Vanishing Dirichlet energy and oscillation on necks}}
\label{subsec:vanishing-necks}

The next lemma proves vanishing Dirichlet energy and oscillation on residual necks, including type~II necks.

\begin{lemma}\label{lem:noneck}
For every residual edge $e$, let $\mathcal N_{n,e}(R)$ be the neck remaining after an end segment of length $R$ in logarithmic coordinates is removed at every artificial end.
Then
\begin{align}\label{noneck}
\lim_{R\to\infty}\limsup_{n\to\infty}
E\bigl(u_n;\mathcal N_{n,e}(R)\bigr)=0,\qquad
\lim_{R\to\infty}\limsup_{n\to\infty}
\operatorname{Osc}_{\mathcal N_{n,e}(R)}u_n=0.
\end{align}
The same conclusion holds on any sequence of subintervals between two consecutive selected windows whose distances from both windows tend to infinity.
In the type~II case, retain the boundary circle $t=0$ and remove an end segment only at the artificial end.
\end{lemma}

\begin{proof}
The conformal factor expansions in Lemma~\ref{lem:general-boundary-removability} and \cite[Theorem~3.1]{kuwert2012} show that the boundary terms on excised neighborhoods of branch points vanish as those neighborhoods shrink.
Proposition~\ref{prop:finite-profile-selection}(iii) verifies the hypothesis on unit strips, including the fixed enlargements used in the proof of Theorem~\ref{deca}.
Apply that theorem with $\epsilon=\sqrt{\vartheta}$; the choice of $\vartheta$ ensures that $\epsilon\leq\epsilon_*$ and $\gamma=1-C_0\epsilon>0$.
It applies to every residual full, longitudinal half-, or type~II cylinder.

Fix $R\geq4$ before taking the limit in $n$.
For a full or longitudinal half-neck, apply Theorem~\ref{deca} to $u_n$ from both artificial ends.
The total neck energy is bounded by $\Lambda$, so
\[
\limsup_{n\to\infty}
E\bigl(u_n;\mathcal N_{n,e}(R)\bigr)
\leq
C\Lambda e^{-\gamma R},
\qquad
\gamma=1-C_0\epsilon>0.
\]
The oscillation estimate in the same theorem gives
\[
\limsup_{n\to\infty}
\operatorname{Osc}_{\mathcal N_{n,e}(R)}u_n
\leq
C\Lambda^{1/2}e^{-\gamma R/2}.
\]
If a type~II residual region still contains the boundary circle $t=0$, apply the one-sided estimate, retaining that boundary and removing an end segment only at the artificial end.
If the core of a terminal disk profile already contains this circle, the connecting residual region has two artificial ends.
In logarithmic coordinates it is a full cylinder, so we apply the estimate from both ends.
These estimates hold for every fixed $R$; sending $R\to\infty$ proves \eqref{noneck}.

Let $I_n$ be a sequence of subintervals between two consecutive selected windows whose distances from both windows tend to infinity.
For every fixed $R$, the corresponding cylindrical subregion lies in $\mathcal N_{n,e}(R)$ for all sufficiently large $n$.
Its limiting Dirichlet energy and oscillation are therefore bounded by $C\Lambda e^{-\gamma R}$ and $C\Lambda^{1/2}e^{-\gamma R/2}$, respectively.
Letting $R\to\infty$ proves the assertion for these subintervals.
Since there are finitely many edges, these estimates give
\[
\lim_{R\to\infty}\limsup_{n\to\infty}
\sum_e E(u_n;\mathcal N_{n,e}(R))=0,
\qquad
\lim_{R\to\infty}\limsup_{n\to\infty}
\max_e\operatorname{Osc}_{\mathcal N_{n,e}(R)}u_n=0.
\]
If there are no residual edges, the oscillation term is omitted.
\end{proof}

\Needspace{6\baselineskip}
\subsection{\texorpdfstring{Compactification and ghost components}{Compactification and ghost components}}
\label{subsec:compactification}

We extend the nonconstant limiting maps across their punctures and recover local weak $W^{2,2}$ convergence at auxiliary marks where the limit is unbranched.
Lemma~\ref{lem:noneck} then lets us assemble the surface of Definition~\ref{def:free-boundary-bubble tree-surface}, retaining the necessary ghosts.
We show that the permitted contractions preserve vanishing neck energy and oscillation.

On a profile with nonconstant limiting map, use the smoothly converging conformal charts from Subsection~\ref{subsec:generic-extraction} on a principal core and the fixed profile coordinates on every other core.
Write $v_n$ for the maps pulled back to these charts and $v$ for their limit, with induced metrics $e^{2\lambda_n}|dz|^2$ and $e^{2\lambda}|dz|^2$ away from the marked points.
On a compact subset avoiding these points, the noncollapsed alternative gives uniform upper and lower bounds for $\lambda_n$, weak $W^{2,2}$ convergence of $v_n$, and strong $W^{1,r}$ convergence for every finite $r$.
After passing to a subsequence, the normal projections $P_n^\perp$ and the functions $\lambda_n$ converge almost everywhere to $P^\perp$ and $\lambda$, respectively.
Multiplication by the uniformly bounded coefficients $e^{-\lambda_n}P_n^\perp$ sends each fixed $L^2$ test field to a strongly convergent $L^2$ sequence.
Weak convergence of $D^2v_n$ therefore gives, for $i,j\in\{1,2\}$,
\[
 e^{-\lambda_n}P_n^\perp\partial_i\partial_jv_n
 \rightharpoonup
 e^{-\lambda}P^\perp\partial_i\partial_jv
 \quad\text{weakly in }L^2.
\]
Since the sum of the squared Euclidean norms of these four fields is the total curvature density,
\[
 |A_{v_n}|^2\,d\mu_{v_n}
 =
 e^{-2\lambda_n}
 \sum_{i,j=1}^2
 |P_n^\perp\partial_i\partial_jv_n|^2\,dx\,dy,
\]
weak lower semicontinuity gives the local total curvature inequality.
A nonnegative partition of unity gives this inequality on each compact core, including the portions inherited from $\partial\Sigma$ away from the excluded marked points.
Summing over the disjoint embedded cores and exhausting the nonconstant components gives
\[
 \sum_{i\in\mathfrak I_*}T_{\R^K}(u_\infty^i)
 \leq\liminf_{n\to\infty}T_{\R^K}(u_n).
\]
The right-hand side is finite by Proposition~\ref{prop:automatic-curvature-branching-bound}.
The integrals on the left are over the nonconstant component domains with their branch points and nodal, attaching, and concentration marks removed.
The same exhaustion and local strong $W^{1,2}$ convergence give finite area.

Fill the interior punctures using \cite[Theorem~3.1]{kuwert2012} and the boundary punctures using Lemma~\ref{lem:general-boundary-removability}.
The free boundary conditions persist under these extensions.
By Lemma~\ref{lem:noneck}, the values at the two ends of every ordinary edge agree, and the one-sided matching condition holds at each type~II special point.

Let $a$ be a subthreshold $\nu$-mark on a nonconstant profile, so that $\mu(\{a\})<4\delta$. 
The small-$\mu$ estimate gives strong $W^{1,2}$ and uniform convergence near $a$.
The punctured limit is a conformal immersion with finite area and total curvature.
The interior or boundary removability result therefore extends the limit across $a$.
If the extension is branched at $a$, then $a$ belongs to the branch set allowed in Theorem~\ref{thm:intro-main-compactness}(ii).
If the extension is unbranched at $a$, we show that the sequence has a uniform local $W^{2,2}$ bound across $a$.

Shrink the coordinate disk centered at $a$, or the half disk centered at $a$ on the boundary inherited from $\partial\Sigma$, so that it contains no other $\nu$-mark and no branch point of the limiting immersion.
Write
\[
 g_{v_n}
 =
 e^{2\lambda_n}|dz|^2.
\]
On an annulus or half-annulus about $a$, the noncollapsed alternative for the conformal factor gives uniform $L^\infty$ and $W^{1,2}$ bounds for $\lambda_n$.
Let
\[
 dw_n=\frac14|\mathbf H_{v_n}|^2\,d\mu_{v_n},\qquad
 d\kappa_n=|A_{v_n}|^2\,d\mu_{v_n},
\]
and let $d\zeta_n=K_{g_{v_n}}\,d\mu_{v_n}$ be the signed Gauss curvature measure.
These measures have weak limits $w$, $\kappa$, and $\zeta$, respectively, along the fixed subsequence.
The Euclidean Gauss equation is the measure identity
\begin{equation}\label{eq:auxiliary-mark-gauss-measures}
    \kappa_n=4w_n-2\zeta_n.
\end{equation}
We determine the atom of the limiting Gauss curvature measure at $a$.
In this chart, $\cB_{v_n}$ denotes the branch set of the coordinate map, with branching orders $m_{v_n}(b)$.
In an interior disk, the distributional conformal factor equation is
\begin{equation}\label{eq:auxiliary-mark-conformal-factor}
 -\Delta\lambda_n
 =
 \zeta_n
 -
 2\pi\sum_{b\in\cB_{v_n}}m_{v_n}(b)\delta_b.
\end{equation}
The bounds on total curvature and weighted branching give a uniform bound for the total variation of the right-hand side.
On a smaller disk, subtract the Green potential of the right-hand side.
Green potentials of bounded measures are precompact in $L^1$.
Split the logarithmic kernel into a bounded continuous part away from the diagonal and a part near the diagonal with uniformly small $L^1$-norm.
The remaining functions are harmonic.
The annular bounds for $\lambda_n$ control their normalization.
The potentials have bounded $L^1$-norm on the annulus.
Fubini's theorem therefore gives a surrounding circle on which the harmonic remainders have uniformly bounded $L^1$-norm.
The Poisson formula then bounds these harmonic remainders on every smaller disk.
Thus $\lambda_n$ is precompact in $L^1$ across $a$, and its limit agrees on the punctured disk with the conformal factor $\lambda$ of the extended limiting map $v$.

Passing to distributions in \eqref{eq:auxiliary-mark-conformal-factor} gives the atom identity
\begin{equation}\label{eq:auxiliary-mark-interior-atom}
 \zeta(\{a\})= 2\pi m_{\mathrm{int}}-2\pi m_v(a),
\end{equation}
where $\zeta$ is the weak limit of $\zeta_n$, $m_{\mathrm{int}}$ is the total order of the interior branch points of $v_n$ converging to $a$, and $m_v(a)$ is the branch order of the extended limiting map $v$ at $a$.
Since the total branching order is bounded, these integer sums are constant after passing to a further subsequence.

At a mark on the boundary inherited from $\partial\Sigma$, reflect $\lambda_n$ evenly.
With the same convention for the inward conormal,
\[
 k_{g_{v_n}}\,ds_{g_{v_n}}
 =
 -\partial_y\lambda_n\,dx
\]
on the flat diameter, where the derivative denotes the weak boundary flux from \eqref{eq:weak-boundary-curvature}.
Hence the reflected equation contains the line measure $2k_{g_{v_n}}\,ds_{g_{v_n}}$, the paired interior branch atoms, and the branch atom on the full disk associated with each boundary branch point.
Every weak limit of these line measures has no point mass at $a$.
Indeed,
\[
 |k_{g_{v_n}}|\,ds_{g_{v_n}}
 \leq
 C_{\cK}|dv_n|\,dx
\]
on the diameter, while the small-$\mu$ $W^{2,p}$ estimate may be taken with $p=3/2$.
The trace theorem then bounds $dv_n|_{\partial\mathbb H}$ in $L^q$ for some $q>1$, so these boundary densities are uniformly integrable.
Apply the same Green potential argument to the even extensions and compare the point masses to obtain
\begin{equation}\label{eq:auxiliary-mark-boundary-atom}
 \zeta(\{a\})
 =
 2\pi m_{\mathrm{int}}
 +
 \pi m_{\partial}
 -
 \pi m_v(a),
\end{equation}
where $m_{\partial}$ is the total order of the boundary branch points converging to $a$.
These signs and coefficients follow from the distributional conformal factor equation.

Assume now that the extended limit is unbranched at $a$.
Then $m_v(a)=0$, and \eqref{eq:auxiliary-mark-interior-atom} or \eqref{eq:auxiliary-mark-boundary-atom} gives $\zeta(\{a\})\geq0$.
By conformality, the Dirichlet measures equal the induced area measures $\mu_{v_n}$.
Strong $W^{1,2}$ convergence through the mark shows that their weak limit has no atom there.
Consequently,
\[
    w(\{a\})=\mu(\{a\})<4\delta.
\]
Taking atoms in \eqref{eq:auxiliary-mark-gauss-measures} gives
\[
 \kappa(\{a\})
 =
 4w(\{a\})-2\zeta(\{a\})\geq0,
 \qquad
 \zeta(\{a\})\leq2w(\{a\})<8\delta<\pi.
\]
The last inequality follows from $64m_0\delta<\vartheta<\pi$.
If any prelimit branch point converged to $a$, then \eqref{eq:auxiliary-mark-interior-atom} would give $\zeta(\{a\})\geq2\pi$, or \eqref{eq:auxiliary-mark-boundary-atom} would give $\zeta(\{a\})\geq\pi$.
Either conclusion contradicts the bound above.
Thus no branch point converges to $a$.
Shrink the chart to exclude the other limiting branch points and marks.
The maps are then unbranched there for all sufficiently large $n$.

Moreover,
\[
 \nu(\{a\})
 =
 w(\{a\})+\kappa(\{a\})
 \leq5w(\{a\})<20\delta.
\]
Choose a still smaller disk or half disk $D$, with boundary of zero measure for the weak limit, such that $\nu_n(D)<24\delta<\vartheta$ for all sufficiently large $n$.
The maps are unbranched on $D$.
Theorem~\ref{sub} in the interior case and Lemma~\ref{lem:boundary-conformal-factor-compactness} in the boundary case give the local collapse or noncollapse alternatives for $v_n$ in the chosen chart $D$.
The maps already converge to a nonconstant immersion on an annulus or half-annulus contained in $D$, so collapse is impossible.
The noncollapsed alternative therefore gives a uniform $W^{2,2}$ bound on every smaller disk or half disk and strong $W^{1,r}$ convergence for every finite $r$, including across $a$.

Thus, on a nonconstant component, each auxiliary mark is either a limiting branch point or a point across which the required $W^{2,2}$ convergence holds.
On a profile core with constant limiting map, every auxiliary curvature mark is retained in the finite exceptional set $\cB_\infty^i$ from Definition~\ref{def:componentwise-free-boundary-bubble tree-convergence}.
Theorem~\ref{sub} in interior charts and the collapsing alternative of Lemma~\ref{lem:boundary-conformal-factor-compactness} give weak local $W^{2,2}$ convergence away from this finite set, together with the strong $W^{1,2}$ and uniform convergence used in the argument on the cores and necks.

Proposition~\ref{prop:finite-profile-selection} gives a finite incidence graph, including all constant profiles, before any contraction.
By Theorem~\ref{gapth}, every nonconstant sphere or free boundary disk also has a compact core carrying a fixed positive amount of $E+W$.

Before contracting constant profiles, apply Lemma~\ref{lem:noneck} to every edge of the finite graph.
Both the Dirichlet energy and oscillation of $u_n$ vanish on each connecting region.
We suppress an unmarked nonprincipal constant profile only when its vertex has exactly two incident edges and replacing these edges by a single edge preserves the node type and all boundary data inherited from the original surface.
In particular, we retain a constant disk component at a boundary attachment when it is needed to attach a sphere at an interior point of a disk.

A concatenated edge may contain a collapsed profile core with positive limiting $\mu_n$-mass.
We therefore do not apply Theorem~\ref{deca} across that core.
Use the common exhaustion index $j$ for all retained and suppressed profiles.
An ordinary concatenated edge corresponds to a path
\[
 v_0\mathbin{-}e_1\mathbin{-}v_1\mathbin{-}\cdots
 \mathbin{-}e_m\mathbin{-}v_m
\]
in the graph before contraction, where $v_0,v_m$ are retained and $v_1,\ldots,v_{m-1}$ are suppressed.
After truncating both retained endpoint cores at index $j$, the region of this edge is
\[
 \Omega_n(j)
 :=\bigcup_{a=1}^m\overline{\mathcal N_{n,e_a}(j)}
 \;\cup\;
 \bigcup_{a=1}^{m-1}\Phi_{n,v_a,j}(K_{v_a,j}).
\]
Neighboring pieces in this common core and neck decomposition share their boundary cross-sections.
For fixed $j$, strong $W^{1,2}$ convergence to a constant and uniform convergence on each suppressed core give
\[
 \lim_{n\to\infty}
 E\bigl(u_n;\Phi_{n,v_a,j}(K_{v_a,j})\bigr)=0,
 \qquad
 \lim_{n\to\infty}
 \operatorname{Osc}_{\Phi_{n,v_a,j}(K_{v_a,j})}u_n=0.
\]
The estimates in the proof of Lemma~\ref{lem:noneck} apply to each $\mathcal N_{n,e_a}(j)$, since all its artificial ends are truncated at the same logarithmic length $j$.
Continuity at the shared cross-sections and the triangle inequality bound the oscillation on $\Omega_n(j)$ by the sum of the oscillations on its constituent cores and necks.
Adding their Dirichlet energies and using the finite number of pieces gives
\[
 \limsup_{n\to\infty}E(u_n;\Omega_n(j))\leq Ce^{-\gamma j},
 \qquad
 \limsup_{n\to\infty}\operatorname{Osc}_{\Omega_n(j)}u_n
 \leq Ce^{-\gamma j/2}.
\]
For a one-sided type~II end, the corresponding union retains the boundary circle and is truncated only at its artificial end; the same argument uses the one-sided neck estimate.
Letting $j\to\infty$ proves vanishing Dirichlet energy and oscillation along every concatenated edge.

\repair{Constant terminal profiles and constant profiles carrying concentration marks are also retained when needed to record attaching nodes.
A retained constant component may record a concentration point, a junction, or a bordered node type.
The finite graph and its permitted contractions give a finite connected bordered stratified surface with the component and nodal incidence relations in Definition~\ref{def:free-boundary-bubble tree-surface}.}

The local graphs are connected and acyclic by Proposition~\ref{prop:finite-profile-selection}(v).
The permitted contractions preserve these properties, the boundary circles of type~II collars, and unresolved one-branch markings, as required in Definition~\ref{def:free-boundary-bubble tree-surface}.
Principal branch vertices are identified only in the global graph, preserving the cycles of the original nodal surface.

For a bordered stratified surface, let $g(\Sigma_\infty)$ denote the sum of the genera of the components of its normalization.
Pinching a separating closed curve preserves this sum, while pinching a nonseparating closed curve decreases it by one.
Pinching a properly embedded boundary arc or a boundary component cannot increase this sum.
Sphere and disk bubbles have genus zero, and attaching them at points does not change the sum.
Hence every bubble tree surface constructed here satisfies $g(\Sigma_\infty)\leq g(\Sigma)$.

\Needspace{6\baselineskip}
\subsection{\texorpdfstring{Area measures and completion of the proof}{Area measures and completion of the proof}}
\label{subsec:general-completion}

We prove the area identity, convergence of pushforward area measures, and Willmore lower semicontinuity, completing Theorem~\ref{thm:intro-main-compactness}.

\begin{prop}\label{prop:area}
 
Consider the finite bubble tree decomposition above, with pairwise disjoint smoothly embedded component cores.
On principal cores, use the smoothly converging metrics $\widetilde h_n$ constructed in Subsection~\ref{subsec:generic-extraction}; on nonprincipal cores, use fixed smooth metrics on the corresponding components $\Sigma_\infty^i$ compatible with their conformal structures.
Assume strong $W^{1,2}$ convergence on every fixed core and vanishing Dirichlet energy on the complementary neck regions in the iterated limit of Definition~\ref{def:no-neck-properties}.
Then
\begin{equation*}
\lim_{n\to\infty}\cA(u_n)
=\cA(u_\infty)
=\sum_{i\in\mathfrak I_*}
\cA(u_\infty^i).
\end{equation*}
For every continuous function $\varphi$ on the fixed compact target set $Q$,
\begin{equation}\label{eq:weighted-area-measure-identity}
\lim_{n\to\infty}
\int_\Sigma\varphi(u_n)\,d\mu_{u_n}
=\sum_{i\in\mathfrak I_*}
\int_{\Sigma_\infty^i}
\varphi(u_\infty^i)\,d\mu_{u_\infty^i}.
\end{equation}
Equivalently,
\[
(u_n)_\#\mu_{u_n}
\rightharpoonup
\sum_{i\in\mathfrak I_*}
(u_\infty^i)_\#
\mu_{u_\infty^i}
\]
as finite Radon measures.
Constant ghost components contribute no area.
 
\end{prop}

\begin{proof}
 
By \eqref{eq:sec2-energy-area},
\begin{equation}\label{eq:energy-equals-area-density}
\frac12|du|_h^2\,d\vol_h=d\mu_u,
\qquad E_h(u;\mathfrak{A})=\cA(u;\mathfrak{A})
\end{equation}
for every measurable $\mathfrak{A}$.
The branch and special points have zero area measure.

Fix compact cores $K_{i,j}$ of the components and their embeddings $\Phi_{n,i,j}$.
For a fixed pair $(i,j)$, write $v_n=u_n\circ\Phi_{n,i,j}$.
On principal cores $K_{i,j}$, use the smoothly converging metrics $\widetilde h_n$ constructed in Subsection~\ref{subsec:generic-extraction}.
On other cores, use a fixed smooth metric on $\Sigma_\infty^i$ compatible with its conformal structure.
Strong $W^{1,2}$ convergence and smooth metric convergence give strong $L^1$ convergence of the energy densities and hence, by \eqref{eq:energy-equals-area-density}, of the area densities.

For $\varphi\in C(Q)$, the functions $\varphi(v_n)$ are uniformly bounded and converge in measure to $\varphi(u_\infty^i)$.
Strong $L^1$ convergence of the densities and bounded convergence in measure against the limiting density give
\[
\int_{K_{i,j}}\varphi(v_n)\,d\mu_{v_n}
\longrightarrow
\int_{K_{i,j}}\varphi(u_\infty^i)\,d\mu_{u_\infty^i}.
\]

Let $\mathcal N_{n,j}$ be the complement of the disjoint embedded cores.
Its weighted area is bounded by
\[
\left|\int_{\mathcal N_{n,j}}\varphi(u_n)\,d\mu_{u_n}\right|
\leq \|\varphi\|_\infty E_{h_n}(u_n;\mathcal N_{n,j}),
\qquad
\lim_{j\to\infty}\limsup_{n\to\infty}
E_{h_n}(u_n;\mathcal N_{n,j})=0.
\]
The limiting tails have weighted area at most $\|\varphi\|_\infty E(u_\infty^i;\Sigma_\infty^i\setminus K_{i,j})$, which tends to zero by finite energy.
Sum over the finitely many components, then let $n\to\infty$ followed by $j\to\infty$.
This proves \eqref{eq:weighted-area-measure-identity}.
Taking $\varphi=1$ proves the area and energy identities.
 
\end{proof}

\begin{proof}[Proof of Theorem~\ref{thm:intro-main-compactness}]
Take the disjoint embeddings of component cores from Proposition~\ref{prop:finite-profile-selection}, with the metrics used in the proof of Proposition~\ref{prop:area}.
Choose the exhaustion indices $j(n)\to\infty$ sufficiently slowly by a diagonal argument to obtain the embeddings and domains in Definition~\ref{def:componentwise-free-boundary-bubble tree-convergence}.
Lemma~\ref{lem:local-convergence}, puncture removability, and the argument at auxiliary marks prove parts (i)-(ii).
On each constant component, the same local lemma gives the weak $W^{2,2}$ convergence required by Definition~\ref{def:componentwise-free-boundary-bubble tree-convergence} away from its finite exceptional set.

Strong $W^{1,2}$ convergence on fixed cores includes the subthreshold marks and collapsed cores.
Together with Lemma~\ref{lem:noneck}, this verifies the hypotheses of Proposition~\ref{prop:area}.
We obtain the Dirichlet energy and area identities and the convergence of pushforward area measures with continuous weights on the target.

Use the charts and weak--strong multiplication argument of Subsection~\ref{subsec:compactification}.
The upper and lower conformal factor bounds on noncollapsed cores give weak $L^2$ convergence of $e^{-\lambda_n}\Delta u_n$ to $e^{-\lambda_\infty}\Delta u_\infty$.
A partition of unity transfers lower semicontinuity to the fixed cores.
Summing over the disjoint cores and then exhausting the components gives
\[
\sum_{i\in\mathfrak I_*}
W_{\R^K}(u_\infty^i)
\leq
\liminf_{n\to\infty}W_{\R^K}(u_n).
\]
Uniform convergence on the finitely many cores and the oscillation estimate in Lemma~\ref{lem:noneck} ensure that all attachment values agree.
The union of the component images is therefore connected.
The original images converge to this union in Hausdorff distance.

If every core collapses, the same finite decomposition into cores and necks gives $\operatorname{Osc}_{\Sigma}u_n\to0$.
By compactness of $Q$, a subsequence then collapses to a point.
\end{proof}

The theorem preserves total area and records the parametrized components and their nodal incidence relations.
The disk and annulus applications in Section~\ref{sec:applications} also exclude splitting and degeneration to obtain minimizers with the prescribed topology and constraints.

\vspace{1cm}

\section{\texorpdfstring{Hamiltonian stationary Lagrangian immersions with free boundary}{Hamiltonian stationary Lagrangian immersions with free boundary}}
\label{sec:HSL-immersions}

We prove the first-variation characterization in Theorem~\ref{thm:1st_vari_natural_support} and the boundary estimates for Theorem~\ref{main theorem}.

\Needspace{6\baselineskip}
\subsection{\texorpdfstring{The HSL system and first variation}{The HSL system and first variation}}
\label{subsec:hsl-system}

Definition~\ref{def:hsl-variation-classes} specifies the variation classes used in Theorem~\ref{thm:1st_vari_natural_support}.
Corollary~\ref{coro:hsl-support-geometries} treats special supports.
We then write the homogeneous HSL system~\eqref{sys1}--\eqref{sys2} used below.

Let $(M,\omega,\bar g,J)$ be a $2n$-dimensional symplectic manifold with a compatible almost complex structure $J$ and Riemannian metric $\bar g$ on $M$, so that
\[
\bar g(X,Y)=\omega(X,JY),\qquad \omega(JX,JY)=\omega(X,Y).
\]
Let $\cK\subset M$ be a smooth embedded real submanifold of dimension $1\leq\mathbf k\leq2n-1$, and let $\Sigma$ be a compact smooth $n$-dimensional manifold with nonempty smooth boundary.
An immersion $u:\Sigma\to M$ is Lagrangian if $u^*\omega=0$; equivalently, $Jdu(T\Sigma)$ is orthogonal to $du(T\Sigma)$.
We consider Lagrangian immersions $u$ satisfying the boundary constraint $u(\partial\Sigma)\subset\cK$.

A smooth ambient field $X$ is Lagrangian if the contracted $1$-form $\omega\contraction X$ is closed.
If $X$ is also tangent to $\cK$ along all of $\cK$, Cartan's formula shows that its flow preserves the Lagrangian condition and the support constraint.
For $f\in C_c^\infty(M)$, let $\nabla f$ denote its gradient with respect to $\bar g$.
The Hamiltonian field $J\nabla f$ is \emph{infinitesimally admissible} for $u$ when
\[
J\nabla f(u(x))\in T_{u(x)}\cK, \qquad x\in\partial\Sigma.
\]
Tangency along the initial boundary image need not preserve the constraint for nonzero time.

\begin{definition}\label{def:hsl-variation-classes}
The immersion $u$ is \emph{ambient Hamiltonian stationary}, or HSL in the ambient sense, if
\[
\left.\frac{d}{dt}\right|_{t=0}
\cA\bigl(\phi_{J\nabla f}(t)\circ u\bigr)=0
\]
for every $f\in C_c^\infty(M)$ for which $J\nabla f$ is infinitesimally admissible.

An \emph{exact infinitesimal variation on the source} is a smooth section $V\in\Gamma(u^*TM)$ satisfying
\begin{equation}\label{eq:hsl-source-exact-field}
\omega(V,du(X))=-d\varphi(X),\qquad X\in T\Sigma,
\end{equation}
for one globally defined $\varphi\in C^\infty(\Sigma,\R)$.
It is admissible when $V(x)\in T_{u(x)}\cK$ on $\partial\Sigma$.
We say that $u$ is \emph{stationary under exact infinitesimal source variations} if $\delta\cA(u)(V)=0$ for every such admissible field.
\end{definition}

Every admissible ambient Hamiltonian test gives an admissible exact source field with $\varphi=f\circ u$.
Here $\cA$ and $\mu_u$ denote the induced volume and its measure in dimension $n$, and area and its measure when $n=2$. We use the trace mean curvature convention and the definition of $\alpha$ in \eqref{eq:intro-geometric-mean-curvature-form} in every dimension.
An exact source potential need not be an ambient pullback.
Both notions in Definition~\ref{def:hsl-variation-classes} concern initial velocities.
We do not assert that every admissible source field is realized by a constrained Lagrangian family.

\begin{theorem}\label{thm:1st_vari_natural_support}
Let $(M^4,\omega,J,\bar g)$ be a K\"ahler surface, and let $\cK\subset M$ be a smooth embedded supporting submanifold.
Let $\Sigma$ be a compact oriented surface with nonempty smooth boundary, and let $u:\Sigma\to M$ be a smooth Lagrangian immersion satisfying $u(\partial\Sigma)\subset\cK$.
Use the induced metric $g_u$ and trace mean curvature vector $\mathbf H_{u,M}$ from Subsection~\ref{sec:fbbci}, and let $\alpha$ be the mean curvature $1$-form defined in \eqref{eq:intro-geometric-mean-curvature-form}.
Let $\nu$ and $\tau$ be the outward unit conormal and unit boundary tangent vectors with respect to the metric $g_u$, chosen so that $(\nu,\tau)$ is positively oriented. By \eqref{eq:induced-boundary-geometry}, the image vectors satisfy $\eta_u=du(\nu)$ and $T=du(\tau)$.
Then $u$ is stationary under exact infinitesimal source variations if and only if there exists a single-valued $q\in C^\infty(\partial\Sigma,\R)$, periodic on every boundary circle, such that
\begin{equation}\label{eq:hsl-coupled-natural-boundary}
\begin{gathered}
d_{g_u}^*\alpha=0\quad\text{in }\Sigma,\\
\eta_u+qJT\perp T_u\cK,\qquad
\alpha(\nu)+\partial_\tau q=0\quad\text{on }\partial\Sigma.
\end{gathered}
\end{equation}
This source equivalence holds without any boundary embedding assumption.
If the entire boundary map $u|_{\partial\Sigma}$ is an embedding, then \eqref{eq:hsl-coupled-natural-boundary} is also equivalent to ambient Hamiltonian stationarity in Definition~\ref{def:hsl-variation-classes}, with interior self-intersections allowed.

More generally, let $(M^{2n},\omega,J,\bar g)$ be a K\"ahler manifold, and let $\cK\subset M$ be a smooth embedded supporting submanifold.
Let $\Sigma$ be a compact smooth $n$-dimensional manifold with nonempty smooth boundary, and let $u:\Sigma\to M$ be a smooth Lagrangian immersion satisfying $u(\partial\Sigma)\subset\cK$.
Let $\nu$ be the outward unit conormal with respect to the induced metric $g_u$ and set $\eta_u=du(\nu)$. Then the homogeneous conditions
\begin{equation}\label{eq:HSL-real-hypersurface-support}
 d_{g_u}^*\alpha=0\quad\text{in }\Sigma,\qquad
 \alpha(\nu)=0,\qquad du(\nu)\perp T_u\cK
 \quad\text{on }\partial\Sigma
\end{equation}
are sufficient for stationarity under both classes of infinitesimal tests.
In particular, the volume of $u$ is stationary under compactly supported Hamiltonian flows preserving the support constraint.
\end{theorem}

\begin{proof}
Since $u$ is Lagrangian, \eqref{eq:hsl-source-exact-field} is equivalent to
\[
V=Jdu(\nabla_{g_u}\varphi)+du(Y)
\]
for a smooth tangent field $Y$ on $\Sigma$.
Equation~\eqref{eq:hsl-source-exact-field} fixes the normal component of $V$ and leaves its tangential component free.
Let $\dd\vol_{g_u|_{T\partial\Sigma}}$ be the induced boundary volume measure, equal to $ds_{g_u}$ when $n=2$. The identity $J\mathbf H_{u,M}=du(\alpha_{g_u}^\sharp)$ and the first variation formula give, in every dimension,
\begin{equation}\label{eq:hsl-exact-first-variation}
\begin{aligned}
\delta\cA(u)(V)
&=-\int_\Sigma\langle\mathbf H_{u,M},V\rangle\,d\mu_u
  +\int_{\partial\Sigma}\langle V,\eta_u\rangle\,\dd\vol_{g_u|_{T\partial\Sigma}}\\
&=\int_\Sigma\langle\alpha,d\varphi\rangle_{g_u}\,d\mu_u
  +\int_{\partial\Sigma}\langle V,\eta_u\rangle\,\dd\vol_{g_u|_{T\partial\Sigma}}\\
&=\int_\Sigma\varphi d_{g_u}^*\alpha\,d\mu_u
  +\int_{\partial\Sigma}
   \bigl(\varphi\alpha(\nu)+\langle V,\eta_u\rangle\bigr)\,\dd\vol_{g_u|_{T\partial\Sigma}}.
\end{aligned}
\end{equation}
For an ambient Hamiltonian field, set $\varphi=f\circ u$ and choose a frame $e_1,\ldots,e_n$ on $\Sigma$ orthonormal with respect to the induced metric $g_u$.
Then
\begin{equation}\label{eq:first-variation-general-natural-support-2}
\langle\mathbf H_{u,M},J\nabla f\rangle
=\sum_{i=1}^n e_i(\varphi)\langle\mathbf H_{u,M},Jdu(e_i)\rangle
=-\langle\alpha,d\varphi\rangle_{g_u}.
\end{equation}
On a surface, the boundary part of \eqref{eq:hsl-source-exact-field} reads
\begin{equation}\label{eq:hsl-boundary-potential-derivative}
\partial_\tau\varphi=\langle V,JT\rangle.
\end{equation}
We first derive the interior equation.
For source variations, taking $\varphi$ compactly supported in the interior and $V=Jdu(\nabla_{g_u}\varphi)$ in \eqref{eq:hsl-exact-first-variation} gives $d_{g_u}^*\alpha=0$.
For ambient stationarity, put $c=d_{g_u}^*\alpha$ and $C=u(\partial\Sigma)$.
Every $f\in C_c^\infty(M\setminus C)$ is admissible, so
\[
\int_\Sigma(f\circ u)c\,d\mu_u=0,
\qquad
u_\#(c\,d\mu_u)=0\quad\text{on }M\setminus C.
\]
Fix $p\in u(\operatorname{int}\Sigma)\setminus C$.
The compact discrete fiber $u^{-1}(p)$ is finite.
Choose disjoint compact coordinate disks $D_1,\ldots,D_m\Subset\operatorname{int}\Sigma$, with these preimages in their interiors, such that $u$ is an embedding on a neighborhood of each disk.
The compact image of $\Sigma\setminus\bigcup_i\operatorname{int}D_i$ misses $p$.
Hence an ambient neighborhood $O\Subset M\setminus C$ of $p$ can be chosen so that
\[
u^{-1}(O)\subset\bigcup_i\operatorname{int}D_i,
\qquad S_i:=O\cap u(D_i)
\]
are embedded surfaces relatively closed in $O$, without artificial disk boundaries in $O$.
For fixed $i$, each $\mathfrak A_{ij}=S_i\cap S_j$ is relatively closed in $S_i$, and its relative boundary has empty interior.
Consequently,
\[
G_i=S_i\setminus\bigcup_{j=1}^m\partial_{S_i}\mathfrak A_{ij}
\]
is open and dense in $S_i$.
Near $y\in G_i$, each other sheet is either absent or contains an open neighborhood of $y$ in $S_i$.
In the latter case, inclusion into $S_j$ is a local diffeomorphism between embedded surfaces of the same dimension.
Thus all nonempty sheets coincide in a smaller ambient neighborhood.
Their induced metrics, trace mean curvatures, mean curvature forms, and scalars $c$ agree under reparametrization.
The restriction of the zero pushforward measure there is thus
\[
k\,c_i\,\mathcal H^2_{\bar g}|_{S_i},
\qquad c_i=c\circ(u|_{D_i})^{-1},\qquad k\geq1.
\]
Since $k\geq1$, we have $c_i=0$ on $G_i$ and, by continuity, on $S_i$.
No zero-area assumption on sheet intersections is needed.
Finally, $C$ has zero two-dimensional Hausdorff measure, so it contains no relatively open subset of any embedded local sheet.
The set $u^{-1}(M\setminus C)$ is therefore dense in the interior, and continuity proves $c=0$ everywhere.

A smooth field $W\in\Gamma(u|_{\partial\Sigma}^*T\cK)$ extends to an admissible exact source field if and only if
\begin{equation}\label{eq:hsl-boundary-moment}
\int_\Gamma\langle W,JT\rangle\,ds_{g_u}=0
\qquad\text{for each boundary component }\Gamma.
\end{equation}
Necessity follows from \eqref{eq:hsl-boundary-potential-derivative}.
Conversely, \eqref{eq:hsl-boundary-moment} gives a periodic function $\varphi_\Gamma$ on each circle with
\[
\partial_\tau\varphi_\Gamma=\langle W,JT\rangle.
\]
Its additive constant may be chosen independently on each component.
Prescribe also $\partial_\nu\varphi=\langle W,J\eta_u\rangle$.
In a geodesic inward collar $(p,r)$, the function
\[
\varphi(p,r)=\varphi_\Gamma(p)-r\langle W(p),J\eta_u(p)\rangle
\]
near $r=0$, multiplied by a cutoff supported in the collar, realizes both prescriptions.
The disjoint collars and extension by zero give a smooth global potential.
In the frame $\eta_u,T,J\eta_u,JT$, its field $Jdu(\nabla_{g_u}\varphi)$ is the normal component of $W$.
Extend $du^{-1}(W-Jdu(\nabla_{g_u}\varphi))$ from the boundary to a smooth tangent field $Y$ on $\Sigma$.
Then $Jdu(\nabla_{g_u}\varphi)+du(Y)$ realizes $W$.

When the entire boundary map is embedded, ambient potentials realize these boundary fields and independent constants.
Regard them as single-valued data on the compact embedded curve $C$ and put
\[
A=-JW,\qquad B=A-\langle A,T\rangle T.
\]
The first jet $f|_C=\varphi_C$, $\nabla f|_C=A$ is compatible because
\[
d\varphi_C(T)=\langle W,JT\rangle=\langle A,T\rangle.
\]
In a tubular neighborhood of $C$, write a point as $\exp_p\xi$, with $\xi\in(T_pC)^{\perp_{\bar g}}$, and define
\[
f(\exp_p\xi)=\chi(|\xi|^2)
\bigl(\varphi_C(p)+\langle B(p),\xi\rangle\bigr),
\]
where $\chi$ equals one near zero and has support strictly inside the tubular neighborhood.
Extension by zero gives $f\in C_c^\infty(M)$ with the prescribed value and gradient, and hence $J\nabla f=W$ along $C$.
Interior sheets meeting its neighborhood contribute nothing, since $d_{g_u}^*\alpha=0$ is already proved.

Assume now either source stationarity or ambient stationarity with embedded entire boundary.
On one boundary component $\Gamma$, write $a=\alpha(\nu)$.
Take $W=0$ and prescribe the potential to be one on $\Gamma$ and zero on the other boundary components.
The preceding constructions and \eqref{eq:hsl-exact-first-variation} give
\[
\int_\Gamma a\,ds_{g_u}=0.
\]
There is therefore a smooth periodic $q_0$ on $\Gamma$ with $\partial_\tau q_0=-a$.
For a tangent support field $W$ on $\Gamma$ satisfying \eqref{eq:hsl-boundary-moment}, extended by zero on the other components, integration by parts yields
\[
0=\int_\Gamma\bigl(\varphi_\Gamma a+\langle W,\eta_u\rangle\bigr)\,ds_{g_u}
 =\int_\Gamma\langle W,\eta_u+q_0JT\rangle\,ds_{g_u}.
\]
Define linear functionals on all smooth sections of $u|_\Gamma^*T\cK$ by
\[
L(W)=\int_\Gamma\langle W,JT\rangle\,ds_{g_u},
\qquad
F(W)=\int_\Gamma\langle W,\eta_u+q_0JT\rangle\,ds_{g_u}.
\]
Then $F$ vanishes on $\ker L$.
If $L\ne0$, choose $W_0$ with $L(W_0)=1$.
Since $W-L(W)W_0\in\ker L$, we have $F(W)=F(W_0)L(W)$ for every $W$.
Set $q=q_0-F(W_0)$.
The integral of $\langle W,\eta_u+qJT\rangle$ is zero for every tangent support field $W$.
Taking $W=\operatorname{pr}_{T_u\cK}(\eta_u+qJT)$ proves $\eta_u+qJT\perp T_u\cK$ pointwise, while $\partial_\tau q=-a$ is unchanged.
If $L=0$, taking $W=\operatorname{pr}_{T_u\cK}(JT)$ shows that $JT\perp T_u\cK$.
Every $W$ now lies in $\ker L$, so $F(W)=\int_\Gamma\langle W,\eta_u\rangle\,ds_{g_u}=0$ for all $W$, and $\eta_u\perp T_u\cK$.
In this case $q=q_0$ satisfies both equations.
Applying this construction on each boundary component proves necessity in \eqref{eq:hsl-coupled-natural-boundary}, including where the coefficient of $q$ vanishes.

Conversely, suppose \eqref{eq:hsl-coupled-natural-boundary} holds.
For every admissible exact field, \eqref{eq:hsl-boundary-potential-derivative} gives
\[
\langle V,\eta_u\rangle=-q\langle V,JT\rangle=-q\partial_\tau\varphi.
\]
Thus \eqref{eq:hsl-exact-first-variation} becomes
\[
\delta\cA(u)(V)
=\int_{\partial\Sigma}\bigl(\varphi\alpha(\nu)-q\partial_\tau\varphi\bigr)\,ds_{g_u}
=\int_{\partial\Sigma}\varphi\bigl(\alpha(\nu)+\partial_\tau q\bigr)\,ds_{g_u}=0.
\]
Every ambient test is an exact source test, so sufficiency needs no boundary embedding assumption.
In arbitrary dimension, \eqref{eq:HSL-real-hypersurface-support} makes all three terms of \eqref{eq:hsl-exact-first-variation} vanish.
\end{proof}

\begin{coro}\label{coro:hsl-support-geometries}
Let $M$, $\cK$, $\Sigma$, and $u$ satisfy the surface hypotheses of Theorem~\ref{thm:1st_vari_natural_support}, and retain the notation $g_u$, $\alpha$, $\nu$, $\tau$, $\eta_u$, and $T$.
Below, stationarity is taken with respect to all admissible exact infinitesimal source variations.
Each equivalence also holds for ambient Hamiltonian stationarity if the entire boundary map $u|_{\partial\Sigma}$ is an embedding; interior self-intersections are allowed.
\begin{enumerate}
\item If $\cK$ is a complex curve, then $u$ is stationary if and only if $d_{g_u}^*\alpha=0$ and $\alpha(\nu)=0$.
Orthogonal contact is automatic, and the multiplier is $q=0$.
\item If $\cK$ is a real symplectic surface, choose locally a unit $v\in T_u\cK$ perpendicular to $T$.
Then $u$ is stationary if and only if $d_{g_u}^*\alpha=0$ and
\begin{equation}\label{eq:hsl-symplectic-support}
q=-\frac{\langle\eta_u,v\rangle}{\omega(T,v)},
\qquad \alpha(\nu)+\partial_\tau q=0.
\end{equation}
The denominator is nonzero, and the quotient is independent of the sign of $v$.
Symplecticity does not in general force $q=0$.
\item If $\cK$ is Lagrangian, then $u$ is stationary if and only if
\begin{equation}\label{eq:hsl-lagrangian-support-flux}
d_{g_u}^*\alpha=0,\qquad \eta_u\perp T_u\cK,\qquad
\int_\Gamma\alpha(\nu)\,ds_{g_u}=0
\end{equation}
for every boundary component $\Gamma$. If in addition $\cK$ is minimal and $\alpha=d\psi$ for a globally defined smooth real function $\psi$ on $\Sigma$, then stationarity, the homogeneous system \eqref{eq:HSL-real-hypersurface-support}, and the conditions $\mathbf H_{u,M}=0$ and $\eta_u\perp T_u\cK$ are equivalent.
If the Lagrangian support is minimal and the target is K\"ahler--Einstein, stationarity implies exactness on every connected genus-zero surface, so no separate exactness assumption is needed in this case.
\item If $\cK$ is a real hypersurface with local unit normal $n_{\cK}$, then $u$ is stationary if and only if $d_{g_u}^*\alpha=0$ and
\begin{equation}\label{eq:hsl-hypersurface-coupled}
\langle n_{\cK},\eta_u\rangle\ne0,\qquad
\langle n_{\cK},J\eta_u\rangle=0,\qquad
\alpha(\nu)+\partial_\tau\!\left(
\frac{\langle n_{\cK},JT\rangle}{\langle n_{\cK},\eta_u\rangle}\right)=0.
\end{equation}
Here $q=\langle n_{\cK},JT\rangle/\langle n_{\cK},\eta_u\rangle$ is independent of the sign of $n_{\cK}$.
If the initial immersion also satisfies $\langle T,Jn_{\cK}\rangle=0$, stationarity is equivalent to \eqref{eq:HSL-real-hypersurface-support}.
When $\ker\bigl(\bar g(Jn_{\cK},\cdot)|_{T\cK}\bigr)$ is a contact distribution, the condition $\langle T,Jn_{\cK}\rangle=0$ means that the boundary curve is Legendrian.
\item More generally, if the initial boundary map satisfies
\begin{equation}\label{eq:hsl-initial-complex-tangency}
Jdu(T\partial\Sigma)\subset T_u\cK,
\end{equation}
then stationarity is equivalent to \eqref{eq:HSL-real-hypersurface-support}.
Both \eqref{eq:hsl-initial-complex-tangency} and the additional condition in item~\textup{(4)} are hypotheses on the initial immersion; admissibility of the test fields is still defined by Definition~\ref{def:hsl-variation-classes} and does not require preservation of these additional conditions.
\end{enumerate}
\end{coro}

\begin{proof}
For a complex curve, $T_u\cK=\operatorname{span}\{T,JT\}$.
Since $\eta_u$ is orthogonal to $T$ and, by the Lagrangian condition, to $JT$, contact is orthogonal.
Pairing $\eta_u+qJT$ with $JT$ gives $q=0$, so Theorem~\ref{thm:1st_vari_natural_support} proves item~\textup{(1)}.
For a symplectic surface, $T,v$ span $T_u\cK$ and $\omega(T,v)=\langle JT,v\rangle\ne0$.
The vector $\eta_u+qJT$ is orthogonal to $T$, and its pairing with $v$ vanishes exactly for the quotient in \eqref{eq:hsl-symplectic-support}.
The quotient is unchanged under $v\mapsto-v$, so it defines a smooth global $q$ and proves item~\textup{(2)}.

If $\cK$ is Lagrangian, $JT\perp T_u\cK$.
The geometric equation in \eqref{eq:hsl-coupled-natural-boundary} therefore reduces to $\eta_u\perp T_u\cK$.
A periodic solution of $\partial_\tau q=-\alpha(\nu)$ exists on each boundary circle exactly when its flux integral is zero, proving \eqref{eq:hsl-lagrangian-support-flux}.
At orthogonal contact we also have
\begin{equation}\label{eq:hsl-lagrangian-boundary-form}
\alpha(X)=\alpha_{\cK}(du(X)),\qquad X\in T\partial\Sigma,
\end{equation}
where $\alpha_{\cK}(Z)=\bar g(J\mathbf H_{\cK},Z)$ for $Z\in T\cK$, and $\mathbf H_{\cK}$ is the trace mean curvature vector of $\cK\subset M$.
Indeed, the tangent frames of the immersion and support along the boundary are $(T,\eta_u)$ and $(T,J\eta_u)$, respectively.
For a Lagrangian submanifold $L\subset M$, the cubic tensor $(X,Y,Z)\mapsto\langle A_{L,M}(X,Y),JZ\rangle$ is symmetric.
Consequently, writing $\widehat X=du(X)$, we have
\[
\alpha(X)
=-\langle\overline\nabla_{\widehat X}T,JT\rangle
 -\langle\overline\nabla_{\widehat X}\eta_u,J\eta_u\rangle.
\]
For $\alpha_{\cK}(\widehat X)$, the first term is the same along the common boundary curve, and the second agrees because $\overline\nabla J=0$ gives
\[
\langle\overline\nabla_{\widehat X}(J\eta_u),J(J\eta_u)\rangle=\langle\overline\nabla_{\widehat X}\eta_u,J\eta_u\rangle.
\]
This proves \eqref{eq:hsl-lagrangian-boundary-form}.
If $\cK$ is minimal, stationarity and this identity imply $\alpha(\tau)=0$.
When $\alpha=d\psi$, the function $\psi$ is therefore constant, say $c_\Gamma$, on each boundary component.
Using coclosedness and the componentwise zero flux, we obtain
\[
\int_\Sigma|\alpha|_{g_u}^2\,d\mu_u=\int_\Sigma\langle\alpha,d\psi\rangle_{g_u}\,d\mu_u=\sum_\Gamma c_\Gamma\int_\Gamma\alpha(\nu)\,ds_{g_u}=0.
\]
Hence $\alpha=0$ and $\mathbf H_{u,M}=-Jdu(\alpha_{g_u}^\sharp)=0$.
Conversely, minimality and orthogonal contact imply stationarity by the first variation.
The same integration by parts shows that the homogeneous system and exactness imply $\alpha=0$, proving all three equivalences.
In the K\"ahler--Einstein case, $d\alpha=u^*\operatorname{Rc}=0$.
For a connected genus-zero surface the boundary cycles generate $H_1(\Sigma;\R)$; \eqref{eq:hsl-lagrangian-boundary-form} makes their periods zero.
Thus $\alpha$ is exact and the preceding argument applies.
For positive genus, handle periods require an additional hypothesis.
When the immersion is minimal, any componentwise constant $q$, including zero, solves the Lagrangian-support boundary system.

For a hypersurface, the geometric equation is equivalent to
\[
n_{\cK}=\sigma\frac{\eta_u+qJT}{\sqrt{1+q^2}},\qquad \sigma\in\{1,-1\}.
\]
It gives $\langle n_{\cK},\eta_u\rangle\ne0$, $\langle n_{\cK},J\eta_u\rangle=0$, and the stated quotient for $q$.
Conversely, $n_{\cK}\perp T$ together with $n_{\cK}\perp J\eta_u$ puts $n_{\cK}$ in $\operatorname{span}\{\eta_u,JT\}$; the nonzero coefficient of $\eta_u$ recovers the vector $\eta_u+qJT$ up to a nonzero scalar factor and proves \eqref{eq:hsl-hypersurface-coupled}.
The condition $\langle T,Jn_{\cK}\rangle=0$ is equivalent to $\langle n_{\cK},JT\rangle=0$, and hence to $q=0$ under these equations.
This proves item~\textup{(4)}.
Finally, under \eqref{eq:hsl-initial-complex-tangency}, pairing $\eta_u+qJT$ with $JT\in T_u\cK$ gives $q=0$.
Theorem~\ref{thm:1st_vari_natural_support} proves item~\textup{(5)} with the full test class stated in the corollary.
\end{proof}

For Lagrangian supports, stationarity gives zero flux on each boundary component.
The additional minimality and exactness hypotheses also give the pointwise conclusion.
Schoen's discussion of special Lagrangian and complex hypersurface supports uses criticality under relative Lagrangian variations \cite[Section~3]{Schoen2005}; a closed deformation $1$-form need not be exact.
For hypersurfaces in $\C^2$, Gaia obtains the homogeneous characterization under the complex-tangent hypothesis in a localized formulation of stationarity \cite[Theorem~7]{Gaia2025}.
In complex dimension $n$, a coisotropic support of codimension $r$ has complex tangent space $T\cK\cap JT\cK$ of complex dimension $n-r$.
Under \eqref{eq:hsl-initial-complex-tangency} it contains the isotropic space $du(T\partial\Sigma)$ of dimension $n-1$.
Thus $n-1\leq n-r$, so $r\leq1$; higher-codimension coisotropic supports give no nonempty additional examples under this condition.

For an ambient stationary immersion with a nonembedded boundary, the proof still gives $d_{g_u}^*\alpha=0$.
Its remaining variational condition is precisely
\[
\int_{\partial\Sigma}\bigl((f\circ u)\alpha(\nu)-df(J\eta_u)\bigr)\,ds_{g_u}=0
\]
for all initially admissible $f\in C_c^\infty(M)$, since $\langle J\nabla f,\eta_u\rangle=-df(J\eta_u)$.
At a shared boundary image, this condition sums the mean curvature flux and conormal contributions over the sheets.

Embeddedness of each boundary component does not suffice for the ambient equivalence.
Consider the cut-torus example
\[
E=\mathbb C/(10\mathbb Z+10i\mathbb Z),\qquad
M=E\times E,\qquad
\cK=\{[1]\}\times E,
\]
with $\Sigma=[0,2\pi]\times S^1$ and
\[
u(s,\theta)=([e^{is}],[e^{i\theta}]).
\]
The vectors $u_s=(ie^{is},0)$ and $u_\theta=(0,ie^{i\theta})$ are orthogonal unit vectors, so $u$ is conformal, and the product symplectic form vanishes on their span.
The two boundary circles have conormals perpendicular to $T_u\cK$ and have the same image in $\cK$.
The trace convention gives
\[
g_u=ds^2+d\theta^2,\qquad \mathbf H_{u,M}=(-e^{is},-e^{i\theta}),\qquad \alpha=-ds-d\theta.
\]
Thus $d_{g_u}^*\alpha=0$, whereas $\alpha(\nu)=1$ at $s=0$ and $\alpha(\nu)=-1$ at $s=2\pi$.
Nevertheless, $u$ is ambient Hamiltonian stationary.
For $\varphi=f\circ u$,
\[
\langle\mathbf H_{u,M},J\nabla f\rangle=\varphi_s+\varphi_\theta.
\]
The $\theta$-integral vanishes by periodicity, and the $s$-integral vanishes because $\varphi(0,\theta)=\varphi(2\pi,\theta)$.
The geometric boundary term vanishes by orthogonality.
Ambient tests see only the sum of the boundary contributions.
A source potential equal to one near $s=0$ and zero near $s=2\pi$ has zero boundary velocity and first variation $2\pi$.
Thus source stationarity and componentwise zero flux fail.

An embedded disk can also have nonzero coupling.
Up to a source reflection and a unitary coordinate change, the disk below is a restriction of the cylindrical model in \cite[Lemma~11]{Gaia2025}, whose boundary cancellation is verified in the proof of that lemma using equation~(24).
We construct a compact support for this disk and compute its boundary terms.
For $0<r<10^{-2}$, let
\[
\Sigma_r=\{x^2+y^2\leq r^2\},\qquad u(x,y)=(e^{ix},y)\in\C^2.
\]
This is an embedded conformal Lagrangian disk, and
\[
g_u=dx^2+dy^2,\qquad \mathbf H_{u,M}=(-e^{ix},0),\qquad \alpha=-dx,\qquad d_{g_u}^*\alpha=0.
\]
On its boundary, set
\[
\eta_u=\frac{xu_x+yu_y}{r},\qquad T=\frac{-yu_x+xu_y}{r},\qquad q=y.
\]
Then $\alpha(\nu)=-x/r$ and $\partial_\tau q=x/r$.
To construct the support, use the injective ambient chart
\[
z_1=(1+a)e^{ix},\qquad z_2=y+ib
\]
near $(1,0)$ and define
\[
\rho=x^2+y^2+a^2+b^2-r^2+2y^2a+2xyb,\qquad R=(x^2+y^2+a^2+b^2)^{1/2}.
\]
Let $\Omega$ be the image of $\{R<2r,\rho<0\}$ in this chart.
Writing $\rho=R^2-r^2+F_3$, we have $|F_3|\leq4R^3$ and
\[
R\partial_R\rho=2R^2+3F_3\geq2R^2-12R^3>0\quad(0<R\leq2r).
\]
Moreover, $\rho<0$ on $R=r/2$ and $\rho>0$ on $R=2r$.
Each ray meets the zero level once and transversely, so $\cK=\partial\Omega$ is a smooth compact embedded hypersurface in the chart.
The chart is injective on this set because $1+a>0$ and $|x|<\pi$.
Since $\rho\circ u=x^2+y^2-r^2$, the disk lies in $\overline\Omega$ with boundary on $\cK$.
At $a=b=0$ the coordinate frame is orthonormal, $Ju_x=-\partial_a$, and $Ju_y=\partial_b$.
Consequently,
\[
JT=\frac{y\partial_a+x\partial_b}{r},\qquad
\nabla\rho=2xu_x+2yu_y+2y^2\partial_a+2xy\partial_b
=2r(\eta_u+yJT).
\]
The outward support normal is therefore $(\eta_u+yJT)/\sqrt{1+y^2}$, and \eqref{eq:hsl-coupled-natural-boundary} holds.
For an admissible exact field with potential $\varphi$, support tangency gives $\langle V,\eta_u\rangle=-y\partial_\tau\varphi$.
Hence
\[
\delta\cA(u)(V)
=\int_{\partial\Sigma_r}\left(-\varphi\frac{x}{r}-y\partial_\tau\varphi\right)\,ds_{g_u}
=\int_{\partial\Sigma_r}\varphi\left(-\frac{x}{r}+\partial_\tau y\right)\,ds_{g_u}=0.
\]
Thus the embedded disk is stationary for both variation classes, although $\alpha(\nu)\not\equiv0$ and contact is nonorthogonal where $y\ne0$.
For sufficiently small $r$, the support is also strictly convex: the function $F=\rho+r^2$, viewed in ambient Euclidean coordinates near $p=(1,0)$, satisfies $dF(p)=0$ and $D^2_{\mathrm{Eucl}}F(p)=2I$.
Choose a small convex ambient ball $B$ about $p$ on which $D^2_{\mathrm{Eucl}}F$ is positive definite, and then shrink $r$ so that the sublevel $\{F\leq r^2\}\cap B=\overline\Omega$ is compactly contained in $B$.
The second fundamental form of its outward level boundary is $D^2_{\mathrm{Eucl}}F/|\nabla F|$ restricted to the tangent space, and is therefore positive definite.

Every Lagrangian immersion in a K\"ahler manifold satisfies $d\alpha=u^*\operatorname{Rc}$; see \cite{Dazord1981,schoen2001}.
For a smooth Lagrangian immersion of a surface, choose a smooth metric $h$ in the conformal class of $g_u$ and write $g_u=e^{2\lambda}h$.
Then
\[
 \operatorname{tr}_h\nabla du=e^{2\lambda}\mathbf H_{u,M}
 =-Jdu(\alpha_h^\sharp),\qquad
 d_{g_u}^*\alpha=e^{-2\lambda}d_h^*\alpha,
 \qquad \nu_h=e^\lambda\nu,\quad \tau_h=e^\lambda\tau.
\]
Thus both boundary equations in \eqref{eq:hsl-coupled-natural-boundary} are multiplied by $e^\lambda$, giving \eqref{eq:intro-hsl-intrinsic-system}.
The form $\alpha$ remains the geometric mean curvature $1$-form defined using $g_u$.

We now restrict to the homogeneous HSL system $q=0$, with orthogonal contact and zero normal trace.
Its Hodge equations, written using the smooth metric $h$ on $\Sigma$, are
\begin{equation}\label{1form}
    \left\{
    \begin{aligned}
        d\alpha &= u^*\operatorname{Rc} &&\text{in }\Sigma, \\
        d_h^*\alpha &= 0 &&\text{in }\Sigma, \\
        \alpha(\nu_h) &= 0 &&\text{on }\partial\Sigma.
    \end{aligned}
    \right.
\end{equation}
\begin{rmk}
Let $M$ be a complete K\"ahler manifold of complex dimension one and $\cK\subset M$ a compact smooth embedded curve.
Every smooth immersed curve in $M$ is Lagrangian.
For an arclength parametrization $\gamma:[0,L]\to M$, write $\nabla_s\gamma'=\kappa_gJ\gamma'$.
Our convention gives
\[
\alpha=-\kappa_g\,ds,
\qquad
d^*\alpha=\frac{d\kappa_g}{ds}.
\]
Thus $d^*\alpha=0$ means that $\kappa_g$ is constant; see \cite[Example~3.2]{BettiolPiccioneSantoro2014}.
For closed curves, this constant need not vanish.
At an endpoint, $\alpha(\nu)=0$ forces $\kappa_g=0$.
Hence the arcs satisfying this homogeneous HSL system are exactly the geodesic arcs with endpoints on $\cK$ and meeting $\cK$ orthogonally.

Under a uniform length bound, such arcs have a subsequence converging in $C^\infty([0,1])$ to a free boundary geodesic arc or a constant map after constant speed parametrization. Indeed, the initial points lie in $\cK$ and the initial speeds are bounded, so the claim follows from smooth dependence for the geodesic equation.
\end{rmk}

Let $(M,\omega,J,\bar g)$ be a complete K\"ahler surface, and let $(\Sigma,h)$ be a Riemann surface with nonempty boundary.
For a smooth branched conformal Lagrangian immersion $u$, write $u^*\bar g=e^{2\lambda}h$ away from its finite branch set $\cB_u$, as in \eqref{eq:sec2-intrinsic-conformality}.
For branched conformal Lagrangian immersions satisfying the homogeneous HSL system, \eqref{1form} holds away from $\cB_u$.
Fix an isometric embedding $M\hookrightarrow\R^K$ and a conformal disk or half disk away from $\cB_u$.
Set $u_i=\partial_{x_i}u$, $|\nabla u|_{\bar g}^2=\sum_i\bar g(u_i,u_i)$, and $\Delta=\partial_{x_1}^2+\partial_{x_2}^2$.
In a half disk, the diameter represents a portion of $\partial\Sigma$ on which the free boundary conditions hold.
The intrinsic mean curvature satisfies
\begin{equation}\label{sysH}
\Delta_{g_u} u-(g_u)^{ij}A^M_{\R^K}(u_i,u_j)=\mathbf H_{u,M},
\end{equation}
where $A^M_{\R^K}$ is the second fundamental form of the fixed embedding, with the notation of Remark~\ref{rem:intro-euclidean-reduction}.
Since $(g_u)_{ij}=\tfrac12|\nabla u|_{\bar g}^2\delta_{ij}$ and the normal bundle in $M$ is spanned by $Ju_1,Ju_2$, write $\mathbf H_{u,M}=H^{x_1}Ju_1+H^{x_2}Ju_2$.
Our sign convention gives
\begin{equation}\label{alphaH}
\alpha_{x_i}=-\tfrac12|\nabla u|_{\bar g}^2H^{x_i},\qquad i=1,2,
\end{equation}
and hence
\begin{equation}\label{intro h}
\tfrac12|\nabla u|_{\bar g}^2\mathbf H_{u,M}=-\sum_{i=1}^2\alpha_{x_i}J(u)u_i.
\end{equation}
Combining \eqref{sysH} and \eqref{alphaH} yields
\begin{equation}\label{lapu}
\Delta u=-\sum_{i=1}^2\alpha_{x_i}J(u)u_i+\sum_{i=1}^2A^M_{\R^K}(u_i,u_i).
\end{equation}
For branched maps, these identities first hold on $\overline\Sigma\setminus\cB_u$.
Using the smooth metric $h$ on $\Sigma$ avoids defining unit vectors for the degenerate tensor $u^*\bar g$ at branch points.
This does not establish variational equivalence at those points.
For a coclosed form $\alpha\in L^2(T^*\Sigma,h)$, the homogeneous weak normal trace $\alpha(\nu_h)=0$ means
\[
\int_\Sigma\langle\alpha,d\psi\rangle_h\,\dd\vol_h=0
\qquad\text{for every }\psi\in C^\infty(\overline\Sigma).
\]
To extend the weak flux across the finite branch set, choose logarithmic cutoffs $\chi_\varepsilon$ in disks and half disks, vanishing near the branch points, tending to one away from them, and satisfying $\|d\chi_\varepsilon\|_{L^2}\to0$.
The homogeneous Green identity away from the branch set gives
\[
0=\int_\Sigma\langle\alpha,d(\chi_\varepsilon\psi)\rangle_h\,\dd\vol_h
\]
for smooth $\psi$.
The cutoff error is bounded by $\|\psi\|_\infty\|\alpha\|_{L^2}\|d\chi_\varepsilon\|_{L^2}$, so passage to the limit gives the weak identity on the full surface.
Together with the Lagrangian condition $u^*\omega=0$, the homogeneous system with respect to the smooth metric $h$ on $\Sigma$ is
\begin{equation}\label{sys1}
\left\{
\begin{aligned}
\Delta_hu&=-J(u)du(\alpha_h^\sharp)
+\operatorname{tr}_h A^M_{\R^K}(u)(du,du)
&&\text{in }\Sigma,\\
u(\partial\Sigma)&\subset\cK,\qquad
du(\nu_h)\perp T_u\cK
&&\text{on }\partial\Sigma.
\end{aligned}
\right.
\end{equation}
\begin{equation}\label{sys2}
\left\{
\begin{aligned}
d\alpha&=u^*\operatorname{Rc}&&\text{in }\Sigma,\\
d_h^*\alpha&=0&&\text{in }\Sigma,\\
\alpha(\nu_h)&=0&&\text{on }\partial\Sigma.
\end{aligned}
\right.
\end{equation}
Here $\Delta_hu$ is the componentwise Laplace--Beltrami operator for the fixed Euclidean embedding of $M$.
These equations initially hold away from $\cB_u$, with boundary traces as above.
We prove their weak extension, regularity, and removability below.
In a conformal coordinate chart, multiplication by the conformal factor of $h$ gives \eqref{lapu}.
In the local estimates below, $\mathbb A*\mathbb B$ denotes a linear combination of contractions.
The source $u^*\operatorname{Rc}$ is quadratic in the first derivatives of $u$.
With the Willmore convention of Subsection~\ref{sec:fbbci}, the Lagrangian condition and conformal invariance of the $L^2$ norm of $1$-forms in dimension two give
\begin{equation}\label{hs:normalization}
\begin{aligned}
W_M(u;U)&=\frac14\int_{U\setminus\cB_u}|\mathbf H_{u,M}|^2\,d\mu_u\\
&=\frac14\int_{U\setminus\cB_u}|\alpha|_{g_u}^2\,d\mu_u\\
&=\frac14\int_{U\setminus\cB_u}|\alpha|_h^2\,\dd\vol_h.
\end{aligned}
\end{equation}
In particular, finite $W_M(u;U)$ makes $\alpha$ an $L^2$ form on $U$, and, if $U$ lies in a conformal coordinate chart,
\[
\int_U\bigl(\alpha_{x_1}^2+\alpha_{x_2}^2\bigr)\,dx=4W_M(u;U).
\]

\Needspace{6\baselineskip}
\subsection{\texorpdfstring{Geodesic reflection for the HSL system}{Geodesic reflection for the HSL system}}
\label{sec:geodesic-reflection}

We reflect the homogeneous HSL system on conformal half disks whose diameters represent portions of $\partial\Sigma$; the formulas are given in Proposition~\ref{extsys}.
Rescaling and finite coverings transfer the estimates to $\Sigma$.

\begin{lemma}\label{lemdis}
Fix an isometric embedding $M\hookrightarrow\R^K$, a compact set $Q\Subset M$, and a compact smooth support $\cK\subset M$.
Let $u:D_1^+\to M$ be a branched conformal HSL immersion with free boundary in the local sense of Definition~\ref{def:fb-branched-conformal}, and suppose that $u(D_1^+)\subset Q$.
Thus its branch set is finite, $u\in W^{2,2}_{\loc}(D_1^+\cup I_1,M)$ across that set, and it has finite area and finite total curvature.
Assume that $u(I_1)\subset\cK$, that the conormal contact is perpendicular, and that the pair $(u,\alpha)$ satisfies
\begin{align*} \begin{cases}
 \Delta u = -\alpha * J(u) * \nabla u
 + A^M_{\R^K}(u)(\nabla u,\nabla u),
 &x\in D_1^+,\\
 u_{x_2}\perp T_u\cK,
 &\text{a.e. on }I_1,
 \end{cases}
 \end{align*}
and
\begin{align*} \begin{cases}
 d\alpha=\operatorname{Rc}(u)*\nabla u*\nabla u,\\
 d^*\alpha=0.
 \end{cases}
 \end{align*}
There are positive constants $C$ and $\epsilon_0$, depending only on the fixed embedding, the geometry on $Q$, and $\cK$, such that if
\[
 E(u;D_1^+)+W_M(u;D_1^+)
 \leq\epsilon_0^2,
 \]
then
\begin{equation}\label{dist}
 \operatorname{dist}(u(x),\cK)
 \leq
 C\left(
 \|\nabla u\|_{L^2(D_1^+)}
 +\|\alpha\|_{L^2(D_1^+)}
 \right)
 \end{equation}
for every $x\in D_{1/2}^+$, and
\begin{equation}\label{osc}
 \operatorname{Osc}^{\R^K}_{D_{1/4}^+}u
 \leq C\left(\|\nabla u\|_{L^2(D_1^+)}
 +\|\alpha\|_{L^2(D_1^+)}\right).
 \end{equation}
 \end{lemma}

The distance in Lemma~\ref{lemdis} is measured in the fixed Euclidean embedding, and $\operatorname{Osc}^{\R^K}$ has the meaning in \eqref{eq:oscillation-definition} after composing $u$ with that embedding.
The finite branch set and finite total curvature in Definition~\ref{def:fb-branched-conformal} are used in the proof.

\begin{proof}
The assumed local regularity gives $u\in W^{2,2}(D_{7/8}^+,M)$.
The comparison \eqref{eq:sec2-willmore-comparison} gives $E+W_{\R^K}\leq C(E+W_M)$.
Choose $\epsilon_0$ so that Theorem~\ref{intsma0} applies with $p=3/2$ on $D_{7/8}^+$.
Sobolev embedding and Morrey's inequality then give
\[
\|\nabla u\|_{W^{1,3/2}(D_{3/4}^+)}+\|\nabla u\|_{L^6(D_{3/4}^+)}+\operatorname{Osc}^{\R^K}_{D_{1/2}^+}u\leq C\|\nabla u\|_{L^2(D_1^+)}.
\]
For $x=(x_1,x_2)\in D_{1/2}^+$, the continuous trace $u(x_1,0)$ lies in $\cK$.
Thus the same oscillation bounds $\operatorname{dist}(u(x),\cK)$, proving \eqref{dist} and \eqref{osc}.
\end{proof}

We extend the tubular reflection to the ambient space so that the image of $u$ need not stay near $\cK$. Choose a Euclidean tubular neighborhood $\mathcal T$ of $\cK$ and a cutoff $\chi\in C_c^\infty(\mathcal T)$ equal to one on a smaller tubular neighborhood.
Define a globally smooth map $S:\R^K\to\R^K$ by
\[
 S(q)=q+2\chi(q)\bigl(\pi_{\cK}(q)-q\bigr)
 \quad\text{on }\mathcal T,
 \qquad
 S(q)=q\quad\text{outside }\mathcal T.
\]
Its first two derivatives are bounded.
On $\cK$, $S$ fixes each point, while $DS$ acts as the identity on $T\cK$ and as minus the identity on $(T\cK)^\perp$.

Let $\rho(x_1,x_2)=(x_1,-x_2)$.
On the lower half disk put $v=u\circ\rho$, and define
\[
 \widehat u(x)=
 \begin{cases}
 u(x),&x_2\geq0,\\
 S(v(x)),&x_2<0,
 \end{cases}
 \qquad
 \widehat\alpha=
 \begin{cases}
 \alpha,&x_2\geq0,\\
 \rho^*\alpha,&x_2<0.
 \end{cases}
\]
Thus
\begin{equation}\label{defhatalpha-pullback}
    \widehat\alpha_{x_1}(x)=\alpha_{x_1}(\rho(x)),
    \qquad
    \widehat\alpha_{x_2}(x)=-\alpha_{x_2}(\rho(x))
    \quad (x_2<0).
\end{equation}
\begin{prop}\label{extsys}
Let $u$ and $\alpha$ be smooth on the closed upper half disk, assume $u(\overline{D_1^+})\subset Q\Subset M$, and suppose distributionally that
\[
\Delta u
=\sum_{j=1}^2A^M_{\R^K}(u)(u_{x_j},u_{x_j})-\sum_{j=1}^2\alpha_{x_j}J(u)u_{x_j},
\qquad
d\alpha=u^*\operatorname{Rc},
\qquad
d^*\alpha=0.
\]
Assume on the diameter that
\[
u\in\cK,\qquad u_{x_2}\perp T_u\cK,\qquad \alpha_{x_2}=0.
\]
Then
\[
\widehat u\in C^1(D_1)\cap W^{2,q}(D_1,\R^K)
\quad\text{for every finite }q,
\]
and $\widehat\alpha$ is continuous and piecewise smooth.
The following distributional equations have no measure-valued terms supported on the reflection line $x_2=0$:
\begin{equation}\label{equation:global-form}
    \Delta\widehat u =\Upsilon_{\widehat u}+\widehat F.
\end{equation}
Here
\begin{align}\label{defE}
\Upsilon_{\widehat u}:=
\begin{cases}
\displaystyle
\sum_{j=1}^2A^M_{\R^K}(u)(u_{x_j},u_{x_j}),
&x_2>0,\\[6pt]
\displaystyle
\sum_{j=1}^2D^2S(v)(v_{x_j},v_{x_j})+DS(v)\sum_{j=1}^2A^M_{\R^K}(v)(v_{x_j},v_{x_j}),
&x_2<0,
\end{cases}
\end{align}
and
\begin{align}\label{defF}
F&:=-\sum_{j=1}^2\alpha_{x_j}J(u)u_{x_j},\\
\label{defF1}
\widehat F
&:=\begin{cases}
\displaystyle
-\sum_{j=1}^2\alpha_{x_j}J(u)u_{x_j},
&x_2>0,\\[6pt]
\displaystyle
-DS(v)\sum_{j=1}^2
\widehat\alpha_{x_j}J(v)v_{x_j},
&x_2<0.
\end{cases}
\end{align}
Moreover,
\begin{align}\label{extalpha}
d\widehat\alpha
&=\begin{cases}
\operatorname{Rc}(u)(u_{x_1},u_{x_2})\,dx_1\wedge dx_2,
&x_2>0,\\
\operatorname{Rc}(v)(v_{x_1},v_{x_2})\,dx_1\wedge dx_2,
&x_2<0,
\end{cases}
&d^*\widehat\alpha&=0.
\end{align}
There is a constant depending only on the fixed geometry such that
\[
 |\Delta\widehat u|
 \leq
 C\bigl(|\nabla v|^2+|\widehat\alpha||\nabla v|\bigr)
 \quad\text{on }D_1^-,
 \qquad
 |d\widehat\alpha|\leq C|\nabla v|^2,
\]
with the analogous estimates on the upper half, and
\[
 \|\nabla\widehat u\|_{L^2(D_1)}
 \leq C\|\nabla u\|_{L^2(D_1^+)},
 \qquad
 \|\widehat\alpha\|_{L^2(D_1)}
 =\sqrt2\,\|\alpha\|_{L^2(D_1^+)}.
\]
\end{prop}

\begin{proof}
The map $S$ is smooth with bounded first and second derivatives.
Differentiating $u(x_1,0)\in\cK$ along the diameter gives $u_{x_1}\in T_u\cK$, while the free boundary condition gives $u_{x_2}\in(T_u\cK)^\perp$.
Therefore
\[
 DS(u)u_{x_1}=u_{x_1},
 \qquad
 DS(u)(-u_{x_2})=u_{x_2}.
\]
The values and first derivatives of the two pieces of $\widehat u$ have matching traces.
Sobolev gluing and the chain rule give the asserted $C^1$ and $W^{2,q}$ regularity.
In particular, its distributional second derivatives contain no interface measure.

On the lower half, $\widehat\alpha$ is given by \eqref{defhatalpha-pullback}.
Its tangential traces agree across the diameter, and both normal traces vanish.
Hence $\widehat\alpha$ is continuous, and neither $d\widehat\alpha$ nor $d^*\widehat\alpha$ contains an interface measure.
Pulling the equations back by the domain reflection gives \eqref{extalpha}.

On the lower half disk, the sign changes in $\widehat\alpha_{x_2}$ and $v_{x_2}$ cancel.
Hence
\[
 \Delta v=\sum_{j=1}^2A^M_{\R^K}(v)(v_{x_j},v_{x_j})-\sum_{j=1}^2\widehat\alpha_{x_j}J(v)v_{x_j}.
\]
The chain rule for $\widehat u=S(v)$ gives
\[
 \Delta\widehat u=\sum_{j=1}^2D^2S(v)(v_{x_j},v_{x_j})+DS(v)\Delta v,
\]
which proves \eqref{equation:global-form}-\eqref{defF1}.
The bounds on the target tensors over $Q$ and on the first two derivatives of $S$ give the pointwise estimates.
The two norm relations follow by reflection and a change of variables.
\end{proof}
 
On a sufficiently small tubular neighborhood of $\cK$ in $M$, the geodesic reflection $\sigma$ is an involution.
If $u(D_1^+)$ lies in this neighborhood, define $U=\sigma(u\circ\rho)$ on the lower half disk.
Then $u\circ\rho=\sigma(U)$, and boundedness of $D\sigma$ gives $|\nabla(u\circ\rho)|\leq C|\nabla U|$.
The chain rule, using a smooth Euclidean extension of $\sigma$, gives
\[
|\Delta U|\leq C\bigl(|\nabla U|^2+|\widehat\alpha||\nabla U|\bigr),
\qquad |d\widehat\alpha|\leq C|\nabla U|^2,
\qquad d^*\widehat\alpha=0.
\]
The same trace matching as in Proposition~\ref{extsys} removes interface measures.
The equations and estimates also hold where $du=0$, since we never divide by $du$.

% Superseded inactive draft removed; retained in the project audit.

We use the global ambient reflection for estimates without restricting the image to this tubular neighborhood and the local involution when inverse derivative bounds are needed.
The global reflected map need not be conformal or $M$-valued.

\Needspace{6\baselineskip}
\subsection{\texorpdfstring{Boundary regularity}{Boundary regularity}}
\label{subsec:hsl-regularity}

Under the homogeneous boundary conditions, Theorem~\ref{higher} gives higher estimates for smooth branched HSL immersions with small Dirichlet energy.
We then prove boundary smoothness for weak solutions with $u\in W^{2,2}$ and $\alpha\in L^2$.

\begin{prop}\label{intsma}
Fix $1<p<2$, $0<r<R$, and $K\in\mathbb N$.
Here $(\widehat u,\widehat\alpha)$ denotes the pair on the full disk obtained from $(u,\alpha)$ on the upper half disk by the global ambient reflection in Proposition~\ref{extsys}.
The image need not lie in a tubular neighborhood of the support.
Suppose that $\widehat u\in W^{2,2}(D_R,\R^K)$ satisfies \eqref{equation:global-form} and \eqref{extalpha}, and that the fixed coefficients in $\Upsilon_{\widehat u}$ and $\widehat F$ are bounded by $C_0$.
There are constants $\varepsilon_1,\varepsilon_2>0$ and $C<\infty$, depending only on $p,r,R,K,C_0$, such that if
\begin{align}\label{intsma1-1}
\|\widehat\alpha\|_{L^2(D_R)}<\varepsilon_1,
\end{align}
and
\begin{align}\label{intsma1-2}
\|\nabla\widehat u\|_{L^2(D_R)}<\varepsilon_2,
\end{align}
then
\begin{align*}
\|\nabla\widehat u\|_{W^{1,p}(D_r)}
\leq C\|\nabla\widehat u\|_{L^2(D_R)}.
\end{align*}
For a conformal Lagrangian map with geometric mean curvature form $\alpha$, the flat coordinate norm in \eqref{intsma1-1} satisfies
\[
\|\widehat\alpha\|_{L^2(D_R)}=\sqrt2\,\|\alpha\|_{L^2(D_R^+)}=2\sqrt{2W_M(u;D_R^+)}.
\]
\end{prop}

\begin{proof}
Choose a radial cutoff $\eta\in C_c^\infty(D_R)$ equal to one on $D_r$.
Put $h=\eta(\widehat u-\overline{\widehat u})$, where $\overline{\widehat u}$ is the average over $D_R$, and set
\[
p^*=\frac{2p}{2-p},\qquad
A=\|\nabla\widehat u\|_{L^2(D_R)}+\|\widehat\alpha\|_{L^2(D_R)}.
\]
Proposition~\ref{extsys}, reflection of the lower-half integrals, and H\"older's inequality give
\[
\|\eta\Delta\widehat u\|_{L^p(D_R)}\leq CA\|\eta\nabla u\|_{L^{p^*}(D_R^+)}.
\]
This uses only upper-half integrals and bounded derivatives of $S$, not invertibility of $DS$.
On the upper half, $\eta\nabla u=\nabla h-(u-\overline{\widehat u})\nabla\eta$.
Sobolev embedding and Poincar\'e's inequality therefore yield
\[
\|\eta\nabla u\|_{L^{p^*}(D_R^+)}
\leq C\|h\|_{W^{2,p}(D_R)}+C\|\nabla\widehat u\|_{L^2(D_R)}.
\]
Since $p<2$, the cutoff terms in $\Delta h$ have $L^p$-norm at most $C\|\nabla\widehat u\|_2$.
The zero-boundary Calder\'on--Zygmund estimate~\cite{GilbargTrudinger2001} gives
\[
\|h\|_{W^{2,p}(D_R)}
\leq CA\|h\|_{W^{2,p}(D_R)}
+C(1+A)\|\nabla\widehat u\|_{L^2(D_R)}.
\]
Choose $\varepsilon_1,\varepsilon_2$ in \eqref{intsma1-1}-\eqref{intsma1-2} so that $CA\leq1/2$.
Absorption and $\eta=1$ on $D_r$ give the assertion.
\end{proof}

On the compact image set $Q$, \eqref{eq:sec2-willmore-comparison} gives
\[
E(u)+W_{\R^K}(u)\leq C_Q\bigl(E(u)+W_M(u)\bigr).
\]
Proposition~\ref{intsma} requires both $\|\alpha\|_{L^2}$ and $\|du\|_{L^2}$ to be small.
For smooth branched HSL immersions, the Hodge equation for $\alpha$ allows Theorem~\ref{higher} to use small Dirichlet energy alone.

\begin{theorem}\label{higher}
Fix an isometric embedding $M\hookrightarrow\R^K$, a compact subset $Q\subset M$, and a compact smooth supporting submanifold $\cK\subset M$.
There exists $\epsilon_H=\epsilon_H(K,Q,\cK)>0$, independent of $p$, with the following property.
Let
\[
u:D_1^+\to M
\]
be a smooth branched conformal HSL immersion with $u(D_1^+)\subset Q$, and let $\alpha$ be its smooth geometric mean curvature $1$-form.
Assume
\[
\Delta u=\sum_{j=1}^2A^M_{\R^K}(u)(u_{x_j},u_{x_j})-\sum_{j=1}^2\alpha_{x_j}J(u)u_{x_j},
\qquad
d\alpha=u^*\operatorname{Rc},
\qquad
d^*\alpha=0
\]
and, on the diameter,
\[
u\in\cK,\qquad u_{x_2}\perp T_u\cK,\qquad \alpha_{x_2}=0.
\]
If
\[
E(u;D_1^+)<\epsilon_H^2,
\]
then for every $p>2$,
\begin{align}\label{holes}
\|\nabla u\|_{W^{1,p}(D_{1/2}^+)}
\leq
C_p\left(
\int_{D_1^+}|\nabla u|^2\,dx
+
\int_{D_1^+}|\mathbf H_{u,M}|^2\,d\mu_u
+1
\right)^2.
\end{align}
The threshold $\epsilon_H$ is the same for every $p>2$; no smallness assumption on $W_M(u;D_1^+)$ is needed.
\end{theorem}

\begin{proof}
Set
\[
a=\|\nabla u\|_{L^2(D_1^+)},
\qquad
b=\|\alpha\|_{L^2(D_1^+)}.
\]
Fix $R_*=7/8$.
By smoothness up to the diameter, Proposition~\ref{extsys} applies on $D_{R_*}^+$ after a fixed rescaling.
The reflected fields on $D_{R_*}$ satisfy
\[
\|\nabla\widehat u\|_{L^2(D_{R_*})}\leq C a,
\qquad
\|\widehat\alpha\|_{L^2(D_{R_*})}\leq\sqrt2\,b,
\]
and
\[
|\Delta\widehat u|
\leq C\bigl(|\nabla u\circ\rho_\pm|^2
+|\alpha\circ\rho_\pm|\,|\nabla u\circ\rho_\pm|\bigr),
\qquad
|d\widehat\alpha|
\leq C|\nabla u\circ\rho_\pm|^2,
\qquad
d^*\widehat\alpha=0,
\]
where $\rho_\pm$ is the identity or the domain reflection.

For $0<t<R_*$, set
\[
X(t) = \|\nabla\widehat u\|_{L^4(D_t)} + \|D^2\widehat u\|_{L^{4/3}(D_t)}, \qquad A_4(t)=\|\widehat\alpha\|_{L^4(D_t)}.
\]
Fix $0<r<R<R_*$, put $s=(r+R)/2$, and set $\delta=R-r$.
A cutoff Calder\'on--Zygmund estimate~\cite{GilbargTrudinger2001} in $L^{4/3}$ for $\widehat u$ and the embedding $W^{1,4/3}\hookrightarrow L^4$ give the estimate for $X(r)$:
\[
X(r) \leq C\delta^{-2}a + C\|\nabla u\|_{L^{8/3}(D_s^+)}^2
+ CaA_4(s).
\]
Using the interpolation inequality
\[
\|\nabla u\|_{L^{8/3}}^2 \leq \|\nabla u\|_{L^2}\|\nabla u\|_{L^4}
\]
yields
\begin{equation}\label{eq:higher-coupled-map}
X(r)\leq C\delta^{-2}a+CaX(s)+CaA_4(s).
\end{equation}
The local Hodge estimate \eqref{eq:hodge-local-bootstrap} in $L^{4/3}$, with a cutoff supported in $D_R$, gives
\[
A_4(s)\leq C\delta^{-1}b+C\|\nabla u\|_{L^{8/3}(D_R^+)}^2 \leq C\delta^{-1}b+CaX(R).
\]
Substituting this estimate into \eqref{eq:higher-coupled-map} gives
\begin{equation}\label{eq:higher-hole-filling}
X(r)\leq C\delta^{-2}a(1+b)+C(a+a^2)X(R).
\end{equation}
Choose $\epsilon_H$ so that $a<\sqrt2\,\epsilon_H$ implies $\theta:=C(a+a^2)<1/2$.
Fix $r_0<R_*$ and set $r_k=R_*-(R_*-r_0)\tau^k$ with $\tau=3/4$.
Then $\theta/\tau^2<8/9$, and iteration of \eqref{eq:higher-hole-filling} gives
\[
X(r_0)\leq C_{R_*,r_0,\tau}a(1+b)\sum_{k=0}^{m-1}(\theta/\tau^2)^k+\theta^mX(r_m).
\]
Smoothness up to the diameter and Proposition~\ref{extsys} give $X(R_*)<\infty$, so $\theta^mX(r_m)\leq\theta^mX(R_*)\to0$.
The geometric series converges.
Taking $r_0=3/4$ yields
\begin{equation}\label{eq:higher-L4}
\|\nabla u\|_{L^4(D_{3/4}^+)}+\|D^2u\|_{L^{4/3}(D_{3/4}^+)}
\leq Ca(1+b).
\end{equation}
Applying the preceding Hodge estimate on smaller disks also gives
\begin{equation}\label{eq:higher-alpha-L4}
\|\alpha\|_{L^4(D_{2/3}^+)}\leq C\bigl(b+a^2(1+b)\bigr)\leq C(1+b).
\end{equation}
The $L^2$ Hodge estimate for $\widehat\alpha$, the bound $|d\widehat\alpha|\leq C|\nabla u\circ\rho_\pm|^2$, and \eqref{eq:higher-L4} give
\[
\|\alpha\|_{W^{1,2}(D_{5/8}^+)}\leq C(1+b)^2.
\]
Thus, for every finite $q\geq2$,
\begin{equation}\label{eq:higher-all-q-alpha}
\|\alpha\|_{L^q(D_{5/8}^+)}\leq C_q(1+b)^2.
\end{equation}
The $L^2$ map estimate, together with \eqref{eq:higher-L4} and \eqref{eq:higher-alpha-L4}, gives
\[
\|D^2u\|_{L^2(D_{5/8}^+)}
\leq Ca(1+b)^2.
\]
Consequently,
\begin{equation}\label{eq:higher-all-q-map}
\|\nabla u\|_{L^q(D_{5/8}^+)}\leq C_qa(1+b)^2
\end{equation}
for every finite $q\geq2$.

For $p>2$, we apply the cutoff Calder\'on--Zygmund estimate~\cite{GilbargTrudinger2001} in $L^p$ to $\widehat u$ on $D_{1/2}$, using a cutoff supported in $D_{5/8}$.
H\"older's inequality and \eqref{eq:higher-all-q-alpha}-\eqref{eq:higher-all-q-map} then bound the first and second derivatives of $u$ by
\[
%\begin{split}
\|\nabla u\|_{W^{1,p}(D_{1/2}^+)}
\leq
C_p\left(
\|\nabla\widehat u\|_{L^p}
+\|\nabla u\|_{L^{2p}}^2
+\|\alpha\|_{L^{2p}}\|\nabla u\|_{L^{2p}}
\right) \leq C_pa(1+b)^4.
%\end{split}
\]
Since $a\leq1$ and
\[
b^2=\int_{D_1^+}|\mathbf H_{u,M}|^2\,d\mu_u=4W_M(u;D_1^+),
\]
the last bound implies \eqref{holes}.
The only absorption step requiring smallness is \eqref{eq:higher-hole-filling}, and its threshold is independent of $p$.
\end{proof}

We next prove qualitative regularity for $u\in W^{2,2}$ and $\alpha\in L^2$ satisfying the same equations and weak boundary traces.
The Sobolev chain rule and the trace matching in Proposition~\ref{extsys} give $\widehat u\in W^{2,2}$.
Hence $\nabla u\in L^q_{\loc}$ on the half disk for every finite $q$.
We first justify the reflected Hodge equations distributionally.
Write $d\alpha=f\,dx_1\wedge dx_2$ on the upper half disk; the Ricci source $f$ belongs locally to $L^2$.
For $\phi\in C_c^\infty(D)$, use $\phi+\phi\circ\rho$ in the weak zero-flux identity to obtain $d^*\widehat\alpha=0$.
For the exterior derivative, the scalar coefficient $\phi-\phi\circ\rho$ has zero trace on the diameter.
Density in $W_0^{1,2}(D^+)$ extends the distributional identity to this test.
It follows that $d\widehat\alpha=f\,dx_1\wedge dx_2$ on the upper half and $d\widehat\alpha=-(f\circ\rho)\,dx_1\wedge dx_2$ on the lower half, with no distribution supported on the diameter.
The reflected source $|d\widehat\alpha|\leq C|\nabla u\circ\rho_\pm|^2$ therefore belongs locally to $L^2$.
For $1<r<\infty$, a $1$-form $\gamma\in L^2(D_R)$ satisfying
\[
d\gamma=f\,dx\wedge dy,\qquad d^*\gamma=0,
\]
with $f\in L^r(D_R)$ obeys the local Hodge estimate
\begin{equation}\label{eq:hodge-local-bootstrap}
\|\gamma\|_{W^{1,r}(D_\rho)}
\leq
C\bigl(\|f\|_{L^r(D_R)}+\|\gamma\|_{L^2(D_R)}\bigr),
\qquad D_\rho\Subset D_R,
\end{equation}
which follows by solving $\Delta v=f$ with zero Dirichlet data and subtracting the form $-v_y\,dx+v_x\,dy$ from $\gamma$.
The Dirichlet estimate in \cite[Theorems~9.9 and~9.15]{GilbargTrudinger2001} controls $v$ in $W^{2,r}$.
The remainder is componentwise harmonic, so its $L^2$ norm controls its derivatives on smaller disks; see \cite[proof of Theorem~2.10 and Corollary~2.20]{Folland1995}.
Applying this estimate to the reflected Hodge system gives $\widehat\alpha\in W^{1,2}_{\loc}$, hence $\widehat\alpha\in L^q_{\loc}$ for every finite $q$.
The right-hand side of the reflected map equation then belongs to $L^p_{\loc}$ for every finite $p$, so $\widehat u\in W^{2,p}_{\loc}$.

On the original half disk write $\alpha=c\,dx+d\,dy$ and $d\alpha=f\,dx\wedge dy$.
The Hodge equations and boundary condition are
\[
d_x-c_y=f,\qquad c_x+d_y=0,\qquad d=0,
\]
so
\[
\Delta d=f_x, \qquad \Delta c=-f_y, \qquad c_y=-f.
\]
The first two identities hold in the half disk; the last holds weakly on its diameter.
Fix $p>2$ and pass to a smaller half disk.
The regularity $u\in W^{2,p}$ and smoothness of the Ricci form imply $f\in W^{1,p}$, while the local Hodge estimate gives $\alpha\in W^{1,p}$.
The scalar Dirichlet estimate~\cite[Theorems~9.9 and~9.15]{GilbargTrudinger2001} applied to $\Delta d=f_x$ with zero Sobolev trace gives $d\in W^{2,p}$.
The identities $c_x=-d_y$ and $c_y=d_x-f$ then give $c\in W^{2,p}$, hence $\alpha\in W^{2,p}$.
In Fermi coordinates for $\cK$, the normal variables of $u$ satisfy homogeneous Dirichlet conditions and the tangential variables satisfy homogeneous Neumann conditions.
The scalar Laplacian remains the principal part, and the right-hand sides lie in $W^{1,p}$.
Tangential difference quotients preserve these boundary conditions.
The Dirichlet estimate and even reflection for the homogeneous Neumann problem give the tangential derivatives in $W^{2,p}$.
The map equation recovers the remaining normal derivatives, so $u\in W^{3,p}$.
Iteration on smaller half disks makes $u$ and $\alpha$ smooth up to the diameter.
This argument requires no smallness hypothesis.

For use in the compactness argument, fix $\Lambda<\infty$ and $0<r<1/2$.
If $E(u;D_1^+)<\epsilon_H^2$ and $W_M(u;D_1^+)\leq\Lambda$, then Theorem~\ref{higher} bounds $u$ in $W^{2,p}$ on smaller half disks for every finite $p>2$, while the Hodge equation bounds $\alpha$ in $W^{1,p}$ there.
A fixed finite cover of $\cK$ by Fermi charts bounds the coefficients and their derivatives.
The $W^{2,p}$ bound with $p>2$ gives a uniform source radius on which the image lies in one chart.
The boundary iteration has uniform constants, so for every integer $m\geq0$,
\begin{equation}\label{eq:hsl-uniform-smooth-estimates}
 \|\nabla u\|_{C^m(\overline{D_r^+})}
 +\|\alpha\|_{C^m(\overline{D_r^+})}
 \leq C_{m,r,\Lambda},
\end{equation}
where the constant also depends on the fixed embedding, $Q$, and $\cK$, but not on $u$.
For a weak $W^{2,2}$ input, the preceding qualitative regularity argument first gives smoothness, after which the same conclusion holds.

There is no positive Willmore gap: the real unit disk $u(x,y)=(x,y)\in\R^2\subset\mathbb C^2$ has $du(\nu)=u$ on $S^3$ and $\mathbf H_{u,M}=\alpha=0$.
It is a nonconstant conformal Lagrangian immersion satisfying the homogeneous HSL system \eqref{sys1}-\eqref{sys2} and stationary under every admissible Hamiltonian variation.

\Needspace{6\baselineskip}
\subsection{\texorpdfstring{Boundary removability}{Boundary removability}}
\label{subsec:hsl-removability}

Theorem~\ref{remo} removes a boundary puncture for the homogeneous HSL system under finite area and finite intrinsic Willmore energy.
It is the boundary counterpart of Chen--Ma's interior removability theorem~\cite{chen2021}.

\begin{theorem}\label{remo}
 
Let $M$ be a K\"ahler surface isometrically embedded in a Euclidean space, let $Q\Subset M$ be compact, and let $\cK\subset M$ be a compact smooth support.
Let
\[
u:D_1^+\setminus\{0\}\longrightarrow Q
\]
be a smooth branched conformal HSL immersion with only finitely many branch points, and let $\alpha$ be its geometric mean curvature $1$-form.
Assume
\[
\cA(u;D_1^+)+W_M(u;D_1^+)<\infty,
\qquad
d\alpha=u^*\operatorname{Rc},
\qquad
d^*\alpha=0.
\]
On the punctured diameter, assume
\[
u\in\cK, \qquad u_y\perp T_u\cK, \qquad\alpha_y=0.
\]
Then $u$ and $\alpha$ extend smoothly across $0$ on the original closed upper half disk, and all three boundary conditions persist.
 
\end{theorem}

\begin{proof}
 
Shrink the half disk so that it contains no branch point away from the puncture.
Write
\[
z=e^{-t+i\theta},\qquad t>T_0,\qquad 0\leq\theta\leq\pi.
\]
Both the Dirichlet energy and $\|\alpha\|_2^2=4W_M(u)$ are conformally invariant.
The embedding comparison \eqref{eq:sec2-willmore-comparison} therefore implies that the Dirichlet energy, Euclidean Willmore energy, and $L^2$-norm of $\alpha$ on the tail tend to zero.

Choose $T_0$ so that every unit strip in the tail, together with its fixed enlargement, lies below the thresholds of Theorem~\ref{intsma0} and the interior estimate proved in Theorem~\ref{deca}.
The tail is smooth and unbranched, with compact image and finite Dirichlet and Euclidean Willmore energies.
Thus the unbranched alternative of Theorem~\ref{deca} applies without a total curvature assumption.
We use Theorem~\ref{intsma0} on boundary charts and the full disk estimate on interior charts.
Morrey's inequality gives
\[
\operatorname{Osc}^{\R^K}_{C_j}u\leq C E(u;C_j^*)^{1/2}
\]
on each unit tail strip.
Each $C_j$ meets a longitudinal side $\theta=0$ or $\theta=\pi$ whose image lies in $\cK$.
Thus the tail eventually lies in any prescribed tubular neighborhood of $\cK$.

Applying the unbranched alternative of Theorem~\ref{deca} to $u$ on the longitudinal half-cylinder $[T_0,2\tau-T_0]\times[0,\pi]$, with center $\tau$ and central half-length $1$, gives exponential decay of the Dirichlet energy on the central strip.
More precisely, for some $\delta=1-C_0\epsilon>0$,
\[
\int_{[\tau-1,\tau+1]\times[0,\pi]}|\nabla u|^2\leq Ce^{-\delta\tau}.
\]
Summing the strip estimates yields
\[
\int_{[T,\infty)\times[0,\pi]}|\nabla u|^2
\leq Ce^{-\delta T}.
\]
Since the strip oscillations are summable, $u$ converges uniformly along the tail to a single point of $\cK$.

Reflect across the two longitudinal edges using the local geodesic reflection $\sigma$ from Subsection~\ref{sec:geodesic-reflection}, and denote the periodic doubled map by $\widetilde U$.
Returning to disk coordinates gives a reflected map $U:D_{r_0}\setminus\{0\}\to\R^K$ that extends continuously to $0$ and satisfies
\begin{equation}\label{morrey_growth}
\int_{D_r}|\nabla U|^2\leq Cr^\delta.
\end{equation}

Choose
\[
1<p_0<\frac{4}{4-\delta}.
\]
Cover each dyadic upper half-annulus $D_{2r}^+\setminus D_r^+$, together with a fixed enlargement, by a uniformly finite family of boundary half disks and interior disks.
The scaled estimates, together with \eqref{morrey_growth}, give
\[
\|D^2u\|_{L^{p_0}(D_{2r}^+\setminus D_r^+)}
\leq
C r^{\,2/p_0-2+\delta/2}.
\]
On the reflected lower half, the chain rule gives terms of the form
\[
D\sigma\,D^2u
\quad\text{and}\quad
D^2\sigma[\nabla u,\nabla u].
\]
To control $D^2\sigma[\nabla u,\nabla u]$, apply the disk and half disk estimates to $u$ with
\[
p_1=\frac{2p_0}{1+p_0}.
\]
The embedding $W^{2,p_1}\hookrightarrow W^{1,2p_0}$ in dimension two then gives
\[
\|\nabla u\|_{L^{2p_0}(D_{2r}^+\setminus D_r^+)}
\leq
C r^{\,1/p_0-1+\delta/2}.
\]
Therefore
\[
\||\nabla u|^2\|_{L^{p_0}(D_{2r}^+\setminus D_r^+)}
\leq
C r^{\,2/p_0-2+\delta},
\]
which is stronger than the preceding $D^2u$ bound.
Consequently,
\[
\|D^2U\|_{L^{p_0}(D_{2r}\setminus D_r)}
\leq
C r^{\,2/p_0-2+\delta/2}.
\]
The $p_0$-th powers sum over the dyadic annuli because
\[
2-2p_0+\frac{p_0\delta}{2}>0.
\]
Define $U(0)$ by continuity.
The finite Dirichlet energy gives an extension of $U$ across the point as a $W^{1,2}$ map.
Each distributional second derivative differs from its classical $L^{p_0}$ value off $0$ by a distribution supported at $0$.
Both terms belong to $W^{-1,2}$: the distributional derivatives of the $L^2$ first derivatives lie in $W^{-1,2}$, while $L^{p_0}\hookrightarrow W^{-1,2}$ because $p_0>1$.
In dimension two, no nonzero distribution supported at one point belongs to $W^{-1,2}$; in particular, the Fourier $W^{-1,2}$-norm of the Dirac mass diverges.
Hence $U\in W^{2,p_0}_{\loc}(D_{r_0})$.
On a sufficiently small disk the reflected equation gives
\[
|\Delta U|
\leq b_0|\nabla U|,
\qquad
b_0=C\bigl(|\nabla U|+|\widehat\alpha|\bigr),
\]
and $\|b_0\|_{L^2}$ tends to zero with the radius.
Where $|\nabla U|\neq0$, define matrix fields
\[
B_j=\frac{\Delta U\otimes U_{x_j}}{|\nabla U|^2}, \qquad j=1,2,
\]
where $(a\otimes b)\xi=a\langle b,\xi\rangle$, and set $B_j=0$ where $\nabla U=0$.
Then
\[
\Delta U=B_1U_x+B_2U_y,
\qquad
\left(\|B_1\|_2^2+\|B_2\|_2^2\right)^{1/2}
\leq\|b_0\|_2.
\]
Fix $p_0<p<2$, choose a cutoff $\eta$, and put $w=\eta U$.
Then
\[
\Delta w=B_1w_x+B_2w_y+g,
\]
where
\[
g=2\nabla\eta\cdot\nabla U+U\Delta\eta-B_1U\eta_x-B_2U\eta_y.
\]
The embedding $W^{2,p_0}\hookrightarrow W^{1,2p_0/(2-p_0)}$, boundedness of $U$, and $B_j\in L^2$ show that $g\in L^p$.
Let $G_q$ be the inverse Dirichlet Laplacian on $L^q$.
For $q=p_0,p$,
\[
\|G_q(B_1v_x+B_2v_y)\|_{W^{2,q}}
\leq
C_q
\left(\|B_1\|_2^2+\|B_2\|_2^2\right)^{1/2}
\|v\|_{W^{2,q}},
\]
because $W^{2,q}\hookrightarrow W^{1,2q/(2-q)}$.
Shrink the disk so that both operator norms are less than $1/2$.
The contraction mapping theorem gives a unique $W^{2,p}$ solution.
It also belongs to $W^{2,p_0}$, so uniqueness in that space identifies it with $w$.
Thus $U\in W^{2,p}_{\loc}$ for every fixed $p<2$.

Since $\widehat\alpha\in L^2$, logarithmic cutoffs $\chi_r$ with $\|d\chi_r\|_{L^2}\to0$ extend both reflected Hodge equations across the puncture.
The cutoff error is bounded by
\[
\|\widehat\alpha\|_{L^2}\|d\chi_r\|_{L^2},
\]
and the Ricci source is in $L^1$.
Hence, distributionally across the origin,
\[
|d\widehat\alpha|\leq C|\nabla U|^2,
\qquad
d^*\widehat\alpha=0.
\]
Choose $p<2$ sufficiently close to $2$ that
\[
W^{2,p}\hookrightarrow W^{1,q}\quad\text{for some }q>6.
\]
The Ricci source then lies in $L^{q/2}$.
The local Hodge estimate \eqref{eq:hodge-local-bootstrap} gives
\[
\widehat\alpha\in W^{1,q/2}_{\loc}\subset L^\infty_{\loc}.
\]
The reflected map equation has right-hand side in $L^{q/2}$, so Calder\'on--Zygmund estimates~\cite{GilbargTrudinger2001} yield
\[
U\in W^{2,q/2}_{\loc}\subset W^{2,2}_{\loc}.
\]
We now return to the original upper half disk.
Adding one point does not change
\[
\|\alpha\|_{L^2}^2=\int|\mathbf H_{u,M}|^2\,d\mu_u.
\]
The signed equations and the weak boundary conditions hold across the origin. \repair{Applying the regularity argument for weak solutions following Theorem~\ref{higher} to $(u,\alpha)$ in the original Fermi half disk coordinates proves that both $u$ and $\alpha$ are smooth up to the boundary.} Continuity gives $u(0)\in\cK$, and the smooth boundary traces extend $u_y\perp T_u\cK$ and $\alpha_y=0$ to $0$.
\end{proof}

\Needspace{6\baselineskip}
\subsection{\texorpdfstring{Branch points and extension under finite total curvature}{Branch points and extension under finite total curvature}}
\label{subsec:hsl-branch-points}

Lemma~\ref{lem:hsl-branch-orders} gives integer orders for the branch points of nonconstant smooth solutions of the homogeneous HSL system.
Under finite area and total curvature, Theorem~\ref{w22} gives a $W^{2,2}$ boundary extension for branched conformal immersions with free boundary and a smooth extension in the homogeneous HSL case.

\begin{lemma}\label{lem:hsl-branch-orders}
Let $u$ and $\alpha$ be smooth on a closed upper half disk.
Assume that $u$ is weakly conformal and Lagrangian on the set $\{z:du(z)\ne0\}$, that
\[
\Delta u=\sum_{j=1}^2A^M_{\R^K}(u)(u_{x_j},u_{x_j})
-\sum_{j=1}^2\alpha_{x_j}J(u)u_{x_j},
\qquad
d\alpha=u^*\operatorname{Rc},
\qquad
d^*\alpha=0,
\]
and that
\[
u(I_1)\subset\cK,\qquad
u_y\perp T_u\cK,\qquad
\alpha_y=0
\quad\text{on }I_1.
\]
Then either $u$ is constant or the zero set of $du$ is finite on every smaller closed half disk.
At every interior or boundary zero $z_0$, there is an integer $m\geq1$ such that
\[
\lambda(z)=m\log|z-z_0|+\omega(z)
\]
on the punctured coordinate disk or half disk, where $u^*\bar g=e^{2\lambda}|dz|^2$ and $\omega\in C^0\cap W^{1,2}$.
\end{lemma}

\begin{proof}
Denote the reflected map and form from Proposition~\ref{extsys} by $U=\widehat u$ and $\widehat\alpha$, respectively, and put $v(x,y)=u(x,-y)$ on the lower half disk.
Near a boundary point, continuity ensures that $u$ lies in the smaller tubular neighborhood where $S$ is the involutive normal reflection.
Thus $DS$ is invertible there, and the derivatives of $v$ on the lower half are bounded linear combinations of the derivatives of $U$.

Since $u$ and $\alpha$ are smooth, one factor in each quadratic term of \eqref{defE} is bounded.
The terms in \eqref{defF1} are already linear in $\nabla U$.
Thus, on both halves of every smaller reflected disk, bounded measurable real $K\times K$ matrices $B_1,B_2$ give
\begin{equation}\label{eq:reflected-branch-linear-system}
\Delta U=B_1U_x+B_2U_y
\end{equation}
distributionally.

Set $w=U_x-iU_y$.
Equation \eqref{eq:reflected-branch-linear-system} gives
\[
\overline\partial w=Aw+D\overline w,
\qquad A=\frac{B_1+iB_2}{4},\quad D=\frac{B_1-iB_2}{4}.
\]
Define $C=0$ on $\{w=0\}$ and $C=(Aw+D\overline w)w^*/|w|^2$ elsewhere.
Then $|C|\leq|A|+|D|$ and $\overline\partial w=Cw$.

Fix $p>2$ and a small disk $D_r(z_0)$ in the reflected domain.
Extend $f$ by zero and set $Nf(z)=(2\pi)^{-1}\int\log|z-\zeta|f(\zeta)\,d\mathcal L^2(\zeta)$.
Then the Cauchy transform $T$ satisfies $Tf=4\partial Nf$, with kernel $1/(\pi z)$; see \cite{Vekua1962} for the general construction.
After subtracting the constant $(2\pi)^{-1}\log r\int f\,d\mathcal L^2$ and scaling to a unit disk, Young's inequality and \cite[Theorem~9.9]{GilbargTrudinger2001} bound $D^2Nf$ in $L^p$ on $D_r$ by $C_p\|f\|_{L^p}$.
H\"older's inequality for the Cauchy kernel gives the supremum bound, and Young's inequality controls the lower-order $L^p$ term for $r\leq1$.
Thus
\[
\overline\partial(Tf)=f,\qquad
\|Tf\|_\infty\leq c r^{1-2/p}\|f\|_{L^p},
\qquad
\|Tf\|_{W^{1,p}(D_r)}\leq c_p\|f\|_{L^p}.
\]
Choose $r$ so that $cr\|C\|_\infty<1/3$.
The map $P\mapsto I+T(CP)$ is then a contraction on bounded matrix fields.
Its fixed point satisfies $\overline\partial P=CP$, $P\in W^{1,p}$, and $\|P-I\|_\infty<1/2$, so $P$ is invertible.
The Sobolev product rule gives $\overline\partial(P^{-1}w)=0$.
Thus $h=P^{-1}w$ is holomorphic.

If $h\not\equiv0$, then $h(z)=(z-z_0)^m h_0(z)$, where $m$ is the least vanishing order of its nonzero components and $h_0(z_0)\neq0$.
These local factorizations and the identity theorem imply that either $U$ is constant or the zeros of $dU$ are isolated.
A smaller closed half disk is compact in the reflected domain, so it contains only finitely many zeros.
On the upper half disk, conformality gives $|w|^2=2e^{2\lambda}$.
Consequently,
\[
\lambda(z)=m\log|z-z_0|+\omega(z),
\qquad
\omega=\log|Ph_0|-\tfrac12\log2\in C^0\cap W^{1,p}
\subset C^0\cap W^{1,2}.
\]
The last assertion follows from the Sobolev chain rule because $|Ph_0|$ is bounded away from zero near $z_0$.
\end{proof}

Lemma~\ref{lem:hsl-branch-orders} applies to the smooth extensions supplied by Theorem~\ref{remo}.
The next removability theorem is the boundary analogue of the local result used in \cite{kuwert2012,chen2014}.
\begin{theorem}\label{w22}
 
Let $Q\Subset M$ be the fixed compact target region used in this section, and let
\[
u\in
W_{\loc}^{2,2}
((D_1^+\cup I_1)\setminus\{0\},M)
\]
be a branched conformal immersion satisfying
\[u((D_1^+\cup I_1)\setminus\{0\})\subset Q\]
with free boundary on $\cK$.
Assume
\[
\int_{D_1^+\setminus\{0\}}
|A_{u,M}|^2\,d\mu_u<\infty,
\qquad
\mu_u(D_1^+\setminus\{0\})<\infty,
\]
where $(g_u)_{ij}=e^{2\lambda}\delta_{ij}$ in a conformal coordinate.
Then
\[
u\in W^{2,2}_{\mathrm{loc}}(D_1^+\cup I_1,M).
\]
In addition, if $u$ is HSL on the punctured half disk, satisfies \eqref{sys1}-\eqref{sys2}, and in particular has the homogeneous normal trace
\[
\alpha_y=0
\]
on the punctured diameter, then $u$ and its mean curvature $1$-form $\alpha$ extend smoothly across $0$.
\end{theorem}

\begin{proof}
Conformality and finite area give
\[
\int_{D_1^+\setminus\{0\}}|\nabla u|^2\,dx
=2\mu_u(D_1^+\setminus\{0\})<\infty.
\]
Since the image lies in the compact set $Q$, the map is bounded and therefore belongs to $W^{1,2}(D_1^+)$; the removed point lies on the boundary of this open half disk.
Using the finite-branch-set convention for the branched immersion in this theorem, shrink to a half disk containing no branch point other than the possible puncture.
The embedding comparison \eqref{eq:sec2-total-curvature-comparison} gives finite Euclidean total curvature there.
Lemma~\ref{lem:general-boundary-removability} gives a $W^{2,2}$ extension near $0$.
Together with the assumed local regularity away from $0$, this proves $u\in W^{2,2}_{\loc}(D_1^+\cup I_1,M)$, independently of Theorem~\ref{remo}.

For the HSL assertion, $|\mathbf H_{u,M}|^2\leq2|A_{u,M}|^2$ gives finite intrinsic Willmore energy.
Local weak regularity following Theorem~\ref{higher} makes $u$ and $\alpha$ smooth away from $0$.
The HSL equations with the stated sign conventions, geometric free boundary conditions, and zero normal trace then meet the hypotheses of Theorem~\ref{remo}, which gives their smooth extension on the original half disk.
\end{proof}

\vspace{1cm}

\section{\texorpdfstring{Proof of the HSL compactness theorem}{Proof of the HSL compactness theorem}}
\label{section5}

We prove Theorem~\ref{main theorem} using the fixed isometric embedding $(M,\bar g)\hookrightarrow\R^K$.
By \eqref{eq:distance bound} and the Hopf--Rinow theorem, there is a compact set $Q\Subset M$ containing every image $u_n(\Sigma)$.
Using the embedding comparison \eqref{eq:sec2-willmore-comparison}, fix a common bound $\Lambda$ for $\cA(u_n)+W_M(u_n)$ and $E(u_n)+W_{\R^K}(u_n)$.
The geometric coefficients determined by $M$ and $\cK$ and all derivatives used below are uniformly bounded on $Q$.
The reflected system in \eqref{equation:global-form}--\eqref{extalpha} supplies the boundary estimates, while compactness, rescaling, and neck arguments use the original conformal domains.
Each input pair $(u_n,\alpha_n)$, with $\alpha_n$ its geometric mean curvature $1$-form, satisfies the homogeneous HSL system \eqref{sys1}--\eqref{sys2}, including orthogonal contact and zero normal trace.
Subsection~\ref{subsec:hsl-limit} passes this system to every nonconstant limiting component, and Subsection~\ref{subsec:hsl-first-variation} proves the first variation formula.

\Needspace{6\baselineskip}
\subsection{\texorpdfstring{Convergence of the HSL system}{Convergence of the HSL system}}
\label{subsec:hsl-limit}

Lemma~\ref{lem:hsl-local-convergence} gives local smooth convergence and identifies the nonconstant limits as smooth branched HSL immersions satisfying the same homogeneous boundary conditions.
We first check that each smooth branched HSL immersion $u_n$ belongs to the class in Definition~\ref{def:fb-branched-conformal}.
Finite intrinsic Willmore energy gives $\alpha_n\in L^2$ in conformal coordinates.
The logarithmic cutoff argument in Subsection~\ref{subsec:hsl-system} extends the Hodge equations and the homogeneous flux identity across the finite branch set.
Since $u_n$ is smooth, the local Hodge estimates give smoothness of $\alpha_n$ at an interior branch point.
At a boundary branch point, reflection first gives local $W^{1,p}$ regularity for every finite $p$, and the Dirichlet--Neumann estimates on the original half disk following Theorem~\ref{higher} give smoothness up to the boundary.
For fixed $n$, the proof of Lemma~\ref{lem:hsl-branch-orders} gives $u_{n,x}-iu_{n,y}=(z-z_0)^m a$ near each branch point, where $a\in W^{1,p}$ for some $p>2$ and $|a|$ is bounded away from zero.
Thus the projective Gauss map $G_n=[u_{n,x}+iu_{n,y}]=[\overline a]$ extends in $W^{1,p}$ across the branch point.
The Gauss map identity gives
\[
\int |A_{u_n}|^2\,d\mu_{u_n}
\leq C\int |DG_n|^2<\infty
\]
on each such disk or half disk; see \cite[Section~2.1]{kuwert2012}.
Smoothness gives finite curvature away from the branch points, so each input belongs to the required class.
The embedding comparison gives uniform Euclidean area and Willmore bounds.
Proposition~\ref{prop:automatic-curvature-branching-bound} then bounds total curvature and branching uniformly.

Applying the Euclidean construction in Section~\ref{sec:general-compactness-proof} to $u_n$, viewed as $\R^K$-valued maps through the fixed isometric embedding, gives a finite bubble tree and componentwise weak $W^{2,2}$ convergence on compact subsets away from the selected exceptional points.
Choose $\vartheta$ also below the interior HSL regularity threshold in \cite[Section~3]{chen2021} and below $\epsilon_H^2$ from Theorem~\ref{higher}.
Choose $\delta$ by \eqref{eq:selection-delta-margin}.
The extraction remains valid with these smaller thresholds.

\begin{lemma}\label{lem:hsl-local-convergence}
Let $u_n$ be the HSL sequence considered above, with geometric mean curvature $1$-forms $\alpha_n$ satisfying \eqref{sys1}--\eqref{sys2}, and assume the local hypotheses of Lemma~\ref{lem:local-convergence}.
With the corresponding component comparison maps $\Phi_n^i$, a further subsequence satisfies
\[
u_n\circ\Phi_n^i\longrightarrow u_\infty^i,
\qquad
(\Phi_n^i)^*\alpha_n\longrightarrow\alpha_\infty^i
\quad\text{locally in }C^\infty
\]
on smaller compact subsets covered by those hypotheses, including their boundary portions inherited from $\partial\Sigma$.
\end{lemma}

\begin{proof}
For a weak $W^{2,2}$ input, $\nabla u\in L^q$ for every finite $q$, so \eqref{eq:hodge-local-bootstrap} gives $\alpha\in W^{1,p}_{\loc}$ for every finite $p$.
The map equation and the qualitative boundary regularity argument following Theorem~\ref{higher} then make the pair smooth, including at its branch points and up to the inherited boundary, before we apply the quantitative estimates.

On a principal core, use the actual conformal charts $\Psi_n\circ\theta_n$ constructed in Subsection~\ref{subsec:generic-extraction}; on the other cores, use the fixed conformal charts.
On an interior disk with $W_{\R^K}(u_n)<\vartheta$, the small Willmore estimate for branched conformal maps in \cite[Section~3]{chen2021}, with $p=12/7$, gives a uniform $W^{1,12/7}$ bound for $\nabla u_n$ on a smaller disk, and hence a uniform $L^{12}$ bound.
The Ricci source $u_n^*\operatorname{Rc}$ is therefore bounded in $L^6$.
Since $\|\alpha_n\|_{L^2}^2\leq4\Lambda$, the local Hodge estimate \eqref{eq:hodge-local-bootstrap} gives a uniform $W^{1,6}$ bound for $\alpha_n$ on a smaller disk.
The map equation then gives a uniform $W^{2,6}$ bound for $u_n$.
Alternating the Hodge and map estimates, using that $W^{1,6}$ is a multiplication algebra in dimension $2$, gives uniform bounds of every order on smaller disks.
These estimates apply to branched inputs and do not require a positive lower bound for $|du_n|$.

On a half disk meeting the boundary inherited from $\partial\Sigma$, joint smallness gives $E(u_n)<\vartheta<\epsilon_H^2$.
Theorem~\ref{higher} and \eqref{eq:hsl-uniform-smooth-estimates}, together with $\|\alpha_n\|_{L^2}^2\leq4\Lambda$, give uniform smooth estimates on a smaller closed half disk.

A retained $\nu$-mark has a fixed neighborhood $U$ with $\mu_n(U)<8\delta$ for all sufficiently large $n$, as in Subsection~\ref{subsec:generic-extraction}.
By \eqref{eq:selection-delta-margin}, $8\delta<\vartheta$, so the same interior or boundary estimates apply across the mark, regardless of its $\nu$-mass.
On a principal core, first choose the neighborhood $U$ and then a smaller limiting conformal chart compactly contained in $U$; the smoothly converging charts $\theta_n$ remain inside $U$ for large $n$.
Their smooth convergence and that of their inverses transfer the estimates to the fixed core.
Apply the estimates on a slightly larger core near its artificial boundary; use half-disk estimates only along the boundary inherited from $\partial\Sigma$.
A finite covering and one diagonal subsequence over the core exhaustions give smooth convergence on every fixed core, including retained marks and branch points.
The constants may depend on the fixed core, $Q$, the support geometry, and $\Lambda$, but not on $n$.
\end{proof}

By Lemma~\ref{lem:hsl-local-convergence}, the pulled-back pairs converge smoothly on each fixed principal or nonprincipal core, including the retained marks.
Write $\alpha_\infty^i$ for the form limit on the corresponding component.
On each principal core in $\widetilde\Sigma_0$, use the smoothly converging metrics $\widetilde h_n$ from Subsection~\ref{subsec:generic-extraction}, which are compatible with $\Psi_n^*j_n$; on every nonprincipal core, use its fixed compatible metric.
Their Hodge stars, codifferentials, and outward unit conormals converge smoothly.
Conformal changes of these metrics preserve coclosedness and zero normal trace.
The map equation, conformality, the Lagrangian identity, and orthogonal contact therefore pass to the limit, together with
\[
d\alpha_\infty^i=(u_\infty^i)^*\operatorname{Rc},
\qquad
d_{h_i}^*\alpha_\infty^i=0,
\qquad
\alpha_\infty^i(\nu_{h_i})=0,
\]
where $h_i$ is a limiting smooth compatible metric on the component.
For a nonconstant component, the limiting map equation identifies $\alpha_\infty^i$ with its geometric mean curvature $1$-form away from the branch points.
For a constant component, $\alpha_\infty^i$ is harmonic with zero normal trace.

Consider a nonconstant component.
The compactification in Subsection~\ref{subsec:compactification} gives finite area and total curvature and a finite branch set.
Each interior special point has a full punctured-disk neighborhood with compact image and finite intrinsic Willmore energy inherited from the construction.
This includes the unique normalized principal preimage of a type~II node, each interior attachment point, and the point at infinity of a sphere bubble.
The interior removability result of Chen--Ma~\cite{chen2021} therefore extends $u_\infty^i$ and $\alpha_\infty^i$ smoothly across these points.
At a boundary special point, Theorem~\ref{remo} applies because orthogonal contact and zero normal trace have been preserved.
Lemma~\ref{lem:hsl-branch-orders} gives the integer orders of the resulting branch points.

If $u_\infty^i$ is constant, the uniform bound and exhaustion give $\alpha_\infty^i$ finite $L^2$ norm.
Logarithmic cutoffs extend its closedness and coclosedness across each interior puncture, and local Hodge regularity gives a smooth extension.
At a boundary puncture, reflect the tangential coefficient evenly and the normal coefficient oddly; the zero normal trace gives the same harmonic extension argument.
Thus $\alpha_\infty^i$ extends smoothly to the compactified constant domain.

Intrinsic Willmore lower semicontinuity follows from weak convergence of the mean curvature forms.
Since the $L^2$-norm of $1$-forms is conformally invariant,
\[
W_M(u_n)
=
\frac14\int_\Sigma |\alpha_n|_{h_n}^2\,\dd\vol_{h_n}
=
\frac14\|\alpha_n\|_{L^2(\Sigma,h_n)}^2.
\]
Weak $L^2$ convergence of the pulled-back forms on pairwise disjoint cores, with the smoothly converging compatible metrics described above, followed by exhaustion, gives
\[
\sum_{i\in\mathfrak I_*}
W_M(u_\infty^i)
\leq
\liminf_{n\to\infty}W_M(u_n).
\]
Every nonconstant component of the limit is therefore a smooth branched conformal HSL immersion; constant components remain ghost components.
The sum includes only nonconstant components; the limiting forms on constant domains need not vanish.

\Needspace{6\baselineskip}
\subsection{\texorpdfstring{First variation and completion of the proof}{First variation and completion of the proof}}
\label{subsec:hsl-first-variation}

Proposition~\ref{prop:stationary} gives the first variation of the limiting varifold and its stationarity under infinitesimally admissible ambient Hamiltonian fields.
Together with the preceding convergence results, this completes the proof of Theorem~\ref{main theorem}.

\begin{prop}\label{prop:stationary}
 
Let $V_\infty$ be the parametrized integral varifold in Theorem~\ref{main theorem}, with its nonconstant components indexed by $\mathfrak I_*$.
We use the induced metrics $g_{u_\infty^i}$ and area measures $d\mu_{u_\infty^i}$ from Subsection~\ref{sec:fbbci}, and the boundary arclength measures $ds_{g_{u_\infty^i}}$ from \eqref{eq:induced-boundary-geometry}.
Each preimage outside the branch sets contributes once to these integrals.
Thus multiple covers and coincident component images retain their multiplicities, while constant components contribute no area.
Let $\mathbf H_{u_\infty^i,M}$ be the trace mean curvature vector in $M$ on $\Sigma_\infty^i\setminus\cB_{u_\infty^i}$.
The boundary integrals use the outward unit image conormal $\eta_{u_\infty^i}$ from \eqref{eq:induced-boundary-geometry}, defined almost everywhere on $\partial\Sigma_\infty^i\setminus\cB_{u_\infty^i}$.
Define
\begin{equation}\label{eq:mean curvature measure}
    \mathcal H_\infty
=
\sum_{i\in\mathfrak I_*}
(u_\infty^i)_\#
(\mathbf H_{u_\infty^i,M}\,d\mu_{u_\infty^i})
\end{equation}
and the vector-valued boundary conormal measure
\[
\mathcal B_\infty(X)
=
\sum_{i\in\mathfrak I_*}
\int_{\partial\Sigma_\infty^i}
\langle X(u_\infty^i),\eta_{u_\infty^i}\rangle\,ds_{g_{u_\infty^i}}.
\]
Thus, for a compactly supported smooth vector field $X$ on $M$, the vector-valued measure in \eqref{eq:mean curvature measure} acts by
\[
\mathcal H_\infty(X)
=\sum_{i\in\mathfrak I_*}
\int_{\Sigma_\infty^i}
\langle\mathbf H_{u_\infty^i,M},X(u_\infty^i)\rangle\,d\mu_{u_\infty^i}.
\]
Then $\mathcal H_\infty\ll\|V_\infty\|$.
Write
\[
\overline{\mathbf H}_\infty
=
\frac{d\mathcal H_\infty}{d\|V_\infty\|}.
\]
For every compactly supported smooth vector field $X$ on $M$,
\begin{equation}\label{eq:limit-varifold-first-variation}
\delta V_\infty(X)
=
-\int_M
\langle\overline{\mathbf H}_\infty,X\rangle
\,d\|V_\infty\|
+\mathcal B_\infty(X).
\end{equation}
The measure $\mathcal B_\infty$ is supported on $\cK$ and is singular with respect to $\|V_\infty\|$.

At $\|V_\infty\|$-almost every image point $y$, all sheets through $y$ have the same approximate tangent plane.
For such $y$, set
\[
\Theta(y)=\sum_{i\in\mathfrak I_*}
\#\bigl((u_\infty^i)^{-1}(y)\setminus\cB_{u_\infty^i}\bigr).
\]
The area formula gives $1\leq\Theta(y)<\infty$ almost everywhere, and
\begin{equation}\label{eq:fiberwise-mean-curvature-average}
\overline{\mathbf H}_\infty(y)
=
\frac1{\Theta(y)}
\sum_{i\in\mathfrak I_*}
\sum_{\substack{x\in (u_\infty^i)^{-1}(y)\\x\notin\cB_{u_\infty^i}}}
\mathbf H_{u_\infty^i,M}(x).
\end{equation}
Moreover,
\begin{equation}\label{eq:limit-mean-curvature-L2-bound}
\int_M|\overline{\mathbf H}_\infty|^2\,d\|V_\infty\|
\leq
\sum_{i\in\mathfrak I_*}
\int_{\Sigma_\infty^i}|\mathbf H_{u_\infty^i,M}|^2\,d\mu_{u_\infty^i}.
\end{equation}
For every $f\in C_c^\infty(M)$, without imposing a boundary tangency condition,
\begin{equation}\label{eq:all-potential-curvature-identity}
\int_M
\langle\overline{\mathbf H}_\infty,J\nabla f\rangle
\,d\|V_\infty\|
=0.
\end{equation}
If $X(u_\infty^i(x))\in T_{u_\infty^i(x)}\cK$ at every boundary point outside $\cB_{u_\infty^i}$ of every nonconstant component, then $\mathcal B_\infty(X)=0$.
In particular, this holds when $X$ is tangent to $\cK$ along all of $\cK$.
Hence $\delta V_\infty(J\nabla f)=0$ for every $f\in C_c^\infty(M)$ whose Hamiltonian field satisfies the preceding boundary tangency condition on every limiting component, that is, for every infinitesimally admissible ambient Hamiltonian field in the sense of Definition~\ref{def:hsl-variation-classes}.
For a general field, the singular boundary term in \eqref{eq:limit-varifold-first-variation} remains present.
The field $\overline{\mathbf H}_\infty$ represents only the absolutely continuous interior part of the first variation.
 
\end{prop}

\begin{proof}
 
Remove pairwise disjoint coordinate disks and half disks of radius $r$ around the finitely many branch points of each component.
Denote the remaining domain by $\Omega_{i,r}$ and its artificial boundary by $\Gamma_{i,r}=\partial\Omega_{i,r}\setminus\partial\Sigma_\infty^i$.
Let $\eta_{i,r}$ be the image under $u_\infty^i$ of the outward unit conormal along $\Gamma_{i,r}$ with respect to the induced metric $g_{u_\infty^i}$.
The classical first variation formula gives
\[
\begin{aligned}
\delta\bigl((u_\infty^i)_\#|\Omega_{i,r}|\bigr)(X)
&=
-\int_{\Omega_{i,r}}
\langle\mathbf H_{u_\infty^i,M},X(u_\infty^i)\rangle\,d\mu_{u_\infty^i}
\\
&\quad+
\int_{\partial\Sigma_\infty^i\cap\partial\Omega_{i,r}}
\langle X(u_\infty^i),\eta_{u_\infty^i}\rangle\,ds_{g_{u_\infty^i}}
+\int_{\Gamma_{i,r}}
\langle X(u_\infty^i),\eta_{i,r}\rangle\,ds_{g_{u_\infty^i}}.
\end{aligned}
\]
Smoothness and vanishing of the derivative at each branch point make the artificial boundary image lengths tend to zero as $r\downarrow0$.
Their integrals are bounded by $\|X\|_\infty$ times these lengths and vanish.
The omitted area tends to zero, and $\mathbf H_{u_\infty^i,M}\in L^2(d\mu_{u_\infty^i})$, so the first variation and interior integral converge to those on the full component.
The integrals over $\partial\Sigma_\infty^i\cap\partial\Omega_{i,r}$ converge because the smooth boundary map has finite length.
Letting $r\downarrow0$ and summing gives $\delta V_\infty(X)=-\mathcal H_\infty(X)+\mathcal B_\infty(X)$.
The absolute continuity proved below yields \eqref{eq:limit-varifold-first-variation}.

Let $\mathfrak{A}\subset M$ be a Borel set. Put
\[
\mu_i=(u_\infty^i)_\#\mu_{u_\infty^i},
\qquad
\|V_\infty\|=\sum_i\mu_i.
\]
If $\|V_\infty\|(\mathfrak{A})=0$, then every $\mu_i(\mathfrak{A})$ vanishes.
Cauchy--Schwarz applied to $(u_\infty^i)_\# (\mathbf H_{u_\infty^i,M}\,d\mu_{u_\infty^i})(\mathfrak{A})$ therefore proves $\mathcal H_\infty\ll\|V_\infty\|$.
The area formula identifies the integer sheet multiplicity $\Theta$ and the common approximate tangent plane at almost every image point; intersections of sheets with different tangent planes are $\|V_\infty\|$-null.
Applying the area formula to the vector measure gives the fiberwise average \eqref{eq:fiberwise-mean-curvature-average}.
Jensen's inequality on each finite fiber gives \eqref{eq:limit-mean-curvature-L2-bound}.

The total variation of $\mathcal B_\infty$ is bounded by the sum of the pushforward boundary length measures.
Its support is contained in the finite union of boundary images $u_\infty^i(\partial\Sigma_\infty^i)\subset\cK$.
This union is countably $1$-rectifiable and has zero $\mathcal H^2$-measure.
Since $\|V_\infty\|$ is absolutely continuous with respect to $\mathcal H^2$ on its countably $2$-rectifiable image, it follows that $\mathcal B_\infty\perp\|V_\infty\|$.
Sheetwise orthogonality gives $\mathcal B_\infty(X)=0$ for fields tangent to $\cK$ along the boundary images, without requiring these images to be embedded.

For an arbitrary $f\in C_c^\infty(M)$, choose on each compactified component $\Sigma_\infty^i$ a smooth metric $h_i$ compatible with its conformal structure.
Let $\alpha_\infty^i$ be the limiting form from Lemma~\ref{lem:hsl-local-convergence}, which is the mean curvature $1$-form of $u_\infty^i$ away from its branch points.
The map $u_\infty^i$ and form $\alpha_\infty^i$ are smooth there by Subsection~\ref{subsec:hsl-limit}, with $d_{h_i}^*\alpha_\infty^i=0$ and $\alpha_\infty^i(\nu_{h_i})=0$.
The identity \eqref{eq:first-variation-general-natural-support-2} and conformal invariance of the inner-product integral of $1$-forms give
\[
\begin{aligned}
 \int_{\Sigma_\infty^i}\langle\mathbf H_{u_\infty^i,M},J\nabla f(u_\infty^i)\rangle\,d\mu_{u_\infty^i} =-\int_{\Sigma_\infty^i}\langle\alpha_\infty^i,d(f\circ u_\infty^i)\rangle_{h_i}\,\dd\vol_{h_i}
=0.
\end{aligned}
\]
The last equality follows by integration by parts on $\Sigma_\infty^i$ with respect to the smooth metric $h_i$.
Summing over $i$ and using the Radon--Nikodym definition of $\overline{\mathbf H}_\infty$ proves \eqref{eq:all-potential-curvature-identity}.
For every $f\in C_c^\infty(M)$ whose Hamiltonian field is tangent to $\cK$ along the limiting boundary images, the preceding identity and $\mathcal B_\infty(J\nabla f)=0$ give $\delta V_\infty(J\nabla f)=0$ by \eqref{eq:limit-varifold-first-variation}.
\end{proof}

The construction in Section~\ref{sec:general-compactness-proof}, the HSL identification in Subsection~\ref{subsec:hsl-limit}, and Propositions~\ref{prop:area} and~\ref{prop:stationary} prove all conclusions of Theorem~\ref{main theorem}.

\Needspace{6\baselineskip}
\subsection{\texorpdfstring{Related compactness statements}{Related compactness statements}}
\label{subsec:hsl-consequences}

Corollary~\ref{u0} gives compactness for weak solutions of the homogeneous HSL system and preserves prescribed image points in the Hausdorff limit.
This includes the free boundary analogues of \cite{chen2014}.

\begin{coro}\label{u0}\label{u01}\label{exif}
Let $(M,\omega,J,\bar g)$ be a complete K\"ahler surface, let $\cK\subset M$ be a compact smooth support, and let $\Sigma$ be a compact connected bordered surface of fixed topological type.
Let $h_k$ be the representatives of constant curvature obtained from the complex doubles, and suppose that
\[
 u_k\in W^{2,2}_{b,c,\mathrm{fb}}\bigl((\Sigma,h_k),M,\cK\bigr)
\]
are Lagrangian immersions whose geometric mean curvature $1$-forms $\alpha_k$ satisfy \eqref{sys1}-\eqref{sys2} distributionally, including the orthogonal boundary condition for $u_k$ and zero normal trace for $\alpha_k$.
Assume that $u_k(\Sigma)\subset Q$ for a fixed compact set $Q\Subset M$ and that
\[
 \sup_k\bigl(\cA(u_k)+W_M(u_k)\bigr)\leq\Lambda<\infty.
\]
After passing to a subsequence, the sets $u_k(\Sigma)$ converge in Hausdorff distance to a compact connected set $L\subset Q$.
Either $L=\{q\}$ for some $q\in Q$, or there are a finite connected bordered stratified surface $\Sigma_\infty$, with $g(\Sigma_\infty)\leq g(\Sigma)$, and a continuous map $u_\infty:\Sigma_\infty\to M$ such that
\[
 L=u_\infty(\Sigma_\infty)
\]
and $u_k$ converges to $u_\infty$ in the sense of Definition~\ref{def:componentwise-free-boundary-bubble tree-convergence}.
In the latter alternative, every nonconstant component of $u_\infty$ is a smooth branched conformal HSL immersion satisfying \eqref{sys1}-\eqref{sys2}, constant components may be retained as ghost components, and
\[
 \lim_{k\to\infty}\cA(u_k)=\cA(u_\infty),
 \qquad
 W_M(u_\infty)\leq\liminf_{k\to\infty}W_M(u_k).
\]
In the noncollapsed alternative, the pushforward area measures converge as in Proposition~\ref{prop:area}.
If $y_1,\ldots,y_m\in u_k(\Sigma)$ for every $k$, then $y_1,\ldots,y_m\in L$.
In particular, two distinct fixed image points rule out collapse to a point.
\end{coro}

\begin{proof}
The qualitative regularity argument following Theorem~\ref{higher}, together with its interior analogue using \eqref{eq:hodge-local-bootstrap}, makes each weak input pair smooth locally, including across its finite branch set and up to the inherited boundary.
Pass to a bordered Deligne--Mumford subsequence as in Subsection~\ref{subsec:generic-extraction}.
The embedding comparison \eqref{eq:sec2-willmore-comparison} and Theorem~\ref{thm:intro-main-compactness} give a bubble tree subsequence or point collapse.
Lemma~\ref{lem:hsl-local-convergence} and Subsection~\ref{subsec:hsl-limit} identify and extend the HSL limits and give intrinsic Willmore lower semicontinuity.
Subsection~\ref{subsec:compactification} gives the genus bound, and Proposition~\ref{prop:area} gives the area identity and convergence of the area measures.

In either alternative, the general compactness theorem gives Hausdorff convergence of the images to $L$.
Since each $u_k(\Sigma)$ is nonempty, compact, and connected, $L\subset Q$ is compact and connected.
For each fixed image point,
\[
\operatorname{dist}(y_j,L)
\leq d_{\mathrm{Haus}}\bigl(u_k(\Sigma),L\bigr)\longrightarrow0.
\]
Since $L$ is closed, it contains every $y_j$; two distinct such points exclude a singleton limit.
\end{proof}

The fixed compact image set gives the uniform bounds on the geometry of $(M,\bar g)$ used above.
Requiring only that each image meet a fixed extrinsic ball does not give these bounds for a general complete target.

\vspace{1cm}

\section{Applications of bubble tree compactness}\label{sec:applications}

In this section, we prove Theorems~\ref{app:small-annuli}, \ref{app:small-disks}, \ref{app:calibration-barrier}, and~\ref{app:hsl-energy}. We begin with three auxiliary lemmas on weak convergence of pushforward boundary arclength measures, approximation by smooth immersions, and local densities of pushforward area measures. We then treat Willmore energy minimization for annuli and disks in Subsections~\ref{sp:section} and~\ref{ga:section}, the quantitative Willmore lower bound under a calibration constraint in Subsection~\ref{bar:section}, and the decomposition of the limiting intrinsic Willmore energy for free boundary HSL sequences in Subsection~\ref{hs:section}.

For a branched immersion $u:\Sigma\to(N,\bar g)$, we use the induced metric $g_u$ and finite branch set $\cB_u$ from Subsection~\ref{sec:fbbci}. The weak branched conformal immersions considered below belong to $W^{2,2}(\Sigma,N)$. On $\overline\Sigma\setminus\cB_u$, their induced metrics are locally bounded above and below by positive multiples of a smooth metric on $\Sigma$. Branched conformal immersions with free boundary are understood in the sense of Definition~\ref{def:fb-branched-conformal}. On a compact domain, an embedding in this weak class means a continuous injective map and smoothness is stated separately.

We use the induced boundary arclength measure $\dd s_{g_u}$ and total boundary length $L(u)$ from \eqref{eq:induced-boundary-geometry}.
For Euclidean targets, $W(u)$ has the normalization in \eqref{eq:intro-energy-normalization}. All integrals and pushforward measures retain the multiplicity of the parametrization.

Let $h$ be a smooth metric on $\Sigma$ with respect to which
$u$ is weakly conformal, and let $\nu_h$ and $\tau_h$ be the outward unit conormal and a unit boundary tangent with respect to $h$. The unit image conormal $\eta_u$ and unit image tangent $T$ along $\partial\Sigma\setminus\cB_u$ are those in \eqref{eq:induced-boundary-geometry}. When the supporting submanifold $\mathcal{K}$ is a two-sided hypersurface $S\subset\R^K$, choose a smooth unit normal $n_S$ and use the convention for $h^S$ in \eqref{ga:support-curvature}, with $\nabla^{\mathbb R^K}$ in place of $\nabla^{\mathbb R^3}$. If $\eta_u=n_S\circ u$, the geodesic curvature of $\partial\Sigma$ with respect to $g_u$, computed with the inward unit conormal, satisfies $k_{g_u}=h^S(T,T)$ by Lemma~\ref{lem:free-boundary-gauss-bonnet}. We retain this conormal sign convention when identifying the positive measure $u_\#(\dd s_{g_u})$ on $S$ from the first variation of area.

The next lemma allows us to retain a prescribed boundary length in the bubble tree limit. Under a common conormal sign convention, the first variation of area identifies the weak limit of the pushforward boundary arclength measures on $S$. Strong convergence of the tangential derivatives of the parametrized boundary maps also identifies the limit of the corresponding measures on the bundle of unoriented tangent lines.

\begin{lemma}\label{lem:app-contact-measures}
Let $S\subset\R^K$ be a smooth compact embedded two-sided supporting hypersurface with a chosen smooth unit normal $n_S$. Let $u_j:(\Sigma,h_j)\to\R^K$ satisfy the hypotheses of Theorem~\ref{thm:intro-main-compactness}, with free boundary $u_j(\partial\Sigma)$ on $S$, and suppose that the noncollapsed alternative of Theorem \ref{thm:intro-main-compactness} holds. Assume
\[
 \eta_{u_j}=n_S\circ u_j
 \quad\text{a.e. on }\partial\Sigma\setminus\cB_{u_j}.
\]
Write $u_\infty^i:\Sigma_\infty^i\to\R^K$ for the nonconstant limiting bubble tree components. Then we have
\begin{equation}\label{eq:contact-scalar}
 (u_j)_\#\dd s_{g_{u_j}}\rightharpoonup
 \sum_i(u_\infty^i)_\#\dd s_{g_{u_\infty^i}},
 \qquad L(u_j)\longrightarrow\sum_i L(u_\infty^i).
\end{equation}
More generally, let $\Psi(x,\ell)$ be continuous on the bundle whose fiber at $x\in S$ consists of the one-dimensional unoriented subspaces $\ell\subset T_xS$. Then
\begin{equation}\label{eq:contact-tangent}
 \int_{\partial\Sigma}\Psi(u_j,[T_j])\dd s_{g_{u_j}}
 \longrightarrow
 \sum_i\int_{\partial\Sigma_\infty^i}
                 \Psi(u_\infty^i,[T_i])\dd s_{g_{u_\infty^i}}.
\end{equation}
Here $T_j$ and $T_i$ are the unit tangent vectors to the parametrized boundary maps $u_j|_{\partial\Sigma}$ and $u_\infty^i|_{\partial\Sigma_\infty^i}$, and $[T_j]$ and $[T_i]$ are the unoriented lines they span. The convergence in \eqref{eq:contact-scalar} is the weak convergence of finite Radon measures on $S$. In particular,
\[
 \int_{\partial\Sigma}k_{g_{u_j}}\,\dd s_{g_{u_j}}
 \longrightarrow
 \sum_i\int_{\partial\Sigma_\infty^i}k_{g_{u_\infty^i}}\,\dd s_{g_{u_\infty^i}},
\]
where $k_{g_{u_j}}$ is the geodesic curvature of $\partial\Sigma$ with respect to $g_{u_j}$ and $k_{g_{u_\infty^i}}$ is the geodesic curvature of $\partial\Sigma_\infty^i$ with respect to $g_{u_\infty^i}$, both computed with the inward unit conormal.
\end{lemma}

\begin{proof}
Let $u$ denote one of the maps $u_j$ or one of the nonconstant limiting maps, and let $h$ be a smooth metric on its domain with respect to which $u$ is weakly conformal. Away from the finite branch set, write $g_u=e^{2\lambda}h$. Conformality gives
\[
 \operatorname{div}_uX\,\dd\mu_u
 =\langle du,d(X\circ u)\rangle_h\,\dd\vol_h,
 \qquad
 \Delta_hu\,\dd\vol_h=\mathbf H_u\,\dd\mu_u
\]
for $X\in C_c^\infty(\R^K,\R^K)$. Using \eqref{eq:induced-boundary-geometry}, we obtain
\begin{equation}\label{eq:app-contact-first-variation}
\begin{aligned}
 \int_\Sigma\operatorname{div}_uX\,\dd\mu_u=-\int_\Sigma\mathbf H_u\cdot X(u)\,\dd\mu_u
   +\int_{\partial\Sigma}X(u)\cdot\eta_u\,\dd s_{g_u}.
\end{aligned}
\end{equation}
For a limiting bubble tree component, $\Sigma$ in this identity is replaced by $\Sigma_\infty^i$. 

We next show that every nonconstant bordered limiting component has the prescribed conormal sign. Choose compact cores $K_{i,r}$ avoiding the nodal, attaching, and exceptional points in Definition~\ref{def:componentwise-free-boundary-bubble tree-convergence}, as well as the branch points of $u_\infty^i$. As $r\downarrow0$, these cores exhaust the complement of those finitely many points. We choose their intersections with $\partial\Sigma_\infty^i$ to be finite unions of compact boundary arcs. For fixed $r$ and sufficiently large $j$, let $\Phi_j^i$ be the comparison embeddings from Definition~\ref{def:componentwise-free-boundary-bubble tree-convergence}, restricted to a neighborhood of $K_{i,r}$. Their images are pairwise disjoint, and they map boundary arcs into $\partial\Sigma$.

Fix a limiting bubble tree component and a compact core, and write $\widetilde u_j=u_j\circ\Phi_j^i$. As in Subsection~\ref{subsec:generic-extraction}, choose smooth metrics $\widetilde h_j$ conformally equivalent to $(\Phi_j^i)^*h_j$ and converging smoothly to a metric $h_i$ representing the limiting conformal structure. On a neighborhood of the core, we have
\[
 \widetilde u_j\rightharpoonup u_\infty^i
 \quad\text{in }W^{2,2},\qquad
 \widetilde u_j\longrightarrow u_\infty^i
 \quad\text{in }W^{1,p}\text{ for every finite }p>1
 \text{ and uniformly}.
\]
On a compact boundary coordinate interval $I$, the trace theorem bounds each coordinate first derivative of $\widetilde u_j$ in $W^{1/2,2}(I)$. The compact embedding $W^{1/2,2}(I)\hookrightarrow L^2(I)$, together with the weak convergence of these traces to those of $u_\infty^i$, gives strong $L^2(I)$ convergence; see \cite{AdamsFournier2003}. Smooth convergence of the metrics $h_j$ therefore gives
\[
\begin{aligned}
 d\widetilde u_j(\nu_{\widetilde h_j})
 \longrightarrow du_\infty^i(\nu_{h_i}),\qquad
 d\widetilde u_j(\tau_{\widetilde h_j}) \longrightarrow du_\infty^i(\tau_{h_i})
\end{aligned}
\qquad\text{strongly in }L^2(I).
\]
The assumed sign condition, pulled back by $\Phi_j^i$, is
\[
 d\widetilde u_j(\nu_{\widetilde h_j})
 =\bigl|d\widetilde u_j(\tau_{\widetilde h_j})\bigr|
       n_S(\widetilde u_j).
\]
Uniform convergence of $\widetilde u_j$ and continuity of $n_S$ allow us to pass to the limit in this equality in $L^2(I)$ and we obtain
\[
 du_\infty^i(\nu_{h_i})
 =\bigl|du_\infty^i(\tau_{h_i})\bigr|n_S(u_\infty^i).
\]
The limiting map is an immersion away from its finite branch set, so this proves $\eta_{u_\infty^i}=n_S\circ u_\infty^i$ almost everywhere on its boundary. 

To prove convergence of the scalar boundary measures, we pass to the limit in the two interior integrals of \eqref{eq:app-contact-first-variation}. In local coordinates on a fixed core, conformality of $\widetilde{u}_j$ gives
\[
\begin{aligned}
 \operatorname{div}_{\widetilde u_j}X\,\dd\mu_{\widetilde u_j}
 =\sum_{a,b=1}^2\widetilde h_j^{ab}
   \bigl\langle DX(\widetilde u_j)\partial_a\widetilde u_j,
                      \partial_b\widetilde u_j\bigr\rangle
                        \,\dd\vol_{\widetilde h_j},\quad
 \mathbf H_{\widetilde u_j}\,\dd\mu_{\widetilde u_j}
 =\Delta_{\widetilde h_j}\widetilde u_j\,\dd\vol_{\widetilde h_j}.
\end{aligned}
\]
The first density converges in $L^1$. For the second identity, the coordinate formula
\[
 \Delta_{\widetilde h_j}\widetilde u_j
 =|\widetilde h_j|^{-1/2}
   \sum_{a,b=1}^2\partial_a
   \left(|\widetilde h_j|^{1/2}\widetilde h_j^{ab}
                         \partial_b\widetilde u_j\right),
 \qquad |\widetilde h_j|=\det(\widetilde h_{j,ab}),
\]
shows weak $L^2$ convergence to $\Delta_{h_i}u_\infty^i$. Pairing with $X(\widetilde u_j)$, which converges strongly in $L^2$, proves convergence of the mean curvature integral on the core.

Then, we make the passage from cores to the full domains. Set
\[
\begin{aligned}
 a_j(r)=\cA\left(u_j;\Sigma\setminus
                     \bigcup_i\Phi_j^i(K_{i,r})\right),\qquad
 a_\infty(r)=\sum_i\cA\left(u_\infty^i;
                           \Sigma_\infty^i\setminus K_{i,r}\right).
\end{aligned}
\]
Strong convergence of the first derivatives of $u_j$ gives convergence of the area on each fixed core. Subtracting these core areas from the area identity in Theorem~\ref{thm:intro-main-compactness} yields
\[
 \lim_{j\to\infty}a_j(r)=a_\infty(r),\qquad
 \lim_{r\downarrow0}a_\infty(r)=0.
\]
Choose $\Lambda<\infty$ such that $W(u_j)\leq\Lambda$ for every $j$. By the lower semicontinuity of Willmore energy, we have  $\sum_iW(u_\infty^i)\leq\Lambda$. For any measurable region $U\subset\Sigma$ of area $a=\cA(u_j;U)$, we have
\[
\begin{aligned}
 \left|\int_U\operatorname{div}_{u_j}X\,\dd\mu_{u_j}\right|
 &\leq2\|DX\|_\infty a,\\
 \left|\int_U\mathbf H_{u_j}\cdot X(u_j)\,\dd\mu_{u_j}\right|
 &\leq\|X\|_\infty
       \left(\int_U|\mathbf H_{u_j}|^2\,\dd\mu_{u_j}\right)^{1/2}
       a^{1/2}
 \leq2\|X\|_\infty\sqrt{\Lambda a}.
\end{aligned}
\]
The corresponding bounds hold for the sum over omitted limiting regions, with $a_\infty(r)$ in place of $a$. Thus both omitted interior integrals tend to zero when $j\to\infty$ and then $r\downarrow0$. Applying \eqref{eq:app-contact-first-variation} to the sequence $u_j$ and to the limiting bubble tree components now gives
\[
 \int_{\partial\Sigma}X(u_j)\cdot n_S(u_j)\,\dd s_{g_{u_j}}
 \longrightarrow
 \sum_i\int_{\partial\Sigma_\infty^i}
       X(u_\infty^i)\cdot n_S(u_\infty^i)\,\dd s_{g_{u_\infty^i}}.
\]

For $\varphi\in C^\infty(S)$, choose a compactly supported smooth ambient vector field $X$ whose restriction to $S$ is $\varphi n_S$.  Since $|n_S|=1$, the preceding limit is exactly the scalar measure convergence in \eqref{eq:contact-scalar} tested against $\varphi$. Taking $\varphi=1$ gives convergence of the total boundary length in \eqref{eq:contact-scalar}. Moreover, \eqref{eq:sec2-boundary-length-area-willmore} bounds $L(u_j)$ uniformly. Since smooth functions are uniformly dense in $C(S)$, the boundary arc length measure convergence therefore holds for every continuous test function on $S$, proving \eqref{eq:contact-scalar}.

We finally retain the unoriented tangent lines. On a compact boundary coordinate interval $I$ of a core, put
\[
 v_j=\partial_x\widetilde u_j|_I,\qquad
 v=\partial_xu_\infty^i|_I,
\]
which are the derivatives of the parametrized boundary curves, and $v_j\to v$ strongly in $L^2(I)$. Their induced arclength elements are $|v_j|\,dx$ and $|v|\,dx$, respectively. Since the boundary curves take values in $S$, their derivatives are tangent to $S$ almost everywhere. Hence the pulled-back integrands in \eqref{eq:contact-tangent} are
\[
 F_j(x)=|v_j(x)|\Psi\bigl(\widetilde u_j(x),[v_j(x)]\bigr),
 \qquad
 F(x)=|v(x)|\Psi\bigl(u_\infty^i(x),[v(x)]\bigr),
\]
with value zero wherever the corresponding derivative vanishes. The tangential line bundle over the compact hypersurface $S$ is compact, so $\Psi$ is bounded. The function $(y,w)\mapsto|w|\Psi(y,[w])$, for $w\in T_yS$, extends continuously by zero at $w=0$. After taking a subsequence, strong $L^2$ convergence permits a further subsequence on which $v_j\to v$ almost everywhere. Thus $F_j\to F$ almost everywhere along that subsequence. For every measurable $A\subset I$,
\[
 \int_A|F_j|\,dx
 \leq\|\Psi\|_\infty\|v_j\|_{L^2(I)}|A|^{1/2}.
\]
This bound gives uniform integrability of $F_j$, so $F_j\to F$ in $L^1(I)$. The subsequence argument shows that the integral converges along the original sequence. A finite boundary coordinate covering gives convergence on the boundary portions of every fixed core. 
In order to passage the convergence from cores to the full domains, Write $B_{i,r}=K_{i,r}\cap\partial\Sigma_\infty^i$, and denote the omitted boundary lengths by
\[
\begin{aligned}
 \ell_j(r)=\int_{\partial\Sigma\setminus
                   \bigcup_i\Phi_j^i(B_{i,r})}\dd s_{g_{u_j}},\qquad
 \ell_\infty(r)=\sum_i\int_{\partial\Sigma_\infty^i\setminus B_{i,r}}
                             \dd s_{g_{u_\infty^i}}.
\end{aligned}
\]
The comparison images are pairwise disjoint. Subtracting the converging core lengths from the total length convergence in \eqref{eq:contact-scalar} gives
\[
 \lim_{j\to\infty}\ell_j(r)=\ell_\infty(r),\qquad
 \lim_{r\downarrow0}\ell_\infty(r)=0.
\]
The integral over the omitted part of $\partial\Sigma$ has absolute value at most $\|\Psi\|_\infty\ell_j(r)$. The sum of the corresponding limiting integrals has absolute value at most $\|\Psi\|_\infty\ell_\infty(r)$. Combining these bounds with above convergence on the compact cores, and then letting $r\downarrow0$, proves \eqref{eq:contact-tangent}.

For a unit vector $T\in T_xS$, the value $h^S_x(T,T)$ depends only on the unoriented line $[T]$ and is continuous on this bundle. Applying \eqref{eq:contact-tangent} with $\Psi(x,[T])=h^S_x(T,T)$, and using Lemma~\ref{lem:free-boundary-gauss-bonnet}, proves convergence of the stated geodesic curvature integrals.
\end{proof}

Kuwert--M\"uller~\cite[Definition~2.2 and Lemma~2.8]{appKM} include the contact-line varifold in their compactness theory of curvature varifold with orthogonal boundary.
Lemma~\ref{lem:app-contact-measures} identifies the limit of this measure with the sum of the measures induced by the boundary maps of the bubble tree components, including for branched immersions and degenerating domains.
In particular, the constraint $L(u_j)=\ell$ passes to $\sum_iL(u_\infty^i)=\ell$. The next lemma provides smooth approximations of $W^{2,2}$ immersions whose induced metrics are bounded above and below by positive multiples of a smooth metric on $\Sigma$, with convergence of area and Willmore energy.

\begin{lemma}\label{lem:app-weak-approximation}
Let $(\Sigma,h)$ be a smooth compact surface, possibly with boundary, and let $(N,\bar g)$ be a compact smooth Riemannian manifold or Euclidean space.
Suppose $v\in W^{2,2}(\Sigma,N)$ and
\[
 ch\le v^*\bar g\le Ch\quad\text{a.e.},\qquad \text{for some } 0<c<C<\infty.
\]
There are smooth immersions $v_k:\Sigma\to N$ converging strongly in $W^{2,2}$ and uniformly to $v$, with their induced metrics bounded above and below by positive multiples of $h$, uniformly in $k$.
In addition, $\cA(v_k)\to\cA(v)$, $W_N(v_k)\to W_N(v)$, and $\int_\Sigma v_k^*\omega\to\int_\Sigma v^*\omega$ for every fixed smooth $2$-form $\omega$ on $N$. In particular, for every fixed smooth boundary arc $I\subset\partial\Sigma$, its induced length under $v_k$ converges to its induced length under $v$.
\end{lemma}

\begin{proof}
When $N$ is compact, fix an isometric embedding $N\hookrightarrow\R^K$ and regard $v$ as an $\R^K$-valued map.
If $\Sigma$ has boundary, compose $v$ with the small positive time flow of a smooth inward vector field on a double of $\Sigma$. This defines $w_\varepsilon$ on a neighborhood of $\Sigma$ using only interior values of $v$.
The maps $w_{\varepsilon}$ converge strongly in $W^{2,2}$ and uniformly to $v$ as $\varepsilon \downarrow 0$, while their upper and lower metric bounds persist uniformly for small $\varepsilon$. If $\partial\Sigma=\varnothing$, set $w_\varepsilon=v$.

Fix $\varepsilon>0$. We construct the smooth approximations in local coordinates, keeping a uniform lower bound for their differentials. Choose a finite smooth coordinate atlas $\{(U_a,\kappa_a)\}$ on a neighborhood of $\Sigma$ in its complex double, smaller sets $V_a\Subset U_a$ covering $\Sigma$, and a nonnegative smooth partition of unity $\{\chi_a\}$ with $\operatorname{supp}\chi_a\Subset U_a$ and $\sum_a\chi_a=1$ near $\Sigma$. In the closed case, use an atlas on $\Sigma$ itself. The coordinate expressions of $w_\varepsilon$ belong to $W^{2,2}\cap W^{1,\infty}$. The metric bounds give constants $c_0,C_0>0$ independent of small $\varepsilon$ such that, in each of these coordinates, there holds
\[
 c_0|\xi|\leq|Dw_\varepsilon(x)\xi|\leq C_0|\xi|
 \quad\text{for every }\xi\in\R^2
 \text{ and almost every }x.
\]

We first control the oscillation of $Dw_\varepsilon$ on small coordinate balls. For a ball $B_r(x)$ contained in a chart, write $(Dw_\varepsilon)_{B_r(x)}$ for the average of $Dw_\varepsilon$ over the ball. Poincar\'e's inequality in dimension two gives
\[
 \left(\fint_{B_r(x)}
   |Dw_\varepsilon-(Dw_\varepsilon)_{B_r(x)}|^2\right)^{1/2}
 \leq C\left(\int_{B_r(x)}|D^2w_\varepsilon|^2\right)^{1/2}.
\]
For fixed $\varepsilon$, the right-hand side tends to zero uniformly in the center $x$ as $r\downarrow0$. Let $\rho\in C_c^\infty(B_1(0))$ be nonnegative with $\int_{\R^2}\rho=1$, and put $\rho_\delta(z)=\delta^{-2}\rho(z/\delta)$. For $\delta>0$ sufficiently small, define
\[
 z_a=\bigl(\rho_\delta*(w_\varepsilon\circ\kappa_a^{-1})\bigr)
                  \circ\kappa_a,
 \qquad
 z_{\varepsilon,\delta}=\sum_a\chi_a z_a.
\]
Only the values of $z_a$ near $\operatorname{supp}\chi_a$ are used, and we choose $\delta$ so that each convolution uses values within the neighborhood where $w_\varepsilon$ is defined. Convergence of convolution operator in $W^{2,2}$ on compact subcharts, followed by multiplication by the fixed functions $\chi_a$, gives
\[
 z_{\varepsilon,\delta}\longrightarrow w_\varepsilon
 \quad\text{strongly in }W^{2,2}(\Sigma,\R^K)
 \text{ and uniformly as }\delta\downarrow0.
\]
The uniform convergence also follows directly from the Lipschitz bound of $w_{\varepsilon}$, which gives $|z_a(x)-w_\varepsilon(x)|\leq C\delta$.

Next, we  verify that $z_{\varepsilon,\delta}$ is an immersion. Use a chart $(U_b,\kappa_b)$, consider the point lies in $V_b$, and write $x,y$ for its coordinate variables. For sufficiently small $\delta$, all convolution regions contributing at this point lie in $U_b$ and in the corresponding overlapping charts. Set $F_{ab}=\kappa_a\circ\kappa_b^{-1}$ on each such overlap. The chain rule and the change of variables from $\kappa_a$ to $\kappa_b$ give
\[
 Dz_a(x)=\int k_{a,\delta,x}(y)\,Dw_\varepsilon(y)
                  DF_{ab}(y)^{-1}DF_{ab}(x)\,dy,
\]
where the maps are written in the $\kappa_b$ coordinates and
\[
 k_{a,\delta,x}(y)
 =\rho_\delta\bigl(F_{ab}(x)-F_{ab}(y)\bigr)
                         |\det DF_{ab}(y)|.
\]
Each convolution kernel is nonnegative, has integral one, is supported in $B_{C\delta}(x)$, and is bounded by $C\delta^{-2}$.  On its support, we have
\[
 |DF_{ab}(y)^{-1}DF_{ab}(x)-I|\leq C\delta.
\]
Moreover, differentiating the partition of unity yields
\[
 Dz_{\varepsilon,\delta}
 =\sum_a\chi_a Dz_a
   +\sum_a(z_a-w_\varepsilon)\otimes d\chi_a,
\]
Here, the second sum is bounded by $C\delta$. Thus, we obtain
\begin{equation}\label{eq:app-approx-averaged-differential}
\begin{gathered}
 Dz_{\varepsilon,\delta}(x)
   =\int K_{\delta,x}(y)Dw_\varepsilon(y)\,dy+R_{\varepsilon,\delta}(x),
 \qquad |R_{\varepsilon,\delta}(x)|\leq C\delta,\qquad K_{\delta,x}:=\sum_a\chi_a(x)k_{a,\delta,x},\\
 K_{\delta,x}\geq0,\qquad \int K_{\delta,x}=1,\qquad
 K_{\delta,x}\leq C\delta^{-2}.
\end{gathered}
\end{equation}
The support of $K_{\delta,x}$ is contained in $B=B_{C\delta}(x)$, by Cauchy--Schwarz inequality, we have
\[
 \left|\int K_{\delta,x}Dw_\varepsilon-(Dw_\varepsilon)_B\right|
 \leq C\left(\fint_B|Dw_\varepsilon-(Dw_\varepsilon)_B|^2\right)^{1/2}.
\]
There is a Lebesgue point $y\in B$ at which the metric bounds hold and
\[
 |Dw_\varepsilon(y)-(Dw_\varepsilon)_B|
 \leq2\left(\fint_B|Dw_\varepsilon-(Dw_\varepsilon)_B|^2\right)^{1/2}.
\]
Combining these two estimates with \eqref{eq:app-approx-averaged-differential} gives
\[
 |Dz_{\varepsilon,\delta}(x)-Dw_\varepsilon(y)|
 \leq C\delta+C\left(\int_{B_{C\delta}(x)}|D^2w_\varepsilon|^2\right)^{1/2}.
\]
For fixed $\varepsilon$, this bound tends to zero uniformly in $x$ as $\delta \downarrow 0$. Choosing $\delta$ so that it is smaller than $c_0/2$, we obtain
\[
 |Dz_{\varepsilon,\delta}(x)\xi|\geq\frac{c_0}{2}|\xi|
\]
for every $x\in\Sigma$ and $\xi\in\R^2$.

When $N$ is compact, let $\Pi$ be the smooth nearest point projection from a fixed tubular neighborhood of $N\subset\R^K$ onto $N$, and write $P_N(q)$ for the orthogonal projection onto $T_qN$. Since $w_\varepsilon$ takes values in $N$ and $z_{\varepsilon,\delta}-w_\varepsilon=O(\delta)$ uniformly, the map $w_{\varepsilon,\delta}=\Pi\circ z_{\varepsilon,\delta}$ is defined for sufficiently small $\delta$. The differential of $w_\varepsilon$ is tangent to $N$ almost everywhere. If $y$ belongs to the support of $K_{\delta,x}$, the Lipschitz bound and smoothness of $P_N$ give
\[
\begin{aligned}
 |(I-P_N(w_\varepsilon(x)))Dw_\varepsilon(y)|
 &=|(P_N(w_\varepsilon(y))-P_N(w_\varepsilon(x)))Dw_\varepsilon(y)|\\
 &\leq C|w_\varepsilon(y)-w_\varepsilon(x)|\\
 &\leq C\delta.
\end{aligned}
\]
It follows from \eqref{eq:app-approx-averaged-differential} that the same bound holds with $Dz_{\varepsilon,\delta}(x)$ in place of $Dw_\varepsilon(y)$. Since $D\Pi(q)=P_N(q)$ for $q\in N$, we also have
\[
 |D\Pi(z_{\varepsilon,\delta}(x))-P_N(w_\varepsilon(x))|
 \leq C|z_{\varepsilon,\delta}(x)-w_\varepsilon(x)|\leq C\delta.
\]
Consequently, in this case, we obtain
\[
 |Dw_{\varepsilon,\delta}-Dz_{\varepsilon,\delta}|\leq C\delta.
\]
After decreasing $\delta$ once more, the lower differential bound of $w_{\varepsilon,\delta}$ is therefore at least $c_0/4$. If $N$ is Euclidean, set $\Pi=\operatorname{id}$ and $w_{\varepsilon,\delta}=z_{\varepsilon,\delta}$. Converting the coordinate bounds to the metric $h$, in either case we obtain
\begin{equation}\label{eq:app-approx-uniform-metric}
 c_1h\leq w_{\varepsilon,\delta}^*\bar g\leq C_1h
 \qquad\text{on }\Sigma,
\end{equation}
where $c_1,C_1>0$ are independent of small $\varepsilon$ and of the sufficiently small $\delta$ chosen for that $\varepsilon$.

Next, we show that the nearest point projection also preserves the strong $W^{2,2}$ convergence to $w_\varepsilon$. In each coordinate chart, we have
\[
 D^2(\Pi\circ z)
 =D\Pi(z)D^2z+D^2\Pi(z)[Dz,Dz].
\]
Apply this identity to $z=z_{\varepsilon,\delta}$ and to $z=w_\varepsilon$, noting that $\Pi\circ w_\varepsilon=w_\varepsilon$. The first terms converge in $L^2$ because $D^2z_{\varepsilon,\delta}\to D^2w_\varepsilon$ in $L^2$ and $D\Pi(z_{\varepsilon,\delta})\to D\Pi(w_\varepsilon)$ uniformly. For the quadratic terms, by the Sobolev embedding $W^{2,2}\hookrightarrow W^{1,4}$, we have
\[
 \|Dz_{\varepsilon,\delta}\otimes Dz_{\varepsilon,\delta}
                 -Dw_\varepsilon\otimes Dw_\varepsilon\|_{L^2}
 \leq\bigl(\|Dz_{\varepsilon,\delta}\|_{L^4}
                  +\|Dw_\varepsilon\|_{L^4}\bigr)
       \|Dz_{\varepsilon,\delta}-Dw_\varepsilon\|_{L^4}\longrightarrow0.
\]
Uniform convergence of $D^2\Pi(z_{\varepsilon,\delta})$ then gives convergence of the quadratic terms in $L^2$. The zeroth- and first-order terms converge as well. Thus $w_{\varepsilon,\delta}\to w_\varepsilon$ strongly in $W^{2,2}$ and uniformly as $\delta\downarrow0$ with $\varepsilon$ fixed.

Choose $\varepsilon_k\downarrow0$ and then $\delta_k\downarrow0$ so that \eqref{eq:app-approx-uniform-metric} holds and
\[
 \|w_{\varepsilon_k,\delta_k}-w_{\varepsilon_k}\|_{W^{2,2}(\Sigma)}
 +\|w_{\varepsilon_k,\delta_k}-w_{\varepsilon_k}\|_{C^0(\Sigma)}
 <\frac1k.
\]
Setting $v_k=w_{\varepsilon_k,\delta_k}$ gives smooth immersions with the required uniform pulled-back metric bounds, strong $W^{2,2}$ convergence, and uniform convergence to $v$.

We next verify convergence of the associated geometric integrals of $v_k$. After passing to a subsequence, $dv_k\to dv$ almost everywhere. The metric bounds \eqref{eq:app-approx-uniform-metric} imply that the matrices $g_{v_k}$ and their inverses are uniformly bounded in each fixed coordinate chart, and they converge almost everywhere to $g_v$ and its inverse. For $q=v_k$ or $q=v$, let $P_q$ denote the orthogonal projection onto $dq(T\Sigma)$, so that
\[
 P_q=\sum_{a,b=1}^2g_q^{ab}\,\partial_aq\otimes\partial_bq.
\]
With $P_N=I$ in the Euclidean case, the second fundamental form of the immersion in $N$ gives the coordinate formula
\[
 \mathbf H_{q,N}
 =\sum_{a,b=1}^2 g_q^{ab}\bigl(P_N(q)-P_q\bigr)\partial_{ab}q.
\]
The coefficients in this formula are uniformly bounded and converge almost everywhere. Strong $L^2$ convergence of $D^2v_k$ therefore implies
\[
 \mathbf H_{v_k,N}\longrightarrow\mathbf H_{v,N}
 \qquad\text{strongly in }L^2(\Sigma,\dd\vol_h).
\]
Since the area densities $\sqrt{\det g_{v_k}/\det h}$ are uniformly bounded and converge almost everywhere, we have
\[
 \cA(v_k)\longrightarrow\cA(v),\qquad
 W_N(v_k)=\frac14\int_\Sigma|\mathbf H_{v_k,N}|^2
                    \sqrt{\frac{\det g_{v_k}}{\det h}}\,\dd\vol_h
          \longrightarrow W_N(v).
\]
For a fixed smooth $2$-form $\omega$ on $N$, use $\Pi^*\omega$ to extend it to a tubular neighborhood when $N$ is compact. Its pullback in an oriented chart is
\[
 v_k^*\omega
 =\omega_{v_k}(\partial_1v_k,\partial_2v_k)\,dx^1\wedge dx^2.
\]
The coefficients of $\omega$ converge uniformly along $v_k\to v$, and products of the first derivatives converge in $L^1$ by strong $L^2$ convergence. Summing over a finite partition of unity proves $\int_\Sigma v_k^*\omega\to\int_\Sigma v^*\omega$. In the Euclidean case, uniform convergence keeps all images of $v_k$ in a fixed compact set, so the same coefficient bounds apply.

Finally, let $I\subset\partial\Sigma$ be a fixed smooth boundary arc and let $\tau_h$ be a unit tangent to $I$ with respect to $h$. The bounded trace map $W^{1,2}(\Sigma)\to W^{1/2,2}(I)$ gives
\[
 \|dv_k(\tau_h)-dv(\tau_h)\|_{W^{1/2,2}(I)}
 \leq C_I\|v_k-v\|_{W^{2,2}(\Sigma)}\longrightarrow0.
\]
Taking Euclidean norms through the fixed isometric embedding, and writing $|I|_h=\int_I\dd s_h$, we obtain
\[
\begin{aligned}
 \left|\int_I|dv_k(\tau_h)|_{\bar g}\,\dd s_h
          -\int_I|dv(\tau_h)|_{\bar g}\,\dd s_h\right| \leq&\int_I|dv_k(\tau_h)-dv(\tau_h)|\,\dd s_h\\
 \leq &|I|_h^{1/2}
          \|dv_k(\tau_h)-dv(\tau_h)\|_{L^2(I)} \longrightarrow0.
\end{aligned}
\]
This completes the proof of Lemma \ref{lem:app-weak-approximation}. 
\end{proof}

To apply the Li--Yau multiplicity inequality, we next compute the contribution of a single preimage to the density of the pushforward area measure, including for the double obtained by flattening the supporting hypersurface and reflecting across the resulting hyperplane.

\begin{lemma}\label{lem:app-local-multiplicity}
Let $u\in W^{2,2}(\Sigma,\R^K)$ be a branched conformal immersion with finite total curvature. Choose a conformal disk or half-disk coordinate $z=x+iy$ centered at an interior or boundary point of $\Sigma$, respectively, and suppose
\[
 u^*g_{\R^K}=|z|^{2m}e^{2\omega(z)}|\dd z|^2,
 \qquad m\in\N\cup\{0\},\qquad
 \omega\in C^0\cap W^{1,2}.
\]
In the boundary case, let $S\subset\R^K$ be a smooth embedded supporting hypersurface, and assume that $u$ maps the boundary interval of the chart into $S$ and that the Sobolev boundary trace of $\partial_yu$ is orthogonal to $T_uS$ on this interval.
For a pushforward area measure $\mu$, its density at an image point $y$ is $\lim_{r\downarrow0}\mu(B_r(y))/(\pi r^2)$ whenever this limit exists.
After shrinking the coordinate neighborhood about $0$ and restricting the area measure to this neighborhood before taking its pushforward, the following statements hold.
\begin{enumerate}[label=\textup{(\roman*)},leftmargin=2.4em]
\item The pushforward area measure has density $m+1$ at $u(0)$ in the interior case and $(m+1)/2$ at $u(0)$ in the boundary case. At a point outside the branch set, $m=0$.
\item In the interior case, if $F$ is a smooth diffeomorphism between neighborhoods in $\R^K$, the restricted pushforward area measure of $F\circ u$ has density $m+1$ at $F(u(0))$.
\item In the free boundary case, suppose a smooth diffeomorphism $F$ between neighborhoods in $\R^K$ sends $S$ to a hyperplane and its normal lines to lines perpendicular to that hyperplane. The reflected double of $F\circ u$ has density $m+1$ at $F(u(0))$. Consequently, in a closed surface obtained by this reflection, each preimage outside the branch set contributes at least $1$ to the image density.
\end{enumerate}
\end{lemma}

\begin{proof}
At a point outside the branch set, the local positive lower bound for the induced metric and continuity of $\omega$ force $m=0$. Translate $u(0)$ to $0$ and consider
\[
 u_r(z)=r^{-(m+1)}u(rz),\qquad
 \lambda_r(z)=m\log|z|+\omega(rz),
\]
where $\lambda_r$ is the conformal factor of $u_r$ in logarithmic form. The metric expansion gives $|\nabla u(z)|\leq C|z|^m$ almost everywhere near $0$. Integrating this bound on almost every radial segment and using continuity gives $|u(z)|\leq C|z|^{m+1}$ throughout the chart after shrinking it suitably.

On any fixed annulus $a<|z|<b$ with $0<a<b<\infty$, or half-annulus in the boundary case, $\lambda_r$ converges uniformly to $m\log|z|+\omega(0)$ as $r \downarrow 0$. With integrals taken over the corresponding half-annuli in the boundary case, we also have
\[
\begin{aligned}
 \int_{a<|z|<b}|\nabla(\omega(rz))|^2\,\dd x\dd y
 &=\int_{ar<|z|<br}|\nabla\omega|^2\,\dd x\dd y\longrightarrow0,\\
 \int_{a<|z|<b}|A_{u_r}|^2\,\dd\mu_{u_r}
 &=\int_{ar<|z|<br}|A_u|^2\,\dd\mu_u\longrightarrow0.
\end{aligned}
\]
The second equality uses invariance of the integral of the squared norm of the second fundamental form under reparametrization and Euclidean dilation. For coordinate indices $i,j\in\{1,2\}$, the conformal Hessian decomposition is
\[
 \partial_{ij}u_r
 =\partial_i\lambda_r\,\partial_ju_r
  +\partial_j\lambda_r\,\partial_iu_r
  -\delta_{ij}\sum_{\ell=1}^2\partial_\ell\lambda_r\,\partial_\ell u_r
  +(A_{u_r})_{ij}.
\]
On the fixed annulus, uniform bounds on $\lambda_r$ give an $L^\infty$ bound for $\nabla u_r$ and make $\dd\mu_{u_r}$ comparable to $\dd x\dd y$. More precisely, $\sum_{i,j}|(A_{u_r})_{ij}|^2=e^{4\lambda_r}|A_{u_r}|^2$, so the second fundamental form terms tend to zero in $L^2(\dd x\dd y)$. Since $\nabla\lambda_r$ is bounded in $L^2$ there, the Hessian identity yields a uniform $W^{2,2}$ bound for $u_r$. From any sequence $r\downarrow0$ we can therefore extract a diagonal subsequence on which $u_r$ converges weakly in $W^{2,2}$, strongly in $W^{1,2}$, and uniformly on compact annuli or half-annuli to a map $v$. Strong convergence of the derivatives of $u_r$ and uniform convergence of $\lambda_r$ give
\[
 v^*g_{\R^K}=e^{2\omega(0)}|z|^{2m}|\dd z|^2.
\]
Tracing the Hessian decomposition cancels the tangential terms and gives $\Delta u_r=(A_{u_r})_{11}+(A_{u_r})_{22}\to0$ in $L^2$ on each such region. Thus $v$ is harmonic away from $0$ in the interior.

In the boundary case, the rescaled supporting hypersurfaces $r^{-(m+1)}S$ converge smoothly on bounded sets to $T_0S$. On compact intervals of the boundary line away from $0$, the $W^{2,2}$ bound gives bounds in $W^{1/2,2}$ for the first derivative traces of $v_r$. The compact embedding of $W^{1/2,2}$ into $L^2$ and the weak $W^{2,2}$ convergence identify their strong $L^2$ limits as the traces of $v$. Uniform convergence of the maps and smooth convergence of the tangent spaces on the rescaled supporting hypersurfaces now give $v(x,0)\in T_0S$ and $\partial_yv(x,0)\perp T_0S$. Reflect the components of $v$ in $T_0S$ evenly and its normal component oddly. The zero normal derivative of the tangential components and the zero trace of the normal component give weak harmonicity of $v$ across the boundary line away from $0$.  This reflection preserves the displayed metric because the reflections in the domain and the target are Euclidean isometries. In both cases, $|v(z)|\leq C|z|^{m+1}$, so the puncture of the harmonic map $v$ is removable and $v(0)=0$.

The holomorphic vector $z^{-m}(v_x-i v_y)$ has squared norm $2e^{2\omega(0)}$ on the punctured plane. Thus, $z^{-m}(v_x-i v_y)$ is constant, since the Laplacian of the squared norm of a holomorphic vector is $4$ times the squared norm of its complex derivative. Integrating and using $v(0)=0$ gives
\[
 v(z)=\frac{\operatorname{Re}(c z^{m+1})}{m+1},
 \qquad c\cdot c=0,\qquad
 |\operatorname{Re}c|=|\operatorname{Im}c|=e^{\omega(0)}.
\]
Here $c\cdot c=0$ follows from conformality of $v$. The real and imaginary parts of $c$ are orthogonal, so
\[
 |v(z)|=\frac{e^{\omega(0)}}{m+1}|z|^{m+1}.
\]

We use this growth estimate to pass from the rescaled maps to the density. Take $F=\operatorname{id}$ for the original map and write
\[
 w_r(z)=r^{-(m+1)}\bigl(F(u(rz))-F(0)\bigr).
\]
Fix a subsequence along which $u_r\to v$ as above. The first order expansion of $F$ at $0$ and the uniform boundedness of $u_r$ on compact sets give $w_r\to DF(0)u_r$ uniformly on each compact annulus or half-annulus. Since $DF(u(rz))\to DF(0)$ uniformly on these sets, $w_r\to DF(0)v$ both uniformly and strongly in $W^{1,2}$ there. The strong $L^2$ convergence of the first derivatives and the Cauchy--Schwarz inequality then give
$J_2w_r\to J_2(DF(0)v)$ in $L^1$, where
$J_2w_r=|\partial_xw_r\wedge\partial_yw_r|$.

The local bi-Lipschitz bound for $F$ and the lower bound for $|u|$ give
$|w_r(z)|\geq C^{-1}|z|^{m+1}$. Thus the preimages of $B_1(0)$ under $w_r$ and $DF(0)v$ lie in a fixed disk. Moreover, the bound $|\nabla w_r(z)|\leq C|z|^m$ gives
\[
 \int_{|z|<a}J_2w_r\,\dd x\dd y
 \leq C\int_{|z|<a}|z|^{2m}\,\dd x\dd y
 \leq Ca^{2(m+1)},
\]
with the integrals taken over the half-disk in a boundary chart. The same estimate holds for $DF(0)v$ by the formula for $v$.

In (iii), the rescaled map $u_r$ on the other half of the doubled chart is obtained by reflecting the domain across its boundary line and reflecting the target across the flattening hyperplane translated through $0$. These reflections preserve the estimates and convergence just proved. Since $DF(0)$ sends $T_0S$ and its normal line to the direction space of the flattening hyperplane and its normal line, respectively, it intertwines the two target reflections. Hence the limit of $w_r$ on the doubled chart is $DF(0)v$ on the full plane, with $v$ extended by the reflection used above.

After reparametrizing the chart by $z\mapsto rz$, the density ratio at image radius $r^{m+1}$ is therefore
\[
 \frac{1}{\pi}\int_{\{|w_r|<1\}}J_2w_r\,\dd x\dd y,
\]
where in (iii) the integral includes both halves of the doubled chart. Away from $0$, $DF(0)v$ is a local diffeomorphism onto its image plane. Its preimage of the unit circle in that plane thus has planar measure zero. On the part of the fixed disk with $|z|\geq a$, uniform convergence of the maps gives almost everywhere convergence of the indicators of the integration sets. Combining this with $L^1$ convergence of the Jacobians gives convergence of the integrals on that region. Letting $a\downarrow0$ after taking the subsequential limit and using the previous area bound on $|z| <a$, we obtain the limiting ratio
\[
 \frac{1}{\pi}\int_{\{|DF(0)v|<1\}}J_2(DF(0)v)\,\dd x\dd y
\]
over the full plane or half-plane, respectively.

On the full plane, the formula for $v$ gives exactly $m+1$ preimages of every nonzero point of its image plane. The invertible linear map $DF(0)$ preserves this count. Since $B_1(0)$ intersects the image plane of $DF(0)v$ in a disk of area $\pi$, the area formula gives density $m+1$ in the interior and for the doubled map in (iii). For the original map in a boundary chart, both $|v|$ and $J_2v$ depend only on $|z|$. The half-plane therefore contributes half of the full integral, giving density $(m+1)/2$.

Every sequence of scales tending to zero has a subsequence to which this calculation applies, and the resulting ratio is independent of the subsequence. This proves the stated density limits. At a preimage outside the branch set, $m=0$, so an interior or doubled chart contributes density $1$. This completes the proof of Lemma \ref{lem:app-local-multiplicity}
\end{proof}

\subsection{Annular minimization under a Steklov constraint}
\label{sp:section}

In this subsection, we prove Theorem~\ref{app:small-annuli}.
The relation between Steklov eigenvalue optimization and free boundary minimal surfaces is well established; see, for example, \cite{Fraser-Schoen-2016}.
Here we use a lower bound for the first normalized Steklov eigenvalue to prevent disk bubbling in a geometric minimization problem, since a nonconstant free boundary disk may have zero Willmore energy.

Let $\Sigma$ be a compact smooth annulus, let $S\subset\R^K$, $K\geq3$, be a smooth compact two-sided embedded hypersurface, and fix a smooth unit normal $n_S$ on $S$.
We consider the branched conformal immersions with free boundary in Theorem~\ref{thm:intro-main-compactness}, allowing their conformal structures to vary.
We require the outward image conormal to satisfy
\begin{equation}\label{sp:conormal}
 \eta_u=n_S(u)
 \qquad\text{a.e. on }\partial\Sigma\setminus\cB_u.
\end{equation}
For example, this sign holds when $S=\partial\Omega$, $n_S$ is its outward normal, and $u(\Sigma)\subset\overline\Omega$.
We use the area and Willmore energy from Subsection~\ref{sec:fbbci} and the boundary length from \eqref{eq:induced-boundary-geometry}, retaining the multiplicity of the parametrization.

Choose a smooth metric $h$ on $\Sigma$ compatible with the conformal structure of $u$. Let $\rho_u$ be the density of the measure $\dd s_{g_u}$ in \eqref{eq:induced-boundary-geometry} with respect to $\dd s_h$, and set
\[
 Q_h(f,g)=\int_\Sigma\langle\dd f,\dd g\rangle_h\,\dd\vol_h.
\]
Lemma~\ref{lem:general-boundary-removability} and the conformal factor expansion at a branch point give
\[
 \rho_u(x)=|x|^m e^{\omega(x)},
 \qquad m\in\mathbb N\cup\{0\},\quad \omega\in C^0,
\]
in a boundary conformal chart centered at that point.
Away from the finite branch set the conformal factor is locally bounded above and below.
Thus $\rho_u\in L^\infty(\partial \Sigma)$ and is positive almost everywhere.
We use the variational definition~\eqref{eq:app-annular-steklov} with the induced boundary arclength measure $\dd s_{g_u}$.

The weighted trace $W^{1,2}(\Sigma)\to L^2(\partial \Sigma,\dd s_{g_u})$ is compact, and the form
\[
 Q_h(f,f)+\int_{\partial \Sigma}f^2\dd s_{g_u}
\]
is coercive on $W^{1,2}(\Sigma)$.
Indeed, a sequence contradicting coercivity would, by Poincar\'e's inequality and compactness of the trace map, converge to a nonzero constant with zero weighted boundary norm.
Since $L(u)>0$, this is impossible.
The usual variational construction therefore gives a discrete Steklov spectrum, with a simple zero eigenvalue on each connected component.
We index the spectrum from zero, so that $\sigma_0=0<\sigma_1$ on a connected surface. On a finite disjoint union of bordered surfaces, we use the direct sum of these forms and repeat zero according to the number of components; thus $\sigma_1=0$ if at least two components are present.

We apply bounds for the normalized first Steklov eigenvalue over smooth metrics on $\Sigma$ in the conformal class $[h]$ to the possibly branched induced metric $g_u$.
To justify this, approximate $\rho_u$ in $L^p(\partial \Sigma)$, $p>1$, by smooth positive densities.
H\"older's inequality and the trace embedding into every finite $L^q$ give convergence of the corresponding boundary forms in the operator norm on $W^{1,2}(\Sigma)$.
Here is a direct verification for the first eigenvalue.
Set
\[
 C_\rho(f)=\inf_{c\in\R}\int_{\partial \Sigma}(f-c)^2\rho\,\dd s_h
 =\int_{\partial \Sigma}f^2\rho\,\dd s_h
 -\frac{\bigl(\int_{\partial \Sigma}f\rho\,\dd s_h\bigr)^2}
 {\int_{\partial \Sigma}\rho\,\dd s_h}.
\]
Subtract the mean $\bigl(\int_\Sigma\dd\vol_h\bigr)^{-1}\int_\Sigma f\,\dd\vol_h$ from $f$, which changes neither $Q_h$ nor $C_\rho$.
Poincar\'e's inequality then gives
\[
 |C_{\rho_j}(f)-C_\rho(f)|
 \leq\varepsilon_jQ_h(f,f),\qquad\varepsilon_j\longrightarrow0.
\]
Comparing $Q_h\geq\sigma_1(\rho)C_\rho$ with the same inequality for $\rho_j$ proves convergence of their first eigenvalues.
In particular,
\begin{equation}\label{sp:conformal-upper-bound}
 \sigma_1(u)L(u)\leq\sigma_1^*(\Sigma,[h]),
\end{equation}
where the right-hand side is the supremum of the normalized first Steklov eigenvalue over smooth metrics on $\Sigma$ in $[h]$.

Upper semicontinuity under weak convergence of boundary measures is an established ingredient of spectral optimization; see \cite{appKokarev,appGKL}.
We give the cutoff argument in the form needed to compare the parametrized bubble tree components.

\begin{lemma}\label{sp:spectral-transplantation}
Let $u_n:(\Sigma,h_n)\to\R^K$ converge to a noncollapsed bubble tree as in Theorem~\ref{thm:intro-main-compactness}. Let
\[
 B=\bigsqcup_{\substack{i\in\mathfrak I_*,\,\,\partial\Sigma_\infty^i\ne\varnothing}}\Sigma_\infty^i\ne\varnothing
\]
be the disjoint union of the domains of its nonconstant bordered components. Define the Steklov spectrum on $B$ by the direct sum of the Dirichlet forms and the induced boundary arclength measures of these component maps.
Then, with eigenvalues indexed from zero,
\begin{equation}\label{sp:upper-semicontinuity}
 \limsup_{n\to\infty}\sigma_k(u_n)\leq\sigma_k(B)
 \qquad(k\geq0).
\end{equation}
\end{lemma}

\begin{proof}
Fix $\varepsilon>0$.
By smooth density and the min--max characterization, choose a smooth $(k+1)$-dimensional test space on $B$ whose Gram matrix for the induced boundary $L^2$ inner product on $B$ is positive definite and whose maximal Rayleigh quotient is at most $\sigma_k(B)+\varepsilon$.
The test functions and their first derivatives are bounded. Around each branch, nodal, and attaching point, choose a conformal disk or half disk and multiply the functions by a cutoff which vanishes on radius $r^2$, equals one outside radius $r$, and interpolates logarithmically between them.
Its Dirichlet energy is at most $C/|\log r|$, so the cutoff error in the Dirichlet energy of each test function tends to zero. The boundary Gram entries also converge, because the induced boundary measures have bounded densities and hence no atoms.
For a sufficiently small fixed $r$, the truncated space still has dimension $k+1$ and maximal quotient at most $\sigma_k(B)+2\varepsilon$.

Transplant this space using the pairwise disjoint comparison embeddings of the bubble tree and extend the functions by zero. Since they vanish near the artificial boundaries of the selected cores, the extensions belong to $W^{1,2}(\Sigma)$.
Smooth convergence of the pulled-back conformal structures gives convergence of the Dirichlet inner products of every pair of transplanted basis functions.
On each compact boundary interval of the selected cores, the Sobolev boundary traces of the tangential coordinate derivatives of the pulled-back maps $u_n$ converge strongly in $L^q$ for every finite $q$. This follows from weak $W^{2,2}$ convergence and the compact embedding of $W^{1/2,2}$ into $L^q$ on a one-dimensional boundary interval.
Consequently every boundary Gram entry converges as well.
For all sufficiently large $n$, the transplanted space has dimension $k+1$ and maximal quotient at most $\sigma_k(B)+3\varepsilon$.
Min--max and then $\varepsilon\downarrow0$ prove \eqref{sp:upper-semicontinuity}.
\end{proof}

We recall also the weighted disk inequality
\begin{equation}\label{sp:disk-bound}
 \sigma_1({D},\rho)\int_{\partial{D}}\rho\,ds\leq2\pi
\end{equation}
for $\rho\in L^\infty(\partial D)$ with $\rho>0$ almost everywhere, where $ds$ is Euclidean arclength on the unit circle $\partial D$. The constant $2\pi$ is sharp: for the Euclidean unit disk and $\rho\equiv c>0$, the first positive eigenvalue is $1/c$ and the total boundary mass is $2\pi c$.
Indeed, a disk automorphism balances its nonatomic boundary measure.
The real and imaginary parts of this automorphism are then admissible test functions, each with Dirichlet energy $\pi$, and their squared boundary values sum to one.
Summing their variational inequalities proves \eqref{sp:disk-bound}.
For the balancing assertion, the boundary averages of disk automorphisms define a continuous map on the parameter disk whose extension to the parameter circle is the radial unit vector; Brouwer's theorem forces a zero.

We first record a consequence of the boundary monotonicity estimate used in \cite[proof of Theorem~4.1]{appNovagaPozzetta}.

\begin{lemma}\label{sp:coercivity}
Let $\Sigma$ be a compact connected surface with $b\geq1$ smooth boundary components, and let $u\in W^{2,2}(\Sigma,\R^K)$ be a branched conformal immersion with finite area and finite total curvature. Assume $W(u)\leq\Lambda<4\pi$.
Writing $\Gamma=u(\partial \Sigma)$ and $L=L(u)$, we have
\begin{align}
 \sup_{x\in u(\Sigma)}\operatorname{dist}(x,\Gamma)
 &\leq\frac{2L}{4\pi-\Lambda},\label{sp:boundary-distance}\\
 \operatorname{diam}u(\Sigma)&\leq D_{\Lambda,b}L,
 \qquad D_{\Lambda,b}=\frac12+\frac{4b}{4\pi-\Lambda},\label{sp:diameter}\\
 \cA(u)&\leq C_{\Lambda,b}L^2.\label{sp:area-bound}
\end{align}
\end{lemma}

\begin{proof}
For an interior image point $p\notin\Gamma$ with integer area density $m\geq1$, first variation with $(x-p)/|x-p|^2$ and completion of the square give
\[
 W(u)+2\int_{\partial \Sigma}
 \frac{(u-p)\cdot\eta_u}{|u-p|^2}\,\dd s_{g_u}
 =4\pi m+\int_\Sigma
 \left|\frac{\mathbf H_u}{2}
 +\frac{2(u-p)^\perp}{|u-p|^2}\right|^2\dd\mu_u
 \geq4\pi.
\]
The identity follows by removing inner and outer cutoffs, using finite total curvature and the branch expansion. The boundary integral is at most $L/\operatorname{dist}(p,\Gamma)$; approximation by these density points and continuity then give \eqref{sp:boundary-distance} on the full image.

Put $R=2L/(4\pi-\Lambda)$ and write $\Gamma_i=u(\partial_i\Sigma)$, with parametrized length $L_i$.
The shorter of the two boundary arcs joining two points gives $\operatorname{diam}\Gamma_i\leq L_i/2$.
The closed sets $\{x\in u(\Sigma):\operatorname{dist}(x,\Gamma_i)\le R\}$ cover the connected image, so their intersection graph is connected. Each edge joins two boundary-image components at distance at most $2R$, and a simple path uses at most $b-1$ edges. Following this path and using the distances of the two endpoints to $\Gamma$ and boundary arcs of total length at most $L/2$, we obtain
\[
 \operatorname{diam}u(\Sigma)\leq\frac L2+2bR,
\]
which is \eqref{sp:diameter}.
Finally choose $q\in\Gamma$ and use the position field $x-q$ in first variation:
\[
 2\cA(u)\leq
 D_{\Lambda,b}L^2+
 2D_{\Lambda,b}L\sqrt{\Lambda\cA(u)}.
\]
Solving this quadratic inequality in $\sqrt{\cA(u)}$ proves \eqref{sp:area-bound}.
\end{proof}

Fix a closed set $Q\subset\R^K$, $\ell>0$, and $s>2\pi$.
Let $\mathcal E(S,Q;\ell,s)$ be the class of branched conformal immersions of the annulus $\Sigma$ with free boundary on $S$ in the sense of Definition~\ref{def:fb-branched-conformal} satisfying \eqref{sp:conormal} and
\begin{equation*}
 u(\Sigma)\subset Q,\qquad
 L(u)=\ell,\qquad
 \sigma_1(u)L(u)\geq s.
\end{equation*}
No separate area bound is needed in this definition, and $Q=\R^K$ is allowed: under a fixed Willmore energy bound below $4\pi$, Lemma~\ref{sp:coercivity} bounds the area and places every image in a fixed compact neighborhood of the compact supporting hypersurface $S$.

\begin{theorem}\label{sp:annular-compactness}
Let $S\subset\R^K$ be a smooth compact embedded two-sided supporting hypersurface with chosen unit normal $n_S$, and let $Q\subset\R^K$ be closed. Fix $\ell>0$ and $0\leq\Lambda<4\pi$. Let $\Sigma$ be a compact smooth annulus, and let $h_n$ be smooth metrics on $\Sigma$. Suppose $u_n:(\Sigma,h_n)\to\R^K$ are branched conformal immersions with free boundary on $S$ in the sense of Definition~\ref{def:fb-branched-conformal}, satisfying \eqref{sp:conormal} and
\[
 u_n(\Sigma)\subset Q,\qquad L(u_n)=\ell,\qquad W(u_n)\leq\Lambda.
\]
Assume that
\begin{equation*}
 \liminf_{n\to\infty}\sigma_1(u_n)L(u_n)>2\pi.
\end{equation*}
The conformal structures $[h_n]$ remain in a compact subset of the moduli space of annuli. After passage to a subsequence and reparametrization, $u_n$ converges in the sense of Definition~\ref{def:componentwise-free-boundary-bubble tree-convergence} to a bubble tree whose only nonconstant component is a principal annular map $u_\infty$.
The limiting annular immersion satisfies \eqref{sp:conormal},
\[
 u_\infty(\Sigma)\subset Q,\qquad L(u_\infty)=\ell,\qquad
 \sigma_1(u_\infty)L(u_\infty)>2\pi,\qquad
 W(u_\infty)\leq\liminf_{n\to\infty}W(u_n).
\]
Moreover, $\cA(u_n)\to\cA(u_\infty)$, the images converge in Hausdorff distance to $u_\infty(\Sigma)$, and the pushforward boundary arclength measures converge as in \eqref{eq:contact-scalar}.
If, in addition, $u_n\in\mathcal E(S,Q;\ell,s)$ for a fixed $s>2\pi$, then $u_\infty$ belongs to the same class.
\end{theorem}

\begin{proof}
Choose $s_0$ such that $2\pi<s_0<\liminf_n\sigma_1(u_n)L(u_n)$.
After discarding finitely many terms, the sequence belongs to $\mathcal E(S,Q;\ell,s_0)$.
Pass first to a subsequence along which the Willmore energies converge to their original lower limit, and relabel this and all further subsequences.

Lemma~\ref{sp:coercivity} supplies the uniform area bound and a fixed compact set containing the images.
For a degenerating sequence of conformal classes on an annulus, Medvedev's degeneration formula gives $\sigma_1^*(\Sigma,[h_n])\longrightarrow2\pi$; see \cite[Theorem~1.4 and Remark~1.4]{appMedvedev}.
Together with \eqref{sp:conformal-upper-bound} and $\sigma_1(u_n)L(u_n)\geq s_0>2\pi$, this excludes conformal degeneration.
After reparametrizing $\Sigma$ by diffeomorphisms and passing to a subsequence, the conformal structures $[h_n]$ therefore converge smoothly to a conformal structure on the annulus $\Sigma$.
Theorem~\ref{thm:intro-main-compactness} gives a subsequence of $u_n$ whose images either collapse to a point or converge to a finite bubble tree. Point collapse is incompatible with the fixed positive boundary length $L(u_n)=\ell$: if the images lie in balls of radii $r_n\to0$ about a point $q$, first variation with the field $x-q$ gives
\[
 2\cA(u_n)
 \leq2r_n\sqrt{\Lambda\cA(u_n)}+r_n\ell.
\]
Thus $\cA(u_n)\to0$.
Proposition~\ref{prop:free-boundary-length-estimate} then yields
\[
 \ell\leq C_S\bigl(\cA(u_n)
 +\sqrt{\cA(u_n)W(u_n)}\bigr)\longrightarrow0,
\]
contrary to $\ell>0$.

Lemma~\ref{lem:app-contact-measures} gives
\begin{equation}\label{sp:length-identity}
 \ell=\sum_i L(u_\infty^i).
\end{equation}
The identity $\eta_{u_\infty^i}=n_S\circ u_\infty^i$ holds on each limiting boundary portion away from branch, nodal, and attaching points, by strong $L^2$ convergence of the Sobolev boundary traces of the tangential and transverse coordinate derivatives of the pulled-back maps.
The positive total boundary length ensures that at least one nonconstant bordered component occurs. If two such components were present, their disjoint union would have $\sigma_1=0$, and Lemma~\ref{sp:spectral-transplantation} would contradict
\[
 \sigma_1(u_n)\geq s_0/\ell>0.
\]
Every nonconstant closed component has Willmore energy at least $4\pi$.
Indeed, Lemma~\ref{lem:app-local-multiplicity} gives image density at least $1$ at a preimage outside the branch set, and the finite-curvature cutoff argument removes the branch points from first variation.
The geometric monotonicity argument for the weak Li--Yau inequality in \cite[Section~2.2.2, equation~(2.12)]{appMS} is independent of ambient codimension and gives $W\geq4\pi$.
Lower semicontinuity and $\Lambda<4\pi$ therefore exclude closed components.

There is consequently exactly one nonconstant component, and it is bordered.
If it were a disk bubble, \eqref{sp:length-identity}, Lemma~\ref{sp:spectral-transplantation}, and \eqref{sp:disk-bound} would give
\[
 s_0\leq\limsup_n\sigma_1(u_n)L(u_n)\leq2\pi,
\]
a contradiction.
The component must therefore be the principal annulus, and \eqref{sp:length-identity} gives $L(u_\infty)=\ell$. For the extracted subsequence, Lemma~\ref{sp:spectral-transplantation} and this fixed boundary length give
\[
 \sigma_1(u_\infty)L(u_\infty)\geq\limsup_n\sigma_1(u_n)L(u_n)\geq s_0>2\pi.
\]
If the sequence lies in $\mathcal E(S,Q;\ell,s)$ for a fixed $s>2\pi$, the same inequality retains $\sigma_1(u_\infty)L(u_\infty)\geq s$, so the limit belongs to that class.
The image constraint, Willmore lower semicontinuity, and the remaining convergence assertions follow from Theorem~\ref{thm:intro-main-compactness} and Lemma~\ref{lem:app-contact-measures}.
\end{proof}

\begin{coro}\label{sp:willmore-minimizer}
If
\[
 I:=\inf_{u\in\mathcal E(S,Q;\ell,s)}W(u)<4\pi,
\]
then there exists $u_\infty\in\mathcal E(S,Q;\ell,s)$ such that $W(u_\infty)=I$. In particular, $u_\infty\in W^{2,2}(\Sigma,\R^K)$ is a branched conformal immersion of the annulus $\Sigma$ with free boundary on the supporting hypersurface $S$.
\end{coro}

\begin{proof}
Choose $\Lambda$ with $I<\Lambda<4\pi$ and a minimizing sequence whose Willmore energies are at most $\Lambda$.
Theorem~\ref{sp:annular-compactness} gives $u_\infty\in\mathcal E(S,Q;\ell,s)$ with $W(u_\infty)\le I$.
The definition of $I$ gives the reverse inequality.
\end{proof}

We next establish the lower bound and embedding conclusion without a spectral assumption.

\begin{prop}\label{prop:app-short-annular-geometry}
Let $\Omega\subset\R^3$ be a smooth bounded domain and $\Lambda<4\pi$.
Write $S=\partial\Omega$ and let $n_S$ be its outward unit normal. There are constants $\ell_*=\ell_*(\Omega,\Lambda)>0$ and $C=C(\Omega,\Lambda)<\infty$ such that every branched conformal immersion $u:(\Sigma,h)\to\R^3$ of a compact annulus with free boundary on $S$ in the sense of Definition~\ref{def:fb-branched-conformal}, satisfying
\[
 u(\Sigma)\subset\overline\Omega,\qquad
 \eta_u=n_S(u),\qquad
 L(u)=\ell<\ell_*,\qquad W(u)\leq\Lambda
\]
is an unbranched $W^{2,2}$ immersion and a continuous embedding, satisfies $u^{-1}(\partial\Omega)=\partial\Sigma$, and obeys
\[
 W(u)\geq\pi^2-C\ell,\qquad
 \operatorname{diam}u(\Sigma)\leq C\ell,\qquad
 \cA(u)\leq C\ell^2.
\]
\end{prop}

\begin{proof}
Lemma~\ref{sp:coercivity}, with $b=2$, gives the diameter and area bounds. Choose $q\in u(\partial \Sigma)$; the diameter bound allows us to flatten $\partial\Omega$ in a single inward normal chart containing the entire image. Compactness makes this choice uniform in $q$ on $\partial\Omega$.
For its inverse $F$, with $F(q)=0$, we have
\[
 DF=DF(q)+O(\ell),\qquad |D^2F|\leq C,
\]
where $DF(q)$ is orthogonal.
Within this coordinate neighborhood, $F$ sends $\overline\Omega$ into $\{x_3\geq0\}$ and $S$ into $\{x_3=0\}$, and sends lines normal to $S$ to vertical lines.
Set $f=F\circ u$.
Since $\chi(\Sigma)=0$, Lemma~\ref{lem:free-boundary-gauss-bonnet} and the Gauss equation give
\[
 \int_\Sigma|A_u|^2\dd\mu_u
 =4W(u)-4\pi\mathfrak b(u)
   +2\int_{\partial \Sigma}k_{g_u}\dd s_{g_u}
 \leq4\Lambda+C\ell.
\]
The induced metrics of $f$ and $u$ differ by a relative $O(\ell)$.
In a frame orthonormal with respect to the induced metric $g_u$, normal projection of the Sobolev chain rule gives
\[
 A_f=P_f^\perp\bigl(DF\,A_u+D^2F(du,du)\bigr).
\]
Here the tangential Christoffel terms vanish after projection, and the traces defining $\mathbf H_f$ and $\mathbf H_u$ use their respective induced metrics.
Since $DF=DF(q)+O(\ell)$, this identity implies
\[
 |\mathbf H_f-DF(q)\mathbf H_u|\leq C\ell|A_u|+C,
 \qquad |A_f|\leq C(|A_u|+1).
\]
The area ratio is $1+O(\ell)$, so Cauchy--Schwarz, together with $\int_\Sigma|A_u|^2\dd\mu_u\leq4\Lambda+C\ell$ and $\cA(u)=O(\ell^2)$, gives
\begin{equation}\label{sp:flattened-bending}
 |W(f)-W(u)|\leq C\ell.
\end{equation}
Reflect $f$ across $\{x_3=0\}$ on a second copy of $\Sigma$ and denote the map of the doubled torus by $\widehat f$.
In a boundary chart with tangential coordinate $s$ and transverse coordinate $t$, the Sobolev boundary trace of $\partial_s f$ is horizontal, while that of $\partial_t f$ is vertical. The reflected map therefore has matching Sobolev traces of its value and of both coordinate first derivatives across the common boundary of the two copies. Hence $\widehat f\in W^{2,2}$ and
\begin{equation}\label{sp:torus-double-energy}
 W(\widehat f)=2W(f)\leq2\Lambda+C\ell<8\pi
\end{equation}
after decreasing $\ell_*$.
The double has finite total curvature, and first variation extends across a branch of order $m$ by cutoffs, since the conormal flux through a coordinate circle of radius $r$ is $O(r^{m+1})$. Along the seam, the matching conormal traces cancel. Thus $\widehat f$ defines a boundaryless integral image varifold with square-integrable generalized mean curvature. Averaging mean curvature over image fibers can only decrease its squared integral.

By Lemma~\ref{lem:app-local-multiplicity}, an interior branch point has local density at least $2$, and a boundary branch point has density at least $2$ in the complex double $\widehat \Sigma$.
Two distinct preimages outside the branch set of the same image point also give density at least $2$.
The weak Li--Yau inequality \cite[Section~2.2.2, equation~(2.12)]{appMS} and \eqref{sp:torus-double-energy} exclude all these possibilities.
Hence the double is injective and unbranched.
An interior point of $u$ mapping to $\partial\Omega$ would give two distinct interior preimages on the reflection plane in the double, so
$u^{-1}(\partial\Omega)=\partial \Sigma$.
Continuity and injectivity on the compact annulus $\Sigma$ make $u$ a homeomorphism onto its image.

After branch exclusion, the conformal factors of $u$ are bounded above and below on a finite cover. The smooth flattening preserves these bounds, so the doubled map is Lipschitz and its induced metric is bounded above and below by positive multiples of a smooth metric on the doubled torus. We may therefore apply Lemma~\ref{lem:app-weak-approximation} to obtain smooth torus immersions $\widehat f_j$ with
\[
 \widehat f_j\longrightarrow\widehat f\quad\text{strongly in }W^{2,2},
 \qquad W(\widehat f_j)\longrightarrow W(\widehat f).
\]
For large $j$, $W(\widehat f_j)<8\pi$, and the Li--Yau inequality \cite[Theorem~6]{Li-Yau1982} makes them embedded.
The Marques--Neves inequality \cite[Theorem~A]{Marques-Neves2014} yields
\[
 2\pi^2\leq\lim_jW(\widehat f_j)
 =2W(f)\leq2W(u)+C\ell.
\]
This proves the lower bound.
\end{proof}

To construct admissible annuli with Willmore energy close to $\pi^2$, we start from the upper half Clifford torus. In the half-space $\{z\geq0\}\subset\R^3$, consider
\[
 X(\varphi,\theta)=
 ((\sqrt2+\cos\varphi)\cos\theta,
  (\sqrt2+\cos\varphi)\sin\theta,\sin\varphi),
 \quad0\leq\varphi\leq\pi.
\]
This upper half Clifford torus meets $\{z=0\}$ orthogonally, with outward conormal $-e_3$ on both boundary circles.
Direct integration gives
\[
 W(X)=\frac\pi2\int_0^\pi
 \frac{(\sqrt2+2\cos\varphi)^2}{\sqrt2+\cos\varphi}\dd\varphi
 =\pi^2.
\]
Its conformal cylinder has length $\pi$, and the two boundary densities are $a_\pm=\sqrt2\pm1$.
For functions constant in the angular variable, separation of variables gives the positive eigenvalue $2\sqrt2/\pi$. For each angular Fourier mode $n\geq1$, it gives
\[
 \sigma_{n,\pm}=n\left(
 \sqrt2\coth(n\pi)\pm\sqrt{2\coth^2(n\pi)-1}
 \right).
\]
To see this, the Dirichlet-to-Neumann matrix on the $n$th Fourier mode is
\[
 n\begin{pmatrix}
 \coth(n\pi)&-\operatorname{csch}(n\pi)\\
 -\operatorname{csch}(n\pi)&\coth(n\pi)
 \end{pmatrix},
\]
and its boundary weight matrix is $\operatorname{diag}(a_+,a_-)$.
The smaller root increases with $n$, since its rationalized expression is
\[
 \frac{n}{\sqrt2\coth(n\pi)+\sqrt{2\coth^2(n\pi)-1}}.
\]
The value $2\sqrt2/\pi$ exceeds $\sigma_{1,-}$.
Consequently, $\sigma_1(X)=\sigma_{1,-}$, and we obtain
\begin{equation}\label{sp:comparison-gap}
 \sigma_1(X)L(X)=s_{\mathrm{Cl}}=4\pi\sqrt2\sigma_{1,-}>2\pi.
\end{equation}
Indeed, $\coth\pi<9/8$ gives
\[
 \sigma_{1,-}^{-1}=\sqrt2\coth\pi+\sqrt{2\coth^2\pi-1}
 <\frac{9\sqrt2}{8}+\frac{7\sqrt2}{8}=2\sqrt2.
\]
Fix $s\in(2\pi,s_{\mathrm{Cl}})$.
Let $\Omega\subset\R^3$ be a smooth bounded domain, choose $p\in S=\partial\Omega$, and let $n_S$ be its outward unit normal.
Choose a boundary chart $\psi:U\subset\R^2\to S$ with $\psi(0)=p$ and $D\psi(0)$ an isometry onto $T_pS$.
Inward normal coordinates are given by
\[
 \Phi(y,t)=\psi(y)-t\,n_S(\psi(y)),\qquad t\geq0.
\]
For $|y|$ and $t$ sufficiently small, $\Phi$ is a diffeomorphism onto a one-sided tubular neighborhood, and $\Phi(y,t)\in\Omega$ when $t>0$.
The differential $D\Phi(0,0)$ is a Euclidean isometry.
For all sufficiently small $r>0$, set
\begin{equation*}
 u_r=\Phi(rX_1,rX_2,rX_3).
\end{equation*}
The map $u_r$ is a smooth embedding of an annulus into $\overline\Omega$.
On both boundary circles, the outward conormal of $X$ is $-e_t$, and
\[
 D\Phi(y,0)(-e_t)=n_S(\psi(y)).
\]
Since the boundary tangent is mapped into $T_{\psi(y)}S$, the immersion $u_r$ meets $S$ orthogonally with the prescribed outward conormal sign.
The induced metric $u_r^*g_{\R^3}$ determines the conformal structure on the annulus $[0,\pi]\times(\R/2\pi\mathbb Z)$.

To estimate the area, Willmore energy, boundary length, and Steklov spectrum of $u_r$, translate and rescale in $\R^3$ and put
\[
 v_r=\frac{u_r-p}{r}.
\]
Taylor expansion of the fixed smooth tubular map gives a smooth extension to $r=0$ and
\[
 v_r=D\Phi(0,0)X+O(r)\quad\text{in }C^2,
 \qquad v_0=D\Phi(0,0)X.
\]
The induced metrics remain uniformly nondegenerate.
The formulas for mean curvature, area, and length consequently give
\begin{align}
 W(u_r)&=W(v_r)=\pi^2+O(r), \quad\cA(u_r)=r^2\cA(X)+O(r^3),\notag\\
 L(u_r)&=rL(X)+O(r^2),\qquad
 \frac{d}{dr}L(u_r)=L(X)+O(r)>0,
 \quad L(X)=4\pi\sqrt2.\label{sp:comparison-length}
\end{align}
The first equality uses scale invariance of $W$, while the $C^2$ estimate controls the Willmore integrand because mean curvature depends smoothly on the two-jets of an immersion with nondegenerate metric.
The derivative in \eqref{sp:comparison-length} follows from the smooth extension of $v_r$ and $L(u_r)=rL(v_r)$.

The spectral limit also follows directly from the metric comparison.
Using the induced metrics $g_{v_r}$ and $g_X$ from Subsection~\ref{sec:fbbci}, we compare the quadratic forms and boundary arclength measures.
For a fixed constant $C$ and all sufficiently small $r$,
\[
 (1-Cr)Q_{g_X}\leq Q_{g_{v_r}}\leq(1+Cr)Q_{g_X},
 \qquad
 (1-Cr)\,\dd s_{g_X}\leq \dd s_{g_{v_r}}\leq(1+Cr)\,\dd s_{g_X}.
\]
Minimizing the boundary integrals of $(f-c)^2$ over constants $c$ gives the same comparison for their variance forms.
The Rayleigh quotient therefore gives
\[
 \sigma_1(v_r)=\sigma_1(X)+O(r).
\]
Normalized Steklov eigenvalues are scale invariant, so
\begin{equation}\label{sp:tubular-spectrum}
 \sigma_1(u_r)L(u_r)=\sigma_1(v_r)L(v_r)
 =s_{\mathrm{Cl}}+O(r).
\end{equation}
Since $s<s_{\mathrm{Cl}}$ is fixed, $\sigma_1(u_r)L(u_r)>s$ for all sufficiently small $r>0$.
Finally, \eqref{sp:comparison-length} can be inverted near zero.
Every sufficiently small prescribed length $\ell>0$ has a comparison annulus with
\[
 r(\ell)=\frac{\ell}{L(X)}+O(\ell^2),\qquad
 L(u_{r(\ell)})=\ell,\qquad
 W(u_{r(\ell)})\leq\pi^2+C\ell<4\pi.
\]
\begin{proof}[Proof of Theorem~\ref{app:small-annuli}]
Fix $s\in(2\pi,s_{\mathrm{Cl}})$.
By \eqref{sp:tubular-spectrum} and \eqref{sp:comparison-length}, the comparison annuli belong to $\mathcal E_\Omega(\ell;s)$ for every sufficiently small $\ell>0$ and give
$I_\Omega(\ell;s)\leq\pi^2+C\ell$.
Fix $\Lambda_*\in(\pi^2,4\pi)$ and choose $\ell_0>0$, depending on $\Omega$ and $s$, so that this upper bound is below $\Lambda_*$ whenever $0<\ell<\ell_0$.
The class in \eqref{eq:app-annular-class} is precisely $\mathcal E_\Omega(\ell;s)=\mathcal E(\partial\Omega,\overline\Omega;\ell,s)$. Corollary~\ref{sp:willmore-minimizer} gives attainment.
After discarding finitely many terms, every minimizing sequence has Willmore energy at most $\Lambda_*$, so Theorem~\ref{sp:annular-compactness} gives all the asserted subsequential convergence statements.
Every minimizer also has $W(u)\leq\Lambda_*$.
After decreasing $\ell_0$ once more, Proposition~\ref{prop:app-short-annular-geometry} shows that every minimizer is a continuous embedding satisfying $u^{-1}(\partial\Omega)=\partial\Sigma$ and the stated diameter and area bounds. It also gives $I_\Omega(\ell;s)\geq\pi^2-C\ell$.
Combining the upper and lower bounds for $I_\Omega(\ell;s)$ proves \eqref{sp:universal-profile}.
\end{proof}

The doubling construction in Proposition~\ref{prop:app-short-annular-geometry} also gives a quantitative estimate for the resulting closed tori after ambient M\"obius normalization.
Let
\[
 F_{\mathrm{Cl}}(\theta,\varphi)
 =\frac1{\sqrt2}(\cos\theta,\sin\theta,\cos\varphi,\sin\varphi)
 \quad\text{on }\mathbb S^1\times\mathbb S^1.
\]
\begin{coro}\label{coro:app-normalized-annular-double}
Fix $s\in(2\pi,s_{\mathrm{Cl}})$. For every sufficiently small $\ell>0$, let $u_\ell$ be any minimizer of $I_\Omega(\ell;s)$ in Theorem~\ref{app:small-annuli}, and let $\widehat f_\ell$ be its flattened reflected double constructed in Proposition~\ref{prop:app-short-annular-geometry}.
After inverse stereographic projection into $\mathbb S^3$ and a suitable ambient M\"obius transformation, the image of $\widehat f_\ell$ admits a parametrization $F_\ell\in W^{2,2}(\mathbb S^1\times\mathbb S^1,\mathbb S^3)$ that is a homeomorphism onto its image and satisfies
\[
 F_\ell^*g_{\mathbb S^3}
 =e^{2\lambda_\ell}
   (a_\ell\dd\theta^2+2b_\ell\dd\theta\dd\varphi
                            +c_\ell\dd\varphi^2)
\]
and
\[
 \|F_\ell-F_{\mathrm{Cl}}\|_{W^{2,2}}
 +\|\lambda_\ell\|_{L^\infty}
 +\left|a_\ell-\frac12\right|+|b_\ell|
 +\left|c_\ell-\frac12\right|
 \leq C\sqrt\ell.
\]
Here $a_\ell,b_\ell,c_\ell$ are real constants defining a positive definite flat metric, the norms are taken on $\mathbb S^1\times\mathbb S^1$ with its standard flat metric, and $C$ depends only on $\Omega$ and $s$.
\end{coro}

\begin{proof}
The proof of Proposition~\ref{prop:app-short-annular-geometry} shows that $\widehat f_\ell$ is an unbranched $W^{2,2}$ immersion of a torus and a homeomorphism onto its image. It therefore defines an integral varifold of multiplicity $1$ with square-integrable mean curvature.
Its image under inverse stereographic projection has the same properties in $\mathbb S^3$ and is a topological torus.
The conformal invariance of closed Willmore energy, justified in the weak class by Lemma~\ref{lem:app-weak-approximation}, and \eqref{sp:flattened-bending} give
\[
 2\pi^2\leq \mathcal W_{\mathbb S^3}(\widehat f_\ell)
 =W_{\R^3}(\widehat f_\ell)
 \leq2\pi^2+C\ell.
\]
Here $\mathcal W_{\mathbb S^3}=\int(1+|\mathbf H_{\mathbb S^3}|^2/4)\dd\mu$ denotes the conformally invariant spherical functional.
For $\ell_0$ sufficiently small, depending only on $\Omega$ and $s$, the preceding Willmore bound satisfies the smallness hypothesis of \cite[Theorem~1.2]{BiZhou2025} with $\delta^2=C\ell$. Applying that theorem to the projected image of $\widehat f_\ell$ gives the parametrization $F_\ell$ and all the stated estimates.
\end{proof}

The spectral constraint is used to retain the annulus in the compactness argument.
Proposition~\ref{prop:app-short-annular-geometry} itself applies to every sufficiently short annular competitor with $W\leq\Lambda<4\pi$.

\subsection{Small-area disks and concentration on the supporting hypersurface}
\label{ga:section}

In this subsection, we prove Theorem~\ref{app:small-disks} and extend its shape and location estimates to disks whose Willmore energies approach the infimum. Alessandroni--Kuwert~\cite{appAK} constructed smooth immersions of disks with free boundary on $S$ by locally perturbing hemispheres. Their construction first prescribes both the area and a barycenter on $S$, then obtains critical points under the area constraint alone by varying the barycenter. The first-order Willmore energy expansion of this constrained family gives the upper bound used below. We use their upper bound for the Willmore energy to prove attainment in the full class $\mathcal C(S)$, and quantitative rigidity for nearly umbilical spheres then gives estimates for the shape and location of every minimizer. Related evolution and stability results start from smooth local configurations; see Metsch~\cite{appMetsch} and Dall'Acqua--Rupp--Sch\"atzle--Schlierf~\cite{appDRSS}.

Let $\Omega\subset\R^3$ be a smooth bounded domain with supporting hypersurface $S=\partial\Omega$ and outward unit normal $n_S$. We use the class $\mathcal C(S)$ and the quantities $\beta_S(a)$, $\lambda(a)$, and $\mathcal D(u)$ from \eqref{ga:beta-definition}, with the conormal condition \eqref{ga:signed-contact}. For an unbranched $u\in\mathcal C(S)$, Lemma~\ref{lem:free-boundary-gauss-bonnet} and the Gauss equation give
\begin{equation}\label{ga:disk-defect}
 W(u)=2\pi+\frac12\mathcal D(u)-\int_{\partial{D}}h^S(T,T)\dd s_{g_u}.
\end{equation}
Here $\eta_u=n_S\circ u$, so the convention using the inward unit conormal in the Gauss--Bonnet formula gives $k_{g_u}=h^S(T,T)$.

Set $\kappa=\|h^S\|_{L^\infty(S)}$.
For any branched disk, the Gauss equation and Lemma~\ref{lem:free-boundary-gauss-bonnet} imply
\begin{equation}\label{ga:disk-cost}
 2\pi\bigl(1+\mathfrak b(u)\bigr)
 \le W(u)+\kappa L(u).
\end{equation}
Indeed $K_{g_u}\le|\mathbf H_u|^2/4$ and $|k_{g_u}|\le\kappa$.
A nonconstant closed spherical component has Willmore energy at least $4\pi$.

We first estimate the distance of the interior image from the boundary image.

\begin{lemma}\label{ga:boundary-distance}
Let $u\in W^{2,2}(\Sigma,\R^K)$ be a branched conformal immersion with finite total curvature, where $\Sigma$ is a compact surface with nonempty smooth boundary, and put $\Gamma=u(\partial\Sigma)$.
If $W(u)<4\pi$, then
\begin{equation}\label{ga:distance-estimate}
 \sup_{x\in u(\Sigma)}\dist(x,\Gamma)
 \le\frac{\sqrt{\cA(u)}}{\sqrt\pi-\sqrt{W(u)}/2}.
\end{equation}
For a disk,
\begin{equation}\label{ga:diameter-estimate}
 \diam u(\overline{D})
 \le\frac{L(u)}2+
 \frac{2\sqrt{\cA(u)}}{\sqrt\pi-\sqrt{W(u)}/2}.
\end{equation}
\end{lemma}

\begin{proof}
Let $\mu=u_\#\mu_u$ be the pushforward area measure, and consider the corresponding integral varifold with generalized mean curvature vector $\mathbf H$.
Away from $\Gamma$, its first variation has no singular boundary term, and its generalized mean curvature is the fiberwise average of the parametrized mean curvature.
This averaging gives $\int|\mathbf H|^2\,\dd\mu\le\int_\Sigma|\mathbf H_u|^2\,\dd\mu_u=4W(u)$.
The parametrized first variation extends across the finite branch set by the finite-curvature cutoff argument.
Simon's interior monotonicity formula therefore applies to balls disjoint from $\Gamma$; the corresponding argument for branched conformal immersions is recalled in~\cite[Section~2.2.2]{appMS}.

Choose a point $x$ with area density at least $1$ and $0<r<\dist(x,\Gamma)$.
With the convention that $\mathbf H$ is the trace of the second fundamental form, this formula gives
\[
 1\le\frac{\mu(B_r(x))}{\pi r^2}
 +\frac1{16\pi}\int_{B_r(x)}|\mathbf H|^2\dd\mu
 +\frac1{2\pi r^2}\int_{B_r(x)}\mathbf H\cdot(y-x)\dd\mu.
\]
Writing $A=\cA(u)$ and $W=W(u)$, Cauchy--Schwarz yields
\[
 \pi\le\frac A{r^2}+\frac W4+\frac{\sqrt{AW}}r
 =\left(\frac{\sqrt A}r+\frac{\sqrt W}2\right)^2.
\]
Let $r\uparrow\dist(x,\Gamma)$.
Density points are dense in the support of the image varifold, and this support is the continuous image of the immersion.
To see the latter assertion, use nondegeneracy in charts away from the branch points and continuity at the finitely many branch points.
Thus \eqref{ga:distance-estimate} holds on the whole image.
For a disk, $\Gamma$ is the image of one rectifiable closed curve, so $\diam\Gamma\le L(u)/2$.
The triangle inequality proves \eqref{ga:diameter-estimate}.
\end{proof}

After conformal reparametrization, the disks constructed in \cite[Proposition~1]{appAK} belong to $\mathcal C(S)$ and satisfy
\begin{equation}\label{ga:AK-comparison}
 \cA(\phi_{p,\lambda})=2\pi\lambda^2,
 \qquad
 \left|W(\phi_{p,\lambda})-2\pi+\pi H_S(p)\lambda\right|
 \le C_S\lambda^2,
\end{equation}
uniformly for $p\in S$ and small $\lambda>0$.

The maximum $H_*$ is positive.
For example, fix $c\in\Omega$ and let $p\in S$ maximize $|p-c|$.
At this point the normal points along $p-c$ and the hypersurface $S$ lies inside the tangent enclosing ball of radius $R=|p-c|$.
The second derivative comparison gives $h^S\ge R^{-1}g_S$ at $p$, hence $H_S(p)\ge2/R$.
Consequently \eqref{ga:AK-comparison} implies
\begin{equation}\label{ga:strict-upper}
 \beta_S(2\pi\lambda^2)
 \le2\pi-\pi H_*\lambda+C_S\lambda^2<2\pi
\end{equation}
for all sufficiently small $\lambda$.

\begin{prop}\label{ga:weak-attainment}
There exists $a_0=a_0(S)>0$ such that, for every $0<a<a_0$, there is an unbranched conformal immersion $u\in W^{2,2}(D,\R^3)\cap\mathcal C(S)$ satisfying $\cA(u)=a$ and $W(u)=\beta_S(a)$.
\end{prop}

\begin{proof}
Fix $a>0$ small and choose a minimizing sequence $u_j$ with $\cA(u_j)=a$.
By \eqref{ga:strict-upper}, after discarding finitely many terms,
\[
 W(u_j)\le w_0<2\pi
\]
for a constant $w_0$ depending on this fixed $a$.
Put $\Gamma_j=u_j(\partial D)\subset S$. The length bound \eqref{eq:sec2-boundary-length-area-willmore} is uniform, and Lemma~\ref{ga:boundary-distance} places all images in one compact neighborhood of $S$.

We next exclude collapse to a point for this fixed-area sequence.
If $u_j(\overline{D})\subset B_{r_j}(q_j)$, with $q_j\in\Gamma_j$ and $r_j\to0$, smoothness of $S$ gives
\[
 |(x-q_j)\cdot n_S(x)|\le C_Sr_j^2
 \quad(x\in\Gamma_j).
\]
First variation with $X(x)=x-q_j$ then yields
\[
 2a\le2r_j\sqrt{aW(u_j)}+C_Sr_j^2L(u_j)\longrightarrow0,
\]
a contradiction.

By Theorem~\ref{thm:intro-main-compactness}, a subsequence converges to a bubble tree with at least one nonconstant component; the area identity holds and the Willmore energy is lower semicontinuous. Every nonconstant bordered component has disk domain because the original domain is $D$.
The identity \eqref{ga:signed-contact} passes to every limiting boundary portion corresponding to $\partial D$ by strong $L^2$ convergence of the Sobolev boundary traces of the tangential and transverse coordinate derivatives of the pulled-back maps, as in Lemma~\ref{lem:app-contact-measures}.
Each limiting disk $u_\infty^i$ has $\cA(u_\infty^i)\le a$ and $W(u_\infty^i)\le w_0<2\pi$.
By \eqref{eq:sec2-boundary-length-area-willmore} and \eqref{ga:disk-cost},
\[
 W(u_\infty^i)
 \ge2\pi-\kappa C_S\bigl(a+\sqrt{2\pi a}\bigr)>\pi
\]
when $a$ is small enough.
Two such disks would have total Willmore energy greater than $2\pi$, contrary to lower semicontinuity and $W(u_j)\le w_0<2\pi$. A nonconstant closed component would have Willmore energy at least $4\pi$ and is excluded for the same reason.
The bubble tree consequently has exactly one nonconstant disk $u$, and the area identity gives $\cA(u)=a$. Since $u\in\mathcal C(S)$ and lower semicontinuity gives $W(u)\le\beta_S(a)$, the definition of $\beta_S(a)$ forces $W(u)=\beta_S(a)$. Thus $u$ minimizes the Willmore functional over all immersions in $\mathcal C(S)$ with prescribed area $a$. If $u$ had a branch point, its weighted total branch order would satisfy $\mathfrak b(u)\ge1/2$, and
\[
 3\pi\le W(u)+\kappa L(u)
 \le2\pi+\kappa C_S\bigl(a+\sqrt{2\pi a}\bigr)<3\pi,
\]
which is impossible.
\end{proof}

The last argument excludes branches in every sufficiently small-area disk with $W\le2\pi$, whether or not it minimizes.
For these disks, \eqref{eq:sec2-boundary-length-area-willmore}, \eqref{ga:diameter-estimate} and \eqref{ga:disk-defect} give, with $a=2\pi\lambda^2$,
\begin{equation}\label{ga:preliminary-bounds}
 L(u)\le C_S\lambda,
 \qquad \diam u(\overline{D})\le C_S\lambda,
 \qquad \mathcal D(u)\le C_S\lambda,
 \qquad\int|A_u|^2\dd\mu_u\le C_S.
\end{equation}
These preliminary bounds will be improved for minimizers.

We next flatten the supporting hypersurface $S$ and reflect the disk across the resulting plane to obtain a nearly umbilical immersion of a sphere.
Fix $q\in u(\partial{D})$, rotate $n_S(q)$ to $e_3$, and set
\[
 F=\lambda^{-1}R_q(u-q).
\]
By \eqref{ga:preliminary-bounds}, the image of $F$ lies in one fixed ball and $\cA(F)=2\pi$.
In this ball the rescaled hypersurface $\lambda^{-1}R_q(S-q)$ is a graph $z=s(y)$ satisfying
\begin{equation*}
 \|s\|_{C^3}\le C_S\lambda.
\end{equation*}
The estimate is uniform in $q$, by compactness and smoothness of $S$.
Let $n_s$ denote the normal agreeing with $e_3$ at the origin.
Use the normal-coordinate map
\begin{equation*}
 \Phi(y,t)=(y,s(y))+t\,n_s(y),
 \qquad G=\Phi^{-1}\circ F.
\end{equation*}
Define the normal-coordinate map on a larger fixed ball and extend it with a cutoff outside that ball.
For small $\lambda$ this is a diffeomorphism with
\[
 \|\Phi-\operatorname{Id}\|_{C^2}\le C_S\lambda.
\]
At $t=0$, $D\Phi(e_3)=n_s$ and the horizontal tangent plane maps onto the tangent plane of the graph.
Thus $G$ has boundary in $\{x_3=0\}$ with exact outward conormal $e_3$.

\begin{lemma}\label{ga:tracefree-transport}
Let $u\in\mathcal C(S)$ satisfy $\cA(u)=2\pi\lambda^2$ and $W(u)\le2\pi$, with $\lambda>0$ sufficiently small. For $q\in u(\partial D)$, let $F=\lambda^{-1}R_q(u-q)$ and $G=\Phi^{-1}\circ F$ be the rescaled and flattened maps constructed above. Then
\begin{align}
 |\cA(G)-2\pi|&\le C_S\lambda,\qquad
 \mathcal D(G) \le C_S\bigl(\mathcal D(u)+\lambda^2\bigr).\label{ga:transport-defect}
\end{align}
Here $\mathcal D(G)$ denotes the same trace-free curvature integral as in \eqref{ga:beta-definition}, computed for the immersion $G$, and $C_S$ is independent of $u$, $q$, and $\lambda$.
\end{lemma}

\begin{proof}
Put $\Psi=\Phi^{-1}$ and choose an orthonormal tangent frame for $F$ at a point away from its branch set.
The metric $G^*g_{\R^3}$ has matrix $g'=I+E$ with $|E|\le C\lambda$ in this frame.
Let $n$ and $n'$ be the unit normals of $F$ and $G$, respectively. In this frame, write $A_{ij}$ for the scalar second fundamental form of $F$. The scalar second fundamental form of $G$ is
\[
 B'_{ij}=bA_{ij}+Q_{ij},
 \qquad b=n'\cdot D\Psi n,
 \qquad |b|\le C,\quad |Q|\le C\lambda.
\]
Indeed, the component of $D^2F$ tangent to the image of $F$ is mapped to a vector tangent to the image of $G$ by $D\Psi$, so its projection onto $n'$ vanishes. The remaining error is $D^2\Psi$ evaluated on two unit tangent vectors.
Write $A=A^\circ+(H/2)I$, where $H$ is the scalar trace mean curvature of $F$.
The tracefree projection with respect to $g'$ sends $I$ to a tensor of size $O(\lambda)$, since $g'=I+O(\lambda)$.
Consequently
\[
 |(B')^\circ|_{g'}
 \le C|A^\circ|+C\lambda|H|+C\lambda.
\]
Squaring and integrating gives
\[
 \mathcal D(G)\le C \mathcal D(F)+C\lambda^2\int_D(1+|H|^2)\dd\mu_F.
\]
The area and Willmore bounds prove \eqref{ga:transport-defect}; $\mathcal D(F)=\mathcal D(u)$ by scaling invariance.
The area estimate follows from $g'=I+O(\lambda)$.
All formulas are valid almost everywhere for weak immersions by the Sobolev chain rule.
\end{proof}

Reflect $G$ in the plane on a second copy of the disk, and denote the resulting immersion of the sphere by $\mathcal G$.
In boundary coordinates $(s,t)$ with $t=0$ on $\partial D$, the Sobolev boundary trace of $G$ takes values in $\{x_3=0\}$, the trace of $\partial_sG$ is horizontal, and the trace of $\partial_tG$ is vertical. Reflection therefore matches the Sobolev traces of the value and of both coordinate first derivatives across the common boundary of the two disks. The resulting double $\mathcal G$ is a weak $W^{2,2}$ immersion with finite total curvature.
Since $u$ has branch order zero, its branch expansions give upper and lower bounds for the conformal factor on a finite cover; see Lemma~\ref{lem:general-boundary-removability} and \cite[Theorem~3.1]{kuwert2012}. Composing with $\Phi^{-1}$ preserves the corresponding metric bounds. Hence the induced metrics of $G$ and $\mathcal G$ are bounded above and below by positive multiples of smooth metrics on their domains. These bounds may depend on the fixed disk $u$.
The doubled area is $2\cA(G)$ and the doubled tracefree curvature is $2\mathcal D(G)$.
Set
\begin{equation*}
 r=\sqrt{\frac{\cA(G)}{2\pi}}=1+O_S(\lambda),
 \qquad
 \delta=\left(\int_{\mathcal G}|A^\circ|^2\dd\mu\right)^{1/2}.
\end{equation*}
Then, we have 
\begin{equation}\label{ga:double-smallness}
 \delta\le C_S\sqrt{\mathcal D(u)+\lambda^2}\longrightarrow0,
 \qquad
 W(\mathcal G)=4\pi+\frac12\delta^2<8\pi.
\end{equation}
Here we use closed Gauss--Bonnet on the doubled unbranched sphere. The preliminary bound $\mathcal D(u)\le C_S\lambda$ allows us to choose $\lambda_0(S)>0$ such that $1/2\le r\le2$ and $\delta$ is below the universal rigidity threshold whenever $0<\lambda<\lambda_0(S)$, for every disk under consideration.

We first extend the quantitative sphere rigidity of De Lellis--M\"uller~\cite{appDLM05,appDLM06} to the present weak class. We then adjust the conformal parametrization to preserve the reflection symmetry, so that the equator parametrizes $G(\partial D)$.

\begin{lemma}\label{ga:weak-dlm}
Let $U\in W^{2,2}(\mathbb S^2,\R^3)$ be an unbranched weak immersion with finite total curvature and induced metric locally bounded above and below by positive multiples of a smooth metric on $\mathbb S^2$.
Suppose $\cA(U)=4\pi r^2$ with $r>0$, and put $\delta=\|A_U^\circ\|_{L^2(\mathbb S^2,\dd\mu_U)}$. If $\delta$ is smaller than a universal positive constant, then $U$ is a continuous embedding. Moreover, there exist $c\in\R^3$ and a conformal homeomorphism $\psi\in W^{2,2}(\mathbb S^2,\R^3)$ onto $U(\mathbb S^2)$ such that
\begin{equation}\label{ga:dlm-estimate}
 \|\psi-(c+r\operatorname{id})\|_{W^{2,2}(\mathbb S^2)}\le Cr\delta,
 \qquad
 \|e^v/r-1\|_{L^\infty(\mathbb S^2)}\le C\delta.
\end{equation}
Here $\psi^*g_{\R^3}=e^{2v}g_{\mathbb S^2}$, $\operatorname{id}$ is the inclusion of the unit sphere in $\R^3$, and the norms in \eqref{ga:dlm-estimate} use the round metric on $\mathbb S^2$.
\end{lemma}

\begin{proof}
By scaling, we may assume $r=1$.
Finite-total-curvature uniformization gives a bilipschitz change of variables making $U$ a weak conformal immersion of class $W^{2,2}\cap W^{1,\infty}$; see \cite[Proposition~I.2]{appMR}.
After this change of variables, $U^*g_{\R^3}$ is bounded above and below by positive multiples of the round metric $g_{\mathbb S^2}$.
The compact-manifold hypothesis in~\cite{appMR} is satisfied by placing the bounded image $U(\mathbb S^2)$ in a sufficiently large flat torus.
Closed weak Gauss--Bonnet gives $W(U)=4\pi+\delta^2/2<8\pi$.
Lemma~\ref{lem:app-local-multiplicity} gives a local area density of $1$ from each point of $U^{-1}(x)$. The pushforward varifold has no singular first-variation term, and its mean curvature is the fiberwise average of the parametrized mean curvature, whose squared integral is at most $4W(U)$. The weak Li--Yau inequality, recalled in \cite[Section~2.2.2, equation~(2.12)]{appMS}, therefore gives $4\pi\#U^{-1}(x)\le W(U)$ for every image point $x$ (first for any finite set of preimages). Since $W(U)<8\pi$, $U$ is injective; continuity on the compact sphere makes it an embedding.

Applying Lemma~\ref{lem:app-weak-approximation} to $U$ on the fixed sphere gives smooth immersions $U_j$ converging strongly in $W^{2,2}$ and uniformly, whose induced metrics have common upper and positive lower bounds and whose areas, Willmore energies, and integrals of $|A^\circ|^2$ converge to those of $U$. These bounds may depend on $U$, which is sufficient for this approximation passage. Rescaling $U_j$ by factors tending to one gives $\cA(U_j)=4\pi$ while preserving the stated convergence.
For large $j$, $W(U_j)<8\pi$, so the Li--Yau inequality~\cite{Li-Yau1982} makes each $U_j$ an embedding.
Applying \cite[Theorem~1.1]{appDLM05} and \cite[Theorem~2]{appDLM06} to the smooth embedded spheres $U_j(\mathbb S^2)$ gives conformal parametrizations satisfying \eqref{ga:dlm-estimate} with $\delta$ replaced by $\delta_j=\|A_{U_j}^\circ\|_{L^2}\to\delta$ and a universal constant.

The images $U_j(\mathbb S^2)$ remain in a bounded set, so the translation vectors are bounded. Thus the parametrizations $\psi_j$ are bounded in $W^{2,2}$, and their conformal factors have a uniform positive lower and upper bound.
After a subsequence they converge weakly in $W^{2,2}$, strongly in every finite $W^{1,p}$, and uniformly to a weak conformal immersion $\psi$.
Both estimates remain valid for the limiting parametrization.
We next check that its conformal structure agrees with the one fixed by the weak uniformization of $U$.
Set
\[
 T_j=U_j^{-1}\circ\psi_j:\mathbb S^2\longrightarrow \mathbb S^2.
\]
These maps are smooth diffeomorphisms.
The metric bounds for $U_j$ and the conformal-factor bounds for $\psi_j$ give uniform upper and lower singular-value bounds for $DT_j$.
Thus $T_j$ and their inverses are uniformly Lipschitz for the round metric.
A subsequence converges uniformly to a bilipschitz homeomorphism $T$, and uniform convergence gives $\psi=U\circ T$.
In local coordinates,
\[
 DT_j=
 \left[(DU_j^tDU_j)^{-1}DU_j^t\right]
 \circ T_j\,D\psi_j.
\]
The bracketed coefficients are uniformly bounded and converge strongly in every finite $L^p$.
Their compositions with $T_j$ also converge strongly: approximate the limiting coefficient in $L^p$ by a continuous function and control the approximation error using the uniform bounds for change of variables.
Together with strong convergence of $D\psi_j$, the displayed identity gives strong convergence of $DT_j$ in finite $L^p$.
The identities $\psi_j^*g_{\R^3}=T_j^*(U_j^*g_{\R^3})$ therefore pass to the limit.
Since $U$ and $\psi$ are conformal, $T$ is weakly conformal, hence is a conformal or anticonformal M\"obius transformation.
This proves the required conformal compatibility.
The image is $U(\mathbb S^2)$, and $\psi=U\circ T$ is injective because $U$ is.
This proves the lemma.
\end{proof}

Let $R(x_1,x_2,x_3)=(x_1,x_2,-x_3)$.
For the double $\mathcal G$, the interchange of the two halves is an isometry of its induced metric and corresponds geometrically to $R$.
After weak uniformization, the interchange map is bilipschitz and anticonformal almost everywhere.
It therefore satisfies the weak Cauchy--Riemann equations and is an antiholomorphic M\"obius involution.

\begin{lemma}\label{ga:equivariant-gauge}
Let $\mathcal G:\mathbb S^2\to\R^3$ be the reflected double of $G$ constructed above, with $r=(\cA(G)/(2\pi))^{1/2}$ and $\delta=\|A_{\mathcal G}^\circ\|_{L^2(\dd\mu_{\mathcal G})}$. There exist a conformal homeomorphism $w\in W^{2,2}(\mathbb S^2,\R^3)$ onto $\mathcal G(\mathbb S^2)$ and a vector $c\in\{x_3=0\}$ such that
\begin{equation}\label{ga:equivariance}
 w\circ R=R\circ w,
 \qquad
 \|w-(c+r\operatorname{id})\|_{W^{2,2}(\mathbb S^2)}\le C\delta.
\end{equation}
Here $R$ acts on $\mathbb S^2$ by restriction of the reflection in $\{x_3=0\}$, and $w$ maps the equator $\mathbb S^2\cap\{x_3=0\}$ onto $G(\partial D)$.
\end{lemma}

\begin{proof}
To obtain a parametrization that commutes with reflection, we start with the conformal parametrization $\psi$ satisfying \eqref{ga:dlm-estimate}.
Sobolev embedding gives $\psi=c+r\operatorname{id}+O_{C^0}(\delta)$, since $r$ is bounded above and below.
The image is invariant under $R$.
Its maximum and minimum heights are $c_3+r+O(\delta)$ and $c_3-r+O(\delta)$; symmetry gives $|c_3|\le C\delta$.
We may therefore project $c$ to the horizontal plane and absorb the change in the preceding estimate.

Let $\tau$ be the induced anticonformal involution, so that $\psi\circ\tau=R\circ\psi$.
The preceding $C^0$ estimate implies
\[
 \|\tau-R\|_{C^0(\mathbb S^2)}\le C\delta.
\]
Put $\varphi=\tau R$.
This is a M\"obius transformation close to the identity.
Evaluation at three distinct points gives smooth local coordinates on the M\"obius group, so $\varphi$ is $O(\delta)$-close to the identity in every fixed $C^k$ norm.
Take the local square root $\psi_0=\exp\bigl(\tfrac12\log\varphi\bigr)$. 
The identity $\tau^2=\operatorname{id}$ gives $R\varphi R=\varphi^{-1}$ and therefore $R\psi_0 R=\psi_0^{-1}$.
It follows that
\[
 \psi_0 R\psi_0^{-1}=\psi_0^2R=\varphi R=\tau,
 \qquad \|\psi_0-\operatorname{id}\|_{C^k}\le C_k\delta.
\]
Set $w=\psi\circ \psi_0$.
The conjugacy gives $w\circ R=R\circ w$.
Composition by $\psi_0$ preserves the $W^{2,2}$ estimate: its first two derivatives and inverse derivatives are uniformly bounded, and $\psi_0-\operatorname{id}=O(\delta)$ in $C^2$.
The fixed circle of $R$ consequently maps to the fixed circle of the interchange involution, which is $G(\partial D)$.
\end{proof}

\begin{lemma}\label{ga:boundary-moment}
Let $G:D\to\R^3$, $r$, and $\delta$ be as in Lemma~\ref{ga:equivariant-gauge}. Let $B$ be a symmetric bilinear form on the plane $\{x_3=0\}$, with Euclidean norm $|B|$ and trace $\operatorname{tr}B$. Write $g_G=G^*g_{\R^3}$, and let $T$ be the unit tangent to the parametrized curve $G|_{\partial D}$. Then
\begin{equation}\label{ga:moment-estimate}
 \left|\int_{\partial D}B(T,T)\dd s_{g_G}
       -\pi r\operatorname{tr} B\right|
 \le C|B|\delta.
\end{equation}
\end{lemma}

\begin{proof}
Restrict the equivariant parametrization $w$ to the equator.
The trace theorem applied to \eqref{ga:equivariance} gives
\[
 \|\partial_\theta w-r\partial_\theta\operatorname{id}\|_{L^2(\mathbb S^1)}
 \le C\delta.
\]
The function
\[
 Q_B(z)=\frac{B(z,z)}{|z|}\quad(z\ne0),\qquad Q_B(0)=0,
\]
is globally Lipschitz, with Lipschitz constant at most $3|B|$.
This follows by differentiating away from zero and extending continuously at zero.
Thus the integral of $Q_B(\partial_\theta w)$ differs by at most $C|B|\delta$ from that of $Q_B(r\partial_\theta\operatorname{id})$.
The latter integral is
\[
 r\int_0^{2\pi}B\bigl((-\sin\theta,\cos\theta),
                         (-\sin\theta,\cos\theta)\bigr)\dd\theta
 =\pi r\operatorname{tr} B.
\]
Since $w$ parametrizes $G(\partial D)$ once along the equator, the first integral equals $\int_{\partial D}B(T,T)\dd s_{g_G}$. This proves \eqref{ga:moment-estimate}.
\end{proof}

We now transfer Lemma~\ref{ga:boundary-moment} back to $S$ to estimate the geodesic curvature integral $\int_{\partial D}k_{g_u}\,\dd s_{g_u}$ of $\partial D$ in the induced metric $g_u$. This is the boundary term in \eqref{ga:disk-defect}. Combined with the comparison bound \eqref{ga:strict-upper}, its estimate will improve \eqref{ga:preliminary-bounds} and prove Theorem~\ref{app:small-disks}.

\begin{prop}\label{ga:boundary-expansion}
There exist $\lambda_0(S)>0$ and $C_S<\infty$ such that every $u\in\mathcal C(S)$ with $\cA(u)=2\pi\lambda^2$, $0<\lambda<\lambda_0(S)$, and $W(u)\le2\pi$ satisfies
\begin{equation}\label{ga:curvature-remainder}
 \left|\int_{\partial{D}}k_{g_u}\dd s_{g_u}-\pi H_S(q)\lambda\right|
 \le C_S\bigl(\lambda\sqrt{\mathcal D(u)}+\lambda^2\bigr)
 \qquad(q\in u(\partial{D})).
\end{equation}
\end{prop}

\begin{proof}
For $q\in u(\partial D)$, use the rescaling and normal coordinates defined above. Identify $T_qS$ with $\{x_3=0\}$ by $R_q$, and write $B=h^S(q)$ in this orthonormal frame.
Extend this tensor to ambient vectors by tangential projection.
The diameter bound and smoothness of $S$ imply
\[
 h^S(u)=h^S(q)+O_S(\lambda)
\]
on the whole boundary.
The boundary length satisfies $L(u)=\lambda L(F)$, and $L(F)$ is uniformly bounded.
Let $T$ denote the unit tangent to $G|_{\partial D}$. The map $\Phi$ changes rescaled boundary lengths and tangent directions by $O_S(\lambda)$.
Consequently
\[
 \int_{\partial{D}}k_{g_u}\dd s_{g_u}
 =\lambda\int_{\partial D}B(T,T)\dd s_{g_G}+O_S(\lambda^2).
\]
In a boundary parameter, this follows from the Lipschitz bound for $Q_B$ in Lemma~\ref{ga:boundary-moment} and $D\Phi=I+O_S(\lambda)$.
That lemma now gives
\[
 \int_{\partial{D}}k_{g_u}\dd s_{g_u}
 =\pi r H_S(q)\lambda+O_S(\lambda\delta+\lambda^2).
\]
Since $r=1+O_S(\lambda)$ and $\delta\le C_S\sqrt{\mathcal D(u)+\lambda^2}$, this is \eqref{ga:curvature-remainder}.
\end{proof}

\begin{proof}[Proof of Theorem~\ref{app:small-disks}]
Fix $a>0$ sufficiently small. Proposition~\ref{ga:weak-attainment} gives an unbranched minimizer $u\in\mathcal C(S)$ with $\cA(u)=a$; put $\lambda=\lambda(a)$.
It has $W(u)<2\pi$.
By \eqref{ga:disk-defect}, \eqref{ga:strict-upper} and \eqref{ga:curvature-remainder},
\[
 \frac12\mathcal D(u)+\pi\bigl(H_*-H_S(q)\bigr)\lambda
 \le C_S\lambda\sqrt{\mathcal D(u)}+C_S\lambda^2.
\]
Young's inequality gives
\[
 C_S\lambda\sqrt{\mathcal D(u)}\le\frac14\mathcal D(u)+C'_S\lambda^2.
\]
After absorbing $\mathcal D(u)/4$, both $\mathcal D(u)/4$ and $\pi(H_*-H_S(q))\lambda$ on the left are nonnegative.
Thus $\mathcal D(u)\le C_S\lambda^2$ and $H_*-H_S(q)\le C_S\lambda$, proving \eqref{ga:tracefree-bound} and \eqref{ga:location-deficit}.
Using these estimates once more in \eqref{ga:curvature-remainder} and \eqref{ga:disk-defect} gives
\[
 2\pi-\pi H_*\lambda-C_S\lambda^2
 \le W(u)
 \le2\pi-\pi H_*\lambda+C_S\lambda^2.
\]
This proves \eqref{ga:energy-expansion}.

The reflected double satisfies \eqref{ga:double-smallness} and is injective by the Li--Yau argument in Lemma~\ref{ga:weak-dlm}.
Therefore $G$, and hence $u$, is injective.
An interior point of $G$ lying in the plane would coincide with its distinct reflected point in the other half of the double.
This contradicts injectivity, so the connected interior image lies entirely in one open half-space.
The conormal $e_3$ selects the lower half-space as follows. In a conformal collar coordinate $(s,t)$ with $t\ge0$ directed inward, for almost every $s$ the function $t\mapsto G_3(s,t)$ belongs to $W^{2,2}$, has value zero at $t=0$, and has strictly negative first derivative there, by the signed contact condition and the metric lower bound. It is negative for sufficiently small positive $t$. Since the connected interior image contains no point of the plane, the whole interior lies in the lower half-space.
The normal-coordinate map identifies the lower side with $\Omega$ in the coordinate neighborhood of $S$ containing $u(\overline D)$.
This proves \eqref{ga:properness}.
Continuity and injectivity on the compact closed disk make $u:\overline D\to\R^3$ a homeomorphism onto its image.

Now $\delta\le C_S\lambda$ in \eqref{ga:equivariance}, while $r=1+O_S(\lambda)$.
The estimate $w=c+r\operatorname{id}+O_{C^0}(\lambda)$ sends the south pole to the lower side. Thus the lower hemisphere parametrizes the image of $G$, and composing the restriction of $w$ with $\Phi$ gives a parametrization of the rescaled disk.
The Sobolev chain rule and $W^{2,2}\hookrightarrow W^{1,4}$ give
\[
 \|\Phi\circ w-(c+\operatorname{id})\|_{W^{2,2}(\mathbb S^2_-)}\le C_S\lambda.
\]
Write $c=(c_1,c_2,0)$. Since the rescaled image lies in a fixed ball, $c$ is uniformly bounded. Let $p_a\in S$ correspond in the rescaled coordinates to $\Phi((c_1,c_2),0)$. Then $|p_a-q|=O_S(\lambda)$ and $n_S(p_a)-n_S(q)=O_S(\lambda)$. Composing $R_q$ with a rotation whose angle is $O_S(\lambda)$ gives $R_a n_S(p_a)=e_3$. The estimates for $\Phi-\operatorname{Id}$ show that this recentering and rotation preserve the $O_S(\lambda)$ bound, proving \eqref{ga:quantitative-hemisphere}.
Sobolev embedding supplies the corresponding $C^0$ estimate and therefore the Hausdorff bound $C_S\lambda^2$ in the original coordinates. Finally, for any sequence $a\to0$, \eqref{ga:location-deficit} and compactness of $S$ place all limits of boundary points in $\operatorname{argmax}_SH_S$; the diameter bound gives the same conclusion for all image points.
\end{proof}

\begin{coro}\label{ga:nondegenerate-maxima}
In the setting of Theorem~\ref{app:small-disks}, suppose the Hessian of $H_S:S\to\R$ is negative definite at every point of $\operatorname{argmax}_SH_S$. Then, for all sufficiently small $a>0$, every minimizer $u$ of $\beta_S(a)$ satisfies
\begin{equation}\label{ga:nondegenerate-location}
 \sup_{x\in u(\overline{D})}\operatorname{dist}\bigl(x,\operatorname{argmax}_SH_S\bigr)
 \le Ca^{1/4}.
\end{equation}
\end{coro}

\begin{proof}
Compactness of $S$ and the negative definite Hessian at each maximum give
\[
 H_*-H_S(q)\ge c_S\dist\bigl(q,\operatorname{argmax}_SH_S\bigr)^2
\]
for $q\in S$ near $\operatorname{argmax}_SH_S$. On the complement of any fixed neighborhood of the maximum set, $H_*-H_S$ has a positive minimum.
Equation~\eqref{ga:location-deficit} places the boundary within $C_S\sqrt\lambda$ of that set.
The diameter estimate adds only $C_S\lambda$ for interior image points.
Since $\sqrt\lambda=(a/(2\pi))^{1/4}$, this proves \eqref{ga:nondegenerate-location}.
\end{proof}

\begin{coro}\label{ga:almost-minimizers}
Let $u_j\in\mathcal C(S)$ satisfy $\cA(u_j)=2\pi\lambda_j^2$, with $\lambda_j>0$ and $\lambda_j\to0$, and suppose
\[
 W(u_j)\le\beta_S(2\pi\lambda_j^2)+\varepsilon_j,
 \qquad \varepsilon_j=o(\lambda_j),\quad\varepsilon_j\ge0.
\]
Then, for all sufficiently large $j$ and every $q_j\in u_j(\partial D)$, we have
\begin{equation}\label{ga:almost-estimates}
 \mathcal D(u_j)+\lambda_j\bigl(H_*-H_S(q_j)\bigr)
 \le C_S(\lambda_j^2+\varepsilon_j).
\end{equation}
For all sufficiently large $j$, $u_j$ is unbranched and is a continuous embedding of $\overline D$, with $u_j(D)\subset\Omega$.
More precisely, there are $p_j\in S$, $R_j\in O(3)$ with $R_jn_S(p_j)=e_3$, and $W^{2,2}$ parametrizations $v_j:\mathbb S^2_-\to\R^3$ satisfying
\[
 v_j(\mathbb S^2_-)=\lambda_j^{-1}R_j\bigl(u_j(\overline D)-p_j\bigr),\qquad
 \|v_j-\operatorname{id}_{\mathbb S^2_-}\|_{W^{2,2}(\mathbb S^2_-)}
 \le C_S(\lambda_j+\sqrt{\varepsilon_j}).
\]
In particular, every accumulation point of the images $v_j(\S^2_{-})$ belongs to $\operatorname{argmax}_SH_S$.
\end{coro}

\begin{proof}
Since $H_*>0$, \eqref{ga:energy-expansion} and $\varepsilon_j=o(\lambda_j)$ imply $W(u_j)<2\pi$ for all sufficiently large $j$. The branch exclusion following Proposition~\ref{ga:weak-attainment} and the bounds \eqref{ga:preliminary-bounds} therefore apply. Combining \eqref{ga:disk-defect}, \eqref{ga:strict-upper}, and \eqref{ga:curvature-remainder}, we obtain
\[
 \frac12\mathcal D(u_j)+\pi\bigl(H_*-H_S(q_j)\bigr)\lambda_j
 \le C_S\lambda_j\sqrt{\mathcal D(u_j)}
       +C_S\lambda_j^2+\varepsilon_j.
\]
Applying Young's inequality to the first term on the right and absorbing the resulting term on the left gives \eqref{ga:almost-estimates}.
The squared $L^2$ norm of the trace-free second fundamental form of the reflected double is bounded by $C_S(\lambda_j^2+\varepsilon_j)$. Lemmas~\ref{ga:weak-dlm} and~\ref{ga:equivariant-gauge}, followed by the recentering in the proof of Theorem~\ref{app:small-disks}, give the stated estimate for $v_j$. Injectivity of the double gives injectivity of $u_j$; its reflection symmetry and the conormal condition \eqref{ga:signed-contact} give $u_j(D)\subset\Omega$ by the half-space argument in that proof.
Lastly, \eqref{ga:almost-estimates}, divided by $\lambda_j$, gives $H_S(q_j)\to H_*$, and the image diameter tends to zero.
\end{proof}

For closed surfaces, global minimizers with small area concentrate at maxima of ambient scalar curvature; see \cite[Theorem~1.1]{appLM}.
Here the support mean curvature $H_S$ determines the term of order $\sqrt a$ in the Willmore energy expansion.
The local construction in~\cite{appAK} supplies the coefficient of the leading correction, while bubble tree compactness gives global attainment and the boundary integral estimate controls the shape and location of every minimizer.

\subsection{Willmore energy lower bound under a calibration constraint}\label{bar:section}

In this subsection, we apply Miura's \cite{appMiura} diameter estimate to prove Theorem~\ref{app:calibration-barrier} for the closed Riemannian manifold $(N,\bar g)$, the supporting submanifold $\cK$, and the calibrated surface $S_0$ specified there. We then interpret the obstruction through the bubble tree limit and its no-neck property.

\begin{proof}[Proof of Theorem~\ref{app:calibration-barrier}]
Choose $f\in C^\infty(N,[0,1])$ equal to $0$ near $S_0$ and to $1$ on a neighborhood of $\cK$, and set
\[
 \varepsilon=1-\max_{\{f\geq1/4\}}\|\omega\|_{*,\bar g}>0,
 \qquad D=\cA(u)-m\geq0.
\]
The area excess gives
\begin{equation}\label{bar:strip-area}
 \cA\bigl(u;\{f\circ u\geq1/4\}\bigr)\leq D/\varepsilon.
\end{equation}
We prove the estimate for $D<\varepsilon m/8$; the remaining range follows by increasing the additive constant.

We excise small disjoint disks about interior branches and caps about boundary branches. In boundary coordinates, we take the cap edges to be graphs $y=\rho\phi(x/\rho)$, where the fixed positive bump $\phi$ is flat at its endpoints. These cuts leave a smooth connected complement $\Sigma_\rho$ containing an arc of $\partial\Sigma$. The branch expansion $\lambda=m_u(p)\log|z|+\omega_p$, with bounded $\omega_p$, gives the following bound for the total image length of the artificial cuts:
\begin{equation}\label{bar:cut-length}
 L_{\mathrm{art}}(u;\rho)\leq C(u)\sum_p\rho^{m_u(p)+1}\longrightarrow0.
\end{equation}
For each fixed $\rho>0$, the induced metric on the compact surface $\Sigma_\rho$ is bounded above and below by positive multiples of a smooth metric there. Apply Lemma~\ref{lem:app-weak-approximation} to obtain smooth immersions $v_{\rho,k}:\Sigma_\rho\to N$, and write
\[
 d_{\rho,k}=\cA(v_{\rho,k})-\int_{\Sigma_\rho}v_{\rho,k}^*\omega\geq0.
\]
For fixed $\rho$, strong $W^{2,2}$ approximation gives $d_{\rho,k}\to d_\rho:=\cA(u;\Sigma_\rho)-\int_{\Sigma_\rho}u^*\omega$ as $k\to\infty$, together with convergence of the areas, Willmore energies, and induced lengths of the artificial boundary curves. Integrability gives $d_\rho\to D$ as $\rho\to0$. Although the approximants need not satisfy a period or boundary condition, uniform convergence keeps the images of $\partial\Sigma_\rho\cap\partial\Sigma$ under $v_{\rho,k}$ in $\{f=1\}$. Choosing small $\rho$ and then large $k$, we obtain
\[
 \cA(v_{\rho,k})\geq m/2,\qquad d_{\rho,k}<\varepsilon m/4.
\]
Consequently $f\circ v_{\rho,k}$ takes a value below $1/4$: otherwise the pointwise deficit would give $d_{\rho,k}\geq\varepsilon\cA(v_{\rho,k})$.

For brevity, put $v=v_{\rho,k}$ and $d=d_{\rho,k}$. Applying coarea to the analogue of \eqref{bar:strip-area}, choose $a\in(1/4,3/8)$, $b\in(5/8,3/4)$
regular for $f\circ v$ and its restriction to the artificial boundary, with total level-set length at most $Cd$. The connectedness of $\Sigma_\rho$ gives a path, interior except at its boundary endpoint, from a point where $f\circ v<a$ to the retained portion of $\partial\Sigma$ where $f\circ v=1$. A segment joins the levels $a$ and $b$ while remaining between them. The component $\Sigma_{a,b}$ of $\{a\leq f\circ v\leq b\}$ containing this segment therefore meets both levels, and
\[
 \cA(v;\Sigma_{a,b})\leq d/\varepsilon,\qquad
 L(v;\partial \Sigma_{a,b})\leq Cd+L_{\mathrm{art}}(v;\rho).
\]
The subsurface $\Sigma_{a,b}$ is disjoint from the retained portions of $\partial\Sigma$, and the parts of its boundary in $\partial\Sigma_\rho\setminus\partial\Sigma$ meet the two level sets transversely at finitely many corners. We round these corners while retaining the interior crossing path. For fixed $(\rho,k)$, the changes in area, boundary length, and the mean-curvature integral tend to zero with the rounding scale.

Fix an isometric embedding $\iota_N:N\hookrightarrow\R^K$. A smooth ambient extension of $f$ is Lipschitz, so $b-a>1/4$ gives a fixed lower bound $d_0>0$ for the Euclidean diameter of $\iota_N\circ v(\Sigma_{a,b})$. We apply Miura's estimate~\cite[Theorem~1.1]{appMiura} to $\iota_N\circ v$ on the connected subsurfaces with rounded corners. Passing to the limit as the rounding scale tends to zero gives
\begin{align*}
 d_0 \leq C\left(\int_{\Sigma_{a,b}}|\mathbf H_{\iota_N\circ v,\R^K}|\dd\mu_v
                   +L(v;\partial \Sigma_{a,b})\right) \leq C\sqrt{W_N(v)d}+Cd+C L_{\mathrm{art}}(v;\rho).
\end{align*}
Here we used the bounded second fundamental form of $\iota_N$ and the Cauchy--Schwarz estimate
\[
 \int_{\Sigma_{a,b}}|\mathbf H_{v,N}|\dd\mu_v
 \leq2\sqrt{W_N(v)\cA(v;\Sigma_{a,b})}.
\]
First let $k\to\infty$ and then $\rho\to0$. Lemma~\ref{lem:app-weak-approximation}, integrability of the area and curvature densities, and \eqref{bar:cut-length} yield
\[
 d_0\leq C\sqrt{W_N(u)D}+CD.
\]
The last inequality excludes $D=0$ and, for small $D>0$, gives $W_N(u)\geq c/D$ by absorption, proving \eqref{bar:inverse} with constants depending only on the fixed geometric data. For smooth immersions, the boundary constraint $u(\partial\Sigma)\subset\cK$ suffices for the argument; orthogonality is not used. For a more general diameter estimate, see Menne--Scharrer~\cite[Corollary~2]{appMenneScharrer}.

We now construct disks with free boundary whose Willmore energies have the same inverse order as the lower bound. On $N=\mathbb S^2\times\mathbb S^1$, let $g_0$ be the product metric with unit sphere factor, put $S_0=\mathbb S^2\times\{0\}$, and let $\omega_0$ be the pullback of the sphere area form. Choose a coordinate ball $U$ disjoint from $S_0$, a coordinate sphere $\cK\subset U$, a primitive $\xi$ of $\omega_0$ on $U$, and $\chi\in C_c^\infty(U)$ equal to $1$ near $\cK$. Then, $\omega=\omega_0-\dd(\chi\xi)$ is closed, vanishes near $\cK$, and equals $\omega_0$ near $S_0$. To modify the metric, start with $\rho_0=2-\cos t$, where $t$ is the circle coordinate, and add a sufficiently large nonnegative bump supported away from $S_0$ and constant on $\operatorname{spt}\chi$. The resulting smooth function $\rho\geq1$ equals $1$ exactly on $S_0$ and satisfies
\[
 \rho^2>\|\omega\|_{*,g_0}\quad\text{on }N\setminus S_0.
\]
Thus $\bar g=\rho^2g_0$ satisfies \eqref{bar:calibration}, with $m=4\pi$, and $S_0$ is minimal. Recall that $\iota:S_0\hookrightarrow N$ is the inclusion and $j_*:H_2(N;\mathbb Z)\to H_2(N,\cK;\mathbb Z)$ is the natural homomorphism. Let $\Gamma=j_*(\iota_*[S_0])$ be the relative class in Theorem~\ref{app:calibration-barrier}. Since $d\omega=0$ and $\omega|_{T\cK}=0$, integration of $\omega$ depends only on the relative class, so every map representing $\Gamma$ satisfies \eqref{bar:period}.

Choose an embedded path from $S_0$ to $\cK$, disjoint from both in its interior and normal at its endpoints. We flow a small coordinate disk in $S_0$ along a vector field normal near the endpoints and rescale by the smooth first hitting time of $\cK$ to obtain a tubular embedding
\[
 F:\overline{D}_\eta\times[0,1]\longrightarrow N,
 \qquad F(x,0)\in S_0,\quad F(x,1)\in\cK,\quad
 \partial_tF(x,1)\perp T_{F(x,1)}\cK.
\]
Its image meets $S_0\cup\cK$ only in the end disks. Remove $F( D_{2r}\times\{0\})$ from $S_0$ and insert the annulus
\[
 F\bigl(rR(s)(\cos\theta,\sin\theta),rZ(s)\bigr),
 \qquad 0\leq s\leq1,
\]
where a fixed smooth embedded meridian $(R,Z)$ runs from $(2,0)$ horizontally to $(1,1)$ vertically, has nonzero velocity, and agrees with the corresponding straight segments near its ends. Attach the cylinder $F(r\mathbb S^1\times[r,1])$. Both joins are smooth, and the resulting surface $\mathcal S_r$ is an embedded disk meeting $\cK$ orthogonally.

To identify the relative class, consider first the cornered surface obtained by removing the radius-$r$ initial disk and attaching the full cylinder. Cap this surface in $\cK$ by $F(\overline{D}_r\times\{1\})$. With compatible orientations, its difference from $S_0$ bounds the solid tube $F(\overline{D}_r\times[0,1])$. Replacing the corner by the displayed meridian changes the surface by a boundary in a coordinate ball, so the oriented inclusion of $\mathcal S_r$ represents $\Gamma$ and $\int_{\mathcal S_r}\omega=4\pi$.

The pullback $F^*\bar g$ is a smooth positive definite metric on $\overline D_\eta\times[0,1]$. The cylinder area density equals $rj_0(\theta,t)+O(r^2)$, where
\[
 j_0(\theta,t)=
 \left|D_xF(0,t)(-\sin\theta,\cos\theta)
             \wedge\partial_tF(0,t)\right|_{\bar g}>0.
\]
The removed disk and meridian annulus have area $O(r^2)$. Therefore, with $a_1=\int_0^1\int_0^{2\pi}j_0\dd\theta\dd t>0$,
\begin{equation}\label{bar:tube-estimates}
 \cA(\mathcal S_r)=4\pi+a_1r+O(r^2),\qquad W_N(\mathcal S_r)\leq C/r.
\end{equation}
Indeed, the cylinder has area $O(r)$ and $|A|\leq C/r$; the meridian annulus has area $O(r^2)$ and the same curvature bound, while the unchanged part of $S_0$ is minimal. After factoring out $r$, the cylinder density is smooth down to $r=0$, and the other area terms are $r^2$ times smooth functions. Thus $\frac{\dd}{\dd r}\cA(\mathcal S_r)=a_1+O(r)>0$. Every sufficiently small positive excess occurs in this family. Uniformization supplies conformal disk parametrizations of the constructed disks, so \eqref{bar:tube-estimates} gives the upper bound in \eqref{bar:sharp-order}. The lower bound follows from \eqref{bar:inverse} over the full class of branched conformal disk immersions with free boundary on $\cK$ and relative homology class $\Gamma$. The calibration inequality and the constructed family yield the asserted area infimum. By \eqref{bar:inverse}, the Willmore energy diverges along every sequence in this relative class whose areas tend to $4\pi$.
\end{proof}

For a fixed compact domain $\Sigma$, Theorem~\ref{thm:intro-main-compactness} also gives a direct interpretation of the qualitative obstruction. Suppose that the maps satisfy \eqref{bar:period} and $\cA(u_j)\to m$ with bounded $W_N(u_j)$. If the images collapse to a point $q$, then $q\in\cK$ because the boundary is nonempty. The strict comass gap near $\cK$ would give, for some fixed $\varepsilon>0$ and all sufficiently large $j$,
\[
m=\int_\Sigma u_j^*\omega
\leq(1-\varepsilon)\cA(u_j)
\longrightarrow(1-\varepsilon)m<m,
\]
a contradiction. Strong convergence on cores away from branch, nodal, and attaching points, the estimate $|\int_Ru_j^*\omega|\leq\cA(u_j;R)$ on omitted regions, and the area identity give
\[
 \sum_i\left(\cA(u_\infty^i)
             -\int(u_\infty^i)^*\omega\right)=0.
\]
Since every summand is nonnegative, each nonconstant component has image in $S_0$; continuity includes its branch, nodal, and attaching points. The stratified domain is connected and the component maps agree at the nodal and attaching points, so every constant component also maps into $S_0$. By the no-neck property, the entire Hausdorff limit lies in $S_0$. But $u_j(\partial\Sigma)\subset\cK$ supplies a limit point in $\cK$, contradicting $S_0\cap\cK=\varnothing$. The contradiction persists when $S_0$ is a sphere and hence an allowed bubble component.

\subsection{Intrinsic Willmore energy decomposition for HSL limits}\label{hs:section}

In this subsection, we prove Theorem~\ref{app:hsl-energy} and determine when the additional terms in \eqref{hs:complete-formula} vanish.
We use the mean curvature $1$-form $\alpha$ from \eqref{eq:intro-geometric-mean-curvature-form} and the Willmore normalization \eqref{hs:normalization}, with $h$ any smooth metric on $\Sigma$ compatible with its conformal structure.  For the HSL immersions in Theorem~\ref{main theorem}, $\alpha$ extends smoothly across the branch points and satisfies \eqref{sys2}.
Since $(M,\omega,J,\bar g)$ is K\"ahler-Einstein, $\operatorname{Rc}=\lambda\omega$ and $u^*\omega=0$ imply $d\alpha=0$.
Thus $\alpha$ is closed and coclosed with zero normal trace: it is the harmonic representative with absolute boundary conditions of its class in $H^1(\Sigma;\R)$. In particular, if $\alpha=d\varphi$, then
\[
 \|\alpha\|_{L^2(\Sigma,h)}^2
 =\int_\Sigma\langle\alpha,d\varphi\rangle_h\,\dd\vol_h=0.
\]
Every HSL disk satisfying these homogeneous boundary conditions is therefore minimal, as is every HSL sphere. For closed spheres, this conclusion also appears in~\cite{chen2021}.

We assume the homogeneous system \eqref{sys2}, including $\alpha_n(\nu_{h_n})=0$ on each parametrized boundary component, as in Theorem~\ref{app:hsl-energy}.
The variational characterization is given in Theorem~\ref{thm:1st_vari_natural_support} and Corollary~\ref{coro:hsl-support-geometries}.

We begin with the identity for harmonic forms on closed surfaces. The cylinder analysis is classical; see \cite[Section~3]{appBMMS}. The calculation below keeps track of both constant Fourier coefficients and the normalization of the collar contribution.

\begin{lemma}\label{hs:collar-identity}
Let $\Sigma_n$ be closed conformal surfaces of fixed topology converging in the Deligne--Mumford sense, and let $\Sigma_{\widehat v}$ be the compactified components of the normalization of the limit. Let $\eta_n$ be real harmonic $1$-forms on $\Sigma_n$ satisfying
\[
 \sup_n\|\eta_n\|_{L^2(\Sigma_n)}^2=\Lambda<\infty.
\]
All $L^2$ norms of $1$-forms are computed with smooth metrics in the indicated conformal classes and are independent of the choice of representative. For each pinched curve, choose a full Deligne--Mumford collar $C_{e,n}$ whose ends correspond to fixed puncture neighborhoods on the normalized limit. Write $\eta_n=A_{e,n}\,\dd s+B_{e,n}\,\dd\theta$ on $C_{e,n}$ and define $a_{e,n},b_{e,n}$ by \eqref{eq:hsl-collar-coefficients}. After passage to a subsequence, the pullbacks of $\eta_n$ by the Deligne--Mumford comparison maps converge smoothly on compact subsets of each punctured component to a harmonic $1$-form $\eta_{\widehat v}$. Each $\eta_{\widehat v}$ extends smoothly across the punctures, and
\begin{equation}\label{hs:closed-defect}
 \lim_n\|\eta_n\|_{L^2(\Sigma_n)}^2
 =\sum_{\widehat v}\|\eta_{\widehat v}\|_{L^2(\Sigma_{\widehat v})}^2
  +\lim_n\sum_e4\pi T_{e,n}(a_{e,n}^2+b_{e,n}^2).
\end{equation}
\end{lemma}

\begin{proof}
Smooth convergence of the component conformal structures and local harmonic estimates give convergence on a countable exhaustion of the punctured normalized components after one diagonal subsequence. The limits have finite $L^2$ norm. On a punctured conformal disk, the holomorphic differential associated with a real harmonic $1$-form has a Laurent expansion. Angular orthogonality shows that every negative power, including $dz/z$, would have infinite squared norm. Thus no negative powers occur, and the limits extend smoothly over the punctures.

On one collar, omit its indices and write $\eta=A\,ds+B\,d\theta$.
Closedness and coclosedness give
\[
 B_s-A_\theta=0,\qquad A_s+B_\theta=0.
\]
Hence the angular averages $a,b$ are constant in $s$, and
$F=A-iB$ is periodic and holomorphic in $s+i\theta$. Its expansion is
\[
 F(s+i\theta)=a-ib+
       \sum_{k\in\mathbb Z\setminus\{0\}}c_k e^{k(s+i\theta)}.
\]
The constant part has squared norm
\begin{equation}\label{hs:zero-mode-norm}
 \|a_{e,n}\,ds+b_{e,n}\,d\theta\|_{L^2(C_{e,n})}^2
       =4\pi T_{e,n}(a_{e,n}^2+b_{e,n}^2).
\end{equation}
All Fourier terms are mutually orthogonal after angular integration,
also on every subcylinder. For a nonzero mode, the ratio of its squared
norm on $(-T+R,T-R)\times S^1$ to that on $(-T,T)\times S^1$ is
\[
 \frac{\sinh(2|k|(T-R))}{\sinh(2|k|T)}
 \le e^{-2|k|R}\le e^{-2R},\qquad 0<R<T.
\]
It follows that the sum of the squared $L^2$ norms of the nonconstant Fourier parts on the subcylinders $(-T_{e,n}+R,T_{e,n}-R)\times S^1$ of all degenerating collars is at most $\Lambda e^{-2R}$.
Moreover, $4\pi T_{e,n}(a_{e,n}^2+b_{e,n}^2)\le \Lambda$, so
\begin{equation}\label{hs:mode-bound}
 a_{e,n}^2+b_{e,n}^2\le\frac{\Lambda}{4\pi T_{e,n}}.
\end{equation}
For fixed $R$, the squared $L^2$ norm of the constant Fourier part on the two end regions $(-T_{e,n},-T_{e,n}+R)\times S^1$ and $(T_{e,n}-R,T_{e,n})\times S^1$ is $4\pi R(a_{e,n}^2+b_{e,n}^2)$ and tends to zero.

Since there are finitely many components and collars, we can pass to a further subsequence on which $\|\eta_n\|_{L^2(\Sigma_n)}^2$ and every $4\pi T_{e,n}(a_{e,n}^2+b_{e,n}^2)$ converge.
For fixed $R$, the complement of the middle collars converges smoothly to compact component regions $K_{\widehat v}(R)$, which exhaust the punctured components as $R\to\infty$. Orthogonality and the preceding bounds then give
\[
 \left|\lim_n\left(\|\eta_n\|_{L^2(\Sigma_n)}^2-\sum_e4\pi T_{e,n}(a_{e,n}^2+b_{e,n}^2)\right)
       -\sum_{\widehat v}
          \|\eta_{\widehat v}\|_{L^2(K_{\widehat v}(R))}^2\right|
 \le \Lambda e^{-2R}.
\]
Taking first $n\to\infty$ with $R$ fixed and then $R\to\infty$ proves
\eqref{hs:closed-defect}. When there are no pinched curves, smooth convergence on the whole closed surface gives \eqref{hs:closed-defect} without a collar term.
\end{proof}

\begin{proof}[Proof of Theorem~\ref{app:hsl-energy}]
Applying Theorem~\ref{main theorem} to $u_n$ and Deligne--Mumford compactness to the complex doubles $\widehat\Sigma_n$, we pass to a subsequence with limiting component maps and normalized principal domains. Every principal domain is retained, including those with constant limiting maps and those that become unstable when the auxiliary marks are forgotten. If the images of $u_n$ collapse to a single point, this point defines the limiting map on each principal domain. Passing to a further subsequence makes the finite index sets for the principal domains and pinched collars independent of $n$.

A closed and coclosed $1$-form satisfying $\alpha_n(\nu_{h_n})=0$ reflects across $\partial\Sigma$ with even tangential coefficient and odd normal coefficient. The reflected form is weakly closed and coclosed across the seam, hence smooth by elliptic regularity. Its squared norm satisfies
\begin{equation}\label{hs:double-norm}
 \|\widehat\alpha_n\|_{L^2(\widehat\Sigma_n)}^2
 =2\|\alpha_n\|_{L^2(\Sigma,h_n)}^2=8W_M(u_n).
\end{equation}
Applying Lemma~\ref{hs:collar-identity} to the harmonic forms $\widehat\alpha_n$ gives smooth limits on the normalized components and the squared $L^2$ norm decomposition \eqref{hs:closed-defect}. Since the component comparison maps preserve the doubling involution, these limiting forms descend to the harmonic $1$-forms $\alpha_\infty^i$ on the principal domains in the statement. An invariant normalized component is the
double of its bordered quotient; two exchanged components correspond
to one closed principal domain. In either case their combined squared
norm is twice $\|\alpha_\infty^i\|_{L^2(\Sigma_\infty^i)}^2$.

Lemma~\ref{lem:hsl-local-convergence} identifies $\alpha_\infty^i$ with the geometric mean curvature form on each core of a nonconstant principal component away from its branch points and nodal, attaching, and concentration marks. Continuity identifies the smooth extensions across these special points. The same lemma gives smooth convergence on collapsed cores, including half disks along $\partial\Sigma_\infty^i$, although the constant limiting map has no geometric mean curvature $1$-form with which to identify $\alpha_\infty^i$.
Dividing \eqref{hs:closed-defect} by eight gives \eqref{hs:complete-formula}.

The mean curvature $1$-forms of nonconstant sphere and disk bubbles vanish by \eqref{sys2} and the argument following \eqref{hs:normalization}. At a bubble attachment point away from the nodes, the conformal structures converge on a coordinate disk in the principal domain; at a boundary attachment, we work on the complex double. Choose smooth metrics $\widehat h_n$ representing the doubled conformal structures whose pullbacks converge smoothly on this disk. Local harmonic estimates give smooth bounds for $\widehat\alpha_n$ there independently of concentration of the maps. The conformally invariant measures $|\widehat\alpha_n|_{\widehat h_n}^2\,\dd\vol_{\widehat h_n}$ therefore have no concentration at the attachment point. Lemma~\ref{hs:collar-identity} already accounts for the nodal ends. This proves parts \textup{(i)}--\textup{(iii)}, including the total image collapse alternative, and nonnegativity of the additional terms gives the stated criterion for the Willmore energy identity. The fixed-class and genus-zero assertions follow from Corollaries~\ref{hs:fixed-class-corollary} and~\ref{hs:planar-identity}, proved below.
\end{proof}

Bounded changes of the collar endpoints leave the last term in \eqref{hs:complete-formula} unchanged by \eqref{hs:mode-bound}. Reversing a collar coordinate changes signs of the coefficients but preserves their squared sum. The following corollary shows that the collar contribution vanishes when the full cohomology class, including its longitudinal periods, is fixed under markings adapted to the degeneration.

\begin{coro}\label{hs:fixed-class-corollary}
Under the assumptions of Theorem~\ref{app:hsl-energy}, suppose that
there are markings $m_n:\widehat\Sigma\to\widehat\Sigma_n$ from a fixed reference double $\widehat\Sigma$, commuting with the doubling involutions and adapted to the degeneration
(that is, their restrictions to each fixed compact core of the punctured normalized components agree with the chosen Deligne--Mumford comparison maps up to a fixed homotopy) such that
\begin{equation}\label{hs:fixed-class}
 m_n^*[\widehat\alpha_n]=\widehat c
       \quad\hbox{in }H^1(\widehat\Sigma;\R)
\end{equation}
for one fixed class. Then the collar term in
\eqref{hs:complete-formula} vanishes, and
\begin{equation}\label{hs:ghost-formula}
 \lim_nW_M(u_n)
 =\sum_{\substack{i\in\mathfrak I_{\mathrm{pr}}\\u_\infty^i\text{ nonconstant}}}W_M(u_\infty^i)
  +\frac14\sum_{\substack{i\in\mathfrak I_{\mathrm{pr}}\\u_\infty^i\text{ constant}}}\|\alpha_\infty^i\|_{L^2(\Sigma_\infty^i)}^2.
\end{equation}
In this case, $$\lim_nW_M(u_n)=\sum_{\substack{i\in\mathfrak I_{\mathrm{pr}}\\u_\infty^i\text{ nonconstant}}}W_M(u_\infty^i)$$ if and only if 
$[\alpha_\infty^i]=0$ in $H^1(\Sigma_\infty^i;\R)$ for every principal domain $\Sigma_\infty^i$ on which $u_\infty^i$ is constant. In particular, it
holds when all such domains have first Betti number zero.
\end{coro}

\begin{proof}
We first use the notation of Lemma~\ref{hs:collar-identity} and consider harmonic forms $\eta_n$ representing a fixed cohomology class $c$ under the chosen markings. Let $P\subset H_1(\Sigma_n;\R)$, identified by these markings, be the span of the pinched curves $\gamma_e$. Write $P^\perp$ for the subspace of $H^1(\Sigma_n;\R)$ consisting of classes that vanish on $P$. Let $G$ be the dual graph, whose vertices are the normalized components and whose edges are the nodes. If
$p_e=\langle c,[\gamma_e]\rangle$, closedness and Cauchy--Schwarz give
\[
 \|\eta_n\|_{L^2(C_{e,n})}^2
 \ge\frac{T_{e,n}}\pi|p_e|^2.
\]
Thus $p_e=0$, and $c\in P^\perp$. Cutting along the pinched curves
and capping the boundary circles gives the exact sequence
\begin{equation*}
 0\longrightarrow H^1(G;\R)\longrightarrow P^\perp
 \longrightarrow\bigoplus_{\widehat v}H^1(\Sigma_{\widehat v};\R)
 \longrightarrow0.
\end{equation*}
Choose homology bases on the normalized components, a basis of the span $P$ of the pinched curves, and closed curves on $\Sigma_n$ crossing the collars along the fundamental cycles determined by a spanning tree of $G$.
Together these curves form a homology basis of $\Sigma_n$.
A class vanishing on the component and pinched-curve cycles is determined by its periods along the chosen crossing curves, which can be realized by longitudinal forms supported in collars.
Component classes, in turn, have closed representatives vanishing near the capped points. Extending them by zero proves surjectivity. Period convergence identifies the component classes of $c$ with those of $\eta_{\widehat v}$.

For $\delta>0$, we truncate each $\eta_{\widehat v}$ near the punctures. In a disk of radius $2\delta$, write $\eta_{\widehat v}=df$, with $f(0)=0$, and replace it by $\eta_{\widehat v}-d(\chi_\delta f)$, where $\chi_\delta=1$ on the disk of radius $\delta$, vanishes outside the disk of radius $2\delta$, and $|d\chi_\delta|\le C/\delta$. Since $f=O(|z|)$, the change in squared norm is $o_\delta(1)$. Transferring these closed forms to the component regions of $\Sigma_n$ and extending by zero gives $\zeta_{n,\delta}$ with
\[
 \lim_n\|\zeta_{n,\delta}\|_{L^2(\Sigma_n)}^2
 =\sum_{\widehat v}\|\eta_{\widehat v}\|_{L^2(\Sigma_{\widehat v})}^2+o_\delta(1).
\]
The discrepancy $c-[\zeta_{n,\delta}]$ lies in $H^1(G;\R)$.
For fixed $\delta$, choose the portions of the crossing curves inside the component cores as fixed paths in the comparison coordinates.
Their integrals against $\zeta_{n,\delta}$ are uniformly bounded, while the periods of $c$ are fixed under the markings adapted to the degeneration.
Changing the collar portions by twists adds only pinched-curve cycles, on which $c$ has zero period; the remaining collar portions do not meet the support of $\zeta_{n,\delta}$.
Thus the periods of the discrepancy along the chosen crossing curves are uniformly bounded.
We can represent the discrepancy by an edge cochain $K_{e,n,\delta}$ that vanishes on the spanning tree and is bounded on the remaining edges.

In the collar portions disjoint from $\operatorname{supp}\zeta_{n,\delta}$,
choose smooth longitudinal functions with
\[
 \int\varphi_{e,n,\delta}(s)\,ds=1,
 \qquad
 \int|\varphi_{e,n,\delta}(s)|^2\,ds
       \le\frac{C_\delta}{T_{e,n}},
\]
supported away from their endpoints. Then
\[
 \xi_{n,\delta}=\zeta_{n,\delta}
       +\sum_eK_{e,n,\delta}\varphi_{e,n,\delta}(s)\,ds
\]
is smooth and closed and represents exactly $c$, including its periods
along the chosen crossing curves. The supports are disjoint and the added squared
norm tends to zero. Since the harmonic representative minimizes the squared $L^2$ norm among closed $1$-forms representing $c$, we obtain
\[
 \limsup_n\|\eta_n\|_{L^2(\Sigma_n)}^2
 \le\sum_{\widehat v}\|\eta_{\widehat v}\|_{L^2(\Sigma_{\widehat v})}^2+o_\delta(1).
\]
Let $\delta\to0$ and use local convergence and exhaustion for the
reverse inequality. Lemma~\ref{hs:collar-identity} now forces the
collar term to vanish. Apply this argument to
$\eta_n=\widehat\alpha_n$ and use \eqref{hs:double-norm}.
Finally, an absolute harmonic form vanishes exactly when its ordinary
cohomology class is zero.
\end{proof}

A class fixed before arbitrary mapping-class changes need not satisfy \eqref{hs:fixed-class}.
Small angular periods alone do not force the collar contribution to vanish: on a cylinder of length $2T_n$, the form $T_n^{-1/2}d\theta$ has angular period $2\pi T_n^{-1/2}\to0$ but squared $L^2$ norm $4\pi$.

On a principal domain with constant limiting map, the local Lagrangian phases of $u_n$ restrict the class $[\alpha_\infty^i]$. We derive the phase identity below with our sign and normalization conventions; see \cite{appCG} for its relation to mean curvature and Maslov classes.

\begin{prop}\label{hs:phase-support}
Under the assumptions of Theorem~\ref{app:hsl-energy}, let $\Sigma_\infty^i$ be a principal domain whose limiting map $u_\infty^i$ is constant, and let $\alpha_\infty^i$ be its limiting harmonic $1$-form. There exists a class
$\kappa_i\in H^1(\Sigma_\infty^i;\mathbb Z)$ such that
\begin{equation}\label{hs:local-phase-integrality}
 [\alpha_\infty^i]=-2\pi\kappa_i
       \quad\hbox{in }H^1(\Sigma_\infty^i;\R).
\end{equation}
If the supporting submanifold $\cK\subset M$ is a Lagrangian surface, then
$\alpha_\infty^i=0$ on every constant principal domain with nonempty boundary.
No vanishing conclusion is asserted for closed constant domains.
\end{prop}

\begin{proof}
Write $u_\infty^i\equiv q$. On a contractible neighborhood of $q\in M$, choose a smooth unit section $\Omega$ of the canonical bundle $\Lambda^{2,0}T^*M$, and define its real connection $1$-form $\mathfrak a$ by $\bar\nabla\Omega=i\mathfrak a\otimes\Omega$.
For an oriented Lagrangian immersion $u$, choose an oriented orthonormal frame $(e_1,e_2)$ for $g_u$. The normalized evaluation
\[
 e^{i\theta}=\frac{\Omega(du(e_1),du(e_2))}{|\Omega(du(e_1),du(e_2))|}
\]
defines its local Lagrangian phase.
At the point of differentiation, choose a tangent frame whose intrinsic derivative in the direction $X$ vanishes. Differentiating the determinant then gives
\[
 i\,d\theta(X)
 =i\mathfrak a(du(X))
   +i\sum_{j=1}^2\bar g(A_{u,M}(X,e_j),Jdu(e_j))
 =i\mathfrak a(du(X))-i\alpha(X).
\]
The second equality uses the symmetry of the Lagrangian cubic tensor,
which identifies the sum with
$\bar g(\mathbf H_{u,M},Jdu(X))=-\alpha(X)$.
Thus, for each local real lift of the phase,
\begin{equation}\label{hs:phase-connection}
 \alpha=u^*\mathfrak a-d\theta.
\end{equation}
Under $\Omega'=e^{i\chi}\Omega$, both terms change by
$d(\chi\circ u)$, so this identity is independent of the frame.

Choose finitely many smooth cycles generating $H_1(\Sigma_\infty^i;\mathbb Z)$ and avoiding all node, branch, and bubble marks, and place them in one fixed compact core. Restricting the comparison embeddings $\Psi_n$ to this core, put $\widetilde u_n=u_n\circ\Psi_n$ and $\widetilde\alpha_n=\Psi_n^*\alpha_n$. The collapsed alternative of Lemma~\ref{lem:hsl-local-convergence} gives smooth convergence $\widetilde u_n\to q$ and $\widetilde\alpha_n\to\alpha_\infty^i$ on the core. The image of the core eventually lies in the chosen canonical frame neighborhood, and for each chosen cycle $\gamma$,
\[
 \left|\int_\gamma\widetilde u_n^*\mathfrak a\right|
 \le C\operatorname{Length}_{h_0}(\gamma)
                  \sup_\gamma|d\widetilde u_n|_{h_0}\longrightarrow0,
\]
where $h_0$ is a fixed smooth metric on $\Sigma_\infty^i$. Integrating
\eqref{hs:phase-connection} gives
\[
 \int_\gamma\widetilde\alpha_n
 =\int_\gamma\widetilde u_n^*\mathfrak a-2\pi k_n(\gamma),
 \qquad k_n(\gamma)\in\mathbb Z.
\]
Convergence of both integrals forces $k_n(\gamma)$ to stabilize for each chosen cycle $\gamma$. The stabilized integers define $\kappa_i$ and yield \eqref{hs:local-phase-integrality}. Since $\alpha_\infty^i$ extends smoothly, its periods around the filled punctures vanish. The factor $2\pi$ comes from using the oriented phase rather than its square.

Suppose now that $\cK$ is Lagrangian, and use its mean curvature $1$-form $\alpha_{\cK}$ defined in the proof of Corollary~\ref{coro:hsl-support-geometries}. At orthogonal contact, \eqref{eq:hsl-lagrangian-boundary-form} holds away from the branch points. This identity is local and does not require $\cK$ to be orientable.

A bordered normalized principal domain contains an open boundary arc inherited from $\partial\Sigma$ and avoiding the finitely many marks; the comparison maps preserve this boundary. The normalized preimage of a type~II node is an interior point and cannot create such an arc. Choose a half-disk core at an inherited boundary arc of a principal domain $\Sigma_\infty^i$ with $u_\infty^i\equiv q$. Smooth collapsed convergence and \eqref{eq:hsl-lagrangian-boundary-form} give $\alpha_\infty^i(\partial_x)=0$ there, while the zero normal trace gives $\alpha_\infty^i(\partial_y)=0$ in a conformal half-disk chart. Write $\alpha_\infty^i=a\,dx+b\,dy$. Closedness and coclosedness make $a-ib$ holomorphic. Reflecting with tangential coefficient $a$ even and normal coefficient $b$ odd extends this holomorphic function across the arc, where it vanishes. The identity theorem and interior continuation give $\alpha_\infty^i=0$ on the connected component.
\end{proof}

\begin{coro}\label{hs:planar-identity}
Assume the hypotheses of Corollary~\ref{hs:fixed-class-corollary}.
If the compact bordered surface $\Sigma$ has genus zero and the supporting submanifold $\cK\subset M$ is a Lagrangian surface, then
\[
 \lim_nW_M(u_n)
 =\sum_{i:\,u_\infty^i\text{ nonconstant}}
                         W_M(u_\infty^i).
\]
In total image collapse this limit is zero.
\end{coro}

\begin{proof}
Cutting a genus-zero bordered surface along the degenerating curves
and arcs and filling punctures creates no handles. Every closed
normalized principal domain is therefore a sphere, whose harmonic
$1$-forms vanish. Proposition~\ref{hs:phase-support} gives $\alpha_\infty^i=0$ on every bordered principal domain with constant map. Corollary~\ref{hs:fixed-class-corollary} makes the collar term vanish under the stated cohomology hypothesis. Formula~\eqref{hs:ghost-formula} proves the claim,
since all bubbles are minimal.
\end{proof}

The following family of free boundary HSL annuli in a compact flat K\"ahler surface collapses to a point while its Willmore energy remains $2\pi^2$. Let
\[
 E=\C/(10\mathbb Z+10i\mathbb Z),\qquad
 M=E\times E,\qquad
 \cK=\{([z],[z]):[z]\in E\},
\]
with its flat product K\"ahler metric. The supporting submanifold $\cK\subset M$ is the compact complex diagonal. On the fixed annulus
$\Sigma=(\R/2\pi\mathbb Z)\times[0,\pi]$, define, for small $r>0$,
\begin{equation*}
 u_r(s,t)=\bigl([r e^{i(s+t)}],[r e^{i(s-t)}]\bigr).
\end{equation*}
The images lie in an embedded flat coordinate neighborhood of the
origin. The derivatives are
\[
 (u_r)_s=ir\bigl(e^{i(s+t)},e^{i(s-t)}\bigr),\qquad
 (u_r)_t=ir\bigl(e^{i(s+t)},-e^{i(s-t)}\bigr).
\]
They are orthogonal, have equal length, and span a Lagrangian plane.
In these coordinates,
\[
 g_{u_r}=2r^2(ds^2+dt^2),\qquad
 \mathbf H_{u_r,M}=-r^{-2}u_r,\qquad \alpha_r=-2\,ds,
\]
where the middle formula uses the local lift of $u_r$.
Thus $d\alpha_r=d^*\alpha_r=0$ and $\alpha_r(\nu)=0$.
At $t=0,\pi$, both coordinates agree, while $(u_r)_t$ is complex
antidiagonal and hence normal to $\cK$. The orthogonal free boundary condition follows.
The two boundary circles of $\Sigma$ have the same image with a phase
shift.
Finally, we see that 
\begin{equation*}
 \cA(u_r)=4\pi^2r^2\longrightarrow0,\qquad
 W_M(u_r)=2\pi^2.
\end{equation*}
The conformal structure and the reflected cohomology class are fixed. Nevertheless $u_r\to([0],[0])$, while
$\alpha_\infty=-2\,ds$ on the constant principal annulus and
$\frac14\|\alpha_\infty\|_{L^2(\Sigma,ds^2+dt^2)}^2=2\pi^2$. This shows why the term on the domain of constant mapping in \eqref{hs:complete-formula} must be retained even without conformal structure degeneration.

\vskip1cm

% Keep every entry listed in the supplied manuscript, including uncited items.
% \nocite{*}
\bibliographystyle{amsalpha-fullnames}
\bibliography{references}

@article{sacks1981existence,
  author = {Sacks, Jonathan and Uhlenbeck, Karen Keskulla},
  title = {The existence of minimal immersions of {$2$}-spheres},
  journal = {Ann. of Math. (2)},
  fjournal = {Annals of Mathematics. Second Series},
  volume = {113},
  year = {1981},
  number = {1},
  pages = {1--24},
  issn = {0003-486X},
  mrclass = {58E12 (53C42 58E20)},
  mrnumber = {604040},
  mrreviewer = {John\ C.\ Wood},
}

@article{Laurain-Petrides,
  author = {Laurain, Paul and Petrides, Romain},
  title = {Existence of min-max free boundary disks realizing the width of a manifold},
  journal = {Adv. Math.},
  fjournal = {Advances in Mathematics},
  volume = {352},
  year = {2019},
  pages = {326--371},
  issn = {0001-8708,1090-2082},
  mrclass = {53C42 (35J20 35R01 49Q05 58E20)},
  mrnumber = {3961741},
  mrreviewer = {Fernando\ Manfio},
}

@article{gruter-hildebrandt-nitsche1986,
  author = {Gr{\"u}ter, Michael and Hildebrandt, Stefan and Nitsche, Johannes C. C.},
  title = {Regularity for stationary surfaces of constant mean curvature with free boundaries},
  journal = {Acta Math.},
  fjournal = {Acta Mathematica},
  volume = {156},
  year = {1986},
  number = {1--2},
  pages = {119--152},
  doi = {10.1007/BF02399202},
  url = {https://doi.org/10.1007/BF02399202},
  issn = {0001-5962,1871-2509},
  mrclass = {49F22 (53A10)},
  mrnumber = {822332},
  mrreviewer = {Michael\ T.\ Anderson},
}

@article{Li-Yau1982,
  author = {Li, Peter and Yau, Shing-Tung},
  title = {A new conformal invariant and its applications to the {W}illmore conjecture and the first eigenvalue of compact surfaces},
  journal = {Invent. Math.},
  fjournal = {Inventiones Mathematicae},
  volume = {69},
  year = {1982},
  number = {2},
  pages = {269--291},
  issn = {0020-9910,1432-1297},
  mrclass = {53C40 (53C21 58G25)},
  mrnumber = {674407},
  mrreviewer = {Yu.\ Burago},
}

@article{Marques-Neves2014,
  author = {Marques, Fernando Cod\'{a} and Neves, Andr\'{e} Arroja},
  title = {Min-max theory and the {W}illmore conjecture},
  journal = {Ann. of Math. (2)},
  fjournal = {Annals of Mathematics. Second Series},
  volume = {179},
  year = {2014},
  number = {2},
  pages = {683--782},
  doi = {10.4007/annals.2014.179.2.6},
  url = {https://doi.org/10.4007/annals.2014.179.2.6},
  issn = {0003-486X,1939-8980},
  mrclass = {53C42 (49Q20)},
  mrnumber = {3152944},
  mrreviewer = {Andrea\ Mondino},
}

@article{JostLiuZhu2019,
  author = {Jost, J{\"u}rgen and Liu, Lei and Zhu, Miaomiao},
  title = {Asymptotic analysis and qualitative behavior at the free boundary for {S}acks-{U}hlenbeck {$\alpha$}-harmonic maps},
  journal = {Adv. Math.},
  fjournal = {Advances in Mathematics},
  volume = {396},
  year = {2022},
  eid = {108105},
  pages = {Paper No.~108105},
  pagetotal = {68},
  issn = {0001-8708,1090-2082},
  mrclass = {58E20 (53C43)},
  mrnumber = {4370467},
  mrreviewer = {Andreas\ Gastel},
}

@article{Colding-Minicozzi08,
  author = {Colding, Tobias H. and Minicozzi, II, William P.},
  title = {Width and finite extinction time of {R}icci flow},
  journal = {Geom. Topol.},
  fjournal = {Geometry \& Topology},
  volume = {12},
  year = {2008},
  number = {5},
  pages = {2537--2586},
  doi = {10.2140/gt.2008.12.2537},
  url = {https://doi.org/10.2140/gt.2008.12.2537},
  issn = {1465-3060,1364-0380},
  mrclass = {53C44},
  mrnumber = {2460871},
  mrreviewer = {Andrea\ Nicole\ Young},
}

@article{chen2014,
  author = {Chen, Jingyi and Li, Yuxiang},
  title = {{Bubble tree of branched conformal immersions and applications to the {W}illmore functional}},
  journal = {Amer. J. Math.},
  volume = {136},
  year = {2014},
  number = {4},
  pages = {1107--1154},
  doi = {10.1353/ajm.2014.0023},
}

@article{chen2021,
  author = {Chen, Jingyi and Ma, John Man Shun},
  title = {{On the compactness of {H}amiltonian stationary {L}agrangian surfaces in {K}\"ahler surfaces}},
  journal = {Calc. Var. Partial Differential Equations},
  volume = {60},
  year = {2021},
  number = {2},
  eid = {75},
  pages = {Paper No.~75},
  pagetotal = {23},
  doi = {10.1007/s00526-020-01901-7},
}

@article{Chen-Ma-2024,
  author = {Chen, Jingyi and Ma, John Man Shun},
  title = {{Compactness of {H}amiltonian stationary {L}agrangian submanifolds in symplectic manifolds}},
  journal = {Trans. Amer. Math. Soc.},
  volume = {377},
  year = {2024},
  number = {9},
  pages = {6151--6170},
}

@article{Chen-Warren,
  author = {Chen, Jingyi and Warren, Micah},
  title = {{Compactification of the space of {H}amiltonian stationary {L}agrangian submanifolds with bounded total extrinsic curvature and volume}},
  journal = {J. Differential Geom.},
  volume = {126},
  year = {2024},
  number = {1},
  pages = {65--97},
}

@book{Courant1950book,
  author = {Courant, Richard},
  title = {{Dirichlet's principle, conformal mapping, and minimal surfaces}},
  year = {1950},
  publisher = {Interscience Publishers, Inc.},
  address = {New York, N.Y.},
  note = {Appendix by Menahem Schiffer},
}

@article{Courant-Davids,
  author = {Courant, Richard and Davids, Norman},
  title = {{Minimal surfaces spanning closed manifolds}},
  journal = {Proc. Nat. Acad. Sci. U.S.A.},
  volume = {26},
  year = {1940},
  pages = {194--199},
}

@article{Dazord1981,
  author = {Dazord, Pierre},
  title = {{Sur la g\'eom\'etrie des sous-fibr\'es et des feuilletages lagrangiens}},
  journal = {Ann. Sci. \'Ecole Norm. Sup. (4)},
  volume = {14},
  year = {1981},
  number = {4},
  pages = {465--480},
}

@article{douglas1931mapping,
  author = {Douglas, Jesse},
  title = {{The mapping theorem of Koebe and the problem of Plateau}},
  journal = {J. Math. Phys.},
  volume = {10},
  year = {1931},
  number = {1--4},
  pages = {106--130},
}

@article{FraserCPAM,
  author = {Fraser, Ailana M.},
  title = {{On the free boundary variational problem for minimal disks}},
  journal = {Comm. Pure Appl. Math.},
  volume = {53},
  year = {2000},
  number = {8},
  pages = {931--971},
}

@article{Fraser-Schoen-2016,
  author = {Fraser, Ailana M. and Schoen, Richard M.},
  title = {{Sharp eigenvalue bounds and minimal surfaces in the ball}},
  journal = {Invent. Math.},
  volume = {203},
  year = {2016},
  number = {3},
  pages = {823--890},
}

@article{Gaia2025,
  author = {Gaia, Filippo},
  title = {{Free boundary Hamiltonian stationary {L}agrangian discs in {$\mathbb C^2$}}},
  journal = {J. Geom. Anal.},
  volume = {35},
  year = {2025},
  number = {5},
  eid = {160},
  pages = {Paper No.~160},
  pagetotal = {17},
}

@article{Gao-Wang-Zhucvpde,
  author = {Gao, Rui and Wang, Zijun and Zhu, Miaomiao},
  title = {{Bubbling analysis at the free boundary for approximate {$H$}-surfaces}},
  journal = {Calc. Var. Partial Differential Equations},
  volume = {64},
  year = {2025},
  number = {1},
  eid = {12},
  pages = {Paper No.~12},
  pagetotal = {17},
}

@article{Gruter-Jost,
  author = {Gr{\"u}ter, Michael and Jost, J{\"u}rgen},
  title = {{On embedded minimal disks in convex bodies}},
  journal = {Ann. Inst. H. Poincar\'e C Anal. Non Lin\'eaire},
  volume = {3},
  year = {1986},
  number = {5},
  pages = {345--390},
}

@book{helein,
  author = {H{\'e}lein, Fr{\'e}d{\'e}ric},
  title = {{Harmonic maps, conservation laws and moving frames}},
  edition = {2nd},
  series = {Cambridge Tracts in Mathematics},
  volume = {150},
  year = {2002},
  publisher = {Cambridge University Press},
  address = {Cambridge},
  doi = {10.1017/CBO9780511543036},
}

@article{HeleinRomon2002,
  author = {H{\'e}lein, Fr{\'e}d{\'e}ric and Romon, Pascal},
  title = {{Hamiltonian stationary {L}agrangian surfaces in {$\mathbb C^2$}}},
  journal = {Comm. Anal. Geom.},
  volume = {10},
  year = {2002},
  number = {1},
  pages = {79--126},
}

@article{Hildebrandt-Nitsche1979,
  author = {Hildebrandt, Stefan and Nitsche, Johannes C. C.},
  title = {{Minimal surfaces with free boundaries}},
  journal = {Acta Math.},
  volume = {143},
  year = {1979},
  number = {3--4},
  pages = {251--272},
}

@book{hummel1997,
  author = {Hummel, Christoph},
  title = {{Gromov's compactness theorem for pseudo-holomorphic curves}},
  series = {Progress in Mathematics},
  volume = {151},
  year = {1997},
  publisher = {Birkh\"auser Verlag},
  address = {Basel},
  doi = {10.1007/978-3-0348-8952-0},
}

@article{Jost-Liu-Zhu2019,
  author = {Jost, J{\"u}rgen and Liu, Lei and Zhu, Miaomiao},
  title = {{The qualitative behavior at the free boundary for approximate harmonic maps from surfaces}},
  journal = {Math. Ann.},
  volume = {374},
  year = {2019},
  number = {1--2},
  pages = {133--177},
}

@article{Katz-Liu2001,
  author = {Katz, Sheldon and Liu, Chiu-Chu Melissa},
  title = {{Enumerative geometry of stable maps with {L}agrangian boundary conditions and multiple covers of the disc}},
  journal = {Adv. Theor. Math. Phys.},
  volume = {5},
  year = {2001},
  number = {1},
  pages = {1--49},
}

@article{kuwert2012,
  author = {Kuwert, Ernst and Li, Yuxiang},
  title = {{{$W^{2,2}$}-conformal immersions of a closed {R}iemann surface into {$\mathbb R^n$}}},
  journal = {Comm. Anal. Geom.},
  volume = {20},
  year = {2012},
  number = {2},
  pages = {313--340},
}

@article{Lewy1951CPAM,
  author = {Lewy, Hans},
  title = {{On minimal surfaces with partially free boundary}},
  journal = {Comm. Pure Appl. Math.},
  volume = {4},
  year = {1951},
  number = {1},
  pages = {1--13},
}

@article{Lewy1951PNAS,
  author = {Lewy, Hans},
  title = {{On the boundary behavior of minimal surfaces}},
  journal = {Proc. Nat. Acad. Sci. U.S.A.},
  volume = {37},
  year = {1951},
  pages = {103--110},
}

@article{LiWangWeng2021,
  author = {Li, Mingyan and Wang, Guofang and Weng, Liangjun},
  title = {{Lagrangian surfaces with {L}egendrian boundary}},
  journal = {Sci. China Math.},
  volume = {64},
  year = {2021},
  number = {7},
  pages = {1589--1598},
}

@article{Liu-Song-Zhu2022,
  author = {Liu, Lei and Song, Chong and Zhu, Miaomiao},
  title = {{Harmonic maps with free boundary from degenerating bordered {R}iemann surfaces}},
  journal = {J. Geom. Anal.},
  volume = {32},
  year = {2022},
  number = {2},
  eid = {49},
  pages = {Paper No.~49},
  pagetotal = {21},
}

@article{LuoSun2021,
  author = {Luo, Yong and Sun, Linlin},
  title = {{Rigidity theorems for minimal {L}agrangian surfaces with {L}egendrian capillary boundary}},
  journal = {Adv. Math.},
  volume = {393},
  year = {2021},
  eid = {108124},
  pages = {Paper No.~108124},
  pagetotal = {15},
}

@article{Oh1993,
  author = {Oh, Yong-Geun},
  title = {{Volume minimization of {L}agrangian submanifolds under {H}amiltonian deformations}},
  journal = {Math. Z.},
  volume = {212},
  year = {1993},
  number = {2},
  pages = {175--192},
}

@article{Oh1994,
  author = {Oh, Yong-Geun},
  title = {{Mean curvature vector and symplectic topology of {L}agrangian submanifolds in {E}instein-{K}\"ahler manifolds}},
  journal = {Math. Z.},
  volume = {216},
  year = {1994},
  number = {3},
  pages = {471--482},
}

@article{rado1930plateau,
  author = {Rad{\'o}, Tibor},
  title = {{On Plateau's problem}},
  journal = {Ann. of Math. (2)},
  volume = {31},
  year = {1930},
  number = {3},
  pages = {457--469},
}

@article{rado1932contributions,
  author = {Rad{\'o}, Tibor},
  title = {{Contributions to the theory of minimal surfaces}},
  journal = {Acta Litt. Sci. Univ. Szeged},
  volume = {6},
  year = {1932},
  number = {1},
  pages = {1--20},
}

@incollection{Schoen2005,
  author = {Schoen, Richard M.},
  title = {{Special {L}agrangian submanifolds}},
  booktitle = {{Global theory of minimal surfaces}},
  series = {Clay Math. Proc.},
  volume = {2},
  year = {2005},
  pages = {655--666},
  publisher = {Amer. Math. Soc.},
  address = {Providence, RI},
}

@article{schoen2001,
  author = {Schoen, Richard M. and Wolfson, Jon G.},
  title = {{Minimizing area among {L}agrangian surfaces: the mapping problem}},
  journal = {J. Differential Geom.},
  volume = {58},
  year = {2001},
  number = {1},
  pages = {1--86},
}

@article{Smyth1984,
  author = {Smyth, Brian},
  title = {{Stationary minimal surfaces with boundary on a simplex}},
  journal = {Invent. Math.},
  volume = {76},
  year = {1984},
  number = {3},
  pages = {411--420},
}

@article{Struwe1984invention,
  author = {Struwe, Michael},
  title = {{On a free boundary problem for minimal surfaces}},
  journal = {Invent. Math.},
  volume = {75},
  year = {1984},
  number = {3},
  pages = {547--560},
}

@article{Ye1991,
  author = {Ye, Rugang},
  title = {{On the existence of area-minimizing surfaces with free boundary}},
  journal = {Math. Z.},
  volume = {206},
  year = {1991},
  number = {3},
  pages = {321--331},
}

@article{Fraser2002,
  author = {Fraser, Ailana M.},
  title = {{Minimal disks and two-convex hypersurfaces}},
  journal = {Amer. J. Math.},
  volume = {124},
  year = {2002},
  number = {3},
  pages = {483--493},
  doi = {10.1353/ajm.2002.0013},
}

@article{Lin-Sun-Zhou2020,
  author = {Lin, Longzhi and Sun, Ao and Zhou, Xin},
  title = {{Min-max minimal disks with free boundary in Riemannian manifolds}},
  journal = {Geom. Topol.},
  volume = {24},
  year = {2020},
  number = {1},
  pages = {471--532},
  doi = {10.2140/gt.2020.24.471},
}

@article{Chen-Fraser-Pang2015,
  author = {Chen, Jingyi and Fraser, Ailana M. and Pang, Chao},
  title = {{Minimal immersions of compact bordered Riemann surfaces with free boundary}},
  journal = {Trans. Amer. Math. Soc.},
  volume = {367},
  year = {2015},
  number = {4},
  pages = {2487--2507},
  doi = {10.1090/S0002-9947-2014-05990-4},
}

@article{appAK,
  author = {Alessandroni, Roberta and Kuwert, Ernst},
  title = {{Local solutions to a free boundary problem for the {W}illmore functional}},
  journal = {Calc. Var. Partial Differential Equations},
  volume = {55},
  year = {2016},
  number = {2},
  eid = {24},
  pages = {Paper No.~24},
  pagetotal = {29},
  doi = {10.1007/s00526-016-0961-3},
  url = {https://doi.org/10.1007/s00526-016-0961-3},
}

@article{appKM,
  author = {Kuwert, Ernst and M{\"u}ller, Marius},
  title = {{Curvature varifolds with orthogonal boundary}},
  journal = {J. Lond. Math. Soc. (2)},
  volume = {110},
  year = {2024},
  number = {3},
  eid = {e12976},
  pages = {Paper No.~e12976},
  doi = {10.1112/jlms.12976},
  url = {https://doi.org/10.1112/jlms.12976},
}

@article{appMedvedev,
  author = {Medvedev, Vladimir},
  title = {{Degenerating sequences of conformal classes and the conformal {S}teklov spectrum}},
  journal = {Canad. J. Math.},
  volume = {74},
  year = {2022},
  number = {4},
  pages = {1093--1136},
  doi = {10.4153/S0008414X21000171},
  url = {https://doi.org/10.4153/S0008414X21000171},
}

@article{appKokarev,
  author = {Kokarev, Gerasim},
  title = {{Variational aspects of {L}aplace eigenvalues on {R}iemannian surfaces}},
  journal = {Adv. Math.},
  volume = {258},
  year = {2014},
  pages = {191--239},
  doi = {10.1016/j.aim.2014.03.006},
  url = {https://doi.org/10.1016/j.aim.2014.03.006},
}

@article{appGKL,
  author = {Girouard, Alexandre and Karpukhin, Mikhail and Lagac{\'e}, Jean},
  title = {{Continuity of eigenvalues and shape optimisation for {L}aplace and {S}teklov problems}},
  journal = {Geom. Funct. Anal.},
  volume = {31},
  year = {2021},
  number = {3},
  pages = {513--561},
  doi = {10.1007/s00039-021-00573-5},
  url = {https://doi.org/10.1007/s00039-021-00573-5},
}

@article{appNovagaPozzetta,
  author = {Novaga, Matteo and Pozzetta, Marco},
  title = {{Connected surfaces with boundary minimizing the {W}illmore energy}},
  journal = {Math. Eng.},
  volume = {2},
  year = {2020},
  number = {3},
  pages = {527--556},
  doi = {10.3934/mine.2020024},
  url = {https://doi.org/10.3934/mine.2020024},
}

@article{appMS,
  author = {Mondino, Andrea and Scharrer, Christian},
  title = {{Existence and regularity of spheres minimising the {C}anham--{H}elfrich energy}},
  journal = {Arch. Ration. Mech. Anal.},
  volume = {236},
  year = {2020},
  number = {3},
  pages = {1455--1485},
  doi = {10.1007/s00205-020-01497-4},
  url = {https://doi.org/10.1007/s00205-020-01497-4},
}

@misc{BiZhou2025,
  author = {Bi, Yuchen and Zhou, Jie},
  title = {{Quantitative stability of the {C}lifford torus as a {W}illmore minimizer}},
  year = {2025},
  url = {https://arxiv.org/abs/2511.19681v1},
  howpublished = {\href{https://arxiv.org/abs/2511.19681v1}{arXiv:2511.19681v1}},
  eprint = {2511.19681},
  archiveprefix = {arXiv},
  primaryclass = {math.DG},
}

@article{appLM,
  author = {Lamm, Tobias and Metzger, Jan},
  title = {{Minimizers of the {W}illmore functional with a small area constraint}},
  journal = {Ann. Inst. H. Poincar{\'e} C Anal. Non Lin{\'e}aire},
  volume = {30},
  year = {2013},
  number = {3},
  pages = {497--518},
  doi = {10.1016/j.anihpc.2012.10.003},
  url = {https://doi.org/10.1016/j.anihpc.2012.10.003},
}

@article{appDLM05,
  author = {De Lellis, Camillo and M{\"u}ller, Stefan},
  title = {{Optimal rigidity estimates for nearly umbilical surfaces}},
  journal = {J. Differential Geom.},
  volume = {69},
  year = {2005},
  number = {1},
  pages = {75--110},
  doi = {10.4310/jdg/1121540340},
  url = {https://doi.org/10.4310/jdg/1121540340},
}

@article{appDLM06,
  author = {De Lellis, Camillo and M{\"u}ller, Stefan},
  title = {{A {$C^0$} estimate for nearly umbilical surfaces}},
  journal = {Calc. Var. Partial Differential Equations},
  volume = {26},
  year = {2006},
  number = {3},
  pages = {283--296},
  note = {see also the \href{https://www.math.ias.edu/delellis/sites/math.ias.edu.delellis/files/Errata-umbilical-2_0.pdf}{author's correction}},
  doi = {10.1007/s00526-006-0005-5},
  url = {https://doi.org/10.1007/s00526-006-0005-5},
}

@article{appMR,
  author = {Mondino, Andrea and Rivi{\`e}re, Tristan},
  title = {{Willmore spheres in compact {R}iemannian manifolds}},
  journal = {Adv. Math.},
  volume = {232},
  year = {2013},
  number = {1},
  pages = {608--676},
  doi = {10.1016/j.aim.2012.09.014},
  url = {https://doi.org/10.1016/j.aim.2012.09.014},
}

@article{appMetsch,
  author = {Metsch, Jan-Henrik},
  title = {{On the area-preserving {W}illmore flow of small bubbles sliding on a domain's boundary}},
  journal = {Adv. Calc. Var.},
  volume = {17},
  year = {2024},
  number = {4},
  pages = {1445--1486},
  doi = {10.1515/acv-2023-0023},
  url = {https://doi.org/10.1515/acv-2023-0023},
}

@misc{appDRSS,
  author = {Dall'Acqua, Anna and Rupp, Fabian and Sch{\"a}tzle, Reiner Michael and Schlierf, Manuel},
  title = {{Stability of the free boundary {W}illmore problem}},
  year = {2026},
  url = {https://arxiv.org/abs/2601.18388v1},
  howpublished = {\href{https://arxiv.org/abs/2601.18388v1}{arXiv:2601.18388v1}},
  eprint = {2601.18388},
  archiveprefix = {arXiv},
  primaryclass = {math.AP},
}

@article{appMiura,
  author = {Miura, Tatsuya},
  title = {{A diameter bound for compact surfaces and the {P}lateau--{D}ouglas problem}},
  journal = {Ann. Sc. Norm. Super. Pisa Cl. Sci. (5)},
  volume = {23},
  year = {2022},
  number = {4},
  pages = {1707--1721},
  note = {Author version: \href{https://arxiv.org/abs/2010.02797v2}{arXiv:2010.02797v2}},
  doi = {10.2422/2036-2145.202011_006},
  url = {https://doi.org/10.2422/2036-2145.202011_006},
}

@article{appMenneScharrer,
  author = {Menne, Ulrich and Scharrer, Christian},
  title = {{A priori bounds for geodesic diameter. {P}art {III}. {A} {S}obolev--{P}oincar{\'e} inequality and applications to a variety of geometric variational problems}},
  journal = {Duke Math. J.},
  volume = {175},
  year = {2026},
  number = {9},
  pages = {1649--1714},
  note = {Author version: \href{https://arxiv.org/abs/1709.05504v3}{arXiv:1709.05504v3}},
  doi = {10.1215/00127094-2025-0051},
  url = {https://doi.org/10.1215/00127094-2025-0051},
}

@article{appBMMS,
  author = {Buser, Peter and Makover, Eran and Muetzel, Bjoern and Silhol, Robert},
  title = {{Energy distribution of harmonic {$1$}-forms and {J}acobians of {R}iemann surfaces with a short closed geodesic}},
  journal = {Math. Z.},
  volume = {297},
  year = {2021},
  number = {3--4},
  pages = {1899--1952},
  doi = {10.1007/s00209-020-02584-8},
  url = {https://doi.org/10.1007/s00209-020-02584-8},
}

@article{appCG,
  author = {Cieliebak, Kai and Goldstein, Edward},
  title = {{A note on mean curvature, {M}aslov class and symplectic area of {L}agrangian immersions}},
  journal = {J. Symplectic Geom.},
  volume = {2},
  year = {2004},
  number = {2},
  pages = {261--266},
  doi = {10.4310/JSG.2004.v2.n2.a4},
  url = {https://arxiv.org/abs/math/0310046v3},
}

@article{BettiolPiccioneSantoro2014,
  author = {Bettiol, Renato G. and Piccione, Paolo and Santoro, Bianca},
  title = {Equivariant deformations of {Hamiltonian} stationary {Lagrangian} submanifolds},
  journal = {Mat. Contemp.},
  volume = {43},
  year = {2014},
  pages = {61--88},
  url = {https://arxiv.org/abs/1302.6970},
  eprint = {1302.6970},
  archiveprefix = {arXiv},
  primaryclass = {math.DG},
}

@article{HarveyLawson1982,
  author = {Harvey, Reese and Lawson, Jr., H. Blaine},
  title = {Calibrated geometries},
  journal = {Acta Math.},
  volume = {148},
  year = {1982},
  pages = {47--157},
  doi = {10.1007/BF02392726},
  url = {https://doi.org/10.1007/BF02392726},
}

@misc{Liu2004,
  author = {Liu, Chiu-Chu Melissa},
  title = {{Moduli of {J}-holomorphic curves with {L}agrangian boundary conditions and open {G}romov--{W}itten invariants for an {$S^1$}-equivariant pair}},
  year = {2004},
  url = {https://arxiv.org/abs/math/0210257v2},
  howpublished = {\href{https://arxiv.org/abs/math/0210257v2}{arXiv:math/0210257v2}},
  eprint = {math/0210257},
  archiveprefix = {arXiv},
  primaryclass = {math.SG},
}

@book{AdamsFournier2003,
  author = {Adams, Robert A. and Fournier, John J. F.},
  title = {{Sobolev spaces}},
  edition = {2nd},
  series = {Pure and Applied Mathematics},
  volume = {140},
  year = {2003},
  publisher = {Elsevier/Academic Press},
  address = {Amsterdam},
  isbn = {978-0-12-044143-3},
}

@book{GilbargTrudinger2001,
  author = {Gilbarg, David and Trudinger, Neil S.},
  title = {{Elliptic partial differential equations of second order}},
  edition = {2nd},
  series = {Classics in Mathematics},
  year = {2001},
  publisher = {Springer-Verlag},
  address = {Berlin},
  note = {reprint of the 1998 edition},
  doi = {10.1007/978-3-642-61798-0},
}

@book{Folland1995,
  author = {Folland, Gerald B.},
  title = {{Introduction to partial differential equations}},
  edition = {2nd},
  year = {1995},
  publisher = {Princeton University Press},
  address = {Princeton, NJ},
  isbn = {978-0-691-04361-6},
}

@book{Vekua1962,
  author = {Vekua, Il'ya Nestorovich},
  title = {{Generalized analytic functions}},
  series = {International Series of Monographs on Pure and Applied Mathematics},
  volume = {25},
  year = {1962},
  publisher = {Pergamon Press},
  address = {Oxford},
  isbn = {978-0-08-009693-3},
}

@incollection{Allard1987,
  author    = {Allard, William K.},
  title     = {Notes on the theory of varifolds},
  booktitle = {Th\'eorie des vari\'et\'es minimales et applications},
  series    = {Ast\'erisque},
  number    = {154--155},
  pages     = {73--93},
  year      = {1987},
  publisher = {Soci\'et\'e Math\'ematique de France},
  url       = {https://www.numdam.org/item/AST_1987__154-155__73_0/}
}

\end{document}